\documentclass[11.5pt]{report}
\usepackage[francais]{babel}
\usepackage{geometry}

\usepackage[latin1]{inputenc}
\usepackage[francais]{babel}
\usepackage{float}
\usepackage{multirow}
\usepackage{fancyhdr}
\usepackage{amsmath}
\usepackage{amssymb}
\usepackage{amsfonts}
\usepackage{amsthm}
\usepackage{mathrsfs}
\usepackage{dsfont}
\usepackage{latexsym}
\usepackage{graphicx}
\usepackage{natbibf}
\usepackage{fancybox}
\usepackage{epsfig}
\usepackage{vmargin}
\usepackage{color}
\usepackage{multicol}
\usepackage{rotating}

\theoremstyle{plain}

\newcommand{\eop}{\hspace*{\fill} \ensuremath{\Box}}

\newtheorem{theoreme}{Theorème}[chapter]
\newtheorem{theorem}{Theorem}[chapter]
\newtheorem{lemma}{Lemma}[chapter]

\newtheorem{corollaire}{Corollaire}[chapter]
\newtheorem{proposition}{Proposition}[chapter]

\newtheorem{e-definition}{Definition}[chapter]

\newcommand \Cov {\mathrm{Cov}}

\newcommand \Var {\mathrm{Var}}

\def\og~{\guillemotleft}

\def\Rit{\mathbb{R}}
\def\prob{\mathbb{P}}

\def\esp{\mathbb{E}}
\def\Var{\hbox{\rm Var}}
\def\Cov{\hbox{\rm Cov}}

\def\Card{\hbox{\rm Card}}

\def\Nit{\mathbb{N}}
\def\Rit{\mathbb{R}}

\def\prob{\mathbb{P}}

\def\and{\rm and}

\def\and{\rm and}
\def\lg{\textquotedblleft}
\def\rg{\textquotedblright}

\begin{document}

\thispagestyle{empty}
 \pagenumbering{roman}
\begin{center}
{\bf {\large TH\`ESE DE DOCTORAT }}
\end{center}

\begin{center}
{\bf {\large DE L'UNIVERSIT\'E PIERRE ET MARIE CURIE}}
\end{center}

\vskip 0.3cm
\begin{center}
{\large{Ecole Doctorale 386}}
\end{center}

\begin{center}
{\large{Sciences Mathématiques de Paris-Centre}}
\end{center}

\vskip 0.3cm

\begin{center}
{\large\it{\textbf{Sp\'ecialit\'e}\rm:} {\rm Math\'ematiques}}
\end{center}

\begin{center}
{\large\it{\textbf{Option}\rm:}
 {\rm Statistique}}
\end{center}

\vspace{0.3cm}

\begin{center}
{\large{Pr\'esent\'ee par Rawane SAMB}}
\end{center}

\vskip 0.3cm

\begin{center}
{\large{Pour obtenir le grade de}}
\end{center}

\vspace{0.3cm}

\begin{center}
 {\large {DOCTEUR DE L'UNIVERSITE PIERRE ET MARIE CURIE}}
 \end{center}

\vspace{0.5cm}

\begin{center}
{\Large{Sujet de la thèse:}}
\end{center}

\vspace{0.3cm}

\begin{center}
{\bf{\large CONTRIBUTION A L'ESTIMATION NONPARAM\'ETRIQUE}}
\end{center}

\begin{center}
{\bf{\large DE LA DENSIT\'E DES ERREURS DE R\'EGRESSION}}
\end{center}

\vspace{0.5cm}

\begin{center}{{\large Soutenue le 30 juin 2010 devant le jury
compos\'e de : }}\end{center}

\vspace{0.3cm}

\begin{center}{{\large \begin{tabular}{lll}
 $\bullet$  M.~~Denis BOSQ  \qquad & Examinateur\\[3mm]
 $\bullet$  M.~~Emmanuel GUERRE\qquad & Directeur de th\`ese\\[3mm]
 $\bullet$  Mme~~Ingrid VAN KEILEGOM\qquad & Rapporteur\\[3mm]
 $\bullet$  M.~~Christian FRANCQ \qquad & Rapporteur\\[3mm]
 $\bullet$  M.~~Benoît CADRE \qquad & Examinateur
  \end{tabular}}
  }\end{center}

\newpage

\chapter*{Remerciements}
Mes premiers remerciements s'adressent à mon directeur de thèse
Emmanuel Guerre. Qu'il soit assuré de ma sincère gratitude pour la
bienveillance et la grande disponibilité qu'il a toujours
manifestées à mon égard durant toutes mes années de thèse. Il m'a
beaucoup apporté par ses connaissances et sa rigueur
scientifiques.

\vskip 0.1cm
Je remercie très chaleureusement M. Denis Bosq,
Professeur Emérite à l'université Pierre et Marie Curie, qui me
fait l'honneur de présider le jury de cette thèse. Merci aussi à
M. Christian Francq, Professeur à l'université Lille 3, et à Mme
Ingrid Van Keilegom, Professeur à l'université catholique de
Louvain en Belgique, d'avoir accepter de juger ce travail et d'en
avoir été les rapporteurs. Leurs précieuses remarques et
suggestions ont permis d'améliorer la qualité de cette thèse. Je
suis également très reconnaissant envers M. Benoît Cadre,
Professeur à l'ENS Cachan, qui m'a fait l'honneur d'accepter
d'être membre du jury de cette thèse en tant qu'examinateur.

\vskip 0.1cm Je remercie aussi très vivement M. Paul Deheuvels,
Professeur à l'université Pierre et Marie Curie, de m'avoir
accueilli dans son DEA, puis dans son laboratoire. Je profite
également de cette occasion pour remercier M. Philippe
Saint-Pierre, Maître de Conférence du LSTA, pour son aide et ses
encouragements, et tous les autres professeurs et maîtres de
conférence du laboratoire.

\vskip 0.1cm Je suis extrêment reconnaissant envers  M. Alain
Chateauneuf et M. Jean Pierre Leca, respectivement Professeur et
Maître de Conférence à l'université Paris 1 Panthéon-Sorbonne,
pour leur sympathie et leurs qualités humaines. Ma reconnaissance
va également à l'endroit de  Mme Marie-Lou Margaria, Mme  Brigitte
Augarde et tout le personnel de l'université Paris 1.

\vskip 0.1cm Je souhaiterais également remercier les docteurs ou
doctorants que j'ai connus à l'université Pierre et Marie Curie,
notamment Mamadou Koné, Salim Bouzebda, Amadou Oury Diallo, Issam,
Boris Labrador, François-Xavier, Hicham, Tarek, Mouhamed Cherfi,
Choukri, Nabil, Layal, Lynda, Kaouthar, Lahcen Douge, Olivier
Faugeras, Camille Sabbah, Olivier Bouaziz, Claire Coiffard, Clara,
Emmanuel Onzon, Aurélie Fischer, Ousmane Bâ, entres autres, de
m'avoir aidé ou montré leur sympathie au cours de ces années.

\vskip 0.1cm\noindent Je n'oublie pas mes amis Abass Sagna, Assane
Diop, Massèye Gaye, Amadou Lamine Fall, Ange Toulougoussou,
Serigne Touba Sall, Birame Diouf, Abdoulaye Sow, Mansor Sall, Gora
Thiam et son épouse Fifi, Babacar Niang et Khourédia Ndiaye,
Edwige Sophie Mendy,  Lamine Bara Cissé et grand Bikèse à qui je
fait part de ma reconnaissance et de ma gratitude. Merci également
à Penda Sow, Cheikh Dia et son épouse Khady Kâne, Mamy Coumba
Sanou Diouf, Astou Mbacké Fall et Mary Wade pour leur soutien et
leur disponibilité.

\vskip 0.1cm Merci aussi à ma tante Khady Dia et à mes cousins
Ndèye Niang, Ngomez, Ass, Khalifa et Assa. Je les remercie pour
leur gentilesse, leur disponibilité et leur soutien qu'ils
manifestent à mon égard depuis que je suis en France.

\clearpage
 {\it J'ai également une pensée très forte à mes
parents Kiné Guèye, Mademba Samb, Bassirou Ndiaye et à toute ma
famille. Je les remercie pour leur amour et leur soutien
indéfectible et sans commune mesure à mon endroit. Je leur dédie
ce travail. A ces remerciements et dédicaces, j'associe celle qui
est devenue mon amie et ma douce moitié: je veux nommer mon épouse
Rokhaya Dièye. Je lui suis reconnaissant pour son soutien moral.}

\vskip 0.1cm Merci enfin à Massamba Wade, Rawane Wade, Mamadou
Gassama, Samba Camara et à toutes celles et ceux qui m'ont aidé et
accompagné pendant cette longue épreuve, et dont la place me
manque ici pour les citer tous.

\clearpage \noindent
A ma mère Kiné Guèye,
\vskip 0.1cm\noindent
 A mon père
Mademba Samb,
\vskip 0.1cm\noindent
A mon oncle Bassirou Ndiaye,
\vskip 0.1cm\noindent
 A mes frères et soeurs,
 \vskip0.1cm\noindent
 A mon épouse Rokhaya Dièye.

 \tableofcontents
 \clearpage
%\nobreakpage
\parindent 1.2cm

\pagenumbering{roman}
\chapter*{Notations générales}
\addcontentsline{toc}{chapter}{Notations}

\noindent
Les notations suivantes seront utilisées  dans
les différents chapitres de cette thèse.

\vskip 0.5cm \noindent{\bf Ensembles, Nombres, Fonctions} \vskip
0.15cm\noindent ${\rm Card}\left(\Omega\right)$: Cardinal de
l'ensemble $\Omega$. \vskip 0.15cm\noindent $\lfloor x\rfloor$:
Partie entière du réel $x$. \vskip 0.15cm\noindent $a\vee b$: Le
maximum des réels $a$ et $b$. \vskip 0.15cm\noindent $a\wedge b$:
Le minimum des réels $a$ et $b$. \vskip 0.15cm\noindent
$\mathds{1}_{A}$: Fonction indicatrice qui vaut $1$ sur l'ensemble
$A$ et $0$ ailleurs.
\vskip 0.15cm\noindent $f^{(k)}$: Dérivée
$k$-i\`eme de la fonction $f$.

\vskip 0.5cm
\noindent{\bf Variables aléatoires}
\vskip 0.15cm\noindent
Soient $X$ et $Y$ deux variables aléatoires.
\vskip 0.15cm\noindent
$\esp(X)$: Espérance mathématique de $X$.
\vskip 0.15cm\noindent
$\Var(X)$: Variance de $X$.
\vskip 0.15cm\noindent
$\Cov(X,Y)$: Covariance de $X$ et $Y$.
\vskip 0.15cm\noindent
$\|X\|_{p}$: Norme $L_p$ $\left(p\in]0,\infty[\right)$
de $X$ définie par $\|X\|_p=\left(\esp\left(|X|^p\right)\right)^{1/p}$,
avec $\esp\left(|X|^p\right)<\infty$.

\vskip 0.5cm
\noindent{\bf Abréviations et Symboles}
\vskip 0.15cm\noindent
$:=$ {\hspace{3mm}Symbole utilisé pour la définition d'une quantité.}
\vskip 0.15cm\noindent
Soient $(a_n)_{n\geq 1}$ et $(b_n)_{n\geq 1}$ deux suites réelles.
\vskip 0.15cm\noindent
$a_n=o(b_n)$, $n\rightarrow\infty$: Pour tout réel $\epsilon>0$, on a
$\left|a_n/b_n\right|\leq \epsilon$ pour $n$ suffisamment grand.
\vskip 0.15cm\noindent
$a_n=O(b_n)$, $n\rightarrow\infty$: Il existe un réel $C>0$ tel que
$\left|a_n/b_n\right|\leq C$ pour $n$ suffisamment grand.
\vskip 0.15cm\noindent
$a_n\asymp b_n$, $n\rightarrow\infty$: $a_n=O(b_n)$ and $b_n=O(a_n)$
 pour $n$ suffisamment grand.

\clearpage

\pagenumbering{arabic}
\renewcommand{\thechapter}{\arabic{chapter}}

\chapter[Introduction G\'en\'erale]{Introduction G\'en\'erale}
\section{Présentation du sujet}
Soit $(X_{1},Y_{1}),\ldots,(X_{n},Y_{n})$ un \'echantillon de
variables al\'eatoires ind\'ependentes et identiquement
distribu\'ees (i.i.d), de même loi que $(X,Y)$. On suppose que $Y$
est une variable univari\'ee \`a valeurs dans $\Rit$, et que  $X$
désigne une variable explicative multivari\'ee prenant ses valeurs
 dans $\Rit^d$, $d\geq 1$.
Soit $m(\cdot)$ l'esp\'erance conditionnelle de $Y$ sachant $X$,
de telle sorte que le mod\`ele de r\'egression relatif \`a $X$ et
$Y$ s'\'ecrit
\begin{eqnarray}
Y_i = m(X_i) + \varepsilon_i,
\quad i=1,\ldots,n,
\label{MR}
\end{eqnarray}
où les erreurs $\varepsilon_i$ sont supposées être des variables
aléatoires i.i.d,
 indépendantes des $X_i$, de même loi que $\varepsilon$ satisfaisant en
 particulier $\esp[\varepsilon]=0$.

 \vskip 0.3cm
 Dans ce mémoire de thèse, nous étudions l'estimation
nonparamétrique de la densité $f$ de l'erreur du modèle
(\ref{MR}). Cette estimation de la densité de l'erreur de
régression est un important outil descriptif  permettant de
comprendre le comportement des résidus, et de faire des tests
d'hypothèses sur la distribution des erreurs  du modèle ou sur la
fonction de régression. On  pourra consulter, par exemple, Ahmad
et Li (1997), Dette et {\it al.} (2002), Neumeyer et {\it al.}
(2005), pour le test de symétrie de la distribution des erreurs de
régression; Akritas et Van Keilegom (2001), Cheng et Sun (2008),
pour des tests d'ajustement sur la loi des résidus; Gozalo et
Linton (2001), Dette et von Lieres und Wilkau (2001), Neumeyer et
Van Keilegom (2010), pour le test sur l'additivité de la fonction
de régression. Notons aussi que l'estimation de $f$ peut trouver
son importance dans la prévision de $Y_{n+1}$ à partir de
$X_{n+1}$. En effet, on peut prédire $Y_{n+1}$ par l'estimateur du
mode conditionnel ${\rm mod}\left(x\right)$ de $Y_{n+1}$ sachant
que $X_{n+1}=x$, puisque ${\rm
mod}\left(x\right)=m\left(x\right)+\arg\max_{\epsilon\in\Rit}
f(\epsilon)$. Le fait d'estimer $f$ est également très important
dans la détermination d'un intervalle de prédiction pour
$Y_{n+1}$,
%si on veut \lg prévoir au mieux \rg $Y_{n+1}$ sachant $X_{n+1}$,
%%alors un estimateur de la régression
%%$m(X_{n+1})=\esp\left(Y_{n+1}|X_{n+1}\right)$ est suffisant. Mais
%il est très important
%avoir un \lg intervalle de prédiction\rg
%pour $Y_{n+1}$.
 ce qui nécessite d'estimer des quantiles de la
loi $f$.  L'estimation de  $f$ peut aussi servir à estimer la loi
de la variable $Y$, comme relaté dans Escanciano et Jacho-Chavez
(2010). Enfin cette estimation de la loi des résidus peut être
utile pour la construction d'estimateurs nonparamétriques de la
densité et de la fonction de hazard de $Y$ sachant $X$. Voir Van
Keilegom et Veraverbeke (2002).

 \vskip 0.5cm
Pour  estimer la densité $f$ des résidus du modèle (\ref{MR}),
 une première
approche consiste à noter que la densité $f$ se déduit de la
densité $\varphi\left(\cdot| x\right)$ de $Y$ sachant que $X=x$.
Plus précisémment, on a la relation
\begin{eqnarray}
f(\epsilon)
=
\varphi\left(\epsilon+m(x)| x\right).
\label{fe}
\end{eqnarray}
Suivant cette idée, on peut donc en théorie déduire un estimateur
de $f(\epsilon)$ à partir d'une estimation de $\varphi\left(y|
x\right)$ et de $m(x)$.
 Cette approche est cependant sujette au \lg fléau de la dimension\rg:
l'estimation de $\varphi\left(y| x\right)$ ne peut se faire
qu'avec une vitesse très lente lorsque la dimension de $x$ est
 élevée. Les approches proposées dans cette thèse visent à
\lg déconditionner\rg  dans l'expression (\ref{fe}) de
$f(\epsilon)$. En effet, la relation (\ref{fe}) entraîne que
\begin{eqnarray}
f(\epsilon)
=
\int
\varphi\left(\epsilon+m(x)| x\right)
g(x)
dx,
\label{fint}
\end{eqnarray}
où $g(x)$ désigne la densité de $X$.
 Cette nouvelle formule suggère que le
\lg fléau de la dimension\rg n'est peut être pas aussi important
que le laissait penser la première approche basée sur les
estimations de $f\left(y|x\right)$ et de $m(x)$. Deux stratégies
sont mises en oeuvre dans cette thèse pour essayer d'éviter le \lg
fléau de la dimension\rg.
 La première consiste à estimer nonparamétriquement
 chaque résidu $\varepsilon_i$
par $\widehat{\varepsilon}_i=Y_i-\widehat{m}_n(X_i)$, où
$\widehat{m}_n(\cdot)$ désigne un estimateur nonparamétrique de la
fonction de régression $m(\cdot)$. La seconde consiste à procéder
comme dans (\ref{fint}), et à étudier l'estimateur
$$
\widehat{f}_n(\epsilon)
=
\int
\widehat{\varphi}_n
\left(\epsilon+\widehat{m}_n(x)| x\right)
\widehat{g}_n(x)
dx,
$$
où $\widehat{\varphi}_n(\cdot|x)$ et $\widehat{g}_n(x)$ désignent
respectivement
 des estimateurs nonparamétriques de $\varphi(\cdot|x)$ et $g(x)$.

\vskip 0.3cm Le problème de l'estimation de la densité des résidus
d'un modèle  régression est un cas particulier d'un problème plus
général: l'estimation d'un paramètre d'intérêt en présence d'un
paramètre de nuisance. Dans notre cadre, qui se focalise sur
l'estimation de la distribution des résidus, la densité des
résidus $f(\cdot)$ est le paramètre d'intérêt, et la fonction de
régression $m(\cdot)$ le paramètre de nuisance. La présence de ce
paramètre de nuisance dans le modèle va influencer l'estimation du
paramètre d'intérêt. Dans le cas paramétrique, considérons, par
exemple, un \'echantillon $Z,Z_1,\ldots,Z_n$ de variables
al\'eatoires ind\'ependantes et identiquement distribu\'ees, de
densité $f\left(z|\theta,\eta\right)$, où $\theta$ est le
paramètre d'intérêt et $\eta$ le paramètre de nuisance. Une
quantité centrale liée à ces deux paramètres est la matrice
d'information de Fischer
$$
I(\eta,\theta)
 =
 \Var
 \left[
 \nabla f(z|\eta,\theta)
 \right],
$$
o\`u $\nabla f(z|\eta,\theta)$ est le gradient de
$f(z|\eta,\theta)$ par rapport \`a $\eta$ et $\theta$ d\'efini par
\begin{eqnarray*}
\nabla f(z|\eta,\theta) =
\begin{bmatrix}
\frac{\partial}{\partial\eta} f(z|\eta,\theta)
\\
\frac{\partial}{\partial\theta} f(z|\eta,\theta)
\end{bmatrix}.
\end{eqnarray*}
La matrice $I(\eta,\theta)$  s'écrit sous la forme  d'une matrice
en blocs
\begin{eqnarray*}
I(\eta,\theta) =
\begin{bmatrix}
I_{\eta\eta} & I_{\eta\theta}
\\
I_{\theta\eta} & I_{\theta\theta}
\end{bmatrix},
\end{eqnarray*}
o\`u
\begin{eqnarray*}
I_{\theta\theta}
 =
 \Var
 \left[
 \frac{\partial}{\partial\theta}
f(z|\eta,\theta) \right], \;\;
 I_{\eta\eta}
 =
 \Var
 \left[
 \frac{\partial}{\partial\eta}
f(z|\eta,\theta) \right].
\end{eqnarray*}
\vskip 0.5cm \noindent L'inégalité de Fréchet-Darmois-Cramer-Rao
(Borovkov 1987, page 156) montre que l'inverse de la matrice
d'information de Fischer, $I^{-1}\left(\eta, \theta\right)$, est,
au sens de l'ordre sur les matrices, la plus petite
 matrice de variance possible pour les estimateurs sans biais de
$\left(\eta, \theta\right)$. Cette borne $I^{-1}\left(\eta,
\theta\right)$ est atteinte par les estimateurs du maximum de
vraisemblance, comme le rappelle le théorème suivant.
%\vskip 0.5cm \noindent
\begin{theoreme}
{\rm (Borovkov 1987, page 229)}
\\
Soit $(\widehat{\eta}_n,\widehat{\theta}_n)$ un estimateur du
maximum de vraisemblance de $\left(\eta,\theta\right)$. Sous
certaines conditions de r\'egularit\'e, on a la convergence
asymptotique suivante:
\begin{eqnarray*}
\sqrt{n}
\begin{pmatrix}
\widehat{\eta}_n-\eta
\\
 \widehat{\theta}_n-\theta
 \end{pmatrix}
 \stackrel{d}{\longrightarrow}
 \mathcal{N}
 \left(0, I^{-1}(\eta,\theta)\right).
\end{eqnarray*}
\end{theoreme}
\vskip 0.3cm\noindent
 La formule du calcul de l'inverse d'une matrice en blocs
appliquée  à $I(\eta,\theta)$ permet de voir que
\begin{eqnarray*}
I^{-1}(\eta,\theta) =
\begin{bmatrix}
I^{\eta\eta}  & I^{\eta\theta}
\\
I^{\theta\eta} & I^{\theta\theta}
\end{bmatrix},
\end{eqnarray*}
 avec
$$
 I^{\theta\theta}
 =
 \left(
  I_{\theta\theta}
  -
  I_{\theta\eta} I^{-1}_{\eta\eta}
 I_{\eta\theta}
 \right)^{-1}.
$$
\vskip 0.3cm\noindent
 Du th\'eor\`eme pr\'ec\'edent, on d\'eduit la
 loi limite de l'estimateur du paramètre d'intérêt $\theta$.

\begin{corollaire}
 Sous les conditions du th\'eor\`eme
pr\'ec\'edent, on a la convergence asymptotique
$$
\sqrt{n} \left(\widehat{\theta}_n-\theta\right)
\stackrel{d}{\longrightarrow}
 \mathcal{N}
 \left(0, I^{\theta\theta}\right).
$$
\end{corollaire}
\vskip 0.3cm \noindent La matrice $I^{\theta\theta}$ s'interprète,
grâce à l'inégalité de Fréchet-Darmois-Cramer-Rao, comme étant la
meilleure variance possible pour un estimateur sans biais de
$\theta$, $\eta$ étant inconnu. Puisque
$I_{\theta\eta}I^{-1}_{\eta\eta} I_{\eta\theta}$ est
semi-positive,  la formule  de $I^{\theta\theta}$  suggère que
$I^{\theta\theta}$ est, au sens de l'ordre sur les matrices
symétriques,
 plus grande que
 $I^{-1}_{\theta\theta}$ sauf si
$I_{\eta\theta}=0$, condition indiquant que les estimateurs du
maximum de vraisemblance de $\theta$ et $\eta$ sont
asymptotiquement indépendants. Comme la variance asymptotique de
l'estimateur de $\theta$ quand $\eta$ est connu est $I_{\theta
\theta}^{-1}$,
 cette différence entre  $I^{\theta\theta}$ et $I_{\theta\theta}^{-1}$
mesure la perte (en terme d'efficacité) du fait que $\eta$ soit
inconnu quand on veut estimer $\theta$.

\vskip 0.5cm Une autre situation proche du problème de
l'estimation de la densité des résidus est l'estimation de la
fonction de répartition lorsque des paramètres sont inconnus.
Considérons, par exemple, un échantillon
 $X_1,\ldots,X_n$ de variables aléatoires i.i.d
 de fonction de répartition commune $F(x, \theta)$, où $\theta\in\Rit$.
Pour un estimateur $\widehat{\theta}_n$ de $\theta$, on définit la
fonction
 empirique associée
$$
\widehat{F}_n(t) = \frac{1}{n} \sum_{i=1}^n
 \mathds{1}\left(
F(X_i,\widehat{\theta}_n) \leq t \right), \quad t\in[0,1].
$$
Cette fonction de répartition empirique joue un rôle important
pour les tests d'adéquation du modèle considéré. En effet,
$\widehat{F}_n (t)$ doit être proche de t si le modèle est
correctement choisi. Considérons, par exemple, le modèle de
translation
$$
X_i = \theta + \varepsilon_i, \quad i=1,\ldots,n,
$$
où les résidus $\varepsilon_i$ sont de distribution commune
$\psi$. On a $F(x,\theta)=\psi(x-\theta)$. Pour ce modèle
paramétrique, on a
$$
F(X_i,\widehat{\theta}_n) = \psi(X_i-\widehat{\theta}_n) =
\psi(\widehat{\varepsilon}_i),
$$
où $\widehat{\varepsilon}_i$ est le résidu estimé
$X_i-\widehat{\theta}_n$.
 En conséquence, on a
$$
\widehat{F}_n(t) = \frac{1}{n} \sum_{i=1}^n \mathds{1} \left(
\psi(\widehat{\varepsilon}_i) \leq t \right) = \frac{1}{n}
\sum_{i=1}^n \mathds{1} \left( \widehat{\varepsilon}_i \leq
\psi^{-1}(t) \right).
$$
La relation ci-dessus  montre donc que $\widehat{F}_n(t)$ est, à
une transformation de $t$ près, la fonction de répartition
empirique des résidus $\widehat{\varepsilon}_i$.
 Le  processus empirique associé à $\widehat{F}_n$ est
$$
\widehat{y}_n(t) = n^{1/2} \{\widehat{F}_n(t)-t\}, \quad
t\in[0,1].
$$
Ce processus  a été étudié par Durbin (1973) qui obtint  le
résultat suivant.

\begin{theoreme}
Soit $\widehat{\theta}_n$ un estimateur  de $\theta$ tel que
$$
n^{1/2}(\widehat{\theta}_n-\theta) = \frac{1}{n^{1/2}}
\sum_{i=1}^n \ell(x_i,\widehat{\theta}_n) +
 o_{\prob}(1),
$$
où $\ell$ est une fonction mesurable telle que $\esp \left[
\ell(X_1,\theta) \right] =0$. Pour tout $t\in[0,1]$, on définit la
fonction $g(t)$ par
$$
g(t) = g(t,\theta) = \frac{\partial F(x,\theta)} {\partial\theta}
\mid_{x=Q(t,\theta)}, \quad Q(t,\theta) = \inf\{z:
F(z,\theta)=t\},
$$
et on pose
\begin{eqnarray*}
h(t) &=& h(t,\theta) = \int_{-\infty}^{Q(t,\theta)} \ell(x,\theta)
dF(x,\theta),
\\
L(\theta) &=& \esp \left[ \ell^2(X_1,\theta) \right].
\end{eqnarray*}
Alors  sous des conditions de régularité, le processus
$\{\widehat{y}_n(t), 0\leq t\leq 1\}$ converge asymptotiquement
 en distribution
vers un processus gaussien $\{y(t), 0\leq t\leq 1\}$, de moyenne
nulle et de fonction de covariance
$$
\Cov\left(y(t_1), y(t_2)\right) = \min(t_1,t_2) - t_1t_2 -
h(t_1)g(t_2) - h(t_2)g(t_1) + g(t_1)L(\theta)g(t_2),
$$
\end{theoreme}
\vskip 0.5cm\noindent On note que cette fonction de covariance
dépend de la fonction de répartition $F(\cdot, \theta)$ inconnue.
Donc la distribution asymptotique obtenue pour le processus
$\widehat{y}_n(t)$ est différente de la loi limite obtenue pour le
processus  empirique usuel (qui suppose $\theta$ connu),
\begin{eqnarray*}
 y_n(t)
 =
 n^{1/2}
 \{F_n(t)-t\},
 \;\;\;
 F_n(t)
 =
 \frac{1}{n}
 \sum_{i=1}^n
 \mathds{1}
 \left(F(X_i,\theta)\leq t\right).
 \end{eqnarray*}
En effet, il a été démontré que le processus $\{y_n(t), 0\leq
t\leq 1\}$ converge asymptotiquement vers un pont Brownien. Voir,
par exemple, le livre de Billinsgley (1968, p.109).

\vskip 0.5cm
 La suite de cette introduction générale donne
 des exemples d'estimation de paramètres dans le cas  d'un modèle de régression
$Y=m(X)+\sigma(X)\varepsilon$. Ces exemples seront donnés  selon
que le paramètre de nuisance, ici la fonction de régression
$m(\cdot)$,  est paramétrique ou non.

\section{Estimation de la fonction de répartition des résidus
d'un modèle linéaire} On considère le modèle linéaire
\begin{eqnarray}
 Y_i
 =
 \theta^{\top}
 X_i
 +
 \varepsilon_i,
 \quad
 i=1,\ldots,n,
 \label{ml}
\end{eqnarray}
o\`u les erreurs $\varepsilon_i$ sont i.i.d de fonction de
répartition commune $F$. Les variables $X_i$ sont supposées
 non aléatoires. Soit $\widehat{\theta}_n$ un M-estimateur de $\theta$
 (Consulter, par exemple, Huber 1964, 1981).
On s'intéresse au comportement asymptotique de la fonction de
répartition empirique $\widehat{F}_n$ des résidus estimés
 $\widehat{\varepsilon}_i=Y_i-X_i^{\top}\widehat{\theta}_n$,
$$
\widehat{F}_n(t)
=
 \frac{1}{n} \sum_{i=1}^n \mathds{1}
\left(\widehat{\varepsilon}_i\leq t\right),
\quad t\in\Rit,
$$
lorsque la dimension $p$ des régresseurs peut dépendre de la
taille $n$ de l'échantillon. Ce problème a été étudié par Portnoy
(1986) et Mammen (1996). Portnoy (1986) obtient le développement
\begin{eqnarray}
n^{1/2} \left(\widehat{F}_n(t)-F_n(t)\right)
=
\frac{f(t)}{n^{1/2}} \sum_{i=1}^n X_i^{\top}
\left(\widehat{\theta}_n-\theta\right)
+ o_{\prob}(1),
\label{das}
\end{eqnarray}
où $F_n(t)$ est la fonction de répartition empirique basée sur les
vrais résidus. Puis il montre que ce développement (\ref{das}) n'a
lieu que si $p^2/n=O(1)$
 lorsque $n$ tend vers l'infini.
 Mammen (1996) s'intéresse au comportement asymptotique de $\widehat{F}_n$
lorsque $p^2/n$ est divergente. Il considère un M-estimateur
$\widehat{\theta}_{\psi}$ tel que
$$
\widehat{\theta}_{\psi} - \theta - \sum_{i=1}^n X_i
G(\varepsilon_i)
= O_{\prob} \left(\frac{p^2}{n}\right)^{1/2},
\quad
 G(t)
=
 \frac{\psi(t)}
 {\esp\psi^{(1)}(\varepsilon_i)},
\;\;t\in\Rit, \quad \esp\left[G(\varepsilon_i)\right]=0,
$$
où $\psi$ est une fonction dérivable et croissante. Sous des
conditions de régularité, Mammen montre que pour tout
$0<C<\infty$,
\begin{eqnarray}
\sup_{|t|\leq C} \left| n^{1/2}
\left(\widehat{F}_n(t)-F_n(t)\right) - \Delta_n(t) \right| =
o_{\prob}(1), \label{dasu}
\end{eqnarray}
où, si $f$ désigne la densité des résidus,
$$
\Delta_n(t) = \frac{f(t)}{n^{1/2}} \sum_{i=1}^n \left[ X_i^{\top}
\left(\widehat{\theta}_n-\theta\right) \right] +
\frac{f(t)p}{n^{1/2}} \left[ G(t) + \frac{f^{(1)}(t)}{2f(t)} \esp
G^2(\varepsilon_1) \right].
$$
\vskip 0.1cm\noindent
 Dans le résultat (\ref{das}) de Portnoy, il
n'y a pas d'influence asymptotique  de l'estimation des résidus
sur l'estimateur de la distribution $F(t)$ lorsque
$$
\frac{1}{n^{1/2}}
 \sum_{i=1}^n
 X_i^{\top}
(\widehat{\theta}_n-\theta) = \frac{1}{n} \sum_{i=1}^n
 X_i^{\top}
\sqrt{n} (\widehat{\theta}_n-\theta) = o_{\prob}(1).
$$
Donc, puisque $\sqrt{n}(\widehat{\theta}_n-\theta)=O_{\prob}(1)$,
sous des hypothèses de régularité usuelles, la condition ci-dessus
est réalisée lorsque  $\esp[X] = 0$, d'après la loi des grands
nombres. Pour le résultat (\ref{dasu}) de Mammen, il y a un effet
de l'estimation des résidus. En effet, le terme $\Delta_n (t)$ ne
peut pas être négligeable  puisque $p^2/n$ diverge.

\vskip 0.3cm
 L'estimation de la distribution des résidus a aussi
été étudiée dans le cadre des modèles autoregressifs linéaires.
 Dans le autorégressif d'ordre 1 AR$(1)$, on observe
 les variables  aléatoires $X_0,X_1,\ldots,X_n$ telles que
\begin{eqnarray*}
X_i = \rho X_{i-1} + \varepsilon_i,
 \quad
 1\leq i\leq n,
\end{eqnarray*}
où $\rho$ désigne un paramètre réel, et les $\varepsilon_i$ des
variables aléatoires indépendantes et identiquement distribuées
(i.i.d) de densité de probabilité $f$ définie sur $\Rit$. Pour
estimer la fonction de répartition F des résidus,
 on estime d'abord les résidus $\varepsilon_i$ par
$\widehat{\varepsilon}_i=X_i-\widehat{\rho}_nX_{i-1}$,
 $\widehat{\rho}_n$ pouvant être obtenu par  la méthode  des
moindres carrées ordinaires. Le théorème suivant obtenu par Koul
(1992) donne une idée sur l'effet de l'estimation des résidus sur
la loi limite de l'estimateur de $F$.
\begin{theoreme}
Soit $\widehat{\rho}_n$ un estimateur de $\rho$ tel que
$n^{1/2}(\widehat{\rho}_n-\rho)=O_{\prob}(1)$. Alors sous une
hypothèse d'ergodicité de la famille $\{\varepsilon_i, 1\leq i\leq
n\}$, et sous d'autres hypothèses convenables, on a
$$
\sup_{x\in\Rit} \left|
 n^{1/2}
 \left[
 F_n\left(x,
\widehat{\rho}_n\right) - F_n\left(x, \rho\right) \right] \right|
= o_{\prob}(1).
$$
\end{theoreme}

\vskip 0.3cm\noindent
 Le résultat de ce théorème montre que
 l'estimation du paramètre $\rho$ n'a pas un effet asymptotique  sur
l'estimation de la fonction de répartition $F$ des résidus du
modèle précédent. Ceci vient de ce que le modèle AR$(1)$ est très
proche du modèle linéaire (\ref{ml}), les variables $X_i$ étant de
moyenne nulle.

\section{Estimation des moments d'une fonctionnelle de l'erreur}
La fonction de répartition correspond à un moment particulier,
 le moment de la fonction $\mathds{1}(\varepsilon\leq t)$.
Müller, Schick et Wefelmeyer (2004) ont étudié le cas plus général
 d'un moment $\esp h(\varepsilon)$, mais en supposant que $h$ est
 différentiable. Leur cadre d'étude est le
 modèle de régression  nonparamétrique
$Y=m(X)+\varepsilon$, où $\varepsilon$ est indépendante de $X$. La
fonction $h$ est supposée connue. Le modèle est basé sur un
échantillon d'observations i.i.d $(X_1,Y_1),\ldots,(X_n,Y_n)$ de
même loi que $(X,Y)$. Les résidus $\varepsilon_i$ sont estim\'es
par $\widehat{\varepsilon}_i=Y_i-\widehat{m}(X_i)$, où
$\widehat{m}$ est un estimateur non paramétrique de $m$. Les
auteurs proposent d'estimer $\esp[h(\varepsilon)]$  par
$\widehat{H}_n=n^{-1}\sum_{i=1}^n h(\widehat{\varepsilon}_i)$.
 Sous des conditions de régularité, ces auteurs montrent que
 $\widehat{H}_n$ est un estimateur efficace de
 $\esp[h(\varepsilon)]$ tel que
$$
\widehat{H}_n =
 \frac{1}{n}
\sum_{i=1}^n \left[ h(\varepsilon_i) - \esp [h^{(1)}(\varepsilon)]
 \varepsilon_i
 \right]
 +
 o_{\prob}(n^{-1/2}).
$$
\vskip 0.2cm \noindent
 En cons\'equence, la quantité
 $n^{1/2}[\widehat{H}_n-\esp h(\varepsilon)]$ converge asymptotiquement vers une
distribution normale de moyenne nulle et de variance
$$
 \tau_{*}^2
= \esp \left[ \left(h(\varepsilon) - \esp h(\varepsilon) -
\esp[h^{(1)}(\varepsilon)]
 \varepsilon\right)^2
 \right].
$$
 \vskip 0.3cm\noindent
 Un aspect surprenant de ce r\'esultat est que, pour certaines fonctions $h$,
 la variance asymptotique  $\tau_{*}^2$ de $\widehat{H}_n$
  est plus petite que  la variance asymptotique  $\tau^2$ de
 l'estimateur $H_n=n^{-1}\sum_{i=1}^n h(\varepsilon_i)$ bas\'e sur
 les vrais r\'esidus. En effet, supposons, par exemple, que les résidus
suivent une loi normale de moyenne nulle et  variance égale à
$\sigma^2$. Pour simplifier, on suppose que $\sigma^2=1$.
 Puisque la variance asymptotique de l'estimateur $H_n$ est égale
$\tau^2=\esp[(h(\varepsilon)-\esp h(\varepsilon))^2]$,
 on a $\tau_{*}^2<\tau^2$ si
et seulement si
\begin{eqnarray}
0<\esp[h^{(1)}(\varepsilon)]
 <2\esp[\varepsilon h(\varepsilon)]
 \;\;{\rm ou}\;\;
2\esp[\varepsilon h(\varepsilon)] <\esp[h^{(1)}(\varepsilon)]<0.
\label{Momh}
\end{eqnarray}
De plus, dans le cas où la variable  $\varepsilon$ suit une loi
normale de variance $\sigma^2=1$, on a, sous des hypothèses
convenables, $\esp[h^{(1)}(\varepsilon)]=\esp[\varepsilon
h(\varepsilon)]$. En conséquence, la première double inéqualité
dans (\ref{Momh}) est vérifiée si
$\esp[h(\varepsilon)\varepsilon]<0$, alors que la seconde double
inéqualité dans (\ref{Momh}) est satisfaite lorsque
$\esp[h(\varepsilon)\varepsilon]>0$. Cette dernière condition est
par exemple vérifiée lorsque
  $h(z)=z^3$, avec $\varepsilon$ suivant une loi normale centrée réduite.
Ce qui, dans un tel cas, entraîne que $\tau_{*}^2<\tau^2$.
 Un tel paradoxe s'explique par le fait que l'estimateur $\widehat{H}_n$
utilise mieux  le fait que les r\'esidus $\varepsilon_i$ sont de
moyenne nulle.

\section{Estimation nonparamétrique de la densité de l'erreur
dans un modèle autorégressif non linéaire}
Fu et Yang (2008)
étudient la distribution asymptotique d'un estimateur à noyau de
la densité de l'erreur dans un modèle  ${\rm AR}(p)$ non linéaire.
Ce modèle est de la forme
\begin{eqnarray*}
X_i =
g_{\theta}(X_{i-1},\ldots,X_{i-p})
+
\varepsilon_i,
\quad i\geq 1,
\end{eqnarray*}
où $\{X_i,i\in\mathbb{Z}\}$ est strictement stationnaire, et
$\theta=(\theta_1,\ldots,\theta_q)^{\top}\in\Rit^q$. Les
$\varepsilon_i$ sont i.i.d, de densité $f$, avec une moyenne nulle
et une variance $\sigma^2$. On suppose également que les résidus
$\varepsilon_i$ sont indépendantes de la famille
$(X_{i-1},\ldots,X_{i-p})$. Pour un estimateur
$\widehat{\theta}=(\widehat{\theta}_1,
\ldots,\widehat{\theta}_q)^{\top}$,
on estime les résidus $\varepsilon_i$ par
$$
\widehat{\varepsilon}_i = X_i -
g_{\widehat{\theta}}(X_{i-1},\ldots,X_{i-p}),
\quad i\geq 1.
$$
En utilisant ces résidus empiriques, Fu et Yang estiment
nonparamétriquement
 la densité $f$  par
$$
\widehat{f}_n(t)
=
\frac{1}{nh_n} \sum_{i=1}^n
K\left(\frac{\widehat{\varepsilon}_i-t}{h_n}\right),
 \quad
t\in\Rit,
$$
où $(h_n)$ est une suite de réels positifs tendant vers zero quand
$n$ tend vers l'infini, et $K$ une fonction noyau définie sur
$\Rit$. En désignant par
$$
f_n(t) = \frac{1}{nh_n} \sum_{i=1}^n
K\left(\frac{\varepsilon_i-t}{h_n}\right), \quad t\in\Rit,
$$
l'estimateur nonparamétrique de $f$ basé sur les vrais résidus, Fu
et Yang obtiennent le résultat suivant.
\begin{theoreme} {\rm  Fu et Yang (2008)}
\\ Supposons qu'il existe un réel $C_1>0$  tel que l'estimateur $\widehat{\theta}$ vérifie,
avec une probabilité égale à $1$,
\begin{eqnarray}
\lim_{n\rightarrow\infty} \sup\sqrt{\frac{n}{\log\log n}}
\|\widehat{\theta}-\theta\| \leq C_1, \label{param}
\end{eqnarray}
où $\|x\|^2=\sum_{j=1}^q x_j^2$ pour tout
$x=(x_1,\ldots,x_q)^{\top}\in\Rit^q$.
 On suppose également
que la fenêtre $h_n$ satisfait
\begin{eqnarray}
h_n\rightarrow 0, \quad \lim_{n\rightarrow\infty}
\frac{n^{1/2}h_n^{5/2}}{\log\log n} = \infty. \label{fenh}
\end{eqnarray}
Alors sous certaines conditions de régularité, on a la convergence
 en distribution suivante:
$$
\frac{1}{\sqrt{\Var f_n(t)}} \left(\widehat{f}_n(t)-\esp
f_n(t)\right)
 \stackrel{d}{\longrightarrow}
\mathcal{N}\left(0,1\right),
$$
où $\mathcal{N}(0,1)$ désigne la loi normale centrée réduite.
\end{theoreme}
\noindent
 La condition (\ref{param}) est satisfaite par
un estimateur du maximum de vraisemblance sous certaines
conditions proposées par Klimko et Nelson (1978).

\vskip 0.3cm
 Il a été
démontré dans la littérature statistique que $n^{-1/5}$ est
l'ordre de la fenêtre optimale  pour  l'estimation nonparamétrique
de la densité d'une variable aléatoire réelle $\zeta$ à partir
d'un échantillon de variables aléatoires i.i.d
$\zeta_1,\zeta_2,\ldots,\zeta_n$. Pour ce résultat, on peut, par
exemple, se référer aux ouvrages de Bosq et Lecoutre (1987), Scott
(1992), Wand  et Jones (1995). On note que dans le cadre du
théorème précédent, la condition (\ref{fenh})
 ne peut  pas vérifiée lorsque  $h_n$ est d'ordre  $n^{-1/5}$,
 mais que tous les ordres $n^{-(1/5)+\epsilon}$, $\epsilon>0$,
 qui s'en approchent sont possibles.

\section{Estimation de la loi des résidus en régression nonparam\'etrique}

L'\'etude de l'estimation  nonparamétrique d'une distribution de
l'erreur dans un mod\`ele de r\'egression nonparam\'etrique occupe
une place importante  dans la litterature statistique. En effet,
plusieurs r\'esultats inh\'erents \`a ce type d'estimation ont
\'et\'e obtenus au d\'ebut de cette d\'ecennie. On peut citer, par
exemple, Akritas et Van Keilegom (2001)
 dans le cadre de l'estimation non param\'etrique de la
fonction de r\'epartition de l'erreur d'un modèle de régression
hétéroscédastique, puis  Efromovich (2005, 2007) et Cheng (2005)
pour l'estimation nonparam\'etrique de la densit\'e des r\'esidus
d'un  modèle de régression homoscédastique. Plus r\'ecemment,
Wang, Brown, Cai et Levine (2008) se sont int\'eress\'es \`a
l'\'etude de l'influence de la fonction moyenne conditionnelle,
supposée inconnue, sur l'estimation de la  variance conditionnelle
des résidus dans le cas d'un mod\`ele de r\'egression
h\'et\'erosc\'edastique.

\subsection{Estimation de la fonction de répartition des résidus
dans un modèle de régression hétéroscédastique}

Akritas et Van Keilegom (2001) proposent un estimateur
nonparam\'etrique de la fonction de r\'epartition $F$ de l'erreur
$\varepsilon$ dans le mod\`ele de r\'egression
h\'et\'eroscedastique $Y=m(X)+\sigma(X)\varepsilon$, o\`u
$\varepsilon$ est ind\'ependante de $X$, et $m$ et $\sigma$ des
fonctions \lg lisses\rg satisfaisant quelques conditions de
r\'egularit\'e. L'estimateur $\widehat{F}_n$ de $F_{\varepsilon}$
est bas\'e sur l'estimation nonparam\'etrique des r\'esidus
$\varepsilon_i=(Y_i-m(X_i))/\sigma(X_i)$, o\`u
 $(X,Y),(X_1,Y_1),\ldots,(X_n,Y_n)$ d\'esignent un \'echantillon d'observations
ind\'ependantes et identiquement distribu\'ees. Pour l'estimation
de ces résidus, Akritas et Van Keilegom écrivent $m(x)$
 sous la forme
\begin{eqnarray}
m(x) = \int _0^1 F^{-1}(s| x)  ds, \label{msig}
\end{eqnarray}
où $F^{-1}(s| x)=\inf\{y\in\Rit: F(y| x)\geq s\}$, $F(y|
x)=\prob(Y\leq y| x)$.  On note que si la fonction $F$ est
continue,
 le changement de variable $s=F(u| x)$ dans (\ref{msig}) entraîne
$$
\int _0^1 F^{-1}(s| x) ds = \int_{\Rit} u dF(u| x) = \esp\left[Y|
X=x\right] = m(x).
$$
Pour l'estimation de $F_{\varepsilon}$, les auteurs estiment dans
un premier temps  $F(y|x)$ par l'estimateur de Stone (1977)
$$
\widetilde{F}(y| x) = \sum_{i=1}^n W_i(x,a_n) \mathds{1}(Y_i\leq
y),
$$
où les $W_i(x,a_n)$ sont les poids de Nadaraya-Watson (1964)
définis par
$$
W_i(x,a_n) = \frac{K\left(\frac{X_i-x}{a_n}\right)}
 {\sum_{j=1}^n
K\left(\frac{X_j-x}{a_n}\right)}\;,
$$
avec $K$ désignant une fonction noyau, et $a_n$ une fenêtre
tendant vers $0$ lorsque $n$ tend vers l'infini. Dans un deuxième
temps, Akritas et Van Keilegom estiment  $m(x)$ et $\sigma^2(x)$
par
$$
\widehat{m}(x) = \int _0^1 \widetilde{F}^{-1}(s| x) ds, \quad
\widehat{\sigma}^2(x) = \int _0^1 \widetilde{F}^{-1}(s| x)^2 ds -
\widehat{m}^2(x).
$$
Il convient de signaler à nouveau que  le changement de variable
 $s=\widetilde{F}_n^{-1}(y| x)$
entraîne
$$
\int _0^1 \widetilde{F}^{-1}(s| x) ds = \sum_{i=1}^n
Y_iW_i(x,a_n),
$$
ce qui correspond à l'estimateur de Nadaraya-Watson (1964)
classique.

\vskip 0.3cm\noindent Avec l'aide de ces estimateurs de $m(x)$ et
$\sigma(x)$,
 on estime chaque résidu $\varepsilon_i$ par
$\widehat{\varepsilon}_i=(Y_i-\widehat{m}(X_i))/\widehat{\sigma}(X_i)$.
L'estimateur de $F_{\varepsilon}(t)$ basé sur les résidus estimés
est alors défini par
 $\widehat{F}_{\varepsilon}(t)=n^{-1}
\sum_{i=1}^n\mathds{1}\left(\widehat{\varepsilon}_i\leq t\right)$.
Pour la détermination de la loi limite de cet estimateur,
 Akritas et Van Keilegom proposent d'abord un développement
asymptotique de $\widehat{F}_{\varepsilon}(t)$. Ce développement
est donné par le théorème suivant.

\begin{theoreme}
On suppose que la fonction de répartition $F_X$
 de $X$ est trois fois dérivable sur le support $\mathcal{X}$ de $X$, et
que et la densité $f_X$ de $X$ vérifie $\inf_{x\in\mathcal{X}}
f_X(x)>0$. On suppose également que les fonctions $m(\cdot)$ et
$\sigma(\cdot)$ sont deux fois continûment dérivables sur
$\mathcal{X}$ et que $\inf_{x\in\mathcal{X}}\sigma(x)>0$. Alors
pour tout $t\in\Rit$,  on a
$$
\widehat{F}_{\varepsilon}(t) = \frac{1}{n} \sum_{i=1}^n \mathds{1}
\left(\frac{Y_i-m(X_i)}{\sigma(X_i)}\leq t\right) -
F_{\varepsilon}(t) + \frac{1}{n} \sum_{i=1}^n \varphi(X_i,Y_i,t) +
\beta_n(t) + o_{\prob}(n^{-1/2}) + o_{\prob}(a_n^2),
$$
où
\begin{eqnarray*}
\varphi\left(x,y,t\right)
 &=&
-\frac{f_{\varepsilon}(t)}{\sigma(x)} \int \left[ \mathds{1}
\left(y\leq v\right) - F\left(v| x\right) \right] \left[
1+t\frac{v-m(x)}{\sigma(x)} \right] dv,
\\
\beta_n(t) &=& \frac{a_n^2\mu_{K}}{2} \int
\frac{\partial^2}{\partial x^2} \esp\left[\varphi(x, Y, t)|
u\right]\mid_{x=u} dF_{X}(u),
\end{eqnarray*}
avec $f_{\varepsilon}$ désignant la densité de $\varepsilon$,
$\mu_{K}$
 une constante
qui dépend de $K$, et $F_{X}$ la fonction de répartition de $X$.
\label{Akw}
\end{theoreme}
\noindent De ce théorème, Akritas et Van Keilegom déduisent le
corollaire suivant qui donne un résultat de convergence
asymptotique du processus
$n^{1/2}(\widehat{F}_{\varepsilon}(t)-F_{\varepsilon}(t))$. Ce
résultat étend les travaux de Durbin (1973) et Loynes (1980)
concernant la loi asymptotique d'un estimateur de la fonction de
répartition des résidus basé sur des paramètres estimés.
\begin{corollaire}
Supposons que le Théorème \ref{Akw} est vérifié. \vskip
0.1cm\noindent {\bf (i)} Si $na_n^4\rightarrow 0$, alors le
processus
$n^{1/2}(\widehat{F}_{\varepsilon}(t)-F_{\varepsilon}(t))$,
$t\in\Rit$,
 converge en distribution vers un processus gaussien $Z(t)$ de moyenne
$$
\esp Z(t) = \esp \left[ \mathds{1}\left(\varepsilon\leq t\right) -
F_{\varepsilon}(t) + \varphi(X,Y,t) \right] =0,
$$
et de fonction covariance
$$
\Cov\left(Z(t_1), Z(t_2)\right) = \esp \left( \biggl[ \mathds{1}
\left( \varepsilon\leq t_1\right) - F_{\varepsilon}(t_1) +
\varphi\left(X, Y, t_1\right) \biggr] \biggl[ \mathds{1} \left(
\varepsilon\leq t_2\right) - F_{\varepsilon}(t_2) +
\varphi\left(X, Y, t_2\right) \biggr] \right).
$$
{\bf (ii)}
 Si $a_n=Cn^{-1/4}$, avec $C>0$, alors le processus
$n^{1/2}(\widehat{F}_{\varepsilon}(t)-F_{\varepsilon}(t))$,
$t\in\Rit$,
 converge en distribution vers un processus gaussien $\widetilde{Z}(t)$
de moyenne
$$
\esp\widetilde{Z}(t) = \frac{C^2\mu_{K}}{2} \int
\frac{\partial^2}{\partial x^2} \esp\left[\varphi(x, Y, t)|
u\right]\mid_{x=u} dF_{X}(u),
$$
et de même fonction de covariance que le processus $Z(t)$.
\end{corollaire}
 \vskip 0.3cm\noindent
Le premier point du corrolaire précédent montre que si $na_n^4$
tend vers $0$, alors pour tout $t\in\Rit$,
\begin{eqnarray}
n^{1/2}(\widehat{F}_{\varepsilon}(t)-F_{\varepsilon}(t))
\stackrel{d}{\longrightarrow} \mathcal{N} \left(0, Var
Z(t)\right). \label{TCLAKV}
\end{eqnarray}
De plus, puisque $\esp\left[\varphi(X,Y,t)\right]=0$,
 un simple calcul montre que
\begin{eqnarray}
\nonumber \Var Z(t) &=& \esp \biggl[ \mathds{1} \left(
\varepsilon\leq t\right) - F_{\varepsilon}(t) + \varphi\left(X, Y,
t\right) \biggr]^2
\\
&=& F_{\varepsilon}(t) \left(1-F_{\varepsilon}(t)\right) + \esp
\left[ \varphi^2(X,Y,t) + 2\mathds{1}\left(\varepsilon\leq
t\right) \varphi(X,Y,t) \right]. \label{VarAKV}
\end{eqnarray}
Mais par le Théorème Central Limite, l'estimateur
$F_n(t)=n^{-1}\sum_{i=1}^n\mathds{1}\left(\varepsilon_i\leq
t\right)$ de $F_{\varepsilon}(t)$ basé sur les vrais résidus
 satisfait
$$
n^{1/2}(F_n(t)-F_{\varepsilon}(t)) \stackrel{d}{\longrightarrow}
\mathcal{N} \left( 0, F_{\varepsilon}(t)
\left(1-F_{\varepsilon}(t)\right) \right).
$$
Ce résultat, (\ref{TCLAKV}) et (\ref{VarAKV}) montrent que
 la variance asymptotique obtenue avec l'estimateur
$\widehat{F}_{\varepsilon}(t)$ est inférieure à  la variance
asymptotique $F_{\varepsilon}(t)
\left(1-F_{\varepsilon}(t)\right)$ obtenue avec $F_n(t)$ lorsque
$$
\esp\left[ \varphi^2(X,Y,t) + 2 \mathds{1} \left(\varepsilon\leq
t\right) \varphi(X,Y,t) \right] \leq 0.
$$
Dans ce cadre, il ya donc un impact positif causé par l'estimation
des résidus sur la loi limite de l'estimateur de
$F_{\varepsilon}(t)$. Notons que ces résultats ne traitent pas le
cas où l'ordre de $a_n$ est $n^{-1/5}$,  l'ordre optimal de la
fenêtre pour l'estimation de $m(\cdot)$.

\vskip 0.1cm\noindent Dans un article plus récent, Neumeyer et Van
Keilegom (2010) ont établi des résultats comparables à ceux
obtenus par Akritas et Van Keilegom (2001) dans le cas du modèle
de régression hétéroscédastique multiple:
$Y=m(X)+\sigma(X)\varepsilon$, $X\in\Rit^d$, $d\geq 1$.

\subsection{Estimation adaptative de la densité des résidus}

Efromovich (2005, 2007) utilise une méthode adaptative pour
estimer la densité $f_{\varepsilon}$ de l'erreur dans le cas des
modèles de régression homoscédastique et  hétéroscédastique. La
méthode est adaptative par rapport à la régularité de
$f_{\varepsilon}$, mesurée par son ordre $\alpha$ de dérivabilité.
Un estimateur est alors dit adaptatif s'il ne dépend pas de
$\alpha$ mais converge vers $f_{\varepsilon}$ avec la même vitesse
que les estimateurs optimaux construits en connaissant $\alpha$ et
basés sur les vrais résidus.

\vskip 0.3cm Les modèles considérés sont de la forme
$Y=m(X)+\varepsilon$ pour le modèle de régression homoscédastique,
ou de la forme $Y=m(X)+\sigma(X)\xi$, pour le modèle de régression
hétéroscédastique.
 Ces modèles sont basés sur un échantillon d'observations
i.i.d $(X_1,Y_1),\ldots,(X_n,Y_n)$ de même loi que $(X,Y)$. Les
variables $\xi$ et $\varepsilon$ sont supposées centrées et
indépendantes de $X$.
 Les fonction $m(\cdot)$ et $\sigma(\cdot)$ sont inconnues et  définies
 sur $[0,1]$. L'étude d'un estimateur de la densité de l'erreur
 par Efromovich s'est faite suivant la nature du support de l'erreur.
 On distinguera
le cas où le terme d'erreur est à support borné $[-1,1]$, et le
cas où le terme
 résiduel est de support non borné $(-\infty, \infty)$.
  Mais dans cette sous-section,
on ne parlera que du dernier cas. Pour le premier cas, le lecteur
pourra se référer au papier d'Efromovich (2005).

\vskip 0.1cm\noindent Dans le cas où le terme d'erreur est de
support $(-\infty, \infty)$, l'étude se fait donc avec le modèle
de régression homoscédastique $Y=m(X)+\varepsilon$, où la fonction
de régression $m$ est
 supposée  inconnue et définie dans $[0,1]$.
 Pour estimer la densité $f_{\varepsilon}$ de l'erreur $\varepsilon$,
Efromovich utilise un estimateur basé sur un développement en
série de cosinus. L'estimation de $f_{\varepsilon}$ nécessite une
subdivision des observations en trois sous-échantillons. Le
premier sous-échantillon de taille $n_1$ est utilisé pour estimer
la densité marginale $p$ de $X$. La deuxième partie de
l'échantillon (de taille $n_1$) est réservée à l'estimation de la
fonction de régression $m$, alors  que le dernier sous-échantillon
(de taille $n_2=n-2n_1$) est réservé à l'estimation de la densité
$f_{\varepsilon}$. On pose, pour tout  $u\in[0,1]$,
$$
\varphi_0(u)=1,
\quad
 \varphi_j(u)
=
\sqrt{2} \cos(\pi j u),
\quad j>0.
$$
Les estimateurs de $\widehat{p}$ et $\widehat{m}$ sont alors
définis par,
 pour $x\in[0,1]$,
\begin{eqnarray}
\nonumber
\widehat{p}(x)
&=&
\max
\left( b_n^{-1}
,
 n_1^{-1}
\sum_{\ell=1}^{n_1}
\sum_{s=0}^S
\varphi_s(X_{\ell})
\varphi_s(x)
\right),
\\
\widehat{m}(x)
&=&
n_1^{-1}
\sum_{\ell=n_1+1}^{2n_1}
\sum_{s=0}^S
\frac{ Y_{\ell} \varphi_s(X_{\ell})
\varphi_s(x)}
{\widehat{p}(X_{\ell})}\;.
 \label{EstimP}
\end{eqnarray}
où $b_n=4+\ln\ln(n+20)$, $n_1=n_1(n)$  désigne le plus petit
entier supérieur ou égal à $n/b_n$, et $S=S_n$ représente le plus
petit entier supérieur ou égal à $n^{1/3}$.

\vskip 0.3cm\noindent
 Avec l'aide de ces estimateurs de $p$ et $m$, Efromovich
estime les résidus $\varepsilon_{\ell}$, $\ell=2n_1+1,\ldots,n$
par
\begin{eqnarray*}
\widehat{\varepsilon}_{\ell}
=
Y_{\ell}-\widehat{m}(X_{\ell}),
\quad \ell=2n_1+1,\ldots,n.
\end{eqnarray*}
Pour $t\in\Rit$, l'estimateur $\widehat{f}_{\varepsilon}$ de
$f_{\varepsilon}(t)$
 est alors défini, suivant la méthode d'estimation de Pinsker (1980), par
\begin{eqnarray*}
\widehat{f}_{\varepsilon}(t)
 =
 \sum_{j=0}^{k_n}
 \widehat{\mu}_j
\widehat{\theta}_j \varphi_j(t),
\quad
\widehat{\theta}_j
=
(n-2n_1)^{-1}
\sum_{\ell=2n_1+1}^n
\varphi_j \left(
\widehat{\varepsilon}_{\ell}
\right),
\end{eqnarray*}
où $k_n$ est le plus petit entier supérieur ou égal à
$n^{1/5}b_n$, et les $\widehat{\mu}_j$ sont les estimateurs des
coefficients de Fourier
 $\theta_j=\int_{0}^1 f_{\varepsilon}(u)\varphi_j(u)du$. Ces
coefficients sont estimés
 selon la procédure suivante. On subdivise l'ensemble $\Nit$ des entiers naturels
en des blocs non imbriqués $B_k$, $k=1,2,\ldots$ et on pose
$t_k=1/\ln(k+2)$. Les $\widehat{\mu}_j$ sont alors définis par
\begin{eqnarray}
\widehat{\mu}_j = \frac{k^{-2} \sum_{s\in B_k}
\widehat{\theta}_s^2-n^{-1}}{k^{-2} \sum_{s\in B_k}
\widehat{\theta}_s^2}
\mathds{1}
 \left( k^{-2}
\sum_{s\in B_k}
\widehat{\theta}_s^2
>
(1+t_k)n^{-1} \right),
\quad j\in B_k.
 \label{mu}
\end{eqnarray}
Pour  évaluer la performance de l'estimateur
$\widehat{f}_{\varepsilon}(t)$, Efromovich considère l'estimateur
$\overline{f}_{\varepsilon}(t)$ de $f_{\varepsilon}$ basé sur les
vrais résidus. Cet estimateur est défini par
$$
\overline{f}_{\varepsilon}(t) = \sum_{j=0}^{k_n} \bar{\mu}_j
\overline\theta_j \varphi_j(t), \quad \overline{\theta}_j =
(n-2n_1)^{-1} \sum_{\ell=2n_1+1}^n \varphi_j \left(
\varepsilon_{\ell} \right),
$$
où les coefficients $\overline{\mu}_j$ sont définis comme dans
(\ref{mu}) en remplaçant seulement les $\widehat{\theta}_j$ par
les pseudos-estimateurs $\overline{\theta}_j$ des coefficients
$\theta_j$. En définissant l'erreur quadratique moyenne intégrée
$$
{\rm MISE} (\widehat{f}_{\varepsilon}, f_{\varepsilon}) = \esp
\int_0^1 (\widehat{f}_{\varepsilon}(t)-f_{\varepsilon}(t))^2 dt ,
$$
Efromovich obtient le résultat suivant.
\begin{theoreme} {\rm Efromovich (2005)}
\\ On suppose que les fonctions $p$ et $m$ sont de classe $C^1$ sur
$[0,1]$. Alors sous certaines conditions de régularité, on a
$$
{\rm MISE}(\widehat{f}_{\varepsilon}, f_{\varepsilon})
\leq
\left(1+\frac{C}{\ln b_n}\right) {\rm
MISE}(\overline{f}_{\varepsilon}, f_{\varepsilon})
+
\frac{Cb_n^3}{n},
$$
où $C$ est une constante  strictement positive. \label{Efro}
\end{theoreme}
\noindent
 Dans un article plus récent, Efromovich (2007) montre
que le résultat du théorème précédent reste valable sans une
procédure de \lg splitting\rg (subdivision)
 des données de l'échantillon.
\vskip 0.1cm
 Dans le cas où la densité $f_{\varepsilon}$
 admet une dérivée généralisée d'ordre $\alpha\geq 2$,
 Efromovich montre que
l'estimateur $\overline{f}_{\varepsilon}$ basé sur les vrais
résidus
 atteint la vitesse
de convergence minimax $n^{-2\alpha/(2\alpha+1)}$ pour le risque
 quadratique moyen intégré.
 Donc le Théorème \ref{Efro} prouve qu'il n'y a pas de perte
 (au sens de la vitesse minimax)
 du fait de ne pas observer
les résidus.  En conséquence, puisque $\overline{f}_{\varepsilon}$
est adaptatif par rapport à la régularité de $f_{\varepsilon}$,
 il en est de même pour l'estimateur
$\widehat{f}_{\varepsilon}$.

 \vskip 0.3cm
 Dans un article récent,
Plancade (2008) présente un estimateur nonparamétrique de la
densité de l'erreur dans un modèle de régression homoscédastique,
basé sur des techniques de sélection de modèle. Avec cette
méthode, Plancade propose une majoration du risque quadratique
intégré, et obtient la même vitesse minimax que celle obtenue par
Efromovich (2005).

\subsection{Estimation de la  fonction variance en régression hétéroscédastique}
Dans cette sous-section, nous donnons un exemple sur l'influence
de l'estimation  la fonction moyenne $m(\cdot)$ sur l'estimation
de la fonction variance $V(\cdot)$ dans le cas du mod\`ele de
r\'egression h\'et\'erosc\'edastique
\begin{eqnarray}
Y_i = m(x_i) + V^{1/2}(x_i)\varepsilon_i, \quad i=1,\ldots,n,
\label{Cai}
\end{eqnarray}
o\`u $x_i=i/n$, et les $\varepsilon_i$ sont des variables
al\'eatoires i.i.d, centr\'ees, de variance \'egale \`a $1$, et
admettant des moments d'ordre $4$ finis. Dans ce mod\`ele, le
param\`etre d'int\'er\^et est la fonction $V$, et on s'int\'eresse
\`a l'\'etude de l'impact de $m$ sur l'estimation de $V$. La
qualité de cette estimation est fortement dépendante de la
régularité de la fonction de régression $m$. On souhaite évaluer
l'impact de l'estimation de $m$ sur un estimateur de $V$. Ce
problème a été étudié par Wang, Brown, Cai et Levine (2008). Ces
auteurs ont  montré qu'il est possible d'évaluer
 explicitement l'impact de  $m$ sur l'estimateur de $V$.
Cet impact se mesure à l'aide des erreurs
 quadratiques moyennes globale  et locale  définies par
$$
R_n = \esp \int_0^1 \left(V_n(x)-V(x)\right)^2 dx,
 \quad R_n(x)
 =
\esp \left(V_n(x)-V(x)\right)^2.
$$
Ici $V_n(x)$ désigne un estimateur nonparamétrique de $V(x)$.
L'estimateur considéré par Wang et al. (2008) est défini comme
suit. On considère d'abord un noyau $K$ à support dans $[-1,1]$.
Ensuite, pour
 $i=2,\ldots,n-2$, on pose $a_i=\left(x_i+x_{i-1}\right)/2$ et
 $b_i=\left(x_i+x_{i+1}\right)/2$. Enfin,  pour
 $i=2,\ldots,n-2$, $0<h<1/2$ et
 $x\in[0,1]$, on définit
\begin{eqnarray*}
K_i^h(x) = \displaystyle{\int_{a_i}^{b_i}}
\frac{1}{h}K\left(\frac{x-u}{h}\right) du,
\end{eqnarray*}
et on prend cette intégrale de $0$ à $\left(x_1+x_2\right)/2$ pour
$i=1$, et de  $\left(x_{n-1}+x_{n-2}\right)/2$ à $1$ pour $i=n-1$.
Sous certaines hypothèses sur le noyaux  $K$, on peut vérifier que
pour tout $x\in[0,1]$, $\sum_{i=1}^{n-1}K_i^h(x)=1$. L'estimateur
$V_n(x)$ de $V(x)$ est alors défini par
\begin{eqnarray}
V_n(x)
 =
 \frac{1}{2}
 \sum_{i=1}^{n-1}
 K_i^h(x)
 \left(Y_i-Y_{i+1}\right)^2.
 \label{Vn}
\end{eqnarray}
\vskip 0.1cm\noindent Pour $\alpha>0$ et $M>0$, considérons  la
classe de fonctions $M$-lipschitziennes
$$
\mathcal{L}^{\alpha}(M)
 =
 \left\lbrace
 g: \forall\;x, y\in[0,1],
\forall\;
 k=0,\ldots,\lfloor\alpha\rfloor-1,
 \;|g^{(k)}|\leq M,\;
\left| g^{(\lfloor\alpha\rfloor)}(x)
-g^{(\lfloor\alpha\rfloor)}(y) \right| \leq
M|x-y|^{\alpha^{\prime}}
 \right\rbrace,
$$
où $\lfloor\alpha\rfloor$ est le plus grand entier naturel
inférieur à $\alpha$, et
$\alpha^{\prime}=\alpha-\lfloor\alpha\rfloor$.
 On a alors le résultat suivant.
\begin{theoreme} {\rm Wang, Brown, Cai et Levine (2008)}
\\
On considère le modèle de régression (\ref{Cai}), o\`u $x_i=i/n$,
et les $\varepsilon_i$ sont des variables al\'eatoires i.i.d,
centr\'ees, de variance \'egale \`a $1$, et admettant des moments
d'ordre $4$ finis. On suppose qu'il existe des constantes
strictement positives $\alpha$,
 $\beta$, $M_1$ et $M_2$
telles que  $m\in\mathcal{L}^{\alpha}(M_1)$ et
$V\in\mathcal{L}^{\beta}(M_2)$. Alors
 sous des hypoth\`eses convenables, la fenêtre optimale $h_n$
pour  l'estimateur $V_n(x)$ de $V(x)$ est de l'ordre de
$n^{-1/(1+2\beta)}$. De plus, pour un tel choix optimal de $h_n$,
la vitesse de convergence mimimax pour les quantit\'es $R_n$ et
$R_n(x)$ est de l'ordre de $\max\{n^{-4\alpha},
n^{-2\beta/(2\beta+1)}\}$.
\end{theoreme}

\vskip 0.3cm\noindent A l'aide de ce thèorème, on peut comparer la
performance (en terme de vitesse minimax)
 de l'estimateur
$V_n(x)$ à celle de l'estimateur $\widehat{V}_n(x)$ basé sur
l'estimation de $m$ par $\widehat{m}_n$. Cet estimateur
$\widehat{V}_n(x)$ est de la forme
\begin{eqnarray}
\widehat{V}_n(x) = \sum_{i=1}^{n-1} w_i(x)
\left(Y_i-\widehat{m}_n(x_i)\right)^2,
 \label{Vnchap}
\end{eqnarray}
où les $w_i(x)$ sont des fonctions poids. On note qu'avec
l'estimateur $\widehat{V}_n(x)$, la vitesse de convergence minimax
$\max\{n^{-4\alpha}, n^{-2\beta/(2\beta+1)}\}$ ne peut être
obtenue que si la fonction moyenne $m$ est estimée par un
estimateur de $\widehat{m}_n$  faiblement biaisé.  C'est ce qui a
incité  Brown, Cai et Levine (2008) à prendre un estimateur
 $\widehat{m}_n$ de $m$  tel que $\widehat{m}_n(x_i)=Y_{i+1}$.
  Ce qui, reporté dans
 (\ref{Vnchap}), conduit à un estimateur du type (\ref{Vn}).
  Un tel estimateur a une variance
 assez élevée et un biais suffisamment petit, pour $n$ suffisamment grand.
 Mais les auteurs
 ont prouvé qu'une grande variance de $\widehat{m}_n$ ne peut pas affecter
 la vitesse
 de convergence
 de $\widehat{V}_n$. Donc finalement, pour l'estimation de la fonction $V$,
  un estimateur
 optimal $\widehat{m}_n$ est celui de biais minimum,
 et non nécessairement celui d'erreur
 quadratique mimimale. Un enseignement important est que le carré du biais
 de $\widehat{m}_n$ joue un rôle plus important que sa variance.
 En conséquence, utiliser un estimateur qui serait optimal pour l'estimation
 de $m$ n'est pas intéressant ici, car un tel estimateur égalise asymptotiquement
 le carré du biais et la variance.

\subsection{Estimation de la densité des résidus basée sur un estimateur
de Nadaraya-Watson de la fonction de régression}

 Le problème de l'estimation nonparamétrique de la densité $f$ des
 résidus a été considéré par Cheng (2005) dans le cadre du
 modèle de régression nonparamétrique $Y=m(X)+\varepsilon$.
 Dans ce mod\`ele, la fonction de r\'egression $m$ est
 d\'efinie sur $[0,1]$, et les
estimateurs propos\'es se construisent en utilisant les
observations $(X_1,Y_1),\ldots,(X_n,Y_n)$. Ces observations sont
scind\'ees en deux parties. La premi\`ere partie est destin\'ee
\`a l'estimation des r\'esidus $\varepsilon_i=Y_i-m(X_i)$, tandis
que la seconde partie des observations est r\'eserv\'ee \`a la
construction de l'estimateur de $f$. Les estimateurs
$\widehat{\varepsilon}_i$ des r\'esidus $\varepsilon_i$
s'obtiennent \`a partir des estimations des quantit\'es $m(X_i)$.
Pour ce faire, Cheng consid\`ere un entier $r_n$ d\'ependant de
$n$, et satisfaisant
$$
0<r_n\leq n/2, \;\;\;
 \lim_{n\rightarrow\infty}r_n=\infty,
 \;\;\;
\lim_{n\rightarrow\infty}(n-r_n)=\infty.
$$
Il utilise les $r_n$ premi\`eres observations
$(X_1,Y_1),\ldots,(X_{r_n},Y_{r_n})$ pour construire l'estimateur
de la fonction $m(x)$. Cet estimateur de $m(x)$ est celui de
Nadaraya-Watson  bas\'e sur les donn\'ees
$(X_1,Y_1),\ldots,(X_{r_n},Y_{r_n})$:
$$
m_n(x)
 =
\frac{\sum_{i=1}^{r_n}Y_i K\left(\frac{X_i-x}{h_n}\right)}
{\sum_{i=1}^{r_n} K\left(\frac{X_i-x}{h_n}\right)}, \quad
x\in[0,1],
$$
o\`u $h_n$ est une fen\^etre strictement positive tendant vers $0$
quand $n$ tend vers l'infini, et $K$ une fonction intégrable sur
$\Rit$ et d'intégrale $1$. \vskip 0.2cm \noindent
 Le reste des observations
 $(X_{r_n+1},Y_{r_n+1}),\ldots,(X_n,Y_n)$ est utilis\'e pour
 estimer les r\'esidus $\varepsilon_i$ par
 $$
 \widehat{\varepsilon}_i
 =
 Y_i-m_n(X_i),
 \quad
r_n+1\leq i\leq n.
 $$
L'estimateur nonparamétrique de la densité des résidus construit
par  Cheng est alors défini par
\begin{eqnarray*}
\widehat{f}_n(t) = \frac{1}{2(n-r_n)a_n} \sum_{i=r_n+1}^n
 \mathds{1}
\left(
 t-a_n<\widehat{\varepsilon_i}\leq t+a_n
\right),
 \quad
t\in\Rit.
\end{eqnarray*}
Avec cet estimateur, Cheng (2005) obtient le résultat suivant.

\begin{theoreme}
Soit $t\in[0,1]$ tel que $f(t)>0$. Supposons que $0\leq r_n\leq
n/2$ tel que
\begin{eqnarray}
\lim_{n\rightarrow\infty} (n-r_n)a_n^3=0, \quad
\lim_{n\rightarrow\infty} (n-r_n)a_n=\infty, \quad
\lim_{n\rightarrow\infty} \frac{(n-r_n)a_n\log r_n} {r_n h_n} =0.
\label{fen}
\end{eqnarray}
On suppose également que la densité $g$ des $X_i$ est localement
lipchitzienne sur $[0,1]$. Alors sous d'autres hypothèses de
régularité,
 on a la convergence en distribution suivante:
$$
\sqrt{2(n-r_n)a_n}
\left(\frac{\widehat{f}_n(t)-f(t)}{\sqrt{f(t)}}\right)
\stackrel{d}{\longrightarrow} N\left(0,1\right),
$$
o\`u $N(0,1)$ d\'esigne la loi normale centr\'ee r\'eduite.
\end{theoreme}

\vskip 0.1cm\noindent
 Il a été démontré dans la littérature
statistique que $n^{-2/5}$ est la vitesse optimale de convergence
obtenue avec l'estimation nonparamétrique de la densité d'une
variable aléatoire réelle $\zeta$ à partir d'un échantillon de
variables aléatoires i.i.d $\zeta_1,\zeta_2,\ldots,\zeta_n$.
 Pour ce résultat, on
peut, par exemple, se référer aux ouvrages de Bosq et Lecoutre
(1987), Scott (1992), Wand et Jones (1995). Mais pour $0\leq
r_n\leq n/2$,
 le résultat du théorème précédent montre que la vitesse $n^{-2/5}$ pour
l'estimateur $\widehat{f}_n(t)$ ne peut-être atteinte que si la
fenêtre $a_n$ est d'ordre $n^{-1/5}$. Mais pour un tel ordre, la
première condition
 dans (\ref{fen})
ne peut pas être satisfaite. Donc sous les conditions du thèorème
précédent, l'estimateur $\widehat{f}_n (t)$ ne peut pas atteindre
la vitesse optimale $n^{-2/5}$, ni même s'en approcher. En effet,
(\ref{fen}) implique que $a_n =o\left(1/n^{1/3}\right)$, et que la
vitesse de convergence de $\widehat{f}_n (t)$ est
$o\left(1/n^{1/3}\right)$.

\vskip 0.1cm Cette thèse améliore les résultats de Cheng (2005).
En effet, nous verrons que sous des hypothèses convenables, les
estimateurs que nous proposerons pour estimater la loi $f$ des
résidus pourront atteindre la vitesse de convergence $n^{-2/5}$
pour ${\rm dim}(X)\leq 2$, où ${\rm dim}(X)$ désigne la dimension
de la variable explicative $X$.

\chapter[
Contribution de la thèse] {Contribution de la thèse}

\renewcommand{\thefootnote}{\arabic{footnote}}
\setcounter{footnote}{1} \setlength{\baselineskip}{.26in}

%\numberwithin{equation}{section}
\setcounter{subsection}{0} \setcounter{equation}{0}
\renewcommand{\theequation}{\thesection.\arabic{equation}}

\section{Introduction}
 La revue de la littérature faite au Chapitre 1 montre que
 la plupart des auteurs cités précédemment ont utilisé les r\'esidus
estim\'es pour construire un estimateur d'une distribution de
l'erreur. Mais  aucun d'entre eux ne s'est attaché à \'etudier
l'impact de la dimension  de la variable explicative sur
l'estimateur de la loi $f$ des erreurs, ni d'\'evaluer l'influence
de la fen\^etre de premi\`ere \'etape (utilisée pour estimer la
fonction de régression)
 sur l'estimateur final de la densit\'e des r\'esidus.
 La thèse s'attachera donc à \'evaluer l'impact de la
 dimension de la variable $X$ sur l'estimation de la densité $f$.
 Nous tenterons également de
 déterminer les vitesses de
 convergence ponctuelle des estimateurs
  nonparamétriques de $f$.
 Un de nos objectifs majeurs sera aussi de caract\'eriser les
 façons
 optimales de choisir les fen\^etres de  premi\`ere et deuxième
 \'etapes utilisées pour estimer $f$.

\vskip 0.1cm
 Nous donnons maintenant une briève
présentation de nos résulats qui seront  établis dans les deux
prochains chapitres de la thèse.

\section {Estimateur conditionnel nonparamétrique de la densité
des résidus}
 Pour mieux illustrer l'effet de la dimension de la
variable explicative $X$ sur l'estimation de la densité $f$ des
résidus du modèle de régression (\ref{MR}), nous considérons
d'abord une méthode naïve d'estimation de $f$ basée sur la
relation
$$
f(\epsilon| x)
=
\varphi\left(m(x)+\epsilon| x\right),
$$
où $f(\cdot|x)$ et $\varphi(\cdot|x)$ désignent respectivement les
densités de
 $\varepsilon$ et $Y$ sachant que $X=x$. En utilisant l'indépendance de $X$ et
$\varepsilon$, on a donc
\begin{eqnarray*}
f(\epsilon)
=
f(\epsilon| x)
=
\varphi\left(m(x)+\epsilon|x\right).
\end{eqnarray*}
Suivant cette idée, on peut donc déduire un estimateur de
$f(\epsilon)$ à partir
 d'une estimation
de $\varphi(y| x)$ et de $m(x)$. Par conséquent, un estimateur
$\widetilde{f}_n(\epsilon| x)$ de $f(\epsilon)$ est défini par
\begin{eqnarray*}
\widetilde{f}_n (\epsilon| x)
 =
\frac{ \frac{1}{nh_0^dh_1} \sum_{i=1}^n
 K_0\left(\frac{X_i-x}{h_0}\right)
 K_1\left(
 \frac{Y_i-\widehat{m}_n(x)-\epsilon}{h_1}
 \right)
 }
 {\frac{1}{nh_0^d}
 \sum_{i=1}^n
 K_0\left(\frac{X_i-x}{h_0}\right)}\;,
\end{eqnarray*}
où $h_0$, $h_1$ et $b_1$ désignent des fenêtres positives, $K_0$
et $K_1$ sont des fonctions noyaux définies respectivement sur
$\Rit^d$ et $\Rit$,
 et  $\widehat{m}_n(x)$  l'estimateur de  Nadaraya-Watson (1964) de
$m(x)$ defini par
\begin{eqnarray*}
\widehat{m}_n(x)
 =
\frac{ \sum_{j=1}^n Y_{j} K_0\left(\frac{X_j-x}{b_0}\right)}
 {
\sum_{j=1}^n K_0\left(\frac{X_j-x}{b_0}\right) }\;,
\end{eqnarray*}
où $b_0$ est une fenêtre positive.  Le théorème suivant, qui sera
démontré dans la suite de cette thèse,
 permet de mieux illustrer
l'effet négatif de la dimension de $X$ sur le comportement
asymptotique
 de l'estimateur $\widetilde{f}_n(\epsilon|x)$.

\begin{theoreme}
Considérons
$$
\mu_1(x,\epsilon) =
 \frac{\partial^2\varphi\left(x,
m(x)+\epsilon\right)}
 {\partial^2 x}
\int z K_0(z) z^{\top} dz, \quad \mu_2(x,\epsilon) =
\frac{\partial^2\varphi\left(x, m(x)+\epsilon\right)}
 {\partial^2 y}
\int v^2 K_1(v)dv,
$$
et supposons que  $b_0$, $h_0$ et $h_1$ décroissent vers $0$ et
satisfont
 $nh_0^{2d}/\ln n\rightarrow\infty$,
 $\ln(1/h_0)/\ln(\ln n)\rightarrow\infty$ et
$$
 nh_0^dh_1\rightarrow\infty,
 \quad
\left(\frac{nh_0^d}{h_1}\right) \left(b_0^4+\frac{\ln
n}{nb_0^d}\right) =o(1),
$$
 lorsque $n\rightarrow\infty$.
Alors sous des conditions de régularité sur $m$, $g$, $\varphi$,
$K_0$ and $K_1$, on a
\begin{eqnarray*}
\sqrt{nh_0^dh_1}
 \left(
\widetilde{f}_n(\epsilon|x)
 -
 \overline{\widetilde{f}}_n(\epsilon|x)
 \right)
\stackrel{d}{\rightarrow} \mathcal{N} \left( 0,
\frac{f(\epsilon|x)}{g(x)}
 \int\int
 K_0^2(z)
 K_1^2(v)
 dz dv
\right),
\end{eqnarray*}
où
\begin{eqnarray*}
\overline{\widetilde{f}}_n(\epsilon|x)
 =
 f(\epsilon| x)
 +
 \frac{h_0^2\mu_1(x,\epsilon)}{2g(x)}
 +
 \frac{h_1^2\mu_2(x,\epsilon)}{2g(x)}
 +
 o\left(h_0^2+h_1^2\right).
\end{eqnarray*}
\end{theoreme}
\vskip 0.1cm\noindent Le résultat de ce théorème suggère que pour
la normalité asymptotique de l'estimateur
$\widetilde{f}_n(\epsilon|x)$, les fenêtres optimales $h_0$ et
$h_1$ sont celles qui minimisent le développement quadratique
moyenne asymptotique
$$
AMSE\left(\widetilde{f}_n(\epsilon| x)\right)
 =
 \left[
 \frac{h_0^2\mu_1(x,\epsilon)}{2g(x)}
 +
 \frac{h_1^2\mu_2(x,\epsilon)}{2g(x)}
\right]^2
+
\frac{ f(\epsilon|x) \int\! K_0^2(z)dz \int\!K_1^2(v)dv}
{n h_0^d h_1g(x)}\;.
$$
Un simple calcul montre que les fenêtres optimales $h_0$ et $h_1$
sont toutes de l'ordre de $n^{-1/(d+5)}$, conduisant  à une
vitesse de convergence optimale $n^{-2/(d+5)}$ pour l'estimateur
$\widetilde{f}_n(\epsilon| x)$. Par conséquent, dans le cas où
$d=1$, cette vitesse de convergence est de l'ordre de $n^{-2/3}$,
ce qui est pire que la vitesse optimale $n^{-2/5}$ atteinte dans
le cadre de l'estimation d'une densité univariée. Pour la vitesse
optimale de l'estimateur d'une densité univariée,
 on pourra consulter, par exemple, les ouvrages de
 Bosq and Lecoutre (1987),
Scott (1992), Wand and Jones (1995). On note également que
l'exposant $2/(d+5)$ décroît vers $0$ lorsque $d$ devient de plus
en plus grand. Cette situation illustre donc l'impact négatif de
la dimension  de $X$ sur la performance (au sens de la vitesse de
convergence optimale) de l'estimateur
$\widetilde{f}_n(\epsilon|x)$. C'est le  problème du \lg fléau de
la dimension\rg.
 Ce problème est dû au conditionnement par $x$  dans l'expression
$f(\epsilon)=f(\epsilon|x)=\varphi\left(m(x)+\epsilon| x\right)$,
où l'on identifie la densité non conditionnelle $f(\epsilon)$ à la
densité conditionnelle $f(\epsilon|x)$ sous l'hypothèse
d'indépendance de $\varepsilon$ et $X$. Il convient également
d'ajouter  que si on voulait utiliser l'estimateur
$\widetilde{f}_n(\epsilon|x)$, il faudrait résoudre le problème du
choix de $x$. En effet, même si la densité $f(\epsilon)$ ne dépend
pas de $x$, l'estimateur $\widetilde{f}_n(\epsilon|x)$ en dépend.

\vskip 0.1cm Pour palier ce problème du \lg fléau de la
dimension\rg, il faut donc \lg déconditionner\rg  dans
l'expression ci-dessus de $f(\epsilon)$. Deux approches sont alors
proposées dans la suite cette thèse. Ces approches sont résumées
dans les deux sections suivantes.

\section{Estimation de la densité de l'erreur par
utilisation des résidus estimés}

 Cette première approche
 consiste, dans un premier temps,  à estimer nonparamétriquement les
 résidus $\varepsilon_i$ du modèle
(\ref{MR}) par
\begin{equation*}
\widehat{\varepsilon}_i
 =
 Y_i-\widehat{m}_{in},
 \quad
i=1,\ldots,n,
%\label{epschap}
\end{equation*}
où $\widehat{m}_{in} = \widehat{m}_{in} (X_i)$ désigne le \lg
leave-one out\rg
 estimateur à noyau de $m(X_i)$ défini par
\begin{equation*}
\widehat{m}_{in}
 =
\frac{ \sum_{j=1\atop j\neq i}^n Y_{j}
K_0\left(\frac{X_i-X_j}{b_0}\right)} { \sum_{j=1\atop j\neq i}^n
K_0\left(\frac{X_i-X_j}{b_0}\right) }\;.
%\label{mchapi}
\end{equation*}
Dans un deuxième temps, on utilise ces résidus estimés, comme si
c'était les vrais, pour construire un estimateur nonparamétrique
de $f(\epsilon)$. Cette construction  tient compte du fait que
 les $\widehat{m}_n(X_i)$ peuvent être des
estimateurs biaisés des $m(X_i)$ lorsque les variables $X_i$ sont
très proches des bords de leur support $\mathcal{X}$. Par
conséquent, l'estimateur de $f(\epsilon)$ est construit en prenant
les observations $X_i$ dans un ensemble ouvert $\mathcal{X}_0$
intérieur à $\mathcal{X}$.  L'estimateur de $f(\epsilon)$ est donc
défini par
\begin{equation*}
\widehat{f}_{1n}(\epsilon)
 =
\frac{1}{b_1\sum_{i=1}^n
\mathds{1}\left(X_i\in\mathcal{X}_0\right)} \sum_{i=1}^n
\mathds{1} \left( X_i \in \mathcal{X}_0 \right)
 K_1\left(\frac{\widehat{\varepsilon}_i-\epsilon}{b_1} \right).
 %\label{fnchap}
\end{equation*}
En principe, on peut supposer que $\mathcal{X}_0$ est suffisamment
proche de $\mathcal{X}$ de telle sorte que
$\widehat{f}_{1n}(\epsilon)$ se rapproche considérablement de
l'estimateur \lg classique\rg $\sum_{i=1}^n
K\left((\widehat{\varepsilon}_i-\epsilon)/b_1\right)/(nb_1)$.
Néanmoins, dans la suite de cette thèse, nous considérerons un
sous-ensemble fixé $\mathcal{X}_0$, pour des raisons de commodité.
Notons aussi que l'estimateur $\widehat{f}_{1n}(\epsilon)$ ne
dépend d'aucun paramètre inconnu, comme désiré dans la pratique.
Ceci contraste avec l'estimateur idéal nonparamétrique
\begin{equation*}
\widetilde{f}_{1n}(\epsilon)
 =
\frac{1}{b_1\sum_{i=1}^n
\mathds{1}\left(X_i\in\mathcal{X}_0\right)}
\sum_{i=1}^n \mathds{1}
\left( X_i \in \mathcal{X}_0 \right)
 K_1\left(\frac{\varepsilon_i-\epsilon}{b_1}\right),
 \end{equation*}
qui dépend en particulier des résidus non observés
$\varepsilon_i$. Cet estimateur $\widetilde{f}_{1n}(\epsilon)$ est
très proche de l'estimateur $\widehat{f}_{1n}(\epsilon)$, comme le
suggère le théorème suivant.

\begin{theoreme}
Supposons que $b_0$ and $b_1$ décroissent vers $0$ telles que
  $\ln(1/b_0)/\ln(\ln n)\rightarrow\infty$,
  $nb_0^{d^*}/\ln n\rightarrow\infty$,
  $d^*=\sup\{d+2,2d\}$, et
  $n^{(d+8)}b_1^{7(d+4)}\rightarrow\infty$ lorsque
  $n\rightarrow\infty$. Alors sous certaines conditions de
  régularité sur $m$, $g$, $f$, $K_0$ et $K_1$, on a
$$
\widehat{f}_{1n}(\epsilon)-\widetilde{f}_{1n}(\epsilon)
 =
 O_{\prob}
\biggl(R_n(b_0, b_1)\biggr)^{1/2},
\quad
\widehat{f}_{1n}(\epsilon)-f(\epsilon)
 =
 O_{\prob}
\biggl(AMSE(b_1)+R_n(b_0, b_1)\biggr)^{1/2},
$$
où
$$
AMSE(b_1)
 =
 \esp_n
 \left[
\left(\widetilde{f}_{1n}(\epsilon)-f(\epsilon)\right)^2
 \right]
 =
 O_{\prob}
\left(b_1^4 + \frac{1}{nb_1}\right),
$$
et
\begin{eqnarray*}
 R_n(b_0, b_1)
 =
 b_0^4
 +
 \left[
 \frac{1}{(nb_1^5)^{1/2}}
 +
 \left(\frac{b_0^d}{b_1^3}\right)^{1/2}
 \right]^2
 \left(
 b_0^4
 +
 \frac{1}{nb_0^d}
 \right)^2
 +
 \left[
 \frac{1}{b_1}
 +
 \left(\frac{b_0^d}{b_1^7}\right)^{1/2}
 \right]^2
 \left(b_0^4+\frac{1}{nb_0^d}\right)^3.
\end{eqnarray*}
\end{theoreme}

\vskip 0.1cm\noindent Les résultats de ce théorème donnent une
première idée de l'impact de l'estimation des résidus sur
l'estimateur nonparamétrique de la densité $f(\epsilon)$.

\vskip 0.1cm Le théorème suivant détermine la façon optimale de
choisir la fenêtre de première étape $b_0$. A notre connaissance,
cet aspect n'a pas encore été étudié dans la littérature
statistique. Dans ce qui suit,  $a_n \asymp b_n$ signifie que $
a_n= O (b_n)$ et $b_n = O(a_n)$, c'est à dire il existe une
constante $C>0$ telle que $|a_n|/C \leq |b_n| \leq C |a_n|$, pour
$n$ suffisamment grand.

\begin{theoreme}
 On considère la fenêtre
$$
b_0^* = b_0^*(b_1)
 =
\arg\min_{b_0}
 R_n (b_0, b_1),
$$
où la minimisation se fait sur l'ensemble des fenêtres $b_0$
satisfaisant les condtions du théorème précédent. Alors la fenêtre
$b_0^*$ vérifie
$$
b_0^*
\asymp
\max
\left\lbrace
\left(\frac{1}{n^2b_1^3}\right)^{\frac{1}{d+4}}
,
\left(\frac{1}{n^3b_1^7}\right)^{\frac{1}{2d+4}}
\right\rbrace,
$$
et on a
$$
R_n(b_0^*, b_1)
\asymp
\max
\left\lbrace
\left(\frac{1}{n^2b_1^3}\right)^{\frac{4}{d+4}}
,
\left(\frac{1}{n^3b_1^7}\right)^{\frac{4}{2d+4}}
\right\rbrace.
$$
\end{theoreme}
De ce théorème, on déduit le résultat suivant qui donne les
conditions pour lesquelles l'estimateur
$\widehat{f}_{1n}(\epsilon)$ atteind la vitesse optimale
$n^{-2/5}$ lorsque $b_0=b_0^*$.

\begin{theoreme}
On considère la fenêtre
$$
b_1^*
 =
\arg\min_{b_1}
\biggl(AMSE(b_1)+R_n(b_0^*,b_1)\biggr),
$$
où  $b_0^*=b_0^*(b_1)$ est definie comme dans le théorème
précédent. Alors
\begin{enumerate}
\item Pour $d\leq 2$, la fenêtre $b_1^*$ satisfait
$$
b_1^* \asymp \left(\frac{1}{n}\right)^{\frac{1}{5}},
$$
et on a
$$
\biggl( AMSE(b_1^*) + R_n(b_0^*, b_1^*) \biggr)^{\frac{1}{2}}
\asymp
\left(\frac{1}{n}\right)^{\frac{2}{5}}.
$$
\item Pour $d\geq 3$, $b_1^*$ satisfait
$$
b_1^*
\asymp
\left(\frac{1}{n}\right)^{\frac{3}{2d+11}},
$$
et on a
$$
\biggl(
 AMSE(b_1^*)
 +
 R_n(b_0^*, b_1^*)
\biggr)^{\frac{1}{2}} \asymp
\left(\frac{1}{n}\right)^{\frac{6}{2d+11}}.
$$
\end{enumerate}
\end{theoreme}
Ces résultats montrent que pour $d\leq 2$, la vitesse de
convergence de la différence
$\widehat{f}_{1n}(\epsilon)-f(\epsilon)$ est d'ordre $n^{-2/5}$,
ce qui correspond à la vitesse de convergence optimale dans le cas
de l'estimation de la densité d'une variable univariée. Donc dans
ce cas,  il ya un impact positif de l'estimation des résidus sur
l'estimateur de $f(\epsilon)$. Mais pour $d\geq 3$, la vitesse le
taux de convergence $n^{-2/5}$ ne peut pas être atteinte avec
l'estimateur $\widehat{f}_{1n}(\epsilon)$.

 Nous obtenons également le résulat de normalité asymptotique suivant.

\begin{theoreme}
Supposons que
$$
 nb_0^{d+4}=O(1), \quad nb_0^4b_1=o(1),
 \quad
 nb_0^{d}b_1^3\rightarrow\infty,
$$
lorsque $n$ tend vers $\infty$. Alors sous des conditions de
régularité, on a
$$
\sqrt{nb_1} \left( \widehat{f}_{1n}(\epsilon)
-
\overline{f}_{1n}(\epsilon)\right)
\stackrel{d}{\rightarrow}
\mathcal{N} \left(
 0,
\frac{f(\epsilon)}
{\prob\left(X\in\mathcal{X}_0\right)}
\int
K_1^2(v) dv \right),
$$
où
$$
\overline{f}_{1n}(\epsilon)
 =
f(\epsilon)
+
\frac{b_1^2}{2}f^{(2)}(\epsilon)
\int v^2 K_1(v) dv
+
 o\left(b_1^2\right).
$$
\end{theoreme}
 La deuxième approche utilisée pour l'estimation de la
densité $f$ est résumée dans la sous-section suivante.

\section{Estimation de la densité de l'erreur par
intégration d'une loi conditionnelle}

Cette approche consiste d'abord à remarquer que
\begin{eqnarray*}
f(\epsilon)
=
\int \varphi\left(\epsilon+m(x)|x\right)
 g(x)dx
=
\int \varphi\left(x,\epsilon+m(x)\right) dx,
\end{eqnarray*}
où $g$ désigne la densité marginale de $X$,  et
$\varphi(\cdot,\cdot)$ la densité conjointe du couple $(X,Y)$.
 Cette formule suggère donc d'estimer,
 dans un second temps, $f(\epsilon)$ par
$$
\widehat{f}_{2n}(\epsilon)
=
\int \widehat{\varphi}_n
\left(x,\epsilon+\widehat{m}_n(x)\right)
 dx,
$$
où $\widehat{m}_n(x)$ désigne l'estimateur à noyau de
Nadaraya-Watson (1964) de $m(x)$, et $\widehat{\varphi}_n$
l'estimateur nonparamétrique de $\varphi$. Ces estimateurs sont
définis comme suit. On considère des fenêtres $b_{0}=b_0(n)$ et
$b_{1}=b_1(n)$  associées à la variable $X$, et une fenêtre
$h=h(n)$ associé à la variable $Y$. On suppose que $K_0$ et $K_1$
 sont des fonctions noyaux définis dans $\Rit^d$, et que
 $K_2$ désigne une fonction noyau défini dans $\Rit$. Pour
 tout $(x,y)\in\Rit^d\times\Rit$, les estimateurs $\widehat{m}_n(x)$
 et $\widehat{\varphi}_n(x,y)$ sont définis par
\begin{eqnarray*}
\widehat{m}_{n}(x)
&=&
\frac{\sum_{j=1}^{n}Y_{j}
K_0\left(\frac{X_{j}-x}{b_{0}}\right)}
 {\sum_{j=1}^{n}
K_0\left(\frac{X_{j}-x}{b_{0}}\right)},
\\
\widehat{\varphi}_{n}\left( x,y\right)
&=&
\frac{1}{nb_{1}^{d}h}
\sum_{i=1}^{n}
K_1\left(\frac{X_{i}-x}{b_{1}}\right)
K_2\left(\frac{ Y_{i}-y}{h}\right).
\end{eqnarray*}
On considère également
$$
\widetilde{f}_{2n}(\epsilon)
 =
  \int \widehat{\varphi}_n
\left(x,\epsilon+m(x)\right) dx,
$$
l'estimateur par de $f$ basé sur la fonction de régression $m$.
Avec l'aide de ces estimateurs, on obtient d'abord le  théorème
suivant.

\begin{theoreme}
On suppose que $b_0$, $b_1$ et $h$  décroissent vers $0$ telles
que $\ln(1/b_0)/\ln(\ln n)\rightarrow\infty$,
  $b_0^{d}/(nb_0^{2d})^p= O(b_0^{2p})$, $p\in[0,6]$,
 $nb_1^{2d}\rightarrow\infty$ et
$n^{(d+8)}h^{7(d+4)}\rightarrow\infty$ lorsque
$n\rightarrow\infty$. Alors, sous des conditions de régularité sur
$g$, $m$, $f$, $\varphi$, et $K_j$, $j=0,1,2$, on a
$$
\widehat{f}_{2n}(\epsilon)-f(\epsilon)
 =
 O_{\prob}
\biggl(
 AMSE(b_1,h)+RT_n(b_0,b_1,h)
 \biggr)^{1/2},
$$
où
$$
AMSE(b_1,h)
 =
 \esp_n
 \left[
\left(\widetilde{f}_{2n}(\epsilon)-f(\epsilon)\right)^2
 \right]
 =
 O_{\prob}
\left(b_1^4 +h^4+ \frac{1}{nb_1}\right),
$$
et
\begin{eqnarray*}
RT_n(b_0,b_1,h)
&=&
b_{0}^{4}
+ \left(b_{0}^d\vee b_{1}^d\right)
\left[
\frac{1}{nb_1^dh^3}
\left( b_0^4 + \frac{1}{nb_0^{d}}
\right)
+
\frac{1}{nb_0^d}
\right]
\\
&&
+
\left(b_{0}^d\vee b_{1}^d\right)
\left[
\frac{1}{nb_1^dh^5}
 \left(b_{0}^4
 +
 \frac{1}{nb_0^d}
 \right)^2
 +
\frac{1}{n^2b_0^{2d}h^3} \right]
 \\
 &&
 +
\frac{1}{h^2}
\left( b_0^4 + \frac{1}{nb_0^d} \right)^3
 +
\frac{b_{0}^d\vee b_{1}^d}{h^7}
\left( b_0^{4} + \frac{1}{nb_0^d}
\right)^3.
 \end{eqnarray*}
\end{theoreme}

\vskip 0.3cm En se basant sur ce théorème, on retrouve des
résultats similaires à ceux obtenus avec l'estimateur
$\widehat{f}_{1n}(\epsilon)$, notamment ceux relatifs aux choix
optimaux des fenêtres de première et deuxième étape pour
l'estimation de $f(\epsilon)$.

\vskip 0.1cm\noindent
 $\bullet$ Choix optimal de la fenêtre $b_0$

\begin{theoreme}
 On pose $b_0=b_1$, puis on considère la fenêtre
 $$
 b_0^*
 =
 b_0^*(h)
 =
\arg\min_{b_0}
 RT_n (b_0, b_0, h),
$$
où la minimisation se fait sur l'ensemble des fenêtres $b_0$
satisfaisant les hypothèses du théorème précédent. Alors $b_0^*$
vérifie
$$
b_0^*
\asymp
\max
\left\lbrace
\left(\frac{1}{n^2h^3}\right)^{\frac{1}{d+4}}
,
\left(\frac{1}{n^3h^7}\right)^{\frac{1}{2d+4}}
\right\rbrace,
$$
et on a
$$
RT_n(b_0^*,b_0^*, h)
\asymp
\frac{1}{n}
+
\max \left\lbrace
\left(\frac{1}{n^2h^3}\right)^{\frac{4}{d+4}}
 ,
\left(\frac{1}{n^3h^7}\right)^{\frac{4}{2d+4}}
\right\rbrace.
$$
%\label{Optimb0}
\end{theoreme}

\vskip 0.1cm\noindent
 $\bullet$ Choix optimal de la fenêtre $h$

\begin{theoreme}
On considère la fenêtre
$$
h^*
 =
\arg\min_{h}
\biggl(AMSE(b_0^*,h)+RT_n(b_0^*,b_0^*,h)\biggr),
$$
où $b_0^*=b_0^*(h)$ est définie comme dans le théorème précédent.
Alors
 \begin{enumerate}
\item Pour $d\leq 2$, la fenêtre $h^*$ vérifie
$$
h^*
\asymp
\left(\frac{1}{n}\right)^{\frac{1}{5}},
$$
et on a
$$
\biggl( AMSE(b_0^*, h^*)
 +
 RT_n\left(b_0^*, h^*, h^*\right)
\biggr)^{\frac{1}{2}}
\asymp
\left(\frac{1}{n}\right)^{\frac{2}{5}}.
$$
\item Pour $d\geq 3$, $h^*$ satisfait
$$
h^*
\asymp
\left(\frac{1}{n}\right)^{\frac{3}{2d+11}},
$$
et on a
$$
\biggl(
 AMSE(b_0^*, h^*)
 +
 RT_n(b_0^*, b_0^*, h^*)
\biggr)^{\frac{1}{2}} \asymp
\left(\frac{1}{n}\right)^{\frac{6}{2d+11}}.
$$
\end{enumerate}
\end{theoreme}
La conclusion des résulats de ce théorème est la même que celle du
théorème similaire obtenu avec l'estimateur
$\widehat{f}_{1n}(\epsilon)$.

\vskip 0.1cm\noindent
 $\bullet$ Normalité asymptotique

\begin{theoreme}
Supposons que
$$
nb_0^{d+4}=O(1),
\quad nb_0^4h=o(1),
\quad
nb_0^dh^3\rightarrow\infty,
$$
 lorsque $n\rightarrow\infty$. Alors sous certaines conditions de
 régularité on a,
\begin{eqnarray*}
\sqrt{nh}
 \left(
\widehat{f}_{2n}(\epsilon)
 -
 \overline{f}_{2n}(\epsilon)
 \right)
\stackrel{d}{\rightarrow}
\mathcal{N}
\left( 0,
f(\epsilon)
\int K_2^2(v) dv
\right),
\end{eqnarray*}
avec
 \begin{eqnarray*}
\overline{f}_{2n}(\epsilon)
&=&
 f(\epsilon) + \frac{b_0^2}{2}
 \int
 \mathds{1}
 \left(x\in\cal{X}\right)
 \frac{\partial^2 \varphi(x,\epsilon+m(x))}
 {\partial^2 x}
 dx
 \int
 z K_1(z)z^{\top}
 dz
 \\
 &&
 +
 \frac{h^2}{2}
 \int
 \mathds{1}
 \left(x\in\cal{X}\right)
 \frac{\partial^2 \varphi(x,\epsilon+m(x))}
 {\partial^2 y}
 dx
 \int
 v^2 K_2(v)
 dv
 +
 o\left(b_0^2+h^2\right).
\end{eqnarray*}
\end{theoreme}

Pour finir la thèse, nous réaliserons des simulations numériques
pour valider et mieux mettre en exergue les résultats obtenus avec
les estimateurs $\widehat{f}_{1n}$ et $\widehat{f}_{2n}$. Nous
comparerons les performances de ces estimateurs  en terme
d'erreurs quadratiques moyennes globales et locales. Nous
présenterons également des perspectives de recherche pour nos
futurs travaux.

\chapter[
Kernel estimation of the p.d.f  of regression errors using
estimated residuals] {Nonparametric kernel estimation of the
probability density function of  regression errors using estimated
residuals}

\renewcommand{\thefootnote}{\arabic{footnote}}
\setcounter{footnote}{1}
\setlength{\baselineskip}{.26in}

%\numberwithin{equation}{section}
\setcounter{subsection}{0} \setcounter{equation}{0}
\renewcommand{\theequation}{\thesection.\arabic{equation}}

{\bf Abstract:} In this chapter we deal with the nonparametric
density estimation of the regression error term assuming its
independence with the covariate. The difference between the
feasible estimator which uses the estimated residuals and the
unfeasible one using the true
 residuals
is studied. An optimal choice of the bandwidth used to estimate
the residuals is  given. We also study the asymptotic normality of
the feasible kernel estimator and its rate-optimality.

\section{Introduction}

 Consider a sample
$(X,Y),(X_{1},Y_{1}),\ldots,(X_{n},Y_{n})$ of independent and
identically distributed (i.i.d) random variables, where Y is the
univariate dependent variable and the covariate X is of dimension
$d$. Let $m(\cdot)$ be the conditional expectation of $Y$ given
$X$ and let
 $\varepsilon$ be
the related regression error term, so that the regression error
model is
\begin{eqnarray}
Y_i = m(X_i)+\varepsilon_i,
\quad i=1,\ldots,n.
\label{Rm}
\end{eqnarray}
 We wish to estimate the probability distribution function
(p.d.f) of the regression error term, $f(\cdot)$, using the
nonparametric
 residuals.
Our potential applications are as follows.
 First, an estimation of the p.d.f of $\varepsilon$ is an
 important tool for understanding the residuals behavior and
 therefore the fit of the regression model (\ref{Rm}).
This estimation of $f(\cdot)$ can be used for
 goodness-of-fit tests of a specified error  distribution in a
 parametric regression setting. Some examples can be founded in
  Loynes (1980), Akritas and Van  Keilegom (2001),  Cheng and Sun
  (2008).
   The estimation of the density of the regression error term can
   also
 be useful for testing the symmetry of the residuals.
  See Ahmad et Li (1997), Dette et
 {\it al.} (2002). Another interest of the estimation of $f$
 is that it can be used for constructing nonparametric estimators
 for the density and hazard function of $Y$ given $X$, as related
 in Van Keilegom and Veraverbeke (2002). This estimation of $f$ is also
 important when are interested in the estimation
 of the p.d.f of the response variable $Y$. See Escanciano and
 Jacho-Chavez (2010).
  Note also that an estimation of the p.d.f of the
 regression
 errors  can be useful for proposing a mode forecast of
 $Y$ given $X=x$. This mode forecast is based on an estimation of
 $m(x)+\arg\min_{\epsilon\in\Rit} f(\epsilon)$.

\vskip 0.3cm
 Relatively little is known about the nonparametric
estimation of the p.d.f and  the cumulative distribution function
(c.d.f) of the regression error. Up to few
 exceptions, the nonparametric literature focuses on studying the
 distribution of $Y$ given $X$. See Roussas
(1967, 1991), Youndj\'e (1996) and
 references therein.
Akritas and Van Keilegom (2001)  estimate the cumulative
distribution function of the regression error in heteroscedastic
model. The estimator proposed by these authors is based on a
nonparametric estimation of the residuals. Their result show the
impact of the estimation of the residuals on the limit
distribution of the
 underlying estimator
of the cumulative distribution function. M\"{u}ller, Schick and
Wefelmeyer (2004) consider the estimation of
 moments of the regression error.
Quite surprisingly, under appropriate conditions, the estimator
based
 on the true errors
is less efficient than the  estimator which uses the nonparametric
estimated residuals.
 The reason is that the latter estimator  better uses the fact that
 the regression error $\varepsilon$
has mean zero. Efromovich (2005) consider adaptive estimation of
the p.d.f of the
 regression error.
 He gives a nonparametric estimator based on the
estimated residuals, for which the Mean Integrated Squared Error
(MISE)  attains the  minimax rate. Fu and Yang (2008) study the
asymptotic normality of the estimators of the regression error
p.d.f  in nonlinear autoregressive models.
 Cheng (2005)  establishes
the asymptotic normality of an estimator of
 $f(\cdot)$ based on the
estimated residuals. This estimator is  constructed by splitting
the
 sample into two parts: the first part is used for the construction of
 estimator of $f(\cdot)$,
 while the second part of the sample is used for the estimation of the
 residuals.

\vskip 0.3cm
 The focus of this chapter is to estimate  the p.d.f
of the regression error using the estimated residuals,  under the
assumption that the covariate $X$ and the regression error
 $\varepsilon$ are independent.
 In a such setup, it would be unwise to use a conditional approach
 based
on the fact that
$f(\epsilon)=f(\epsilon|x)=\varphi\left(m(x)+\epsilon| x\right)$,
where $\varphi(\cdot| x)$  is the p.d.f of $Y$ given $X=x$.
Indeed,
 the estimation of $m(\cdot)$ and $\varphi(\cdot| x)$ are affected by
the curse of dimensionality, so that the resulting estimator of
$f(\cdot)$ would have considerably a slow rate of convergence if
the
 dimension of $X$ is
high. The  approach proposed here uses a two-steps procedure
which, in
 a
first step, replaces the unobserved regression error terms by some
nonparametric estimator $\widehat{\varepsilon}_i$. In a second
step,
 the
estimated $\widehat{\varepsilon}_i$'s are used to estimate
nonparametrically $f(\cdot)$, as if they were the true
 $\varepsilon_i$'s. If
proceeding so can circumvent the curse of dimensionality, a
challenging
 issue is to evaluate
the impact of the estimated  residuals on the final estimator of
$f(\cdot)$. Hence one of the contributions of our study is to
analyze the effect of the estimation of the residuals on the
regression errors p.d.f. Kernel estimators.
 Next, an optimal
choice of
 the bandwidth used to
 estimate the residuals is  given. Finally, we study the asymptotic
 normality of the
feasible Kernel estimator and its rate-optimality.

\vskip 0.3cm
 The rest of this chapter is organized as follows.
Section 3.2 presents ours estimators
 and proposes an  asymptotic normality of the (naive) conditional
  estimator of the
 density of the regression error. Sections 3.3 and 3.4 group our
assumptions and  main results. The conclusion of this chapter is
given in Section 3.5, while the proofs of our results are gathered
in section 3.6 and in an appendix.

\setcounter{subsection}{0} \setcounter{equation}{0}
\renewcommand{\theequation}{\thesection.\arabic{equation}}

\section{Some nonparametric estimator of the density of the regression error}

 To illustrate the potential impact of the dimension $d$ of the
$X_i$'s,
 let us first consider a naive conditional estimator of the p.d.f
 $f(\cdot)$ of the regression error term $\varepsilon$. Let $\varphi
 (\cdot|x)$ and $f(\cdot|x)$ be respectively the p.d.f. of $Y$ and $\varepsilon$
 given $X=x$. Since $f(\epsilon|x) = \varphi (m(x)+\epsilon|x)$, using
 the independence of $X$ and $\varepsilon$ gives
\begin{equation}
f(\epsilon)
 =
 f(\epsilon|x)
 =
 \varphi\left(m(x)+\epsilon|x\right).
\label{Naive}
\end{equation}
Consider some Kernel functions $K_0$, $K_1$ and some bandwidths
$b_0$, $h_0$
 and $h_1$. The expression (\ref{Naive}) of $f$ suggests to use the
 Kernel nonparametric estimator
\begin{eqnarray*}
\widetilde{f}_n(\epsilon|x)
 =
\frac{ \frac{1}{nh_0^dh_1} \sum_{i=1}^n
 K_0\left(\frac{X_i-x}{h_0}\right)
 K_1\left(\frac{Y_i-\widehat{m}_n(x)-\epsilon}{h_1}\right)
 }{\frac{1}{nh_0^d}
 \sum_{i=1}^n
 K_0\left(\frac{X_i-x}{h_0}\right)}\;,
\end{eqnarray*}
where $\widehat{m}_n(x)$ is the Nadaraya-Watson (1964) estimator
of $m(x)$ defined as
\begin{eqnarray}
\widehat{m}_n(x)
 =
 \frac{
\sum_{j=1}^n Y_{j}
 K_0\left(\frac{X_j-x}{b_0}\right)}
 {
 \sum_{j=1}^n
 K_0\left(\frac{X_j-x}{b_0}\right)}\;.
 \label{mchap}
\end{eqnarray}
The first result presented in this chapter is the following
proposition.
\begin{proposition}
Define
$$
\mu_1(x,\epsilon)
=
 \frac{\partial^2\varphi\left(x,
m(x)+\epsilon\right)}
 {\partial^2 x}
\int z K_0(z) z^{\top} dz, \quad \mu_2(x,\epsilon) =
\frac{\partial^2\varphi\left(x, m(x)+\epsilon\right)}
 {\partial^2 y}
\int v^2 K_1(v)dv,
$$
and suppose that  $h_0$ decrease to $0$ such that
 $nh_0^{2d}/\ln n\rightarrow\infty$,
 $\ln(1/h_0)/\ln(\ln n)\rightarrow\infty$ and
$$
{(\rm\bf{A}_0):} \quad nh_0^dh_1\rightarrow\infty, \quad
\left(\frac{nh_0^d}{h_1}\right) \left(b_0^4+\frac{\ln
n}{nb_0^d}\right) =o(1),
$$
 when $n\rightarrow\infty$. Then under Assumptions  $(A_1)-(A_{10})$
  given in the next section, we have
\begin{eqnarray*}
\sqrt{nh_0^dh_1}
 \left(
\widetilde{f}_n(\epsilon|x)
 -
 \overline{\widetilde{f}}_n(\epsilon|x)
 \right)
\stackrel{d}{\rightarrow} \mathcal{N} \left( 0,
\frac{f(\epsilon|x)}{g(x)}
 \int\int
 K_0^2(z)
 K_1^2(v)
 dz dv
\right),
\end{eqnarray*}
where $g(\cdot)$ is the marginal density of $X$ and
\begin{eqnarray*}
\overline{\widetilde{f}}_n(\epsilon|x)
 =
 f(\epsilon| x)
 +
 \frac{h_0^2\mu_1(x,\epsilon)}{2g(x)}
 +
 \frac{h_1^2\mu_2(x,\epsilon)}{2g(x)}
 +
 o\left(h_0^2+h_1^2\right).
\end{eqnarray*}
\label{Prop}
\end{proposition}
\noindent
 This results suggests that an optimal choice of the
 bandwidths $h_0$ and $h_1$ should achieve the minimum of
the  asymptotic mean square expansion first order terms
$$
AMSE \left( \widetilde{f}_n(\epsilon| x) \right)
 =
 \left[
 \frac{h_0^2\mu_1(x,\epsilon)}{2g(x)}
 +
 \frac{h_1^2\mu_2(x,\epsilon)}{2g(x)}
\right]^2
+
 \frac{ f(\epsilon|x)
 \int\! K_0^2(z)dz \int\!
K_1^2(v)dv} {n h_0^d h_1g(x)}\;.
$$
Elementary calculations yield that the resulting optimal
bandwidths $h_0$ and $h_1$
 are all proportional to
$n^{-1/(d+5)}$, leading to the exact consistency rate
$n^{-2/(d+5)}$ for $\widetilde{f}_n(x|\epsilon)$. In the case
$d=1$, this rate is $n^{-1/3}$, which is worst than the
 rate $n^{-2/5}$ achieved  by the optimal Kernel estimator of an univariate
 density. See  Bosq and Lecoutre (1987),
Scott (1992), Wand and Jones (1995). Note also that the exponent
$2/(d+5)$ decreases to $0$ with
 the dimension $d$. This indicates a negative impact of the dimension $d$
 on the performance of the estimator, the so-called curse of
 dimensionality.
The fact that $\widetilde{f}_n(\epsilon|x)$ is affected by the
curse
 of dimensionality is a consequence of conditioning. Indeed, (\ref{Naive})
 identifies the unconditional $f(\epsilon)$ with the conditional
 distribution
 of the regression error given the covariate.

\vskip 0.3cm To avoid this curse of dimensionality in the
nonparametric kernel estimation of $f(\epsilon)$, our approach
proposed here builds, in
 a first step,  the estimated residuals
\begin{equation}
\widehat{\varepsilon}_i
=
Y_i-\widehat{m}_{in},
\quad
i=1,\ldots,n,
\label{epschap}
\end{equation}
where $\widehat{m}_{in} =\widehat{m}_{in}(X_i)$ is a leave-one out
 version of the Kernel regression estimator (\ref{mchap}),
\begin{equation}
\widehat{m}_{in} = \frac{\sum_{j=1\atop j\neq i}^nY_j
K_0\left(\frac{X_j-X_i}{b_0}\right)} {\sum_{j=1\atop j\neq i}^n
K_0\left(\frac{X_j-X_i}{b_0}\right)}\;. \label{mchapi}
\end{equation}
It is tempting to use, in a second step, the estimated
 $\widehat{\varepsilon}_i$ as if they were the true residuals $\varepsilon_i$. This
 would ignore that the $\widehat{m}_{n}(X_i)$'s can deliver  severely
 biased estimations of the $m(X_i)$'s for those $X_i$ which are close to the
 boundaries of the support $\mathcal{X}$ of the covariate distribution.
 To that aim, our proposed estimator trims the observations $X_i$
 outside an inner subset $\mathcal{X}_0$ of
$\mathcal{X}$,
\begin{equation}
\widehat{f}_{1n}(\epsilon)
 =
\frac{1} {b_1\sum_{i=1}^n \mathds{1}
\left(X_i\in\mathcal{X}_0\right)}
\sum_{i=1}^n \mathds{1} \left(
X_i\in \mathcal{X}_0\right)
K_1\left(\frac{\widehat{\varepsilon}_i-\epsilon}{b_1}\right).
 \label{fnchap}
\end{equation}
This estimator is the so-called two-steps Kernel estimator of
$f(\epsilon)$. In principle, it would be possible to assume that
$\mathcal{X}_0$ grows
 to $\mathcal{X}$ with a negligible rate compared to the bandwidth
 $b_1$.
This would give an estimator close to the more natural Kernel
estimator $\sum_{i=1}^n
 K\left((\widehat{\varepsilon}_i-\epsilon)/b_1\right)/(nb_1)$.
However, in the rest of the paper, a fixed subset $\mathcal{X}_0$
will
 be considered for the sake of simplicity.

Observe that the two steps Kernel estimator
$\widehat{f}_{1n}(\epsilon)$ is a feasible estimator in the sense
that it does not depend on any unknown
 quantity, as desirable in practice. This contrasts with the unfeasible
 ideal Kernel estimator
\begin{equation}
 \widetilde{f}_{1n}(\epsilon)
 =
\frac{1} {b_1\sum_{i=1}^n \mathds{1}
\left(X_i\in\mathcal{X}_0\right)}
\sum_{i=1}^n \mathds{1}\left(
X_i \in \mathcal{X}_0\right)
 K_1\left(\frac{\varepsilon_i-\epsilon}{b_1}\right),
 \label{fn}
\end{equation}
 which depends in particular on the unknown regression error
terms. It
 is however intuitively clear that $\widehat{f}_{1n}(\epsilon)$ and
 $\widetilde{f}_{1n}(\epsilon)$ should be closed, as illustrated by the results of the
next section.

\section{Assumptions}

The following assumptions are used in our mains results.
 \vskip 0.3cm \noindent
 {$\bf(A_1)$}
 {\it
 The support $\mathcal{X}$ of $X$  is a  compact subset of
 $\Rit^d$
and $\mathcal{X}_0$
 is an inner closed subset of $\mathcal{X}$ with non empty
 interior,
}
\medskip
\\{$\bf(A_2)$}
{\it the p.d.f. $g(\cdot)$ of the i.i.d. covariates $X, X_i$  is
strictly positive over $\mathcal{X}_0$, and has continuous second
order partial derivatives  over $\mathcal{X}$, }
\medskip
\\{$\bf(A_3)$}
{\it
 the regression function $m(\cdot)$ has continuous second order partial
derivatives  over  $\mathcal{X}$, }
\medskip
\\{$\bf(A_4)$}
{\it the i.i.d. centered error regression terms $\varepsilon,
 \varepsilon_i$'s, have finite 6th moments, and are independent of the
 covariates $X,X_i$'s,
}
\\{$\bf(A_5)$}
{\it
 the  probability  density function $f(\cdot)$ has  bounded continuous
 second order
 derivatives over  $\Rit$ and satisfies, for $h_p (e) = e^p f(e)$,
 $\sup_{e\in\Rit}|h_p^{(k)}(e)|<\infty$,
 $p\in[0,2]$, $k\in[0,2]$,
  }
\medskip
\\{$\bf(A_6)$}
{\it the p.d.f  $\varphi$ of $(X, Y)$  has bounded continuous
second  order partial derivatives over $\Rit^d\times\Rit$, }
\\{$\bf(A_7)$}
{\it the Kernel $K_0$  is  symmetric, continuous over $\Rit^d$
with support  contained in $[-1/2, 1/2]^d$ and $\int\! K_0 (z) dz
= 1$, }
\medskip
\\{$\bf(A_8)$}
{\it
 the Kernel $K_{1}$ has a compact support,  is three times
 continuously differentiable over
 $\Rit$, and satisfies $\int\! K_1 (v) dv = 1$ and
$\int\! v K_1 (v) dv = 0$, }
\medskip
\\{$\bf(A_{9})$}
{\it
 the bandwidth $b_0$ decreases to $0$ and satisfies,
 for $d^*=\sup\{d+2,2d\}$, $nb_0^{d^*}/\ln n\rightarrow\infty$ and
 $\ln(1/b_0)/\ln(\ln n)\rightarrow\infty$
 when $n\rightarrow\infty$,
}
\\
{$\bf(A_{10})$}
 {\it
the bandwidth $b_1$ decreases to $0$ and satisfies
$n^{(d+8)}b_1^{7(d+4)}\rightarrow\infty$ when
$n\rightarrow\infty$. }

\vskip 0.3cm \noindent
 Assumptions $(A_2)$, $(A_3)$, $(A_5)$ and $(A_6)$ impose that all the
functions
 to be estimated nonparametrically have two bounded derivatives.
 Consequently the conditions  $\int\! z K_0 (z) dz = 0$ and
$\int\!v K_1 (v) dv = 0$, as assumed in $(A_7)$ and $(A_8)$,
represent
 standard conditions ensuring  that the bias of the resulting
 nonparametric estimators (\ref{mchap})
and (\ref{fn}) are of order $b_0^2$ and $b_1^2$.
 Assumption $(A_4)$ states independence between
the regression error terms
 and the covariates, which is the main condition for (\ref{Naive}) to
 hold.
The differentiability of $K_1$ imposed in $(A_8)$ is more specific
to our
 two-steps estimation method. Assumption $(A_8)$ is used to expand the
 two-steps Kernel estimator
$\widehat{f}_{1n}$ in (\ref{fnchap}) around the unfeasible one
$\widetilde{f}_{1n}$ from
 (\ref{fn}), using the residual error estimation
 $\widehat{\varepsilon}_i - \varepsilon_i$'s and the derivatives of $K_1$
 up to third order.
Assumption $(A_9)$ is useful for obtaining the uniform convergence
 of the regression estimator $\widehat{m}_n$ defined in (\ref{mchap}) (see for
 instance Einmahl and Mason, 2005), and also gives a similar consistency
 result for the leave-one-out estimator $\widehat{m}_{in}$ in
 (\ref{mchapi}). Assumption $(A_{10})$ is needed in the study of
  the difference between the
 feasible estimator $\widehat{f}_{1n}$ and the unfeasible
  estimator $\widetilde{f}_{1n}$.

\section{Main results}
This section is devoted to our main results. The first result we
give here concerns the pointwise consistency of the nonparamatric
Kernel estimator $\widehat{f}_{1n}$ of the density $f$. Next, the
optimal first-step and second-step bandwidths used to estimated
$f$ are proposed. We finish this section by establishing an
asymptotic normality for the estimator $\widehat{f}_{1n}$.

\subsection{Pointwise weak consistency}

The next  result gives the order of the difference between the
feasible estimator and the theoretical density of the regression
error at a fixed point $\epsilon$.

\begin{theorem}
Under $(A_1)-(A_5)$ and $(A_7)-(A_{10})$, we have, when $b_0$ and
$b_1$ go to $0$,
$$
\widehat{f}_{1n}(\epsilon)-f(\epsilon)
 =
 O_{\prob}
\biggl(AMSE(b_1)+R_n(b_0, b_1)\biggr)^{1/2},
$$
where
$$
AMSE(b_1)
 =
 \esp_n
 \left[
\left(\widetilde{f}_{1n}(\epsilon)-f(\epsilon)\right)^2
 \right]
 =
 O_{\prob}
\left(b_1^4 + \frac{1}{nb_1}\right),
$$
and
\begin{eqnarray*}
 R_n(b_0, b_1)
 =
 b_0^4
 +
 \left[
 \frac{1}{(nb_1^5)^{1/2}}
 +
 \left(\frac{b_0^d}{b_1^3}\right)^{1/2}
 \right]^2
 \left(
 b_0^4
 +
 \frac{1}{nb_0^d}
 \right)^2
 +
 \left[
 \frac{1}{b_1}
 +
 \left(\frac{b_0^d}{b_1^7}\right)^{1/2}
 \right]^2
 \left(b_0^4+\frac{1}{nb_0^d}\right)^3.
\end{eqnarray*}
 \label{thm1}
\end{theorem}

\noindent The result of Theorem \ref{thm1} is based on the
evaluation of the difference between  $\widehat{f}_{1n}(\epsilon)$
and $\widetilde{f}_{1n}(\epsilon)$.  This evaluation  gives an
indication about the impact of the estimation of the residuals on
the nonparametric estimation of the regression error density.

\subsection{Optimal first-step and second-step bandwidths
for the pointwise weak consistency}

 As shown in the next result, Theorem \ref{Optimbandw1}
gives some guidelines for the choice of the optimal bandwidth
$b_0$ used in the nonparametric regression errors estimation. As
far as we know, the choice of an optimal $b_0$ has not been
addressed before. In what follows, $a_n \asymp b_n$ means that $
a_n= O (b_n)$ and $b_n = O(a_n)$, i.e. that there is a constant
$C>0$ such that $|a_n|/C \leq |b_n| \leq C |a_n|$ for $n$ large
enough.

\begin{theorem}
 Suppose that
$(A_1)-(A_5)$ and $(A_7)-(A_{10})$ are satisfied, and  define
$$
b_0^* = b_0^*(b_1)
 =
\arg\min_{b_0}
 R_n (b_0, b_1).
$$
where the minimization is performed over bandwidth $b_0$
fulfilling $(A_9)$. Then the bandwidth $b_0^*$ satisfies
$$
b_0^*
\asymp
\max
\left\lbrace
\left(\frac{1}{n^2b_1^3}\right)^{\frac{1}{d+4}} ,
\left(\frac{1}{n^3b_1^7}\right)^{\frac{1}{2d+4}}
\right\rbrace,
$$
and we have
$$
R_n(b_0^*, b_1)
\asymp
\max
\left\lbrace
\left(\frac{1}{n^2b_1^3}\right)^{\frac{4}{d+4}} ,
\left(\frac{1}{n^3b_1^7}\right)^{\frac{4}{2d+4}}
 \right\rbrace.
$$
\label{Optimbandw1}
\end{theorem}
Our next theorem gives the conditions for which the estimator
$\widehat{f}_{1n}(\epsilon)$ reaches the optimal rate $n^{-2/5}$
when $b_0$ takes the value $b_0^*$.
 We prove that for
$d\leq 2$, the bandwidth that minimizes the term
$AMSE(b_1)+R_n(b_0^*, b_1)$ has the same order as $n^{-1/5}$,
yielding the optimal order $n^{-2/5}$ for
$\left(AMSE(b_1)+R_n(b_0^*, b_1)\right)^{1/2}$.

\begin{theorem}
Assume that $(A_1)-(A_5)$ and $(A_7)-(A_{10})$ are satisfied, and
set
$$
b_1^*
 =
\arg\min_{b_1} \biggl(AMSE(b_1)+R_n(b_0^*,b_1)\biggr),
$$
where $b_0^*=b_0^*(b_1)$ is defined as in Theorem
\ref{Optimbandw1}.
 Then
\begin{enumerate}
\item For $d\leq 2$, the bandwidth $b_1^*$ satisfies
$$
b_1^* \asymp \left(\frac{1}{n}\right)^{\frac{1}{5}},
$$
and we have
$$
\biggl( AMSE(b_1^*) +
 R_n(b_0^*, b_1^*)
\biggr)^{\frac{1}{2}}
\asymp
\left(\frac{1}{n}\right)^{\frac{2}{5}}.
$$
\item For $d\geq 3$, $b_1^*$ satisfies
$$
b_1^*
\asymp
\left(\frac{1}{n}\right)^{\frac{3}{2d+11}},
$$
and we have
$$
\biggl(
 AMSE(b_1^*)
 +
 R_n(b_0^*, b_1^*)
\biggr)^{\frac{1}{2}}
\asymp
\left(\frac{1}{n}\right)^{\frac{6}{2d+11}}.
$$
\end{enumerate}
\label{Optimbandw2}
\end{theorem}

\vskip 0.3cm\noindent The results of Theorem \ref{Optimbandw2}
show that the rate $n^{-2/5}$ is reachable if and only when $d\leq
2$. These results are derived  from Theorem \ref{Optimbandw1}.
This latter indicates that if $b_1$ is proportional to $n^{-1/5}$,
the bandwidth $b_0^*$ has the same order as
$$
\max \left\lbrace
\left(\frac{1}{n}\right)^{\frac{7}{5(d+4)}}
,
\left(\frac{1}{n}\right)^{\frac{8}{5(2d+4)}}
\right\rbrace
=
\left(\frac{1}{n}\right)^{\frac{8}{5(2d+4)}}.
$$
 For $d\leq 2$, this order of $b_0^*$ is smaller
than the one of the optimal bandwidth $b_{0*}$ obtained  for
pointwise or mean square estimation of $m(\cdot)$ using a Kernel
estimator. In fact, it has been shown in  Nadaraya (1989, Chapter
4) that the optimal bandwidth $b_{0*}$ for estimating $m(\cdot)$
is obtained by minimizing the order of  the risk function
$$
r_n(b_0)
=
\esp\left[
\int \mathds{1}
\left(x\in\mathcal{X}\right)
\left(\widehat{m}_n(x)-m(x)\right)^2
 \widehat{g}_n^2(x)w(x)
 dx
 \right],
$$
where $\widehat{g}_n(x)$ is a nonparametric Kernel estimator of
$g(x)$, and $w(\cdot)$ is a nonnegative weight function, which is
bounded and squared integrable on $\mathcal{X}$. If $g(\cdot)$ and
$m(\cdot)$ have  continuous second order partial derivatives over
their supports, Nadaraya (1989, Chapter 4) shows that $r_n(b_0)$
has the same order as $b_0^4+\left(1/(nb_0^d)\right)$, leading to
the optimal bandwidth $\widehat{b}_{0}=n^{-1/(d+4)}$  for the
convergence of the estimator $\widehat{m}_n(\cdot)$ of $m(\cdot)$
in the set of the square integrable functions on $\mathcal{X}$.

\vskip 0.2cm\noindent For d=1, the optimal order of $b_0^*$ is
$n^{-(1/5)\times(4/3)}$ which goes to 0 slightly faster than
$n^{-1/5}$, the optimal order of the bandwidth $\widehat{b}_{0}$
for the mean square nonparametric estimation of $m(\cdot)$.

\vskip 0.2cm\noindent For $d=2$, the optimal order of $b_0^*$ is
$n^{-1/5}$. Again this order goes to 0 faster than the order
$n^{-1/6}$ of the optimal bandwidth for the nonparametric
estimation of the regression function with two covariates.

\vskip 0.2cm\noindent However, for $d\geq 3$, we note that the
order of $b_0^*$ goes to $0$ slowly than $\widehat{b}_{0}$. Hence
our results show that optimal $\widehat{m}_n(\cdot)$ for
estimating $f(\cdot)$ should use a very small bandwidth $b_0$.
This suggests that $\widehat{m}_n(\cdot)$ should be less biased
and should have a higher variance than the optimal Kernel
regression estimator of the estimation setup. Such a finding
parallels Wang, Cai, Brown and Levine (2008) who show that a
similar result hold when estimating the conditional variance of a
heteroscedastic regression error term. However Wang et {\it al.}
(2008) do not give the order of the optimal bandwidth to be used
for estimating the regression function in their heteroscedastic
setup. These results show that estimators of $m(\cdot)$ with
smaller bias should be preferred in our framework, compared to the
case where the regression function $m(\cdot)$ is the parameter of
interest.

\subsection{Asymptotic normality}
 We give  now  an asymptotic normality of the
estimator $\widehat{f}_{1n}(\epsilon)$.

\begin{theorem}

Assume that
$$
{(\rm\bf{A}_{11}):}
\quad nb_0^{d+4}=O(1),
\quad nb_0^4b_1=o(1),
\quad
 nb_0^{d}b_1^3\rightarrow\infty,
$$
when $n$ goes to $\infty$. Then under $(A_1)-(A_5)$,
$(A_7)-(A_{10})$, we have
$$
\sqrt{nb_1}
\left(
\widehat{f}_{1n}(\epsilon)-\overline{f}_{1n}(\epsilon)
\right)
\stackrel{d}{\rightarrow}
\mathcal{N}
\left(
 0,
\frac{f(\epsilon)}
{\prob\left(X\in\mathcal{X}_0\right)} \int
K_1^2(v) dv \right),
$$
where
$$
\overline{f}_{1n}(\epsilon)
 =
f(\epsilon) + \frac{b_1^2}{2}f^{(2)}(\epsilon)
\int v^2 K_1(v) dv
+ o\left(b_1^2\right).
$$
\label{normalite}
\end{theorem}

\noindent The result of this  theorem shows that the best choice
$b_1^*$ for the bandwidth $b_1$ should achieve the minimum of the
Asymptotic Mean Integrated Square Error
$$
{\rm AMISE}
 =
 \frac{b_1^4}{4}
 \int
 \left(f^{(2)}(\epsilon)\right)^2
 d\epsilon
 \left(
 \int
 v^2K_1(v)dv
 \right)^2
 +
 \frac{1}{nb_1\prob\left(X\in\mathcal{X}_0\right)}
 \int K_1^2(v)dv,
$$
leading to the  optimal bandwidth
$$
b_1^* = \left[ \frac{\frac{1}
 {\prob\left(X\in\mathcal{X}_0\right)}
\displaystyle{\int}K_1^2(v)dv} {\displaystyle{\int}
(f^{(2)}(\epsilon))^2d\epsilon \left(
 \displaystyle{\int}
 v^2K_1(v)dv\right)^2}
 \right]^{1/5}
 n^{-1/5}.
$$
We also note that for $d\leq 2$, $b_1=b_1^*$ and $b_0=b_0^*$,
Theorems \ref{Optimbandw2} and  \ref{Optimbandw1} give
$$
b_1 \asymp \left(\frac{1}{n}\right)^{\frac{1}{5}},
\quad b_0
\asymp \left(\frac{1}{n}\right)^{\frac{8}{5(2d+4)}},
$$
which yields that
$$
nb_0^{d+4} \asymp
\left(\frac{1}{n}\right)^{\frac{12-2d}{5(2d+4)}},
\quad nb_0^4b_1
\asymp \left(\frac{1}{n}\right)^{\frac{16-8d}{5(2d+4)}},
\quad
nb_0^db_1^3 \asymp
\left(\frac{1}{n}\right)^{\frac{4d-8}{5(2d+4)}}.
$$
This shows that for $d=1$, the Assumption $(\rm\bf{A}_{11})$ is
realizable with the optimal bandwidths $b_0^*$ and $b_1^*$. But
with these bandwidths, the last constraint of $(\rm\bf{A}_{11})$
is not satisfied for $d=2$, since $nb_0^db_1^3$ is bounded when
$n\rightarrow\infty$.

\newpage
\section{ Conclusion}
 The aim of this chapter was  to study
  the nonparametric Kernel estimation of the
probability density function of the regression error using the
estimated residuals. The difference between the feasible estimator
which uses the estimated residuals and the unfeasible one using
the true residuals are studied. An optimal choice of the
first-step bandwidth used to estimate the residuals is also
proposed. Again, an asymptotic normality of the feasible Kernel
estimator and its rate-optimality are established. One of the
contributions of this paper is the analysis of  the impact of the
estimated residuals on the regression errors p.d.f. Kernel
estimator.

\vskip 0.3cm
 In our setup, the strategy was to use an approach
based on a two-steps procedure which, in a first step, replaces
the unobserved residuals terms by some nonparametric estimators
$\widehat{\varepsilon}_i$. In a second step, the \lg
pseudo-observations\rg  $\widehat{\varepsilon}_i$ are used to
estimate the p.d.f $f(\cdot)$, as if they were the true
$\varepsilon_i$'s. If proceeding so can remedy the curse of
dimensionality, a challenging issue was to measure the impact of
the estimated residuals on the final estimator of $f(\cdot)$ in
the first nonparametric step, and to find the order of the optimal
first-step bandwidth $b_0$. For this choice of $b_0$, our results
indicates that the optimal bandwidth to be used for estimating the
regression function $m(\cdot)$ should be smaller than the optimal
bandwidth for the mean square estimation of $m(\cdot)$. That is to
say, the best estimator $\widehat{m}_n(\cdot)$ of the regression
function $m(\cdot)$ needed for estimating $f(\cdot)$  should have
a lower bias and a higher variance than the optimal Kernel
regression of the estimation setup. With this appropriate choice
of $b_0$, it has been seen that for $d\leq 2$, the nonparametric
estimator $\widehat{f}_{1n}(\epsilon)$ of $f$ can reach the
optimal rate $n^{-2/5}$, which
 corresponds to the exact consistency rate reached for the
Kernel density estimator of real-valued variable. Hence our main
conclusion is that for $d\leq 2$, the  estimator
$\widehat{f}_{1n}(\epsilon)$  used for estimating $f(\epsilon)$ is
not affected by the curse of dimensionality, since there is no
negative effect coming from the estimation of the residuals on the
final estimator of $f(\epsilon)$.

\newpage

\setcounter{subsection}{0} \setcounter{equation}{0}
\renewcommand{\theequation}{\thesection.\arabic{equation}}

\section{ Proofs section}

\subsection*{Intermediate Lemmas for Proposition
\ref{Prop} and Theorem \ref{thm1}}

\begin{lemma}
Define, for $x\in\mathcal{X}_0$,
$$
\widehat{g}_n(x)
 =
 \frac{1}{nb_0^d}
 \sum_{i=1}^n
 K_0
 \left(
 \frac{X_i-x}{b_0}
 \right),
 \quad
 \overline{g}_n(x)
 =
 \esp\left[\widehat{g}_n(x)\right].
$$
Then under $(A_1)-(A_2)$, $(A_4)$, $(A_7)$ and $(A_9)$,
 we have, when $b_0$ goes to $0$,
$$
\sup_{x\in\mathcal{X}_0} \left|\overline{g}_{n}(x)-g(x)\right|
 =
O\left(b_0^{2}\right), \quad \sup_{x\in\mathcal{X}_0}
 \left|
\widehat{g}_{n}(x)
 -
 \overline{g}_n(x)
\right|
 = O_{\prob}
 \left(
 b_0^4 +
 \frac{ \ln n}{nb_0^d}
 \right)^{1/2},
$$
 and
$$
\sup_{x\in\mathcal{X}_0}
 \left|
 \frac{1}{\widehat{g}_{n}(x)}
   -
\frac{1}{g(x)} \right| =
 O_{\prob}
 \left(
 b_0^4
 +
 \frac{\ln n}{nb_0^d}
\right)^{1/2}.
$$
\label{Estig}
\end{lemma}

\begin{lemma}
Under $(A_1)-(A_4)$, $(A_7)$ and $(A_9)$, we have
$$
\sup_{x\in\mathcal{X}_0} \left| \widehat{m}_n(x)-m(x)
 \right|
 =
O_{\prob}
 \left(
 b_0^4+\frac{\ln n}{nb_0^d}
 \right)^{1/2}.
$$
\label{Estim}
\end{lemma}

\begin{lemma}
Define for $(x,y)\in\Rit^d\times\Rit$,
\begin{eqnarray*}
f_{n}(\epsilon|x)
 =
 \frac{\frac{1}{nh_0^dh_1}
 \sum_{i=1}^n
 K_0\left(\frac{X_i-x}{h_0}\right)
 K_1\left(\frac{Y_i-m(x)-\epsilon}{h_1}\right)
 }{\frac{1}{nh_0^d}
 \sum_{i=1}^n
 K_0\left(\frac{X_i-x}{h_0}\right)}\;,
\end{eqnarray*}
 Then under $(A_1)-(A_3)$, $(A_6)-(A_9)$,
  we have, when $n$ goes to infinity,
$$
\widetilde{f}_n(\epsilon|x)
 -
f_{n}(\epsilon|x)
 =
o_{\prob} \left( \frac{1}{nh_0^dh_1} \right)^{1/2}.
$$
\label{Reste}
\end{lemma}

\begin{lemma}
 Set, for $(x,y)\in\Rit^d\times\Rit$,
$$
\widetilde{\varphi}_{in}(x, y)
 =
 \frac{1}{h_0^dh_1}
 K_0\left(\frac{X_i-x}{h_0}\right)
 K_1\left(\frac{Y_i-y}{h_1} \right).
$$
Then, under  $(A_6)-(A_8)$, we have, for $x$ in
 $\mathcal{X}_0$ and $y$ in $\Rit$,
  $h_0$ and $h_1$ going to $0$, and for some constant $C>0$,
\begin{eqnarray*}
\esp \left[ \widetilde{\varphi}_{in} \left(x, y\right) \right] -
\varphi\left(x, y\right) &=& \frac{h_0^2}{2}
\frac{\partial^2\varphi (x, y)} {\partial^2 x} \int z
K_0(z)z^{\top} dz
 +
\frac{h_1^2}{2} \frac{\partial^2\varphi(x, y)} {\partial^2 y} \int
v^2 K_1(v) dv
\\
&& +\;
 o\left(h_0^{2}+h_1^2\right),
\\
\Var \left[ \widetilde{\varphi}_{in} \left(x, y\right) \right]
 &=&
\frac{\varphi\left(x, y\right)} {h_0^dh_1}
 \int\int
 K_0^2(z)
 K_1^2(v)
 dv dz
+ o\left(\frac{1}{h_0^dh_1} \right),
\\
\esp \left[ \left| \widetilde{\varphi}_{in}\left(x,y\right) - \esp
\widetilde{\varphi}_{in}\left(x, y\right) \right|^3 \right] &\leq&
 \frac
{C\varphi\left(x,y\right)} {h_0^{2d}h_1^2} \int \int \left| K_0(z)
K_1\left(v\right) \right|^3 dz dv +
o\left(\frac{1}{h_0^{2d}h_1^2}\right).
\end{eqnarray*}
\label{Momvarphi}
\end{lemma}

\begin{lemma}
Set
$$
f_{in}(\epsilon)
= \frac{\mathds{1}\left(X_i \in
\mathcal{X}_0\right)} {b_1\prob\left(X\in\mathcal{X}_0\right)}
K_1\left( \frac{\varepsilon_i - \epsilon}{b_1}\right).
$$
Then under  $(A_4)$, $(A_5)$ and $(A_8)$, we have, for $b_1$ going
to $0$, and for some constant $C>0$,
\begin{eqnarray*}
\esp
 f_{in}(\epsilon)
&=& f(\epsilon) + \frac{b_1^2}{2} f^{(2)}(\epsilon) \int v^2
K_1(v) dv + o\left(b_1^2\right),
\\
\Var \left(f_{in}(\epsilon)\right) &=& \frac {f(\epsilon)}
{b_1\prob\left(X\in\mathcal{X}_0\right)} \int K_1^2(v) dv +
o\left(\frac{1}{b_1}\right),
\\
\esp \left| f_{in}(\epsilon) - \esp
 f_{in}(\epsilon)
\right|^3 &\leq&
 \frac{
C f(\epsilon)}
 {b_1^2\prob^2\left(X\in\mathcal{X}_0\right)}
 \int
 \left|
 K_1(v)
\right|^3 dv + o\left(\frac{1}{b_1^2}\right).
\end{eqnarray*}
\label{Momfin1}
\end{lemma}

\begin{lemma}
Define
\begin{eqnarray*}
S_n
 &=&
 \sum_{i=1}^n
 \mathds{1}
 \left(X_{i} \in \mathcal{X}_0\right)
 \left(\widehat{m}_{in}-m(X_{i})\right)
 K_1^{(1)}
 \left(\frac{\varepsilon_{i}-\epsilon}{b_1}\right),
 \\
 T_n
 &=&
 \sum_{i=1}^n
 \mathds{1}
 \left(X_{i} \in \mathcal{X}_0\right)
 \left(\widehat{m}_{in}-m(X_{i})\right)^2
 K_1^{(2)}
 \left(\frac{\varepsilon_{i}-\epsilon}{b_1}\right),
 \\
 R_n
 &=&
 \sum_{i=1}^n
 \mathds{1}
 \left( X_{i}
 \in\mathcal{X}_0\right)
 \left(\widehat{m}_{in}-m(X_{i})\right)^3
 \int_{0}^{1}
 (1-t)^2
 K_1^{(3)}
 \left(
 \frac{
 \varepsilon_i-t(\widehat{m}_{in} -m(X_i))-\epsilon}{b_1}
 \right)
 dt.
\end{eqnarray*}
Then under $(A_1)-(A_5)$ and $(A_7)-(A_{10})$, we have, for $b_0$
and $b_1$ small enough,
\begin{eqnarray*}
S_n
 &=&
 O_{\prob}
\left[ b_0^2\left(nb_1^2+(nb_1)^{1/2}\right)
 +
\left( nb_1^4+\frac{b_1}{b_0^d} \right)^{1/2}
 \right],
 \\
 T_{n}
 &=&
 O_{\prob}
\left[ \left( nb_1^3 + \left(nb_1\right)^{1/2} +
\left(n^2b_0^db_1^{3}\right)^{1/2} \right)
\left(b_0^4+\frac{1}{nb_0^d}\right) \right],
\\
R_n &=& O_{\prob} \left[ \left( nb_1^3 +
\left(n^2b_0^db_1\right)^{1/2} \right)
\left(b_0^4+\frac{1}{nb_0^d}\right)^{3/2} \right].
\end{eqnarray*}
\label{STR}
\end{lemma}

\begin{lemma}
Under  $(A_5)$ and $(A_8)$ we have, for some constant $C>0$, and
for any $\epsilon$ in $\Rit$ and $p\in[0,2]$,
\begin{eqnarray}
\left|
 \int
 K_1^{(1)}
 \left( \frac{e-\epsilon}{b_1} \right)^2
 e^pf(e) de
 \right|
 \leq C b_1,
 &&
 \left|
  \int
  K_1^{(1)}
  \left(
 \frac{e-\epsilon}{b_1}
 \right)
 e^pf(e)de
  \right|
 \leq C b_1^2,
\label{MomderK1}
\\
\left| \int
 K_1^{(2)}
 \left( \frac{e-\epsilon}{b_1} \right)^2
  e^pf(e)
  de
\right|
 \leq
  C b_1,
  &&
 \left|
\int K_1^{(2)}
 \left(
 \frac{e-\epsilon}{b_1}
 \right)
 e^pf(e)de
\right|
 \leq C b_1^3,
 \label{MomderK2}
\\
\left| \int
 K_1^{(3)}
 \left( \frac{e-\epsilon}{b_1} \right)^2
  e^pf(e)
  de
\right|
 \leq
  C b_1,
  &&
 \left|
\int K_1^{(3)}
 \left(
 \frac{e-\epsilon}{b_1}
 \right)
 e^pf(e)de
\right|
 \leq
 C b_1^3.
 \label{MomderK3}
\end{eqnarray}
\label{MomderK}
\end{lemma}

\begin{lemma}
Set
$$
\beta_{in}
 =
\frac{\mathds{1} \left(X_i\in\mathcal{X}_0\right)}
{nb_0^d\widehat{g}_{in}} \sum_{j=1, j\neq i}^n
\left(m(X_j)-m(X_i)\right) K_0\left(\frac{X_j-X_i}{b_0}\right).
$$
Then, under $(A_1)-(A_5)$ and $(A_7)-(A_{10})$, we have, when
$b_0$ and $b_1$ go to $0$,
\begin{eqnarray*}
\sum_{i=1}^n \beta_{in}
 K_1^{(1)}
 \left(
 \frac{\varepsilon_i-\epsilon}{b_1}
 \right)
 =
 O_{\prob}
 \left(b_0^2\right)
 \left(nb_1^2+(nb_1)^{1/2}\right).
\end{eqnarray*}
\label{Betasum}
\end{lemma}

\begin{lemma}
Set
$$
\Sigma_{in}
 =
\frac{\mathds{1}\left(X_i\in\mathcal{X}_0\right)}
{nb_0^d\widehat{g}_{in}} \sum_{j=1, j\neq i}^n
 \varepsilon_j
K_0\left(\frac{X_j-X_i}{b_0}\right).
$$
Then, under $(A_1)-(A_5)$ and $(A_7)-(A_{10})$, we have
\begin{eqnarray*}
\sum_{i=1}^n \Sigma_{in}
 K_1^{(1)}
 \left(
 \frac{\varepsilon_i-\epsilon}{b_1}
 \right)
 =
 O_{\prob}
 \left(
 nb_1^4
 +
 \frac{b_1}{b_0^d}
\right)^{1/2}.
\end{eqnarray*}
\label{Sigsum}
\end{lemma}

\begin{lemma}
Let $\esp_n[\cdot]$ be the conditional mean given
$X_1,\ldots,X_n$. Then under $(A_1)-(A_5)$ and $(A_7)-(A_{9})$,
 we have, for $b_0$ going to $0$,
\begin{eqnarray*}
\sup_{1\leq i\leq n}
 \esp_{n}
  \biggl[
  \mathds{1}
  \left(X_i \in \mathcal{X}_0\right)
 (\widehat{m}_{in} - m(X_i))^4
 \biggr]
 &=&
 O_{\prob}
 \left(b_0^4+\frac{1}{nb_0^d}\right)^2,
 \\
 \sup_{1\leq i\leq n}
 \esp_{n}
 \biggl[
 \mathds{1}
 \left(X_i \in \mathcal{X}_0\right)
 (\widehat{m}_{in} - m(X_i))^6
 \biggr]
 &=&
 O_{\prob}
 \left(b_0^4+\frac{1}{nb_0^d}\right)^3.
\end{eqnarray*}
\label{BoundEspmchap}
\end{lemma}

\begin{lemma}
Assume that  $(A_4)$ and $(A_7)$ hold. Then, for any $1 \leq i
\neq j \leq n$, and for any $\epsilon$ in $\Rit$,
$$
\left(\widehat{m}_{in} - m(X_i), \varepsilon_i\right)
 \mbox{\it and }
 \left( \widehat{m}_{jn} - m(X_j), \varepsilon_j\right)
$$
are independent given $X_1, \ldots, X_n$,
 provided that $\| X_i - X_j \| \geq C b_0$,
for some constant $C>0$. \label{Indep}
\end{lemma}

\begin{lemma}
Let $\Var_n(\cdot)$ and $Cov_n(\cdot)$ be respectively the
conditional variance and the conditional covariance given
$X_1,\ldots,X_n$, and set
\begin{eqnarray*}
\zeta_{in}
 =
\mathds{1} \left(X_i \in \mathcal{X}_0\right) (\widehat{m}_{in} -
m(X_i))^2 K_1^{(2)} \left( \frac{\varepsilon_i-\epsilon}{b_1}
\right).
\end{eqnarray*}
Then under $(A_1)-(A_5)$ and $(A_7)-(A_{9})$, we have, for $n$
going to infinity,
\begin{eqnarray*}
\sum_{i=1}^n \Var_n\left( \zeta_{in}\right) &=&
O_{\prob}\left(nb_1\right)
 \left(b_0^4+\frac{1}{nb_0^d}\right)^2,
 \\
 \sum_{i=1}^n
 \sum_{j=1\atop j\neq i}^n
  \Cov_n
 \left(
 \zeta_{in},\zeta_{jn}
 \right)
&=& O_{\prob}\left(n^2b_0^db_1^{7/2}\right)
\left(b_0^4+\frac{1}{nb_0^d}\right)^2.
\end{eqnarray*}
\label{sumzeta}
\end{lemma}

\noindent
All these lemmas are proved in Appendix A.

\newpage
\subsection*{Proof of Proposition \ref{Prop}}
Define $f_{n}(\epsilon|x)$ as in Lemma \ref{Reste}, and note that
by this lemma, we have
 \begin{eqnarray}
 \widetilde{f}_n(\epsilon|x)
 =
 f_{n}(\epsilon|x)
 +
 o_{\prob}
\left( \frac{1}{nh_0^dh_1} \right)^{1/2}.
 \label{fncond1}
\end{eqnarray}
The asymptotic distribution of the first term in (\ref{fncond1})
is derived by applying the Lyapounov Central Limit Theorem for
triangular arrays (see e.g Billingsley 1968, Theorem 7.3). Define
for $x\in\mathcal{X}_0$ and $y\in\Rit$,
\begin{eqnarray*}
 \widetilde{\varphi}_n(x,y)
 =
 \frac{1}{nh_0^dh_1}
\sum_{i=1}^n
 K_0\left(\frac{X_i-x}{h_0}\right)
 K_1\left(\frac{Y_i-y}{h_1}\right),
 \quad
 \widetilde{g}_n(x)
 =
 \frac{1}{nh_0^d}
 \sum_{i=1}^n
 K_0\left(\frac{X_i-x}{h_0}\right),
\end{eqnarray*}
and observe that
\begin{eqnarray}
f_{n}(\epsilon|x)
 =
\frac{\widetilde{\varphi}_n\left(x,m(x)+\epsilon\right)}
{\widetilde{g}_n(x)}\;.
 \label{fncond2}
\end{eqnarray}
Let now $\widetilde{\varphi}_{in}(x,y)$ be as in Lemma
\ref{Momvarphi}, and note that
\begin{eqnarray}
\widetilde{\varphi}_n(x,y)
 =
 \frac{1}{n}
 \sum_{i=1}^n
 \biggl(
 \widetilde{\varphi}_{in}(x,y)
 -
 \esp
 \left[
 \widetilde{\varphi}_{in}(x,y)
 \right]
 \biggr)
 +
 \esp
 \left[
 \widetilde{\varphi}_{1n}(x,y)
 \right].
 \label{fnoncond}
\end{eqnarray}
The second and third inequalities in Lemma \ref{Momvarphi} give,
since  $h_0^d h_1$ goes to $0$,
\begin{eqnarray*}
\frac{
\sum_{i=1}^n
\esp
\left|
 \widetilde{\varphi}_{in}(x,y)
 -
 \esp
 \widetilde{\varphi}_{in}(x,y)
\right|^3 } { \left( \sum_{i=1}^n
 \Var
 \left[
 \widetilde{\varphi}_{in}(x,y)
 \right]
\right)^3 }
\leq
\frac { \frac { Cn\varphi(x,y) } {h_0^{2d}h_1^2}
\displaystyle{\int}
\displaystyle{\int} \left| K_0(z) K_1(v)
\right|^3 dz dv
+
o \left( \frac{n}{h_0^{2d}h_1^2} \right) } {
\left( \frac { n\varphi(x,y) }{h_0^{d}h_1}
 \displaystyle{\int}
\displaystyle{\int}
 K_0^2(z) K_1^2(v) dv dz + o \left(
\frac{n}{h_0^{d}h_1} \right)
 \right)^3
 }
 =O(h_0^d h_1) = o(1).
\end{eqnarray*}
 Hence  the Lyapounov Central Limit
Theorem gives, since $nh_0^dh_1$ diverges under $(\rm\bf{A}_0)$,
$$
\frac { \sum_{i=1}^n
 \left\lbrace
 \widetilde{\varphi}_{in}(x,y)
 -
\esp
\left[
\widetilde{\varphi}_{in}(x,y)
\right]
\right\rbrace}{
\left( \sum_{i=1}^n \Var
\left[
\widetilde{\varphi}_{in}(x,y)
\right]
\right)^{1/2}}
\stackrel{d}{\rightarrow}\mathcal{N}(0,1),
$$
so that
\begin{eqnarray}
\frac{\sqrt{nh_0^dh_1}}{n}
 \sum_{i=1}^n
 \biggl(
 \widetilde{\varphi}_{in}(x,y)
 -
 \esp
 \left[
 \widetilde{\varphi}_{in}(x,y)
 \right]
 \biggr)
 \stackrel{d}{\rightarrow}
 \mathcal{N}
 \left(
  0,
  \varphi(x,y)
  \int\int
  K_0^2(z)
  K_1^2(v)
  dz dv
 \right).
\label{normal}
\end{eqnarray}
Further, a similar proof as the one of Lemma \ref{Estig} gives
\begin{eqnarray}
\frac{1}{\widetilde{g}_n(x)}
 =
 \frac{1}{g(x)}
 +
 O_{\prob}
 \left(
 h_0^4+\frac{\ln n}{nh_0^d}
 \right)^{1/2}.
 \label{gtilde}
\end{eqnarray}
Hence by this equality, it follows that, taking $y=m(x)+\epsilon$
in (\ref{normal}), and by (\ref{fncond1})-(\ref{fnoncond}),
\begin{eqnarray*}
\sqrt{nh_0^dh_1}
 \left(
\widetilde{f}_n(\epsilon|x)
 -
\overline{f}_n(\epsilon|x) \right)
\stackrel{d}{\rightarrow}
\mathcal{N}
\left( 0,
\frac{f(\epsilon|x)}{g(x)}
\int\int
 K_0^2(z)
 K_1^2(v)
 dz dv
\right),
\end{eqnarray*}
where
$$
\overline{f}_n(\epsilon|x)
=
\frac{
\esp\left[
\widetilde{\varphi}_{1n}
\left(x,m(x)+\epsilon\right)
\right]}
{\widetilde{g}_n(x)}\;.
$$
This yields  the result of Proposition \ref{Prop}, since the first
equality of Lemma \ref{Momvarphi} and (\ref{gtilde}) yield, for
$h_0$ and $h_1$ small enough,
\begin{eqnarray*}
\overline{f}_n(\epsilon|x)
 &=&
 f(\epsilon| x)
+ \frac{h_0^2}{2g(x)} \frac{\partial^2\varphi\left(x,
m(x)+\epsilon\right)} {\partial^2 x} \int z K_0(z) z^{\top} dz
 \\
 &&
 +\;
 \frac{h_1^2}{2g(x)}
 \frac{\partial^2\varphi (x,m(x)+\epsilon)}
 {\partial^2 y}
 \int v^2 K_1 (v) dv
 +
 o\left(h_0^2+h_1^2\right).
 \eop
\end{eqnarray*}

\subsection*{Proof of Theorem \ref{thm1}}
The proof of the theorem is based upon the following equalities:
\begin{eqnarray}
\nonumber
 \widehat{f}_{1n}(\epsilon)-\widetilde{f}_{1n}(\epsilon)
&=& O_{\prob}
 \left[
 b_0^2
 +
 \left(
 \frac{1}{n}
 +
 \frac{1}{n^2b_0^db_1^3}
 \right)^{1/2}
 \right]
 +
 O_{\prob}
 \left[
 \frac{1}{(nb_1^5)^{1/2}}
 +
 \left(\frac{b_0^d}{b_1^3}\right)^{1/2}
 \right]
 \left(
 b_0^4
 +
 \frac{1}{nb_0^d}
 \right)
 \\
 &&
 +\;
 O_{\prob}
 \left[
 \frac{1}{b_1}
 +
 \left(\frac{b_0^d}{b_1^7}\right)^{1/2}
 \right]
 \left(b_0^4+\frac{1}{nb_0^d}\right)^{3/2},
 \label{fnchapfn1}
 \end{eqnarray}
 and
\begin{eqnarray}
\widetilde{f}_{1n}(\epsilon)-f(\epsilon)
 =
 O_{\prob}
 \left(b_1^4+\frac{1}{nb_1}\right)^{1/2}.
 \label{fnchapfn2}
 \end{eqnarray}
Indeed, since $\widehat{f}_{1n}(\epsilon)-f(\epsilon)
=\left(\widetilde{f}_{1n}(\epsilon)-f(\epsilon)\right)
+\widehat{f}_{1n}(\epsilon)-\widetilde{f}_{1n}(\epsilon)$, it then
follows by (\ref{fnchapfn2}) and (\ref{fnchapfn1}) that
\begin{eqnarray*}
\widehat{f}_{1n}(\epsilon)-f(\epsilon)
&=&
O_{\prob}
\left[
 b_1^4
 +
 \frac{1}{nb_1}
 +
 b_0^4
 +
 \frac{1}{n}
 +
 \frac{1}{n^2b_0^db_1^3}
 +
 \left(
 \frac{1}{(nb_1^5)^{1/2}}
 +
 \left(\frac{b_0^d}{b_1^3}\right)^{1/2}
 \right)^2
 \left(
 b_0^4
 +
 \frac{1}{nb_0^d}
 \right)^2
 \right]^{1/2}
 \\
 &&
 +\;
 O_{\prob}
 \left[
 \left(
 \frac{1}{b_1}
 +
 \left(\frac{b_0^d}{b_1^7}\right)^{1/2}
 \right)^2
 \left(b_0^4+\frac{1}{nb_0^d}\right)^3
 \right]^{1/2}.
 \end{eqnarray*}
This yields the result of the Theorem,  since  under $(A_9)$ and
$(A_{10})$, we have
\begin{eqnarray*}
 \frac{1}{n}
  =
 O\left(
 \frac{1}{nb_1}
 \right),
 \quad
\frac{1}{n^2b_0^db_1^3}
= O\left(\frac{b_0^d}{b_1^3}\right)
\left(b_0^4+\frac{1}{nb_0^d}\right)^2.
 \end{eqnarray*}
 Hence, it remains
to prove (\ref{fnchapfn1}) and (\ref{fnchapfn2}). For this,
 define $S_n$, $R_n$ and $T_n$ as in Lemma \ref{STR}.
Since $\widehat{\varepsilon}_i-\varepsilon_i =
-\left(\widehat{m}_{in}-m(X_{i})\right)$ and that $K_1$ is three
times continuously differentiable under $(A_8)$,  the third-order
Taylor expansion
 with integral remainder gives
\begin{eqnarray*}
 \widehat{f}_{1n} (\epsilon)- \widetilde{f}_{1n}(\epsilon)
 &=&
\frac{1}{b_1\sum_{i=1}^n \mathds{1}(X_i \in\mathcal{X}_0)}
\sum_{i=1}^n \mathds{1} \left(X_i\in \mathcal{X}_0\right)
\left[
 K_1\left(
 \frac{\widehat{\varepsilon}_i-\epsilon}{b_1}
 \right)
 -
 K_1\left(
 \frac{\varepsilon_i-\epsilon}{b_1}
 \right)
 \right]
 \\
& = &
 -\frac{1}{b_1\sum_{i=1}^n\mathds{1}(X_i \in\mathcal{X}_0)}
\left( \frac{S_n}{b_1} - \frac{T_n}{2b_1^2} + \frac{R_n}{2b_1^3}
\right).
\end{eqnarray*}
Therefore, since
$$
\sum_{i=1}^n \mathds{1} \left(X_i \in \mathcal{X}_0\right)
 =
n\left( \prob \left( X \in \mathcal{X}_0\right)
 +
 o_{\prob}(1)
 \right),
$$
 by the Law of large numbers, Lemma \ref{STR} then gives
\begin{eqnarray*}
\lefteqn{ \widehat{f}_{1n}(\epsilon)-\widetilde{f}_{1n}(\epsilon)
 =
O_{\prob}\left(\frac{1}{nb_1^2}\right)S_n + O_{\prob}
\left(\frac{1}{nb_1^3}\right) T_n
 +
O_{\prob} \left(\frac{1}{nb_1^4}\right) R_n }
\\
&=&
 O_{\prob}
 \left[
 b_0^2
 \left(1+\frac{1}{(nb_1^3)^{1/2}}\right)
 +
 \left(
 \frac{1}{n}
 +
 \frac{1}{n^2b_0^db_1^3}
 \right)^{1/2}
 \right]
 \\
 &&
 +\;
 O_{\prob}
 \left[
 1
 +
 \frac{1}{(nb_1^5)^{1/2}}
 +
 \left(\frac{b_0^d}{b_1^3}\right)^{1/2}
 \right]
 \left(
 b_0^4
 +
 \frac{1}{nb_0^d}
 \right)
 +
 O_{\prob}
 \left[
 \frac{1}{b_1}
 +
 \left(\frac{b_0^d}{b_1^7}\right)^{1/2}
 \right]
 \left(b_0^4+\frac{1}{nb_0^d}\right)^{3/2}.
 \end{eqnarray*}
This yields (\ref{fnchapfn1}), since  under $(A_9)$ and
$(A_{10})$,  we have $b_0\rightarrow 0$,
$nb_0^{d+2}\rightarrow\infty$ and  $nb_1^3\rightarrow\infty$, so
that
\begin{eqnarray*}
b_0^2
 \left(1+\frac{1}{(nb_1^3)^{1/2}}\right)
&\asymp&
 O\left(b_0^2\right),
\quad
 \left(
 b_0^4
 +
 \frac{1}{nb_0^d}
 \right)
 =
 O\left(b_0^2\right),
\\
\left[
 1
 +
 \frac{1}{(nb_1^5)^{1/2}}
 +
 \left(\frac{b_0^d}{b_1^3}\right)^{1/2}
 \right]
 \left(
 b_0^4
 +
 \frac{1}{nb_0^d}
 \right)
 &=&
 O\left(b_0^2\right)
 +
 \left[
 \frac{1}{(nb_1^5)^{1/2}}
 +
 \left(\frac{b_0^d}{b_1^3}\right)^{1/2}
 \right]
 \left(
 b_0^4
 +
 \frac{1}{nb_0^d}
 \right).
\end{eqnarray*}

\vskip 0.1cm
 For (\ref{fnchapfn2}), note that
\begin{eqnarray}
\esp_n \left[
\left(\widetilde{f}_{1n}(\epsilon)-f(\epsilon)\right)^2 \right] =
\Var_n\left(\widetilde{f}_{1n}(\epsilon)\right) + \biggl( \esp_n
\left[ \widetilde{f}_{1n}(\epsilon)\right] - f(\epsilon)
\biggr)^2, \label{Quadbias}
\end{eqnarray}
with, using $(A_4)$,
\begin{eqnarray*}
\Var_n\left(\widetilde{f}_{1n}(\epsilon)\right)
 =
\frac{1}{\left(b_1\sum_{i=1}^n
\mathds{1}\left(X_i\in\mathcal{X}_0\right)\right)^2} \sum_{i=1}^n
\mathds{1} \left(X_i\in\mathcal{X}_0\right) \Var \left[
 K_1\left(\frac{\varepsilon-\epsilon}{b_1}\right)
 \right].
\end{eqnarray*}
Therefore, since the Cauchy-Schwarz inequality gives
\begin{eqnarray*}
\Var \left[
 K_1\left(\frac{\varepsilon-\epsilon}{b_1}\right)
 \right]
\leq
 \esp
 \left[
 K_1^2\left(\frac{\varepsilon-\epsilon}{b_1}\right)
 \right]
 \leq
 b_1\int
 K_1^2(v)f(\epsilon+b_1v)
 dv,
\end{eqnarray*}
this bound and the equality above yield, under $(A_5)$ and
$(A_8)$,
\begin{eqnarray}
 \Var_n\left(\widetilde{f}_{1n}(\epsilon)\right)
 \leq
 \frac{C}{b_1\sum_{i=1}^n
 \mathds{1}\left(X_i\in\mathcal{X}_0\right)}
 =
 O_{\prob}\left(\frac{1}{nb_1}\right).
 \label{Quadbias1}
\end{eqnarray}
For the second term in (\ref{Quadbias}),  we have
\begin{eqnarray}
\esp_n\left[\widetilde{f}_{1n}(\epsilon)\right] =
\frac{1}{b_1\sum_{i=1}^n
 \mathds{1}\left(X_i\in\mathcal{X}_0\right)}
\sum_{i=1}^n \mathds{1} \left(X_i\in\mathcal{X}_0\right) \esp
\left[
 K_1\left(\frac{\varepsilon-\epsilon}{b_1}\right)
 \right].
 \label{Quadbias2}
\end{eqnarray}
By $(A_8)$,  $K_1$  is symmetric, has  a compact support, with
 $\int\! v K_1(v)=0$  and $\int \! K_1(v)dv=1$.
 Therefore, since under $(A_5)$ $f$ has bounded continuous second
order derivatives, this yields for some
$\theta=\theta(\epsilon,b_1v)$,
\begin{eqnarray*}
\lefteqn{
 \esp
 \left[
 K_1\left(\frac{\varepsilon-\epsilon}{b_1}\right)
 \right]
 =
  b_1
\int
 K_1(v)
 f(\epsilon+b_1v)
 dv
 }
\\
&=&
 b_1
 \int
 K_1(v)
\left[
 f(\epsilon)
+
 b_1 v f^{(1)}(\epsilon)
 +
 \frac{b_1^2v^2}{2}
 f^{(2)}(\epsilon+\theta b_1v)
 \right]
 dv
 \\
  &=&
 b_1f(\epsilon)
 +
 \frac{b_1^3}{2}
 \int v^2 K_1(v)
 f^{(2)}(\epsilon+\theta b_1v)
 dv.
\end{eqnarray*}
Hence  this equality and (\ref{Quadbias2}) give
$$
\esp_n \left[\widetilde{f}_{1n}(\epsilon)\right]
 =
 f(\epsilon)
 +
 \frac{b_1^2}{2}
 \int v^2 K_1(v)
 f^{(2)}(\epsilon+\theta b_1v)
 dv,
$$
so that
$$
\biggl( \esp_n \left[\widetilde{f}_{1n}(\epsilon)\right] -
f(\epsilon) \biggr)^2
 =
 O_{\prob}\left(b_1^4\right).
$$
Combining this result with  (\ref{Quadbias1}) and
(\ref{Quadbias}), we obtain, by the Tchebychev inequality,
$$
\widetilde{f}_{1n}(\epsilon)-f(\epsilon) = O_{\prob}
\left(b_1^4+\frac{1}{nb_1}\right)^{1/2}.
$$
This proves (\ref{fnchapfn2}), and then achieves the proof of the
theorem.\eop

\subsection*{Proof of Theorem \ref{Optimbandw1}}
Recall that
\begin{eqnarray*}
 R_n(b_0, b_1)
=
 b_0^4
 +
 \left[
 \frac{1}{(nb_1^5)^{1/2}}
 +
 \left(\frac{b_0^d}{b_1^3}\right)^{1/2}
 \right]^2
 \left(
 b_0^4
 +
 \frac{1}{nb_0^d}
 \right)^2
 +
 \left[
 \frac{1}{b_1}
 +
 \left(\frac{b_0^d}{b_1^7}\right)^{1/2}
 \right]^2
 \left(b_0^4+\frac{1}{nb_0^d}\right)^3,
\end{eqnarray*}
and note that
$$
\left(\frac{1}{n^2b_1^3}\right)^{\frac{1}{d+4}} = \max
\left\lbrace
 \left(\frac{1}{n^2b_1^3}\right)^{\frac{1}{d+4}}
, \left(\frac{1}{n^3b_1^7}\right)^{\frac{1}{2d+4}} \right\rbrace
$$
if and only if $n^{4-d}b_1^{d+16}\rightarrow\infty$. To find the
order of $b_0^*$, we shall deal with the cases
$nb_0^{d+4}\rightarrow\infty$ and
 $nb_0^{d+4}=O(1)$.
\vskip 0.1cm\noindent First assume that
$nb_0^{d+4}\rightarrow\infty$. More precisely, we  suppose that
$b_0$ is in $\left[(u_n/n)^{1/(d+4)},+\infty\right)$, where $u_n
\rightarrow\infty$. Since $1/(nb_0^d) = O(b_0^4)$ for all these
$b_0$, we have
$$
\left(b_0^4+\frac{1}{nb_0^d}\right)^2 \asymp \left(b_0^4\right)^2,
\quad \left(b_0^4+\frac{1}{nb_0^d}\right)^3 \asymp
 \left(b_0^4\right)^3.
$$
Hence the order of $b_0^*$ is computed  by minimizing the function
\begin{eqnarray*}
 b_0\rightarrow
  b_0^4
 +
 \left[
 \frac{1}{(nb_1^5)^{1/2}}
 +
 \left(\frac{b_0^d}{b_1^3}\right)^{1/2}
 \right]^2
 \left( b_0^4\right)^2
 +
 \left[
 \frac{1}{b_1}
 +
 \left(\frac{b_0^d}{b_1^7}\right)^{1/2}
 \right]^2
 \left(b_0^4\right)^3.
\end{eqnarray*}
Since this function is increasing with $b_0$, the minimum of $R_n
(\cdot,b_1)$ is achieved for $b_{0*}=(u_n/n)^{1/(d+4)}$. We shall
prove later on that this choice of $b_{0*}$ is irrelevant compared
to the one arising when $nb_0^{d+4}=O(1)$.

\vskip 0.1cm Consider now the case $nb_0^{d+4}=O(1)$ i.e
$b_0^4=O\left(1/(nb_0^d)\right)$. This  gives
\begin{eqnarray*}
 \left[
 \frac{1}{(nb_1^5)^{1/2}}
 +
 \left(\frac{b_0^d}{b_1^3}\right)^{1/2}
 \right]^2
 \left(
 b_0^4
 +
 \frac{1}{nb_0^d}
 \right)^2
 &\asymp&
  \left(
 \frac{1}{nb_1^5}
 +
 \frac{b_0^d}{b_1^3}
 \right)
 \left(
 \frac{1}{n^2b_0^{2d}}
 \right),
 \\
 \left[
 \frac{1}{b_1}
 +
 \left(\frac{b_0^d}{b_1^7}\right)^{1/2}
 \right]^2
 \left(b_0^4+\frac{1}{nb_0^d}\right)^3
 &\asymp&
 \left(
 \frac{1}{b_1^2}
 +
 \frac{b_0^d}{b_1^7}
 \right)
 \left(\frac{1}{n^3b_0^{3d}}\right).
 \end{eqnarray*}
Moreover if $nb_0^db_1^4\rightarrow\infty$, we have, since
$nb_0^{2d}\rightarrow\infty$ under $(A_9)$,
$$
 \left(
 \frac{1}{nb_1^5}
 +
 \frac{b_0^d}{b_1^3}
 \right)
 \left(
 \frac{1}{n^2b_0^{2d}}
 \right)
 \asymp
 \frac{b_0^d}{b_1^3}
 \left(
 \frac{1}{n^2b_0^{2d}}
 \right),
 \quad
 \left(
 \frac{1}{b_1^2}
 +
 \frac{b_0^d}{b_1^7}
 \right)
 \left(\frac{1}{n^3b_0^{3d}}\right)
 =
 O\left(
 \frac{b_0^d}{b_1^3}
 \right)
 \left(
 \frac{1}{n^2b_0^{2d}}
 \right).
 $$
Hence the order of $b_0^*$ is obtained by finding the minimum of
the function $b_0^4+\left(1/n^2b_0^db_1^3\right)$. The
 minimization of this function  gives  a solution $b_0$ such that
$$
b_0
\asymp
\left(\frac{1}{n^2b_1^3}\right)^{\frac{1}{d+4}},
\quad
R_n(b_0,b_1)
\asymp
\left(\frac{1}{n^2b_1^3}\right)^{\frac{4}{d+4}}.
$$
This value  satisfies the constraints $nb_0^{d+4}=O(1)$ and
$nb_0^db_1^4\rightarrow\infty$ when
$n^{4-d}b_1^{d+16}\rightarrow\infty$.

\vskip 0.1cm\noindent If now $nb_0^{d+4}=O(1)$ but
$nb_0^db_1^4=O(1)$, we have, since $nb_0^{2d}\rightarrow\infty$,
\begin{eqnarray*}
\frac{1}{nb_1^5}
 \left(
\frac{1}{n^2b_0^{2d}} \right)
  =
 O\left(
 \frac{b_0^d}{b_1^7}
 \right)
 \left(\frac{1}{n^3b_0^{3d}}\right),
 \quad
 \frac{1}{b_1^2}
 \left(\frac{1}{n^3b_0^{3d}}\right)
 =
O\left(
 \frac{b_0^d}{b_1^3}
 \right)
 \left(
\frac{1}{n^2b_0^{2d}} \right)
 =
 O\left(
 \frac{b_0^d}{b_1^7}
 \right)
 \left(\frac{1}{n^3b_0^{3d}}\right).
\end{eqnarray*}
 In this case, $b_0^*$ is obtained by minimizing the function
$b_0^4+\left(1/n^3b_0^{2d}b_1^7\right)$, for which the solution
 $b_0$ verifies
$$
b_0
\asymp
\left(\frac{1}{n^3b_1^7}\right)^{\frac{1}{2d+4}},
\quad
R_n(b_0,b_1)
\asymp
\left(\frac{1}{n^3b_1^7}\right)^{\frac{4}{2d+4}}.
$$
This solution fulfills  the constraint $nb_0^db_1^4=O(1)$ when
$n^{4-d}b_1^{d+16}=O(1)$. Hence we can conclude that for
$b_0^4=O\left(1/(nb_0^d)\right)$, the bandwidth $b_0^*$ satisfies
$$
b_0^*
\asymp
\max
\left\lbrace
\left(\frac{1}{n^2b_1^3}\right)^{\frac{1}{d+4}} ,
\left(\frac{1}{n^3b_1^7}\right)^{\frac{1}{2d+4}}
\right\rbrace,
$$
which leads to
$$
R_n\left(b_0^*, b_1\right)
\asymp
\max \left\lbrace
\left(\frac{1}{n^2b_1^3}\right)^{\frac{4}{d+4}} ,
\left(\frac{1}{n^3b_1^7}\right)^{\frac{4}{2d+4}}
 \right\rbrace.
$$
We need now to compare the solution $b_0^*$ to the candidate
 $b_{0*}=(u_n/n)^{1/(d+4)}$ obtained when $nb_0^{d+4}\rightarrow\infty$.
  For this, we must  do a comparison between the
 orders of $R_n(b_0^*,b_1)$  and $R_n(b_{0*},b_1)$. Since
 $R_n(b_0,b_1)\geq b_0^4$, we have
 $R_n(b_{0*},b_1)\geq(u_n/n)^{4/(d+4)}$, so that, for $n$ large
 enough,
\begin{eqnarray*}
\frac{R_n(b_0^*, b_1)} {R_n(b_{0*},b_1)}
&\leq&
 C\left[
\left(\frac{1}{n^2b_1^3}\right)^{\frac{1}{d+4}}
 +
\left(\frac{1}{n^3b_1^7}\right)^{\frac{4}{2d+4}}
\right]
\left(\frac{n}{u_n}\right)^{\frac{4}{d+4}}
\\
&=&
o(1)
+
O\left(\frac{1}{u_n}\right)^{\frac{4}{d+4}}
\left(
\frac{1}{nb_1^{\frac{7(d+4)}{d+8}}}
\right)^{\frac{4(d+8)}{(2d+4)(d+4)}} = o(1),
\end{eqnarray*}
using $u_n\rightarrow\infty$ and that
$n^{(d+8)}b_1^{7(d+4)}\rightarrow\infty$ by $(A_{10})$. This shows
that $R_n(b_0^*,b_1)\leq R_n(b_{0*},b_1)$ for $n$ large enough.
 Hence the Theorem is proved, since $b_0^*$ is the best candidate for
 the minimization of $R_n(\cdot, b_1)$. \eop

\subsection*{Proof of Theorem \ref{Optimbandw2}}
Recall that  Theorem \ref{Optimbandw1} gives
\begin{eqnarray*}
AMSE(b_1)+R_n(b_0^*, b_1)
\asymp r_1(b_1)+r_2(b_1)+r_3(b_1)
=
F(b_1),
\end{eqnarray*}
where
\begin{eqnarray*}
r_1(h) &=& h^4 + \frac{1}{nh}, \quad \arg\min r_1(h) \asymp
n^{-1/5}=h_1^*, \quad \min r_1(h) \asymp (h_1^*)^4=n^{-4/5},
\\
r_2(h) &=& h^4 + \frac{1}{n^{\frac{8}{d+4}} h^{\frac{12}{d+4}}},
\quad \arg\min r_2(h) \asymp n^{-\frac{2}{d+7}} = h_2^*, \quad
\min r_3(h) \asymp(h_2^*)^4 = n^{-\frac{8}{d+7}},
\\
r_3(h) &=& h^4 + \frac{1}{n^{\frac{12}{2d+4}}
h^{\frac{28}{2d+4}}}, \quad \arg\min r_3(h) \asymp
n^{-\frac{3}{2d+11}} = h_3^*,
\quad \min r_3(h) \asymp(h_3^*)^4 =
n^{-\frac{12}{2d+11}}.
\end{eqnarray*}
Each $r_j(h)$ decreases on $\left[0,\arg\min r_j(h)\right]$ and
increases on $\left(\arg\min r_j(h),\infty\right)$ and that
$r_j(h)\asymp h^4$ on $\left(\arg\min r_j(h),\infty\right)$.
Moreover $\min r_2(h)=o\left(r_3(h)\right)$ and
$h_2^*=o\left(h_3^*\right)$ for all possible dimension $d$, so
that $\min\{r_2(h)+r_3(h)\}\asymp(h_3^*)^4=n^{-\frac{12}{2d+11}}$
and $\arg\min\{r_2(h)+r_3(h)\}\asymp h_3^*=n^{-\frac{3}{2d+11}}$.

\vskip 0.1cm Observe now that $\min\{r_2(h)+r_3(h)\}=O\left(\min
r_1(h)\right)$ is equivalent to
$n^{-\frac{12}{2d+11}}=O\left(n^{-4/5}\right)$ which holds if and
only if $d\leq 2$. Hence assume that $d\leq 2$. Since
$n^{-\frac{12}{2d+11}}=O\left(n^{-4/5}\right)$ also gives
$\arg\min\{r_2(h)+r_3(h)\}\asymp h_3^*=O\left(h_1^*\right)$, we
have
$$
\min F(b_1) \asymp n^{-4/5}
\;\;{\rm and}\; \arg\min F(b_1) \asymp
n^{-1/5}.
$$
The case $d>2$ is symmetric with
$$
\min F(b_1) \asymp n^{-\frac{12}{2d+11}}
\;\;{\rm and}\; \arg\min
F(b_1) \asymp n^{-\frac{3}{2d+11}}.
$$
This ends the proof of the Theorem. \eop

\subsection*{Proof of Theorem \ref{normalite}}
Observe that the Tchebychev inequality gives
$$
\sum_{i=1}^n \mathds{1} \left(X_i\in\mathcal{X}_0\right) =
n\prob\left(X\in\mathcal{X}_0\right) \left[ 1 + O_{\prob}
\left(\frac{1}{\sqrt{n}}\right) \right],
$$
so that
$$
\widetilde{f}_{1n}(\epsilon) = \left[ 1 + O_{\prob}
\left(\frac{1}{\sqrt{n}}\right) \right] f_n(\epsilon),
$$
where
$$
f_n(\epsilon)
 =
\frac{1}{nb_1\prob\left(X\in\mathcal{X}_0\right)} \sum_{i=1}^n
\mathds{1} \left(X_i\in\mathcal{X}_0\right)
 K_1\left(\frac{\varepsilon_i-\epsilon}{b_1}\right).
$$
Therefore
\begin{eqnarray}
\widehat{f}_{1n}(\epsilon)-\esp f_n(\epsilon) = \left(
f_n(\epsilon)-\esp f_n(\epsilon) \right) + \left(
\widehat{f}_{1n}(\epsilon)-\widetilde{f}_{1n}(\epsilon) \right) +
O_{\prob} \left(\frac{1}{\sqrt{n}}\right)
 f_n(\epsilon).
\label{develop}
\end{eqnarray}
Let now $f_{in}(\epsilon)$ be as in Lemma \ref{Momfin1}, and note
that $f_n(\epsilon) =(1/n)\sum_{i=1}^nf_{in}(\epsilon)$. The
second and the third claims in Lemma \ref{Momfin1} yield, since
$b_1$ goes to $0$ under $(A_{10})$,
\begin{eqnarray*}
\frac { \sum_{i=1}^n \esp \left| f_{in}(\epsilon) - \esp
 f_{in}(\epsilon)
\right|^3 } { \left( \sum_{i=1}^n
 \Var f_{in}(\epsilon)
\right)^3 } \leq \frac{ \frac{C nf(\epsilon)}
{\prob\left(X\in\mathcal{X}_0\right)^2b_1^2}
 \displaystyle{\int}
\left| K_1(v) \right|^3 dv + o\left( \frac{n}{b_1^2} \right) } {
\left( \frac{nf(\epsilon)}{\prob\left(
X\in\mathcal{X}_0\right)b_1} \displaystyle{\int}K_1^2(v) dv +
 o\left(\frac{n}{b_1}\right)
\right)^3
 }
 =O(b_1)= o(1).
\end{eqnarray*}
 Hence the Lyapounov Central Limit Theorem gives, since $nb_1$
 diverges under $(A_{10})$,
$$
\frac{f_n(\epsilon)-\esp f_{n}(\epsilon)} {\sqrt{\Var
f_{n}(\epsilon)}} =
 \frac{f_n(\epsilon)-\esp f_{n}(\epsilon)}
 {\sqrt{\frac{\Var f_{in}(\epsilon)}{n}}}
 \stackrel{d}{\rightarrow}
 \mathcal{N}\left(0,1\right),
$$
which yields, using the second equality in Lemma \ref{Momfin1},
\begin{eqnarray}
\sqrt{nb_1} \left( f_n(\epsilon) - \esp f_{n}(\epsilon) \right)
\stackrel{d}{\rightarrow} \mathcal{N} \left( 0, \frac{f(\epsilon)}
{\prob\left(X\in\mathcal{X}_0\right)} \int K_1^2(v) dv \right).
\label{develop1}
\end{eqnarray}
Moreover, note that for $nb_0^db_1^3\rightarrow\infty$ and
$nb_0^{2d}\rightarrow\infty$,
$$
 \frac{1}{nb_1^5}
 \left(
 \frac{1}{nb_0^d}
 \right)^2
 +
 \left(
 \frac{1}{b_1^2}
 +
 \frac{b_0^d}{b_1^7}
 \right)^2
 \left(\frac{1}{nb_0^d}\right)^3
 =
 O\left(
\frac{1}{n^2b_0^db_1^3}\right).
$$
Therefore, since by Assumptions $(\rm{A}_{11})$ and $(A_9)$, we
have $b_0^4=O\left(1/(nb_0^d)\right)$,
$nb_0^db_1^3\rightarrow\infty$ and that
$nb_0^{2d}\rightarrow\infty$, the equality above and
(\ref{fnchapfn1}) then give
\begin{eqnarray*}
\widehat{f}_{1n}(\epsilon)-\widetilde{f}_{1n}(\epsilon)
 &\asymp&
 O_{\prob}
 \left[
 b_0^4
 +
 \frac{1}{n}
 +
 \frac{1}{n^2b_0^db_1^3}
 +
 \left(
 \frac{1}{nb_1^5}
 +
 \frac{b_0^d}{b_1^3}
 \right)
 \left(
 \frac{1}{nb_0^d}
 \right)^2
 +
 \left(
 \frac{1}{b_1^2}
 +
 \frac{b_0^d}{b_1^7}
 \right)
 \left(\frac{1}{nb_0^d}\right)^3
 \right]^{1/2}
 \\
 &\asymp&
 O_{\prob}
 \left(
 b_0^4
 +
 \frac{1}{n}
 +
 \frac{1}{n^2b_0^db_1^3}
\right)^{1/2}.
\end{eqnarray*}
 Hence for $b_1$ going to $0$, we have
$$
\sqrt{nb_1} \left(
\widehat{f}_{1n}(\epsilon)-\widetilde{f}_{1n}(\epsilon) \right)
 =
O_{\prob}
 \left[
 nb_1
 \left(
 b_0^4
 +
 \frac{1}{n}
 +
 \frac{1}{n^2b_0^db_1^3}
 \right)
\right]^{1/2}
 =
 o_{\prob}(1),
$$
since $nb_0^4b_1=o(1)$ and that $nb_0^db_1^2\rightarrow\infty$
under Assumption $(\rm{A}_{11})$. Combining the above result with
(\ref{develop1}) and (\ref{develop}), we obtain
$$
\sqrt{nb_1}
 \left(
 \widehat{f}_{1n}(\epsilon)
 -
 \esp f_n(\epsilon)
\right)
\stackrel{d}{\rightarrow}
\mathcal{N}
 \left(
 0,
\frac{f(\epsilon)}
{\prob\left(X\in\mathcal{X}_0\right)}
\int
K_1^2(v) dv
 \right).
$$
This ends the proof the Theorem, since the first result of Lemma
\ref{Momfin1} gives
$$
\esp f_n(\epsilon)
=
\esp f_{1n}(\epsilon)
=
f(\epsilon)
+
\frac{b_1^2}{2}f^{(2)}(\epsilon)
\int v^2K_1(v) dv +
o\left(b_1^2\right)
:=
\overline{f}_{1n}(\epsilon).
\eop
$$

\newpage
 \setcounter{equation}{0} \setcounter{subsection}{0}
\setcounter{lemma}{0}
\renewcommand{\theequation}{A.\arabic{equation}}
\renewcommand{\thesubsection}{A.\arabic{subsection}}
\begin{center}
\section*{Appendix A:  Proof of the  intermediate results}
\end{center}

\subsection*{Proof of Lemma \ref{Estig}}
First note that  by $(A_7)$, we have
 $
 \int\!
z K_0(z) dz
 =0
 $
 and
 $
 \int\!
K_0(z) dz
 =1$.  Therefore, since $K_0$ is continuous and has a compact support,
$(A_1)$, $(A_2)$ and a second-order Taylor expansion,
 yield,
for  $b_0$ small enough and any $x$ in $\mathcal{X}_0$,
\begin{eqnarray*}
 \lefteqn{
 \left|
 \overline{g}_{n}(x)-g(x)
 \right|
  =
 \left|
 \frac{1}{b_0^d}
 \int
 K_0
 \left(\frac{z-x}{b_0}\right)
 g(z)
 dz
 -
 g(x)
 \right|
  =
 \left|
\int K_0(z) \left[g(x+b_0z)-g(x)\right]
 dz
 \right|
 }
 &&
\\
&=&
 \left|
 \int K_0(z)
 \left[
 b_0 g^{(1)}(x) z
 +
 \frac{b_0^2}{2}
 z g^{(2)} (x + \theta b_0 z)z^{\top}
\right] dz
 \right|,
 \;\theta = \theta (x,b_0 z)\in [0,1]
\\
&=& \left|
 b_0 g^{(1)}(x)
 \int
 z K_0(z) dz
 +
 \frac{b_0^2}{2}
  \int
z g^{(2)}(x + \theta b_0 z) z^{\top} K_0(z) dz \right|
\\
& = & \frac{b_0^2}{2}
 \left|
 \int z g^{(2)}(x+\theta b_0z)
 z^{\top} K_0(z) dz
 \right|
 \leq C b_0^2,
\end{eqnarray*}
so that
$$
 \sup_{x\in\mathcal{X}_0}
 \left|
 \overline{g}_n (x) - g(x)
 \right|
 =
 O\left(b_0^2\right).
$$
This gives the first equality of the lemma. To prove the two last
equalities in the Lemma, note that it is
 sufficient to show that
$$
 \sup_{x\in\mathcal{X}_0}
 \left|
 \widehat{g}_{n}(x)
  -
 \overline{g}_n(x)
 \right|
 =
 O_{\prob}
 \left(
  \frac{\ln n}{nb_0^d}
 \right)^{1/2},
$$
 since $\bar{g}_n (x)$ is asymptotically bounded away from
$0$ over
 $\mathcal{X}_0$
 and that
$|\overline{g}_n (x) - g(x) | = O (b_0^2) $ uniformly for $x$ in
$\mathcal{X}_0$. This follows from Theorem 1 in  Einmahl and Mason
(2005). \eop

\subsection*{Proof of Lemma \ref{Estim}}

For the first equality in the lemma, set
$$
\widehat{r}_n(x)
 =
 \frac{1}{nb_0^d}
 \sum_{j=1}^n
 Y_j K_0
  \left(
 \frac{X_j-x}{b_0}
 \right),
 \quad
 \overline{r}_n(x)
  =
  \esp
 \left[
\widehat{r}_n(x)
 \right] \;,
$$
and observe that
\begin{equation}
\sup_{x\in\mathcal{X}_0}
 \left|\widehat{m}_n(x)-m(x)\right|
 \leq
 \sup_{x\in\mathcal{X}_0}
 \left|
 \widehat{m}_n(x)
 -
 \frac{\overline{r}_n(x)}{\overline{g}_n(x)}
 \right|
 +
 \sup_{x\in\mathcal{X}_0}
 \frac{1}{\left|\overline{g}_n(x)\right|}
 \left|
 \overline{r}_n(x)-\overline{g}_n(x)m(x)
 \right|.
 \label{mchd}
\end{equation}
Consider the first term of (\ref{mchd}). Note that
$\esp^{1/4}\left[Y^4|X=x\right] \leq |m(x)|
+\esp^{1/4}\left[\varepsilon^4\right]$. The compactness of
$\mathcal{X}$
 from $(A_1)$, the continuity of $m(\cdot)$ from $(A_3)$
and $(A_4)$ then give that $\esp\left[Y^4|X=x\right]<\infty$
uniformly for $x\in\mathcal{X}_0$. Hence under $(A_9)$, Theorem 2
in Einmahl and Mason (2005) gives
$$
\sup_{x\in\mathcal{X}_0} \left|
 \widehat{m}_n(x)
 -
\frac{\overline{r}_n(x)}{\overline{g}_n(x)}
 \right|
 =
 O_{\prob}
 \left(\frac{\ln n}{nb_0^d}\right)^{1/2}.
 $$
 For the second term in (\ref{mchd}), a
second-order Taylor expansion
 gives, as in the proof of Lemma \ref{Estig},
$$
 \sup_{x\in\mathcal{X}_0}
 \left|
 \overline{r}_n(x)
 -
 \overline{g}_n(x)m(x)
\right|
 =
 O (b_0^2).
 $$
  This gives the result of lemma since
 Lemma \ref{Estig}
 implies that $\overline{g}_n (x)$ is bounded away from $0$ over
 $\mathcal{X}_0$ uniformly in $x$ and for $b_0$ small enough.
 \eop

\subsection*{Proof of Lemma \ref{Reste}}
Note that under $(A_8)$, the Taylor expansion with integral
remainder gives, for any $x\in\mathcal{X}_0$ and any integer
$i\in[1,n]$,
\begin{eqnarray*}
 K_1\left(
 \frac{Y_i-\widehat{m}_n(x)-\epsilon}{h_1}
 \right)
=
 K_1\left(\frac{Y_i-m(x)-\epsilon}{h_1}\right)
 -
 \frac{1}{h_1}
 \left(\widehat{m}_n(x)-m(x)\right)
 \int_{0}^1
 K_1^{(1)}
 \left(\frac{Y_i-\theta_n(x, t)}{h_1}
 \right)
  dt,
\end{eqnarray*}
where
 $\theta_n(x,t)=m(x)+\epsilon+t\left(\widehat{m}_n(x)-m(x)\right)$.
 Therefore
 \begin{eqnarray}
 \nonumber
 \widetilde{f}_n(\epsilon|x)
 =
 f_{n}(\epsilon|x)
 -
 \frac{\widehat{m}_n(x)-m(x)}{\widetilde{g}_n(x)}
 \left[
 \frac{1}{nh_0^dh_1^2}
 \sum_{i=1}^n
 K_0\left(\frac{X_i-x}{h_0}\right)
 \int_0^1
 K_1^{(1)}
 \left(\frac{Y_i-\theta_n(x,t)}{h_1}\right)
 dt
 \right].
 \\
 \label{fnt}
 \end{eqnarray}
Now, observe that if $X_i = z$ and $y\in\Rit$, the change of
variable $e=y-m(z)+h_1v$ gives,  under  $(A_1)-(A_5)$ and $(A_7)$,
\begin{eqnarray*}
\lefteqn{ \esp_n
 \left|
 K_1^{(1)}
 \left(\frac{Y_i-y}{h_1}\right)
 \right|
 =
 \esp
 \left|
 K_1^{(1)}
 \left(\frac{\varepsilon_i+m(z)-y}{h_1}\right)
\right|
 }
 \\
 &=&
 \int
 \left|
 K_1^{(1)}
 \left(\frac{e+m(z)-y}{h_1}\right)
\right|
 f(e) de
 \\
 &=&
 h_1
 \int
|K_1^{(1)}(v)|
 f\left((y-m(z)+h_1v\right))
 dv
\leq Ch_1.
\end{eqnarray*}
Hence
$$
\sup_{1\leq i\leq n}
 \int_0^1
 \esp_n
 \left|
 K_1^{(1)}
 \left(\frac{Y_i-\theta_n(x,t)}{h_1}\right)
 \right|
 dt
\leq Ch_1.
$$
With the help of this result and Lemma \ref{Estig}, we have
\begin{eqnarray*}
\lefteqn{ \esp_n \left| \frac{1}{nh_0^dh_1}
 \sum_{i=1}^n
 K_0\left(\frac{X_i-x}{h_0}\right)
 \int_0^1
 K_1^{(1)}
 \left(\frac{Y_i-\theta_n(x,t)}{h_1}\right)
 dt
\right| }
\\
&\leq& \frac{1}{nh_0^dh_1}
 \sum_{i=1}^n
 \left|K_0\left(\frac{X_i-x}{h_0}\right)\right|
 \times
\sup_{1\leq i\leq n} \int_0^1
 \esp_n
 \left|
 K_1^{(1)}
 \left(\frac{Y_i-\theta_n(x,t)}{h_1}\right)
 \right|
 dt
 \\
 &\leq&
\frac{C}{nh_0^d}
 \sum_{i=1}^n
 \left|
 K_0\left(\frac{X_i-x}{h_0}\right)
 \right|
 =
 O_{\prob}(1),
\end{eqnarray*}
so that
$$
\frac{1}{nh_0^dh_1^2}
 \sum_{i=1}^n
 K_0\left(\frac{X_i-x}{h_0}\right)
 \int_0^1
 K_1^{(1)}
 \left(\frac{Y_i-\theta_n(x,t)}{h_1}\right)
 dt
=O_{\prob}\left(\frac{1}{h_1}\right).
$$
Hence from (\ref{fnt}),  (\ref{gtilde}),  Lemma \ref{Estim} and
Assumption $(\rm\bf{A}_0)$, we deduce
$$
 \widetilde{f}_n(\epsilon|x)
 =
 f_{n}(\epsilon|x)
 +
 O_{\prob}
 \left(\frac{1}{h_1}\right)
 \left(b_0^4+\frac{\ln n}{nb_0^d}\right)^{1/2}
 =
 f_{n}(\epsilon|x)
 +
 o\left(\frac{1}{nh_0^dh_1}\right)^{1/2}.
 \eop
$$

\subsection*{Proof of Lemma \ref{Momvarphi} and Lemma \ref{Momfin1}}
We just give the proof of Lemma \ref{Momvarphi}, the proof of
Lemma \ref{Momfin1} being very similar.
 For the first equality of  Lemma \ref{Momvarphi}, note that
\begin{eqnarray*}
\esp
\left[ \widetilde{\varphi}_{in}(x,y) \right]
&=&
\nonumber
\frac{1} {h_0^dh_1}
\int \int K_0 \left( \frac{x_1-x}{h_0} \right)
K_1 \left( \frac{y_1-y}{h_1} \right)
\varphi (x_1, y_1) dx_1 dy_1
\\
&=&
\int \int
K_0(z) K_1(v)
\varphi\left(x+h_0z, y+h_1 v\right)
dz dv.
 \label{espphi}
\end{eqnarray*}
A second-order Taylor expansion gives under $(A_6)$, for $z$ in
the support of $K_0$, $v$ in the support of $K_1$, and $h_0$,
$h_1$ small enough,
\begin{eqnarray*}
\lefteqn{ \varphi\left(x+h_0z, y+h_1 v\right) - \varphi (x, y) }
\\
&=& h_0 \frac{\partial \varphi (x, y)}{\partial x} z^{\top} + h_1
\frac{\partial \varphi (x, y)}{\partial y} v
\\
&& + \frac{h_0^2}{2} z \frac{\partial^2 \varphi (x+\theta h_0 z,
y+\theta h_1 v)}{\partial^2
 x} z^{\top}
+ h_1 h_0 v \frac{\partial^2 \varphi (x+\theta h_0 z, y+\theta h_1
v)}{\partial x\partial y} z^{\top}
\\
&& + \frac{h_1^2}{2} \frac{\partial^2 \varphi (x+\theta h_0 z,
y+\theta h_1 v)}{\partial^2
 y}
v^2,
\end{eqnarray*}
for some $\theta = \theta (x,y,h_0 z, h_1 v)$ in $[0,1]$. This
gives, since $\int\! K_0 (z) dz = \int\! K_1 (v) dv =1$, $\int\!z
K_0(z) dz$ and $\int\! v K_1 (v) dv$ vanish under $(A_7)-(A_8)$,
and by the Lebesgue Dominated Convergence Theorem,
\begin{eqnarray*}
\lefteqn{ \esp \left[ \widetilde{\varphi}_{in}(x,y) \right] -
\varphi(x, y) - \frac{h_0^2}{2} \frac{\partial^2 \varphi(x,
y)}{\partial^2 x} \int z K_0(z)z^{\top} dz - \frac{h_1^2}{2}
\frac{\partial^2 \varphi (x, y)}{\partial^2 y} \int v^2 K_1 (v) dv
} &&
\\
& = & \frac{h_0^2}{2} \int \int z \left( \frac{\partial^2 \varphi
(x+\theta h_0 z, y+\theta h_1 v)}{\partial^2
 x}
- \frac{\partial^2 \varphi (x, y)}{\partial^2 x} \right) z^{\top}
K_0 (z) K_1 (v) dz dv
\\
&& + h_1 h_0 \int \int v \left( \frac{\partial^2 \varphi (x+\theta
h_0 z, y+\theta h_1 v)} {\partial x \partial y} - \frac{\partial^2
\varphi (x, y)}{\partial x \partial y} \right) z^{\top} K_0 (z)
K_1 (v) dz dv
\\
&& + \frac{h_1^2}{2} \int \int \left( \frac{\partial^2 \varphi
(x+\theta h_0 z, y+\theta h_1 v)} {\partial^2 y} -
\frac{\partial^2 \varphi (x, y)}{\partial^2 y} \right) v^2 K_0 (z)
K_1 (v) dz dv
\\
&=& o ( h_0^2 + h_1^2).
\end{eqnarray*}
This proves the first equality of Lemma \ref{Momvarphi}. The
second equality in Lemma  follows similarly, since
\begin{eqnarray*}
\lefteqn{ \Var [\widetilde{\varphi}_{in}(x, y)] = \esp \left[
\widetilde{\varphi}_{in}^2(x, y) \right] - \left( \esp \left[
\widetilde{\varphi}_{in}(x, y) \right] \right)^2 }
\\
& = & \frac{1}{h_0^dh_1} \int \int \varphi\left(x+h_0z,
y+h_1v\right) K_0^2(z) K_1^2(v) dz dv + O(1)
\\
&=& \frac{\varphi(x, y)} {h_0^dh_1} \int \int K_0^2(z) K_1^2(v) dz
dv + o\left(\frac{1}{h_0^dh_1}\right).
\end{eqnarray*}

The last statement of Lemma \ref{Momvarphi} is immediate, since
the Triangular and Convex inequalities give
\begin{eqnarray*}
\esp \left| \widetilde{\varphi}_{in}(x, y) - \esp
\widetilde{\varphi}_{in}(x, y) \right|^3 &\leq& C\esp \left|
\widetilde{\varphi}_{in}(x, y) \right|^3
\\
&\leq& \frac {C\varphi(x, y)} {h_0^{2d}h_1^2 } \int \int \left|
K_0(z) K_1(v) \right|^3 dz dv +
o\left(\frac{1}{h_0^{2d}h_1^2}\right). \eop
\end{eqnarray*}

\subsection*{Proof of Lemma \ref{STR}}
The order of $S_n$ follows from Lemma \ref{Betasum} and Lemma
\ref{Sigsum}. In fact,  since
\begin{eqnarray*}
\mathds{1}(X_i \in \mathcal{X}_0)
 \left(
 \widehat{m}_{in} - m(X_i)
\right)
 &=&
 \frac{\mathds{1}(X_i \in \mathcal{X}_0)}
 {nb_0^d\widehat{g}_{in}}
 \sum_{j=1,j\neq i}^n
 \left( m(X_j)+ \varepsilon_j - m(X_i)\right)
 K_0\left(\frac{X_j-X_i}{b_0}\right)
 \\
 & =&
 \beta_{in} + \Sigma_{in},
\end{eqnarray*}
Lemma \ref{Betasum} and Lemma \ref{Sigsum} give
\begin{eqnarray*}
S_n
 =
 O_{\prob}
\left[ b_0^2 \left(nb_1^2+(nb_1)^{1/2}\right)
 +
\left( nb_1^4 + \frac{b_1}{b_0^d} \right)^{1/2}
 \right],
\end{eqnarray*}
which gives the result for $S_n$.

\vskip 0.3cm For $T_n$, define for any $1\leq i\leq n$,
$$
\esp_{in} [\cdot] = \esp_n \left[
 X_1,\ldots,X_n,\varepsilon_k,
 k\neq i
 \right].
$$
Therefore, since $(\widehat{m}_{in} - m(X_i))$ depends only upon
$\left(X_1,\ldots,X_n,\varepsilon_k,k\neq i\right)$, we have
\begin{eqnarray*}
 \esp_n [T_{n}]
 &=&
 \esp_{n}
 \left[
 \sum_{i=1}^n
 \esp_{in}
 \left[
 \mathds{1}
\left( X_i \in \mathcal{X}_0 \right)
 (\widehat{m}_{in} - m(X_i))^2
  K_1^{(2)}
 \left(\frac{\varepsilon_{i}-\epsilon}{b_1}\right)
 \right]
 \right]
 \\
 &=&
 \esp_{n}
  \left[
  \sum_{i=1}^n
  \mathds{1}
 \left( X_i \in \mathcal{X}_0\right)
(\widehat{m}_{in} - m(X_i))^2
 \esp_{in}
 \left[
 K_1^{(2)}
 \left(\frac{\varepsilon_{i}-\epsilon}{b_1}\right)
 \right]
 \right],
\end{eqnarray*}
with, using $(A_4)$ and Lemma \ref{MomderK}-(\ref{MomderK2}),
\begin{eqnarray*}
 \left|
\esp_{in}
 \left[
 K_1^{(2)}
 \left(\frac{\varepsilon_{i}-\epsilon}{b_1}\right)
 \right]
 \right|
 =
 \left|
 \int
  K_1^{(2)}
 \left(\frac{e-\epsilon}{b_1}
 \right)
 f(e)
 de
 \right|
\leq
 Cb_1^3.
\end{eqnarray*}
 Hence this bound, the equality above, the Cauchy-Schwarz inequality
  and  Lemma \ref{BoundEspmchap} yield that
\begin{eqnarray}
\nonumber
 \left|
 \esp_{n}\left[T_{n}\right]
 \right|
&\leq&
 Cb_1^3
 \sum_{i=1}^n
  \esp_{n}
  \biggl[
  \mathds{1}\left( X_i \in \mathcal{X}_0\right)
 (\widehat{m}_{in} - m(X_i))^2
 \biggr]
 \\\nonumber
 &\leq&
 C n b_1^3
 \left(
 \sup_{1\leq i\leq n}
 \esp_{n}
  \biggl[
  \mathds{1}\left( X_i \in \mathcal{X}_0\right)
 (\widehat{m}_{in} - m(X_i))^4
 \biggr]
 \right)^{1/2}
 \\
 &\leq&
 O_{\prob}
 \left(nb_1^3\right)
 \left(b_0^4+\frac{1}{nb_0^d}\right).
 \label{Boundmean}
\end{eqnarray}

For the conditional variance of $T_n$,  Lemma \ref{sumzeta} gives
\begin{eqnarray*}
\Var_n (T_{n})
 &=&
 \sum_{i=1}^n
 \Var_n\left( \zeta_{in}\right)
 +
 \sum_{i=1}^n
 \sum_{j=1\atop j\neq i}^n
\Cov_n \left( \zeta_{in} , \zeta_{jn}
 \right)
\\
&=& O_{\prob} \left(nb_1\right)
\left(b_0^4+\frac{b_1}{nb_0^d}\right)^2 + O_{\prob}
\left(n^2b_0^db_1^{7/2}\right)
\left(b_0^4+\frac{1}{nb_0^d}\right)^2.
\end{eqnarray*}
Therefore, since $b_1$ goes to $0$ under $(A_{10})$,  this order
and (\ref{Boundmean}) yield,
 applying the  Tchebychev inequality,
\begin{eqnarray*}
T_{n}
 &=&
O_{\prob} \left[ \left(nb_1^3\right)
\left(b_0^4+\frac{1}{nb_0^d}\right)
 +
\left(nb_1\right)^{1/2} \left(b_0^4+\frac{b_1}{nb_0^d}\right) +
\left(n^2b_0^db_1^{7/2}\right)^{1/2}
\left(b_0^4+\frac{1}{nb_0^d}\right) \right]
\\
&=& O_{\prob} \left[ \left( nb_1^3 + \left(nb_1\right)^{1/2} +
\left(n^2b_0^db_1^{3}\right)^{1/2} \right)
\left(b_0^4+\frac{1}{nb_0^d}\right) \right].
\end{eqnarray*}
which gives the result for $T_n$.

\vskip 0.2cm We now compute the order of $R_n$. For this, define
\begin{eqnarray*}
I_{in} &= & \int_{0}^{1}
 (1-t)^2
 K_1^{(3)}
 \left(
 \frac{
 \varepsilon_i-t(\widehat{m}_{in} -m(X_i))-\epsilon}{b_1}
 \right)
 dt,
 \\
R_{in} & =&
 \mathds{1}
 \left(X_i\in\mathcal{X}_0\right)
 \left(\widehat{m}_{in}-m(X_i)\right)^3
 I_{in},
\end{eqnarray*}
and note that $R_n=\sum_{i=1}^n R_{in}$. The order of $R_n$ is
derived by computing its conditional mean and its conditional
variance. For the conditional mean,  observe that
\begin{eqnarray*}
 \esp_n [R_{n}]
&=&
 \esp_n
 \left[
 \sum_{i=1}^n
 \esp_{in}\left[R_{in}\right]
\right]
\\
 &=&
 \esp_{n}
 \left[
 \sum_{i=1}^n
 \mathds{1}
 \left( X_i \in \mathcal{X}_0\right)
 (\widehat{m}_{in} - m(X_i))^3
 \esp_{in}
 \left[ I_{in}\right]\right],
\end{eqnarray*}
with, using $(A_4)$ and Lemma \ref{MomderK}-(\ref{MomderK3}),
\begin{eqnarray*}
 \left|
 \esp_{in}
 \left[ I_{in}\right]
 \right|
& =&
 \left|
 \int_{0}^{1}
 (1-t)^2
 \left[
 \int
  K_1^{(3)}
 \left(
 \frac{e-t(\widehat{m}_{in} -m(X_i))-\epsilon}{b_1}
 \right)f(e)
 de
 \right]
 dt
 \right|
 \\
&\leq&
 Cb_1^3.
\end{eqnarray*}
Therefore the Holder inequality and Lemma \ref{BoundEspmchap}
yield
\begin{eqnarray}
\nonumber
 \left|
 \esp_n\left[R_n\right]
 \right|
 &\leq&
 Cb_1^3
 \sum_{i=1}^n
 \esp_n
 \left[
\left|
\mathds{1}
\left(X_i\in\mathcal{X}_0\right)
\left(\widehat{m}_{in}-m(X_i)\right)
\right|^3
\right]
\\\nonumber
&\leq&
Cb_1^3
\sum_{i=1}^n
\esp_n^{3/4}
\left[ \mathds{1}
\left(X_i\in\mathcal{X}_0\right)
\left(\widehat{m}_{in}-m(X_i)\right)^4
\right]
\\
&\leq&
 O_{\prob}
 \left(nb_1^3\right)
 \left(b_0^4+\frac{1}{nb_0^d}\right)^{3/2}.
 \label{Rn1}
\end{eqnarray}
For the conditional covariance of $R_n$, note that Lemma
\ref{Indep} allows to write
\begin{eqnarray}
\Var_n\left(R_n\right)
=
\sum_{i=1}^n
\Var_n\left(R_{in}\right)
+
\sum_{i=1}^n \sum_{j=1\atop j\neq i}^n
\biggl(\|X_i-X_j\|\leq
Cb_0\biggr)
\Cov_n\left(R_{in}, R_{jn}\right),
\label{Rn2}
\end{eqnarray}
and consider  the  first term in (\ref{Rn2}). We have
\begin{eqnarray*}
\Var_n\left(R_{in}\right)
 \leq
\esp_n\left[R_{in}^2\right] \leq
 \esp_{n}
 \biggl[
 \mathds{1}
 \left( X_i \in \mathcal{X}_0\right)
 (\widehat{m}_{in} - m(X_i))^6
 \esp_{in}
 \left[ I_{in}^2\right]
 \biggr],
\end{eqnarray*}
with, using  $(A_4)$, the Cauchy-Schwarz inequality and Lemma
\ref{MomderK}-(\ref{MomderK3}),
\begin{eqnarray*}
 \esp_{in}
 \left[ I_{in}^2\right]
 &\leq&
 C\esp_{in}
 \left[
\int_{0}^{1}
 K_1^{(3)}
 \left(
 \frac{\varepsilon_i-t(\widehat{m}_{in} -m(X_i))-\epsilon}{b_1}
 \right)^2
 dt
 \right]
\\
&\leq&
 C\int_{0}^{1}
 \left[
 \int
 K_1^{(3)}
 \left(
 \frac{e-t(\widehat{m}_{in} -m(X_i))-\epsilon}{b_1}
 \right)^2
 f(e)
 de
 \right]
 dt
 \\
&\leq&
 Cb_1,
\end{eqnarray*}
so that
\begin{eqnarray*}
 \Var_n\left(R_{in}\right)
\leq
 Cb_1
 \esp_{n}
 \left[
 \mathds{1}
 \left( X_i \in \mathcal{X}_0\right)
 (\widehat{m}_{in} - m(X_i))^6
 \right].
 \end{eqnarray*}
Therefore form Lemma \ref{BoundEspmchap}, we deduce
\begin{eqnarray}
\nonumber
 \sum_{i=1}^n
 \Var_n\left(R_{in}\right)
 &\leq&
 Cnb_1
 \sup_{1\leq i\leq n}
 \esp_{n}
 \left[
 \mathds{1}
 \left( X_i \in \mathcal{X}_0\right)
 (\widehat{m}_{in} - m(X_i))^6
 \right]
 \\
&\leq&
 O_{\prob}
 \left(nb_1\right)
 \left(b_0^4+\frac{1}{nb_0^d}\right)^3.
 \label{Rn4}
\end{eqnarray}
For the second  term in (\ref{Rn2}),  the Cauchy-Schwarz
inequality gives, with the help of the above result for
$\Var_n\left(R_{in}\right)$,
\begin{eqnarray*}
\left|
\Cov_n\left(R_{in}, R_{jn}\right)
 \right|
 &\leq&
 \left(
 \Var_n\left(R_{in}\right)
 \Var_n\left(R_{jn}\right)
 \right)^{1/2}
 \\
 &\leq&
 Cb_1
 \sup_{1\leq i\leq n}
 \esp_{n}
 \left[
 \mathds{1}
 \left( X_i \in \mathcal{X}_0\right)
 (\widehat{m}_{in} - m(X_i))^6
 \right].
\end{eqnarray*}
Hence by Lemma \ref{BoundEspmchap}  and the Markov inequality, we
have
\begin{eqnarray*}
\lefteqn{ \sum_{i=1}^n
\sum_{j=1\atop j\neq i}^n
\biggl(\|X_i-X_j\|\leq Cb_0\biggr)
\left|
\Cov_n\left(R_{in},
R_{jn}\right) \right| }
\\
&\leq&
 O_{\prob}\left(b_1\right)
 \left(b_0^4+\frac{1}{nb_0^d}\right)^3
 \sum_{i=1}^n
 \sum_{j=1\atop j\neq i}^n
 \biggl(\|X_i-X_j\|\leq Cb_0\biggr)
 \\
 &\leq&
 O_{\prob}\left(b_1\right)
 \left(b_0^4+\frac{1}{nb_0^d}\right)^3
 \left(n^2b_0^d\right).
\end{eqnarray*}
This order, (\ref{Rn4}) and (\ref{Rn2}) give, since $nb_0^d$
diverges  under $(A_9)$,
$$
\Var\left(R_n\right) =
 O_{\prob}
 \left(b_0^4+\frac{1}{nb_0^d}\right)^3
 \left(n^2b_0^db_1\right).
$$
Finally, with the help of this result, (\ref{Rn1}) and the
Tchebychev inequality, we arrive at
\begin{eqnarray*}
R_n &=& O_{\prob} \left[ \left(nb_1^3\right)
 \left(b_0^4+\frac{1}{nb_0^d}\right)^{3/2}
 +
\left(n^2b_0^db_1\right)^{1/2}
\left(b_0^4+\frac{1}{nb_0^d}\right)^{3/2} \right]
\\
&=& O_{\prob} \left[ \left( nb_1^3 +
\left(n^2b_0^db_1\right)^{1/2} \right)
\left(b_0^4+\frac{1}{nb_0^d}\right)^{3/2} \right].
 \eop
\end{eqnarray*}

\subsection*{Proof of Lemma \ref{MomderK}}
Set $h_p(e)=e^p f(e)$, $p\in[0,2]$. For  the first inequality of
(\ref{MomderK1}), note that  under $(A_5)$ and $(A_8)$, the change
of variable $e=\epsilon+b_1 v$ give, for any integer $\ell\in[1,
3]$,
\begin{eqnarray}
\nonumber
 \left|
 \int
 K_1^{(\ell)}
 \left(\frac{e-\epsilon}{b_1}\right)^2
 e^pf(e) de
 \right|
 &=&
 \left|
 b_1
 \int
 K_1^{(\ell)}(v)^2 h_p(\epsilon+b_1v)
 dv
 \right|
 \\\nonumber
 &\leq&
 b_1
 \sup_{t\in\Rit}
 |h_p(t)|
 \int
 | K_1^{(\ell)}(v)^2|
 dv
 \\
 &\leq&
 Cb_1,
 \label{Ineg1}
\end{eqnarray}
which yields the first inequality in (\ref{MomderK1}). For the
second inequality in (\ref{MomderK1}), observe that $f(\cdot)$ has
a bounded  continuous derivative under $(A_5)$, and that $\int \!
K_1^{(\ell)}(v)dv =0$ under $(A_8)$. Therefore, since $h_p(\cdot)$
has bounded second order derivatives under $(A_7)$, the Taylor
inequality yields that
\begin{eqnarray*}
\left| \int
 K_1^{(\ell)}
\left(\frac{e-\epsilon}{b_1}\right)
 e^pf(e)de
 \right|
 &=&
  b_1
 \left|
 \int
 K_1^{(\ell)}(v)
 \left[
 h_p(\epsilon+b_1v)-h_p(\epsilon)
 \right]
 \right|
  dv
 \\
 &\leq&
 b_1^2
 \sup_{t\in\Rit}|h_p^{(1)}(t)|
 \int
 |vK_1^{(\ell)}(v)|
  dv
 \leq
 Cb_1^2.
\end{eqnarray*}
which  completes the proof of (\ref{MomderK1}).

 The first inequalities of  (\ref{MomderK2}) and
 (\ref{MomderK3}) follow directly from (\ref{Ineg1}). The second bounds in
(\ref{MomderK2}) and (\ref{MomderK3}) are proved simultaneously.
For this, note that for any integer $\ell\in\{2,3\}$,
$$
\int
 K_1^{(\ell)}
 \left(\frac{e-\epsilon}{b_1}\right)
 h_p(e) de
 =
b_1 \int
 K_1^{(\ell)}(v)
h_p(\epsilon+b_1v) dv.
$$
Under $(A_8)$,  $K_1(\cdot)$ is symmetric, has  a compact support
and two
 continuous derivatives, with
$\int \! K_1^{(\ell)}(v)dv=0$ and
 $\int\! v K_1^{(\ell)}(v) dv=0$ for $\ell\in\{2,3\}$.
 Hence, since by $(A_5)$ $h_p$ has bounded continuous second
order derivatives, this gives for some $\theta=\theta
(\epsilon,b_1 v)$,
\begin{eqnarray*}
 \lefteqn{
 \left|
 \int
 K_1^{(\ell)}
\left(
 \frac{e-\epsilon}{b_1}
 \right)
  h_p(e) de
 \right|
 =
 \left|
  b_1
 \int
 K_1^{(\ell)}(v)
\left[
 h_p(\epsilon+b_1v) - h_p(\epsilon)
 \right]
 dv
 \right|
 }
 &&
\\
&=& \left| b_1 \int
 K_1^{(\ell)}(v)
\left[
 b_1 v h_p^{(1)}(\epsilon)
 +
  \frac{b_1^2v^2}{2}
h_p^{(2)}(\epsilon+\theta b_1v) \right]
 dv
 \right|
 \\
  &=&
 \left|
\frac{b_1^3}{2} \int v^2 K_1^{(\ell)}(v)
 h_p^{(2)}(\epsilon+\theta b_1v)
 dv
\right|
\\
&\leq&
 \frac{b_1^3}{2}
 \sup_{t\in \Rit}|h_p^{(2)}(t)|
 \int
\left| v^2K_1^{(\ell)}(v)
 \right|
 dv
 \leq
 Cb_1^3.
 \eop
\end{eqnarray*}

\subsection*{Proof of Lemma \ref{Betasum}}
Assumption $(A_4)$ and Lemma \ref{MomderK}-(\ref{MomderK1}) give
\begin{eqnarray*}
\left| \esp_n \left[
 \sum_{i=1}^n
 \beta_{in}
 K_1^{(1)}
 \left(
\frac{\varepsilon_i-\epsilon}{b_1} \right)
 \right]
 \right|
 & = &
\left| \esp \left[
 K_1^{(1)}
  \left(
  \frac{\varepsilon -\epsilon}{b_1}
  \right)
  \right]
  \sum_{i=1}^n
  \beta_{in}
   \right|
  \leq C n b_1^2
 \max_{1 \leq i \leq n}
 \left| \beta_{in} \right|,
\\
\Var_n
 \left[
 \sum_{i=1}^n
  \beta_{in}
  K_1^{(1)}
  \left(
\frac{\varepsilon_i-\epsilon}{b_1}
 \right)
 \right]
&\leq&
  \sum_{i=1}^n
 \beta_{in}^2
 \esp
 \left[
 K_1^{(1)}
  \left(
   \frac{\varepsilon - \epsilon}{b_1}
\right)^2
 \right]
  \leq
 C n b_1
  \max_{1 \leq i \leq n}
  \left|
  \beta_{in}
\right|^2 .
\end{eqnarray*}
Hence the (conditional) Markov inequality gives
$$
\sum_{i=1}^n \beta_{in}
 K_1^{(1)}
 \left(
\frac{\varepsilon_i-\epsilon}{b_1} \right) = O_{\prob}
 \left( n
b_1^2 + (nb_1)^{1/2}
 \right)
 \max_{1 \leq i \leq n}
  \left|
\beta_{in}
 \right|,
$$
so that the lemma follows if we can prove that
\begin{equation}
\sup_{1 \leq i \leq n}
 \left|\beta_{in}\right|
 =
 O_{\prob}
\left( b_0^2\right),
 \label{BetasumTBP}
\end{equation}
as established now. For this, define
$$
\zeta_j (x)
 =
 \mathds{1}
 \left( x \in \mathcal{X}_0\right)
\left(m(X_j)-m(x)\right)
 K_0\left( \frac{X_j-x}{b_0}\right),
 \;\;
\nu_{in}(x)
 =
 \frac{1}{(n-1) b_0^d}
 \sum_{j=1, j\neq i}^n
 \left(
\zeta_j (x)-\esp[\zeta_j (x)]
 \right),
$$
and $\bar{\nu}_{n} (x)= \esp[\zeta_j(x)] / b_0^d $, so that
$$
\beta_{in} =
 \frac{n-1}{n} \frac{\nu_{in} (X_i)
 +
 \bar{\nu}_n (X_i)}{\widehat{g}_{in}}
\;.
$$
For  $\max_{1 \leq i \leq n} | \bar{\nu}_n (X_i) |$, first observe
that  a second-order Taylor expansion applied successively to
$g(\cdot)$ and $m(\cdot)$ give, for $b_0$ small enough, and for
any $x$, $z$ in $\mathcal{X}$,
\begin{eqnarray*}
 \lefteqn
 {
 \left[ m(x+b_0z)-m(x)\right]
 g(x+b_0z)
 }
\\
& =& \left[
 b_0 m^{(1)}(x) z
 +
 \frac{b_0^2}{2}
z m^{(2)}(x +\zeta_1 b_0 z)z^{\top}
 \right]
 \left[
 g(x)
 +
 b_0 g^{(1)}(x) z
 +
 \frac{b_0^2}{2}
 z g^{(2)}(x +\zeta_2 b_0z)z^{\top}
 \right],
  \end{eqnarray*}
for some $\zeta_1 =\zeta_1 (x,b_0 z)$ and $\zeta_2 =\zeta_2 (x,b_0
z)$ in $[0,1]$. Therefore, since $\int\!z
 K(z)dz=0$ under $(A_7)$, it follows that, by $(A_1)$, $(A_2)$ and
 $(A_3)$,
\begin{eqnarray}
\nonumber
 \max_{1 \leq i \leq n}
|\bar{\nu}_n (X_i)|
&\leq&
\sup_{x\in\mathcal{X}_0} |\bar{\nu}_n
(x)|
 =
 \sup_{x \in
 \mathcal{X}_0}
 \left|
 \int
 \left(m ( x + b_0 z) - m(x)\right)
 K_0 (z) g(x+b_0z)
 dz
\right|
\\
&\leq&
 Cb_0^2.
 \label{Betasum1}
\end{eqnarray}
 Consider now the term  $\max_{1\leq i \leq n}|\nu_{in}(X_i)|$.
  The Bernstein inequality (see e.g. Serfling (2002)) and
$(A_4)$ give, for any $t>0$,
\begin{eqnarray*}
\prob \left( \max_{1 \leq i \leq n} | \nu_{in} (X_i)|
 \geq t
\right) &\leq & \sum_{i=1}^n \prob
 \left(
 | \nu_{in} (X_i) |
 \geq
t \right)
 \leq
 \sum_{i=1}^n
 \int
 \prob
 \left(
  | \nu_{in} (x) |
\geq t
 \left| X_i = x \right.
 \right) g (x)
  dx
\\
& \leq &
 2n \exp
 \left(
  - \frac{ (n-1) t^2 }
  { 2\sup_{x \in\mathcal{X}_0}
\Var (\zeta_j (x)/b_0^d) + \frac{4M}{3b_0^d} t}
 \right),
\end{eqnarray*}
where $M$ is such that $\sup_{x \in \mathcal{X}_0}|\zeta_j (x)|
\leq M$. The definition of $\mathcal{X}_0$ given in $(A_2)$,
$(A_3)$, $(A_7)$ and the standard Taylor expansion yield, for
$b_0$ small enough,
$$
\sup_{x \in \mathcal{X}_0}
 | \zeta_j (x) |
 \leq C b_0,
 \;\;\;
\sup_{x \in \mathcal{X}_0}
 \Var (\zeta_j (x)/b_0^d)
\leq \frac{1}{b_0^d}
 \sup_{x \in \mathcal{X}_0}
 \int
 \left( m(x +b_0 z) - m(x) \right)^2
 K_0^2 (z) g(x+b_0z) dz
  \leq
  \frac{C
b_0^2}{b_0^d}\;,
$$
so that, for any $t \geq 0$,
$$
\prob
 \left(
 \max_{1 \leq i \leq n}
 | \nu_{in} (X_i) |
 \geq t
\right) \leq 2n
 \exp
 \left( - \frac{(n-1) b_0^d t^2 /b_0^2}{C + C
t/b_0} \right).
$$
This gives
$$
\prob \left( \max_{1 \leq i \leq n} |\nu_{in} (X_i)| \geq \left(
\frac{b_0^2 \ln n}{ (n-1) b_0^d} \right)^{1/2} t\right) \leq 2n
\exp \left(
 - \frac{ t^2 \ln n }
 {C + C t \left( \frac{\ln n}{
(n-1) b_0^d} \right)^{1/2} } \right) = o(1),
$$
provided that $t$ is large enough and under $(A_9)$. It then
follows that
$$
\max_{1 \leq i \leq n} | \nu_{in} (X_i) |
 =
 O_{\prob}
 \left(
 \frac{b_0^2 \ln n}{ n b_0^d} \right)^{1/2}.
$$
This bound, (\ref{Betasum1}) and Lemma \ref{Estig}  show that
(\ref{BetasumTBP}) is proved,  since $b_0^2\ln
n/(nb_0^d)=O\left(b_0^4\right)$ under $(A_9)$, and that
$$
 \beta_{in}
  =
 \frac{n-1}{n} \frac{\nu_{in} (X_i)
 +
 \bar{\nu}_n (X_i)}{\widehat{g}_{in}}\;.
 \eop
$$

\subsection*{Proof of Lemma \ref{Sigsum}}
Note that $(A_4)$ gives that  $\Sigma_{in}$ is independent of
$\varepsilon_i$, and that $\esp_n[\Sigma_{in}]=0$. This yields
\begin{eqnarray}
 \esp_n
 \left[
 \sum_{i=1}^n
 \Sigma_{in}
 K_1^{(1)}
 \left(
 \frac{\varepsilon_i-\epsilon}{b_1}
 \right)
 \right]
  = 0.
\label{EspSigmai}
\end{eqnarray}
Moreover,  observe that
\begin{eqnarray}
\nonumber
 \lefteqn{
 \Var_n
  \left[
 \sum_{i=1}^n
 \Sigma_{in}
 K_1^{(1)}
  \left(
\frac{\varepsilon_i-\epsilon}{b_1}
 \right)
 \right]
 }
\\\nonumber
&=&
 \sum_{i=1}^n
 \Var_n
  \left[
 \Sigma_{in}
 K_1^{(1)}
 \left(
 \frac{\varepsilon_i-\epsilon}{b_1}
 \right)
 \right]
 +
 \sum_{i=1}^n
 \sum_{j=1\atop j\neq i}^n
 \Cov_n
 \left[
 \Sigma_{in}
 K_1^{(1)}
 \left(
 \frac{\varepsilon_{i}-\epsilon}{b_1}
 \right)
 ,
 \Sigma_{jn}
 K_1^{(1)}
 \left(
 \frac{\varepsilon_{j}-\epsilon}{b_1}
 \right)
 \right].
 \\
\label{VarSigm}
\end{eqnarray}
For the sum of variances in (\ref{VarSigm}), Lemma
\ref{MomderK}-(\ref{MomderK1}) and $(A_4)$ give
\begin{eqnarray}
\nonumber \sum_{i=1}^n
 \Var_n
 \left[
 \Sigma_{in} K_1^{(1)}
 \left(
 \frac{\varepsilon_i-\epsilon}{b_1}
 \right)
 \right]
 &\leq&
 \sum_{i=1}^n
 \esp_n
\left[ \Sigma_{in}^2
 \right]
 \esp
 \left[
 K_1^{(1)}
  \left(
\frac{\varepsilon_i-\epsilon}{b_1} \right)^2
 \right]
\\\nonumber
&\leq&
 \frac{C b_1\sigma^2}{(nb_0^d)^2}
 \sum_{i=1}^n
 \sum_{j=1\atop j \neq i}^{n}
 \frac{\mathds{1}
 (X_i \in\mathcal{X}_0)}
{\widehat{g}_{in}^2}
 K_0^2
\left(\frac{X_j-X_i}{b_0}\right)
\\
 &\leq&
 \frac{C b_1\sigma^2}{nb_0^d}
 \sum_{i=1}^n
 \frac{\mathds{1}
 (X_i\in \mathcal{X}_0)
 \widetilde{g}_{in}}
{\widehat{g}_{in}^2}\;,
 \label{VarSigmai}
\end{eqnarray}
where $\sigma^2=\Var(\varepsilon)$ and
$$
\widetilde{g}_{in}
 =
 \frac{1}{n b_0^d}
 \sum_{j=1,j \neq i}^{n}
 K_0^2\left( \frac{X_j-X_i}{b_0}\right).
$$
For the sum of  conditional covariances in (\ref{VarSigm}),
observe that by $(A_4)$ we have
\begin{eqnarray*}
 \lefteqn{
 \sum_{i=1}^n
 \sum_{j=1\atop j\neq i}^n
\Cov_n \left[
 \Sigma_{in}
 K_1^{(1)}
 \left(
 \frac{\varepsilon_{i}-\epsilon}{b_1}
 \right)
 ,
 \Sigma_{jn}
 K_1^{(1)}
 \left(
  \frac{\varepsilon_{j}-\epsilon}{b_1}
 \right)
\right]
 }
\\
&=&
 \sum_{i=1}^n
 \sum_{j=1\atop j\neq i}^n
\esp_n
 \left[
 \Sigma_{in}
 \Sigma_{jn}
 K_1^{(1)}
 \left(
\frac{\varepsilon_{i}-\epsilon}{b_1}
 \right)
 K_1^{(1)}
  \left(
\frac{\varepsilon_{j}-\epsilon}{b_1}
 \right)
  \right]
\\
& = &
 \sum_{i=1}^n
 \sum_{j=1\atop j\neq i}^n
 \frac{\mathds{1}(X_{i}\in\mathcal{X}_0)
\mathds{1}(X_{j}\in\mathcal{X}_0)}
 {(n b_0^d)^2 \widehat{g}_{in} \widehat{g}_{jn}}
 \sum_{k=1\atop k\neq i}^n
\sum_{\ell=1\atop \ell\neq j}^n
 K_0
 \left(
\frac{X_{k}-X_{i}}{b_0}
 \right)
 K_0
 \left(
\frac{X_{\ell}-X_{j}}{b_0} \right) \esp
 \left[
 \xi_{ki}
 \xi_{\ell j}
 \right],
\end{eqnarray*}
where
$$
\xi_{ki}
 =
 \varepsilon_k
 K_1^{(1)}
  \left(
\frac{\varepsilon_{i}-\epsilon}{b_1}
 \right).
$$
Moreover, under $(A_4)$, it is seen that for $k\neq\ell$,
$\esp[\xi_{ki}\xi_{\ell j}] =0 $ when $\Card\{i, j, k,
\ell\}\geq3$. Therefore
 the symmetry of $K_0$  yields that
\begin{eqnarray*}
 \lefteqn{
 \sum_{i=1}^n
 \sum_{j=1\atop j\neq i}^n
\Cov_n \left[ \Sigma_{in}
 K_1^{(1)}
  \left(
\frac{\varepsilon_{i}-\epsilon}{b_1}
 \right),
 \Sigma_{jn}
K_1^{(1)}
 \left(
 \frac{\varepsilon_{j}-\epsilon}{b_1}
 \right)
 \right]
 }
 &&
\\
& = & \sum_{i=1}^n
 \sum_{j=1\atop j\neq i}^n
\frac{\mathds{1}(X_{i} \in \mathcal{X}_0) \mathds{1} (X_{j}
 \in
\mathcal{X}_0)}
 {(n b_0^d)^2 \widehat{g}_{in} \widehat{g}_{jn}}
K_0^2 \left( \frac{X_{j}-X_{i}}{b_0} \right)
 \esp^2
 \left[
\varepsilon K_1^{(1)} \left( \frac{\varepsilon-\epsilon}{b_1}
\right) \right]
\\
&& + \sum_{i=1}^n
 \sum_{j=1\atop j\neq i}^n
 \frac{ \mathds{1} (X_{i}
\in \mathcal{X}_0) \mathds{1} (X_{j} \in \mathcal{X}_0)}
 {(nb_0^d)^2\widehat{g}_{in} \widehat{g}_{jn}}
 \sum_{k=1\atop k\neq  i, j}^n
 K_0 \left(\frac{X_{k}-X_{i}}{b_0}\right)
 K_0\left(\frac{X_{k}-X_{j}}{b_0}\right)
 \esp[\varepsilon^2]
\esp^2 \left[ K_1^{(1)}
 \left( \frac{\varepsilon-\epsilon}{b_1}
\right) \right].
\\
 %\label{CovSigmai}
\end{eqnarray*}
Therefore, since
$$
\sup_{1\leq j\leq n}
 \left(
\frac{\mathds{1}\left(X_j\in\mathcal{X}_0\right)}
{|\widehat{g}_{jn}|} \right)
 =O_{\prob}(1)
$$
by Lemma \ref{Estig},  Lemma \ref{MomderK}-(\ref{MomderK1}) and
$(A_4)$ then give
\begin{eqnarray}
\nonumber
 \lefteqn{
  \left|
 \sum_{i=1}^n
 \sum_{j=1\atop j\neq i}^n
\Cov_n
 \left[
  \Sigma_{in}
 K_1^{(1)}
 \left(
\frac{\varepsilon_{i}-\epsilon}{b_1}
 \right)
 ,
  \Sigma_{jn}
 K_1^{(1)}
 \left(
 \frac{\varepsilon_{j}-\epsilon}{b_1}
\right) \right] \right|
 }
\\
  &=&
 O_{\prob}
 \left(
 \frac{b_1^4}{n b_0^d}
 \right)
 \sum_{i=1}^n
 \frac{\mathds{1}(X_i\in \mathcal{X}_0)\widetilde{g}_{in}}
 {|\widehat{g}_{in}|}
 +
 O_{\prob} (b_1^4)
 \sum_{i=1}^n
 \frac{\mathds{1}(X_i\in \mathcal{X}_0)|g_{in}|}
 {|\widehat{g}_{in}|}\;,
 \label{CovSigmaib}
\end{eqnarray}
where  $\widetilde{g}_{in}$ is defined as  in (\ref{VarSigmai})
and
$$
 g_{in}
 =
 \frac{1}{(n b_0^d)^2}
 \sum_{j=1\atop j\neq i}^n
 \sum_{k=1\atop k\neq j, i}^n
 K_0\left(\frac{X_k-X_i}{b_0}\right)
 K_0\left(\frac{X_k-X_j}{b_0}\right).
$$
The order of the first term in (\ref{CovSigmaib}) follows from
Lemma \ref{Estig}, which  gives
 \begin{eqnarray}
 \sum_{i=1}^n
 \frac{\mathds{1}(X_i\in \mathcal{X}_0)\widetilde{g}_{in}}
 {|\widehat{g}_{in}|}
 =
 O_{\prob}(n).
 \label{CovSigmaib1}
 \end{eqnarray}
Again, by Lemma \ref{Estig}, we have
\begin{eqnarray*}
 \sum_{i=1}^n
 \frac{\mathds{1}(X_i\in \mathcal{X}_0)|g_{in}|}
 {|\widehat{g}_{in}|}
 =
 O_{\prob}(1)
 \sum_{i=1}^n
 \mathds{1}
\left(X_i\in \mathcal{X}_0\right)|g_{in}|,
%\label{CovSigmaib2}
\end{eqnarray*}
with, using  the changes of variables $x_1=x_3+b_0z_1$,
$x_2=x_3+b_0z_2$,
\begin{eqnarray*}
 \esp\left[ \sum_{i=1}^n \mathds{1}
\left(X_i\in\mathcal{X}_0\right)|g_{in}|
\right]
&\leq&
\frac{Cn^3}{(nb_0^d)^2}
\esp \left|
 K_1\left(\frac{X_3-X_1}{h}\right)
K_1\left(\frac{X_3-X_2}{h}\right)
 \right|
\\
&=&
\frac{Cn^3}{n^2h^2}
\int_{\mathcal{X}_0^3}
\left|K_1\left(\frac{x_3-x_1}{h}\right)
K_1\left(\frac{x_3-x_2}{h}\right)
\right|
\prod_{k=1}^3 g(x_k)dx_k
 \\
 &\leq&
 \frac{Cn^3b_0^{2d}}{(nb_0^d)^2}\;.
\end{eqnarray*}
These bounds and the equality above, give under  $(A_2)$ and
$(A_7)$,
$$
 \sum_{i=1}^n
 \frac{\mathds{1}(X_i\in \mathcal{X}_0)|g_{in}|}
 {|\widehat{g}_{in}|}
 =
 O_{\prob}(n).
$$
Hence from  (\ref{CovSigmaib1}), (\ref{CovSigmaib}),
(\ref{VarSigmai}), (\ref{VarSigm}) and Lemma \ref{Estig}, we
deduce, for $b_1$ small enough,
\begin{eqnarray*}
\lefteqn{
 \Var_n
 \left[
 \sum_{i=1}^n
 \Sigma_{in}
 K_1^{(1)}
 \left(
 \frac{\varepsilon_i-\epsilon}{b_1}
 \right)
 \right]
 }
 &&
 \\
 & =&
 O_{\prob}
 \left(
 \frac{b_1}{n b_0^d}
 \right)
 \sum_{i=1}^{n}
 \frac{\mathds{1}(X_i\in \mathcal{X}_0)\widetilde{g}_{in}}
 {\widehat{g}_{in}^2}
 +
 O_{\prob}
 \left(
 \frac{b_1^4}{n b_0^d}
 \right)
\sum_{i=1}^n
 \frac{\mathds{1}(X_i \in \mathcal{X}_0)
\widetilde{g}_{in}} {|\widehat{g}_{in}|}
 +
 O_{\prob}(b_1^4)
 \sum_{i=1}^n
\frac{\mathds{1}(X_i \in\mathcal{X}_0)|g_{in}|}
{|\widehat{g}_{in}|}
\\
&=&
 O_{\prob}
\left( \frac{b_1}{b_0^d}
 +
 \frac{b_1^4}{b_0^d}
 +
 nb_1^4
 \right)
 =
 O_{\prob}
 \left(
 \frac{b_1}{b_0^d}
 +
 nb_1^4
 \right).
\end{eqnarray*}
Finally, this order, (\ref{EspSigmai}) and the Tchebychev
inequality give
$$
\sum_{i=1}^n
 \Sigma_{in}
 K_1^{(1)}
 \left(
 \frac{\varepsilon_i-\epsilon}{b_1}
 \right)
 =
 O_{\prob}
 \left(
 \frac{b_1}{b_0^d}+nb_1^4
 \right)^{1/2}. \eop
$$

\subsection*{Proof of Lemma \ref{BoundEspmchap}}

Define $\beta_{in}$  as in Lemma \ref{Betasum} and set
\begin{eqnarray*}
g_{in}
 =
 \frac{1}{nb_0^d}
 \sum_{j=1, j\neq i}^n
 K_0^4\left(\frac{X_j-X_i}{b_0}\right),
 \quad
\widetilde{g}_{in}
 =
 \frac{1}{nb_0^d}
 \sum_{j=1, j\neq i}^n
 K_0^2\left(\frac{X_j-X_i}{b_0}\right).
\end{eqnarray*}
The proof of the lemma  is based on the following bound:
\begin{eqnarray}
\esp_n
 \biggl[
 \mathds{1}
 \left(X_i\in\mathcal{X}_0\right)
 \left(\widehat{m}_{in} - m(X_i)\right)^k
\biggr]
 \leq
 C
 \left[
 \beta_{in}^k
 +
 \frac{
 \mathds{1}
 \left(X_i\in\mathcal{X}_0\right)
 \widetilde{g}_{in}^{k/2}}
 {(nb_0^d)^{(k/2)}\widehat{g}_{in}^k}
 \right],
 \quad
 k\in\{4,6\}.
 \label{Espm}
 \end{eqnarray}
 Indeed, taking successively $k=4$ and $k=6$ in (\ref{Espm}),
 we have, by (\ref{BetasumTBP}), Lemma \ref{Estig} and
 $(A_9)$,
\begin{eqnarray*}
 \sup_{1\leq i\leq n}
 \esp_n
 \biggl[
 \mathds{1}\left(X_i\in\mathcal{X}_0\right)
 \left(\widehat{m}_{in} - m(X_i)\right)^4
\biggr]
 &=&
 O_{\prob}
 \left(
 b_0^8
  +
 \frac{1}{(nb_0^d)^2}
 \right)
 =
 O_{\prob}
 \left(b_0^4+\frac{1}{nb_0^d}\right)^2,
 \\
 \sup_{1\leq i\leq n}
 \esp_n
 \biggl[
 \mathds{1}
 \left(X_i\in\mathcal{X}_0\right)
 \left(\widehat{m}_{in} - m(X_i)\right)^6
\biggr]
 &=&
 O_{\prob}
 \left(
 b_0^{12}
  +
 \frac{1}{(nb_0^d)^3}
 \right)
 =
 O_{\prob}
 \left(b_0^4+\frac{1}{nb_0^d}\right)^3,
\end{eqnarray*}
which gives the results of the Lemma.
 Hence it remains to prove (\ref{Espm}).
 For this,  define $\beta_{in}$ and  $\Sigma_{in}$  respectively as in
 Lemma \ref{Betasum} and Lemma \ref{Sigsum}.
 Since
 $\mathds{1}(X_i \in\mathcal{X}_0)
 \left(\widehat{m}_{in}- m(X_i)\right)
 =\beta_{in}+\Sigma_{in}$, and that
 $\beta_{in}$ depends only on
 $\left(X_1,\ldots,X_n\right)$, this gives, for $k\in\{4,6\}$
\begin{eqnarray}
\esp_n
 \biggl[
 \mathds{1}
 (X_i \in \mathcal{X}_0)
 \left(\widehat{m}_{in}- m(X_i)\right)^k
 \biggr]
  \leq
 C\beta_{in}^k
 +
 C\esp_n\left[\Sigma_{in}^k\right].
 \label{Espm5}
\end{eqnarray}
The order  of the second term of bound (\ref{Espm5}) is computed
by applying Theorem 2 in Whittle (1960) or the
Marcinkiewicz-Zygmund inequality (see e.g Chow and Teicher, 2003,
p. 386). These inequalities show that for linear form
$L=\sum_{j=1}^n a_j\zeta_j$ with independent mean-zero random
variables
 $\zeta_1,\ldots,\zeta_n$, it holds that, for any $k\geq 1$,
 $$
 \esp
 \left|L^k\right|
 \leq
  C(k)
  \left[
 \sum_{j=1}^n
  a_j^2
\esp^{2/k}
\left|\zeta_j^k \right|
 \right]^{k/2},
 $$
where $C(k)$ is a positive real depending only on $k$. Now,
observe that for any $i\in[1,n]$,
$$
\Sigma_{in}
 =
 \sum_{j=1, j\neq i}^n
 \sigma_{jin},
 \quad
 \sigma_{jin}
 =
 \frac{\mathds{1}
 \left( X_i \in \mathcal{X}_0\right)}
 {nb_0^d\widehat{g}_{in}}
 \varepsilon_j
 K_0\left(\frac{X_j-X_i}{b_0}\right).
$$
Since under $(A_4)$, the $\sigma_{jin}$'s, $j\in[1,n]$, are
centered independent variables given $X_1,\ldots,X_n$,  this
yields, for any $k\in\{4,6\}$,
\begin{eqnarray*}
 \esp_n
 \left[
 \Sigma_{in}^k\right]
 \leq
 C\esp\left[\varepsilon^k\right]
 \left[
 \frac{\mathds{1}
 \left( X_i \in \mathcal{X}_0\right)}
 {(nb_0^d)^2\widehat{g}_{in}^2}
 \sum_{j=1}^n
  K_0^2
 \left(
 \frac{X_j-X_i}{b_0}
 \right)
 \right]^{k/2}
 \leq
 \frac{C\mathds{1}
 \left(X_i\in\mathcal{X}_0\right)\widetilde{g}_{in}^{k/2}}
 {(nb_0^d)^{(k/2)}\widehat{g}_{in}^k}\;.
\end{eqnarray*}
 Hence this  bound and
(\ref{Espm5})  give
$$
 \esp_n
 \biggl[
 \mathds{1}(X_i \in \mathcal{X}_0)
 \left(\widehat{m}_{in}- m(X_i)\right)^k
 \biggr]
 \leq C\left[
 \beta_{in}^k
  +
 \frac{
 \mathds{1}
 \left(X_i\in\mathcal{X}_0\right)
 \widetilde{g}_{in}^{k/2}}
 {(nb_0^d)^{(k/2)}\widehat{g}_{in}^k}
 \right],
$$
which proves (\ref{Espm}), and then completes the proof of the
lemma. \eop

\subsection*{Proof of Lemma \ref{Indep}}

Since $K_0(\cdot)$ has a compact support under $(A_7)$, there is a
$C>0$ such that $\| X_i - X_j \| \geq C b_0$ implies  that for any
integer number $k$ of $[1,n]$, $K_0 ( (X_k - X_i)/b_0) = 0$ if
$K_0 ( (X_j - X_k)/b_0) \neq 0$. Let $D_j \subset [1,n]$ be such
that an integer number $k$ of $[1,n]$ is in $D_j$ if and only if
$K_0 ( (X_j - X_k)/b_0) \neq 0$. Abbreviate $\prob (\cdot| X_1,
\ldots,X_n)$ into $\prob_n$ and assume that $\| X_i - X_j \| \geq
C b_0$ so that $D_i$ and $D_j$ have an empty intersection. Note
also that taking $C$ large enough ensures that $i$ is not in $D_j$
and $j$ is not in $D_i$. It then follows, under $(A_4)$ and since
$D_i$ and $D_j$ only depend upon $X_1,\ldots,X_n$,
\begin{eqnarray*}
\lefteqn{
 \prob_n
 \biggl(
 \left(
 \widehat{m}_{in} - m(X_i),\varepsilon_i
 \right)
 \in A \mbox{ \rm and }
 \left(
\widehat{m}_{jn} - m (X_j),\varepsilon_j \right)
 \in B \biggr)
 }
&&
\\
& = & \prob_n \left( \left( \frac{\sum_{k \in D_i \setminus\{i\}}
 \left( m(X_{k}) - m(X_i) + \varepsilon_{k} \right)
  K_0
\left( (X_{k} - X_i)/b_0 \right)} {\sum_{k \in D_i \setminus\{i\}}
 K_0 \left( (X_{k} - X_i)/b_0 \right)}
 ,
 \varepsilon_i\right) \in A \right.
\\
&& \;\;\;\;\;\;\;\;\;\;\;\;\;\;\;\;\;\;\;\; \left. \mbox{ \rm and}
\left( \frac{ \sum_{\ell \in D_j \setminus \{j \}}
 \left(
m(X_{\ell}) - m(X_j) + \varepsilon_{\ell} \right)
 K_0 \left(
(X_{\ell} - X_j)/b_0 \right) } { \sum_{\ell \in D_j \setminus \{j
\}} K_0 \left( (X_{\ell} - X_j)/b_0 \right)} ,
\varepsilon_j\right) \in B \right)
\\
& = & \prob_n \left( \left( \frac{ \sum_{k \in D_i \setminus \{i
\}} \left( m(X_{k}) - m(X_i) + \varepsilon_{k}\right)
 K_0 \left(
(X_{k} - X_i)/b_0 \right)} {\sum_{k \in D_i \setminus \{i \}}
K_0\left( (X_{k} - X_i)/b_0 \right)} ,
 \varepsilon_i \right) \in
A \right)
\\
&& \;\;\;\;\;\;\;\;\;\;\;\;\;\;\;\;\;\;\;\;
 \times\;
 \prob_n
 \left(
\left( \frac{ \sum_{\ell \in D_j \setminus \{j \}} \left(
m(X_{\ell}) - m(X_j) + \varepsilon_{\ell} \right) K_0 \left(
(X_{\ell} - X_j)/b_0 \right) } { \sum_{\ell \in D_j \setminus \{j
\}} K_0\left( (X_{\ell} - X_j)/b_0\right) } , \varepsilon_j
\right) \in B \right)
\\
& = & \prob_n
 \left(
 \left(\widehat{m}_{in} - m(X_i), \varepsilon_i \right)
 \in A
 \right)
\times \prob_n \left( \left(\widehat{m}_{jn} - m
(X_j),\varepsilon_j \right)
 \in B
\right).
\end{eqnarray*}
This gives the result of Lemma \ref{Indep}, since both
$\left(\widehat{m}_{in} - m (X_i), \varepsilon_i\right)$ and
$\left(\widehat{m}_{jn} - m (X_j), \varepsilon_j\right)$ are
independent given $X_1, \ldots, X_n$. \eop

\subsection*{Proof of Lemma \ref{sumzeta}}
Since $\widehat{m}_{in} - m(X_i)$ depends only upon
 $\left(X_1,\ldots,X_n,\varepsilon_k, k\neq i\right)$,
 we have
\begin{eqnarray*}
\sum_{i=1}^n
 \Var_n
 \left(\zeta_{in}\right)
  \leq
 \sum_{i=1}^n
 \esp_n
 \left[
 \zeta_{in}^2
 \right]
 =
 \sum_{i=1}^n
 \esp_n
 \left[
 \mathds{1}
 \left(X_i\in\mathcal{X}_0\right)
 \left(\widehat{m}_{in} - m(X_i)\right)^4
  \esp_{in}
  \left[
  K_1^{(2)}
  \left(\frac{\varepsilon_{i}-\epsilon}{b_1}\right)^2
  \right]
 \right],
 \end{eqnarray*}
with, using Lemma \ref{MomderK}-(\ref{MomderK2}),
\begin{eqnarray*}
 \esp_{in}
 \left[
 K_1^{(2)}
 \left(\frac{\varepsilon_{i}-\epsilon}{b_1}\right)^2
 \right]
 =
 \int
  K_1^{(2)}
 \left(
 \frac{e-\epsilon}{b_1}
 \right)^2
 f(e) de
 \leq
 C b_1.
\end{eqnarray*}
 Therefore these bounds and Lemma \ref{BoundEspmchap} give
\begin{eqnarray*}
\sum_{i=1}^n
 \Var_n\left(\zeta_{in}\right)
 &\leq&
 Cb_1
 \sum_{i=1}^n
 \esp_n
 \biggl[
 \mathds{1}
 \left(X_i\in\mathcal{X}_0\right)
 (\widehat{m}_{in} - m(X_i))^4
 \biggr]
 \\
 &\leq&
 Cnb_1
 \sup_{1\leq i\leq n}
 \esp_n
 \biggl[
 \mathds{1}
 \left(X_i\in\mathcal{X}_0\right)
 (\widehat{m}_{in} - m(X_i))^4
 \biggr]
 \\
 &\leq&
 O_{\prob}\left(nb_1\right)
 \left(b_0^4+\frac{1}{nb_0^d}\right)^2.
\end{eqnarray*}
which yields the desired result for the conditional variance.

\vskip 0.3cm
 We now prepare to compute the order of the conditional covariance.
 To that aim, observe that Lemma \ref{Indep} gives
\begin{eqnarray*}
\sum_{i=1}^n \sum_{j=1\atop j\neq i}^n \Cov_n \left(
\zeta_{in},\zeta_{jn}
 \right)
 =
  \sum_{i=1}^n
  \sum_{j=1\atop j\neq i}^n
  \mathds{1}
  \biggl(\left\|X_i - X_j\right\|<C b_0\biggr)
  \biggl(
  \esp_n
  \left[
  \zeta_{in}
  \zeta_{jn}
  \right]
  -
 \esp_n\left[\zeta_{in}\right]
 \esp_n\left[\zeta_{jn}\right]
 \biggr).
\end{eqnarray*}
 The order of the term above  is derived from the following
equalities:
\begin{eqnarray}
 \sum_{i=1}^n
 \sum_{j=1\atop j\neq i}^n
 \mathds{1}
 \biggl(
 \left\|X_i - X_j \right\|<C b_0
 \biggr)
 \esp_n\left[\zeta_{in}\right]
 \esp_n\left[\zeta_{jn}\right]
 &=&
 O_{\prob}
 \left(n^2b_0^db_1^6\right)
\left(b_0^4+ \frac{1}{nb_0^d}\right)^2,
  \label{Covzeta2}
  \\
  \sum_{i=1}^n
  \sum_{j=1\atop j\neq i}^n
  \mathds{1}
  \biggl(
  \left\|X_i - X_j\right\|<C b_0
  \biggr)
  \esp_n
  \left[
  \zeta_{in}
  \zeta_{jn}
  \right]
  &=&
  O_{\prob}
 \left(n^2b_0^db_1^{7/2}\right)
\left(b_0^4+ \frac{1}{nb_0^d}\right)^2.
 \label{Covzeta1}
\end{eqnarray}
 Indeed, since $b_1$ goes to $0$ under $(A_{10})$,
 (\ref{Covzeta2}) and (\ref{Covzeta1}) yield that
\begin{eqnarray*}
 \sum_{i=1}^n
 \sum_{j=1\atop j\neq i}^n
 \Cov_n
 \left(\zeta_{in},\zeta_{jn}\right)
 &=&
O_{\prob}
 \left[
 \left(n^2b_0^db_1^6\right)
 \left(b_0^4+ \frac{1}{nb_0^d}\right)^2
 +
 \left(n^2b_0^db_1^{7/2}\right)
 \left(b_0^4+ \frac{1}{nb_0^d}\right)^2
\right]
\\
&=& O_{\prob}
 \left(n^2b_0^db_1^{7/2}\right)
 \left(b_0^4+ \frac{1}{nb_0^d}\right)^2,
\end{eqnarray*}
 which gives the result for the conditional
covariance. Hence,  it remains to prove (\ref{Covzeta2}) and
(\ref{Covzeta1}).  For (\ref{Covzeta2}), note that by $(A_4)$ and
Lemma \ref{MomderK}-(\ref{MomderK2}), we have
 \begin{eqnarray*}
 \left|
 \esp_{n}
 \left[\zeta_{in}\right]
 \right|
 &=&
 \left|
 \esp_n
 \left[
 \mathds{1}
 \left(X_i\in\mathcal{X}_0\right)
 (\widehat{m}_{in} - m(X_i))^2
 \esp_{in}
 \left[
 K_1^{(2)}
 \left(
 \frac{\varepsilon_i-\epsilon}{b_1}
 \right)
 \right]
 \right]
 \right|
 \\
 &\leq&
 Cb_1^3
 \biggl(
  \esp_n
 \biggl[
 \mathds{1}
 \left(X_i\in\mathcal{X}_0\right)
 (\widehat{m}_{in} - m(X_i))^4
 \biggr]
 \biggr)^{1/2}.
 \end{eqnarray*}
Hence from this bound and Lemma \ref{BoundEspmchap} we deduce
\begin{eqnarray*}
\sup_{1\leq i, j\leq n}
 \left|
 \esp_{n}
 \left[\zeta_{in}\right]
 \esp_n
\left[\zeta_{jn}\right]
 \right|
 &\leq&
 Cb_1^6
 \sup_{1\leq i \leq n}
 \esp_n
 \biggl[
 \mathds{1}
 \left(X_i\in\mathcal{X}_0\right)
 (\widehat{m}_{in} - m(X_i))^4
 \biggr]
 \\
 &\leq&
 O_{\prob}
 \left(b_1^6\right)
\left(b_0^4+ \frac{1}{nb_0^d}\right)^2.
\end{eqnarray*}
Therefore, since the Markov inequality gives
\begin{eqnarray}
\sum_{i=1}^n
 \sum_{j=1\atop j\neq i}^n
 \mathds{1}
 \biggl(
 \| X_i - X_j \|< C b_0
 \biggr)
 =
 O_{\prob}(n^2 b_0^d),
 \label{Markov}
 \end{eqnarray}
 it then follows that
\begin{eqnarray*}
  \sum_{i=1}^n
 \sum_{j=1\atop j\neq i}^n
 \mathds{1}
 \biggl(
 \| X_i - X_j \|< C b_0
 \biggr)
 \esp_{n}
 \left[\zeta_{in}\right]
 \esp_{n}
 \left[\zeta_{jn}\right]
  =
 O_{\prob}
 \left(n^2b_0^db_1^6\right)
\left(b_0^4+ \frac{1}{nb_0^d}\right)^2,
\end{eqnarray*}
which   proves (\ref{Covzeta2}).

\vskip 0.3cm
 For (\ref{Covzeta1}), set $Z_{in}= \mathds{1}
\left(X_i\in\mathcal{X}_0\right)\left(\widehat{m}_{in} -
m(X_i)\right)^2$, and note that for $i\neq j$, we have
 \begin{eqnarray}
 \esp_n\left[\zeta_{in}\zeta_{jn}\right]
 =
 \esp_n
 \left[
 Z_{in}
  K_1^{(2)}
 \left(
 \frac{\varepsilon_j-\epsilon}{b_1}
 \right)
  \esp_{in}
 \left[
 Z_{jn}
  K_1^{(2)}
 \left(
 \frac{\varepsilon_i-\epsilon}{b_1}
 \right)
 \right]
 \right],
 \label{Prodzeta}
 \end{eqnarray}
where
\begin{eqnarray}
\nonumber
 \lefteqn{
 \esp_{in}
 \left[
 Z_{jn}
  K_1^{(2)}
 \left(
 \frac{\varepsilon_i-\epsilon}{b_1}
 \right)
 \right]
}
\\\nonumber
&=&
 \beta_{jn}^2
 \esp_{in}
 \left[
 K_1^{(2)}
 \left(
 \frac{\varepsilon_i-\epsilon}{b_1}
 \right)
 \right]
 +
 2\beta_{jn}
 \esp_{in}
 \left[
 \Sigma_{jn}
 K_1^{(2)}
 \left(
 \frac{\varepsilon_i-\epsilon}{b_1}
 \right)
 \right]
 +
 \esp_{in}
 \left[
 \Sigma_{jn}^2
 K_1^{(2)}
 \left(
 \frac{\varepsilon_i-\epsilon}{b_1}
 \right)
 \right].
 \\
\label{Covzeta3}
\end{eqnarray}
The first term of Equality (\ref{Covzeta3}) is treated by using
Lemma \ref{MomderK}-(\ref{MomderK2}). This gives
\begin{eqnarray}
\left| \beta_{jn}^2
 \esp_{in}
 \left[
 K_1^{(2)}
 \left(
 \frac{\varepsilon_i-\epsilon}{b_1}
 \right)
 \right]
 \right|
 \leq
 Cb_1^3
 \beta_{jn}^2.
 \label{Covzeta4}
\end{eqnarray}
Since  under $(A_4)$, the $\varepsilon_j$'s are independent
centered variables, and  are independent of the  $X_j$'s, the
second term in (\ref{Covzeta3}) gives
\begin{eqnarray*}
\esp_{in}
 \left[
 \Sigma_{jn}
   K_1^{(2)}
 \left(
 \frac{\varepsilon_i-\epsilon}{b_1}
 \right)
 \right]
&=&
 \frac{\mathds{1}\left(X_j\in\mathcal{X}_0\right)}
 {nb_0^d\widehat{g}_{jn}}
\sum_{k=1, k\neq j}^n
 K_0
 \left(
 \frac{X_k-X_j}{b_0}
 \right)
\esp_{in}
 \left[
 \varepsilon_k
   K_1^{(2)}
 \left(
 \frac{\varepsilon_i-\epsilon}{b_1}
 \right)
 \right]
 \\
 &=&
 \frac{\mathds{1}\left(X_j\in\mathcal{X}_0\right)}
 {nb_0^d\widehat{g}_{jn}}
 K_0
 \left(
 \frac{X_i-X_j}{b_0}
 \right)
\esp_{in}
 \left[
 \varepsilon_i
   K_1^{(2)}
 \left(
 \frac{\varepsilon_i-\epsilon}{b_1}
 \right)
 \right].
\end{eqnarray*}
Therefore, by $(A_7)$ which ensures that $K_0$ is bounded, the
equality above and Lemma \ref{MomderK}-(\ref{MomderK2}) yield that
\begin{eqnarray}
\left| \beta_{jn} \esp_{in}
 \left[
 \Sigma_{jn}
   K_1^{(2)}
 \left(
 \frac{\varepsilon_i-\epsilon}{b_1}
 \right)
 \right]
\right|
 \leq
 Cb_1^3
\left| \beta_{jn}
\frac{\mathds{1}\left(X_j\in\mathcal{X}_0\right)}
 {nb_0^d\widehat{g}_{jn}}
\right|.
 \label{Covzeta5}
\end{eqnarray}
For the  last term in (\ref{Covzeta3}), we have
\begin{eqnarray*}
\lefteqn{
 \esp_{in}
 \left[
 \Sigma_{jn}^2(x)
   K_1^{(2)}
 \left(\frac{\varepsilon_i-\epsilon}{b_1}\right)
 \right]
}
\\
&=& \frac{1}{(nb_0^d\widehat{g}_{jn})^2}
 \sum_{k=1\atop k\neq j}^n
\sum_{\ell=1\atop\ell\neq j}^n
 K_0\left(\frac{X_k-X_j}{b_0}\right)
 K_0\left(\frac{X_{\ell}-X_j}{b_0}\right)
 \esp_{in}
 \left[
 \varepsilon_k
 \varepsilon_{\ell}
  K_1^{(2)}
 \left(\frac{\varepsilon_i-\epsilon}{b_1}\right)
 \right]
 \\
 &=&
 \frac{1}{(nb_0^d\widehat{g}_{jn})^2}
 \sum_{k=1, k\neq j}^n
 K_0^2\left(\frac{X_k-X_j}{b_0}\right)
 \esp_{in}
 \left[
 \varepsilon_k^2
  K_1^{(2)}
 \left(\frac{\varepsilon_i-\epsilon}{b_1}\right)
 \right],
\end{eqnarray*}
with, using Lemma \ref{MomderK}-(\ref{MomderK2}),
\begin{eqnarray*}
\lefteqn
 {
 \left|
 \esp_{in}
 \left[
 \varepsilon_k^2
  K_1^{(2)}
 \left(\frac{\varepsilon_i-\epsilon}{b_1}\right)
 \right]
 \right|
}
\\
&\leq& \max
 \left\lbrace
\sup_{e\in\Rit} \left|
 \esp_{in}
 \left[
 \varepsilon^2
  K_1^{(2)}
 \left(\frac{\varepsilon-e}{b_1}\right)
 \right]
 \right|,
 \;
 \esp[\varepsilon^2]
\sup_{e\in\Rit}
 \left|
 \esp_{in}
 \left[
 K_1^{(2)}
 \left(\frac{\varepsilon-e}{b_1}\right)
 \right]
 \right|
\right\rbrace
\\
&\leq&
 C b_1^3.
\end{eqnarray*}
Therefore
$$
 \left|
 \esp_{in}
 \left[
 \Sigma_{jn}^2
   K_1^{(2)}
 \left(
 \frac{\varepsilon_i-\epsilon}{b_1}
 \right)
 \right]
 \right|
 \leq
\frac{Cb_1^3}{(nb_0^d\widehat{g}_{jn})^2} \sum_{k=1, k\neq j}^n
K_0^2\left(\frac{X_k-X_j}{b_0}\right).
$$
Substituting this bound, (\ref{Covzeta5}) and (\ref{Covzeta4}) in
(\ref{Covzeta3}), we obtain
$$
\left|
 \esp_{in}
 \left[
 Z_{jn}
  K_1^{(2)}
 \left(
 \frac{\varepsilon_i-\epsilon}{b_1}
 \right)
 \right]
\right|
 \leq
 Cb_1^3 M_n,
$$
where
$$
M_n
=
\sup_{1\leq j\leq n}
\left[
\beta_{jn}^2
+
\left|
\beta_{jn}
\frac{\mathds{1}\left(X_j\in\mathcal{X}_0\right)}
{nb_0^d\widehat{g}_{jn}}
\right|
+
\frac{1}{(nb_0^d\widehat{g}_{jn})^2}
\sum_{k=1, k\neq j}^n
K_0^2\left(\frac{X_k-X_j}{b_0}\right)
\right].
$$
Hence from (\ref{Prodzeta}), the Cauchy-Schwarz inequality, Lemma
\ref{BoundEspmchap} and Lemma \ref{MomderK}-(\ref{MomderK2}), we
deduce
\begin{eqnarray*}
\lefteqn{ \sum_{i=1}^n \sum_{j=1\atop j\neq i}^n
 \mathds{1}
 \biggl(
 \| X_i - X_j \|< C b_0
 \biggr)
\left|
\esp_n
\left[ \zeta_{in}\zeta_{jn} \right]
\right| }
\\
&\leq&
C M_nb_1^3
\sum_{i=1}^n \sum_{j=1\atop j\neq i}^n
 \mathds{1}
 \biggl(
 \| X_i - X_j \|< C b_0
 \biggr)
\esp_n
 \left|
 Z_{in} K_1^{(2)}
\left(\frac{\varepsilon_j-\epsilon}{b_1}\right)
\right|
\\
&\leq&
CM_nb_1^3
\sum_{i=1}^n \sum_{j=1\atop j\neq i}^n
 \mathds{1}
 \biggl(
 \| X_i - X_j \|< C b_0
 \biggr)
\esp_n^{1/2} \left[ Z_{in}^2 \right]
 \esp_n^{1/2}
 \left[
 K_1^{(2)}
 \left(
 \frac{\varepsilon_j-\epsilon}{b_1}
 \right)^2
 \right]
 \\
 &\leq&
 M_n b_1^3
 O_{\prob}
 \left(b_0^4+\frac{1}{nb_0^d}\right)
 (b_1)^{1/2}
\sum_{i=1}^n \sum_{j=1\atop j\neq i}^n
\biggl( \mathds{1}
\left(\|X_i-X_j\|\leq Cb_0\right)
\biggr).
\end{eqnarray*}
Moreover, (\ref{BetasumTBP}) and Lemma \ref{Estig} give, under
$(A_1)$, $(A_7)$ and $(A_9)$,
$$
M_n =
 O_{\prob}
 \left( b_0^4 + \frac{b_0^2}{nb_0^d}
  +
\frac{1}{nb_0^d}
 \right)
=
O_{\prob} \left(b_0^4+
\frac{1}{nb_0^d}\right).
$$
Finally, substituting this order in the bound above, and using
(\ref{Markov}), we arrive at
\begin{eqnarray*}
\sum_{i=1}^n \sum_{j=1\atop j\neq i}^n
 \mathds{1}
 \biggl(
 \| X_i - X_j \|< C b_0
 \biggr)
\esp_n \left[\zeta_{in}\zeta_{jn}\right]
=
 O_{\prob}
\left(n^2b_0^db_1^{7/2}\right)
\left(b_0^4+\frac{1}{nb_0^d}\right)^2.
\end{eqnarray*}
This proves (\ref{Covzeta1}), and then completes the proof of the
theorem. \eop

\chapter[ An integral nonparametric kernel estimator of the p.d.f of
 regression errors] {An integral
nonparametric kernel estimator of the probability density function
of regression errors}

\setcounter{subsection}{0}
\setcounter{equation}{0}
\renewcommand{\theequation}{\thesection.\arabic{equation}}

\renewcommand{\thefootnote}{\arabic{footnote}} \setcounter{footnote}{1}
\setlength{\baselineskip}{.26in}

\setcounter{subsection}{0} \setcounter{equation}{0}

{\bf Abstract:}
 This chapter is devoted to the nonparametric density
estimation of the regression error using an integral method.
 The difference
between the feasible estimator which uses the estimated regression
function and the unfeasible one using the true regression function
is investigated. An optimal choice of the first-step bandwidth
used for estimating this  regression function is proposed. We also
study the asymptotic normality of the feasible integral kernel
estimator and its rate-optimality.

\section{Introduction}
 Consider a sample
$(X,Y),(X_{1},Y_{1}),\ldots,(X_{n},Y_{n})$ of independent and
identically distributed (i.i.d) random variables, where Y is the
univariate dependent variable and the covariate X is of dimension
$d$. Let $m(\cdot)$ be the conditional expectation of $Y$ given
$X$ and let
 $\varepsilon$ be
the related regression error term, so that the regression error
model is
\begin{eqnarray}
Y_i = m(X_i)+\varepsilon_i,
\quad i=1,\ldots,n.
\label{RM}
\end{eqnarray}
The aim of this chapter is to estimate  the p.d.f of the
regression
 error
under the assumption that the covariate $X$ and  the regression
error
 $\varepsilon$ are independent. Indeed, under this assumption, we
 have
\begin{equation}
f(\epsilon)
 =
 f(\epsilon|x)
 =
 \varphi\left(m(x)+\epsilon|x\right).
 \label{Naive2}
\end{equation}
Hence, the  approach proposed here is based on a two-steps
procedure, which, in a first step, uses (\ref{Naive2}) and writes
$f(\epsilon)$ in the integral form
\begin{eqnarray*}
f(\epsilon)
 =
\int
 \mathds{1}
 \left(x\in\mathcal{X}\right)
 \varphi\left(\epsilon+m(x)\mid x\right)
  g(x)
  dx
 =
 \int
 \mathds{1}
 \left(x\in\mathcal{X}\right)
 \varphi\left(x, \epsilon+m(x)\right)
 dx.
\end{eqnarray*}
where $\mathcal{X}$ is the support of the p.d.f $g(\cdot)$ of $X$,
and $\varphi(\cdot,\cdot)$ the joint density of $(X,Y)$. This
formula suggests to estimate $f(\epsilon)$, in a second-step, by
$$
\widehat{f}_{2n}(\epsilon)
=
\int \mathds{1}
 \left(x\in\mathcal{X}\right)
\widehat{\varphi}_n
\left(x,\epsilon+\widehat{m}_n(x)\right)
dx,
$$
where $\widehat{\varphi}$ and  $\widehat{m}_n$  define
respectively some nonparametric estimators of $\varphi$ and $m$.
As in Chapter 2, a challenging
 issue is first to evaluate the impact of the estimated  regression function on
 the final estimator of $f(\cdot)$.
 Next, an optimal choice of the bandwidth used to
 estimate the residuals is proposed. Finally, we study the asymptotic
 normality of the estimator $\widehat{f}_{2n}(\epsilon)$
 and its rate-optimality.

\vskip 0.3cm
 The rest of this chapter is organized as follows.
Section 4.2 is devoted to presentation of ours estimators.
  Sections 4.3 and 4.4  group our
assumptions and  main results. The conclusion of this paper is
given in Section 4.5, while the proofs of our results are gathered
in section 4.6 and in two appendixes.

\section{Presentation of the estimators}

 In what follows, the bandwidths $b_{0}$
and $b_{1}$ are associated with $X$ and $h$ with $Y$, and $K_0$,
$K_1$ and $K_2$ represent some Kernels functions. Then for
$(x,y)\in\Rit^d\times\Rit$, the nonparametric estimators of
$\varphi(x,y)$ and $g(x)$ are respectively defined as
\begin{eqnarray*}
\widehat{\varphi }_{n}\left( x,y\right)
&=&
\frac{1}{nb_{1}^{d}h}
\sum_{i=1}^{n}
K_{1}\left(\frac{X_{i}-x}{b_{1}}\right)
K_{2}\left(\frac{ Y_{i}-y}{h}\right),
\\
\widehat{g}_{n}\left( x\right)
&=&
\frac{1}{nb_{0}^{d}}
\sum_{i=1}^{n}
K_0\left( \frac{X_{i}-x}{b_{0}}\right).
\end{eqnarray*}
The estimation of the  regression function $m(\cdot)$ is given by
the Nadaraya-Watson estimator (1964)
\begin{equation}
\widehat{m}_{n}(x)
=
 \frac{ \sum_{j=1}^{n}Y_{j}
K_0\left(\frac{X_{j}-x}{b_{0}}\right) }{\sum_{j=1}^{n}
K_0\left(\frac{X_{j}-x}{b_{0}}\right)}.
\label{mhat}
\end{equation}
 Since $Y=m(X)+\varepsilon$, we have
\begin{eqnarray*}
\prob
\left(
\varepsilon
\leq
\epsilon\mid X=x
 \right)
 =
 \prob
 \left(
 Y
 \leq
 \epsilon+m(x)\mid X=x
 \right).
\end{eqnarray*}
Then if $f$ represents the probability density function of
$\varepsilon$, and $\varphi$ the joint density of $(X,Y)$, it
follows
\begin{eqnarray}
f(\epsilon)
 =
\int
 \mathds{1}
 \left(x\in\mathcal{X}\right)
 \varphi \left(\epsilon+m(x)|x\right)
  g(x)
  dx
 =
 \int
 \mathds{1}
 \left(x\in\mathcal{X}\right)
 \varphi\left(x, \epsilon+m(x)\right)
 dx,
\label{fth}
\end{eqnarray}
where $\mathcal{X}$ is the support of the p.d.f $g$ of the
covariates. Therefore an estimator of $f(\epsilon)$ is the
so-called \lg Two-steps estimator\rg, defined as
\begin{equation}
\widehat{f}_{2n}(\epsilon)
 =
 \int
 \mathds{1}
 \left(x\in\mathcal{X}\right)
 \widehat{\varphi}_n
 \left(x,\epsilon+\widehat{m}_n(x)\right)
 dx.
\label{fnhat}
\end{equation}
 This  estimator is a feasible estimator
 in the sense that it does not depend on
any unknown
 quantity, as desirable in practice. This contrasts with the unfeasible
 ideal Kernel estimator
\begin{equation}
\widetilde{f}_{2n}(\epsilon)
 =
 \int
 \mathds{1}
 \left(x\in\mathcal{X}\right)
 \widehat{\varphi}_n
 \left(x,\epsilon+m(x)\right)
 dx,
 \label{f2n}
\end{equation}
which depends  in particular on the unknown regression function
$m(\cdot)$. It is however intuitively clear that
$\widehat{f}_{2n}(\epsilon)$ and
 $\widetilde{f}_{2n}(\epsilon)$ should be closed, as illustrated by the results of the
next section.

\section{Assumptions}
\noindent
 {$\bf(H_1)$}
 {\it
 The support $\mathcal{X}$ of $X$  is a known compact subset of
 $\Rit^d$,
  }
\medskip
\\{$\bf(H_2)$}
{\it the p.d.f. $g(\cdot)$ of the i.i.d. covariates $X, X_i$  has
continuous second order partial derivatives  over $\mathcal{X}$.
Moreover, there exists $\alpha>0$ such that $g(x)>\alpha$ for all
$x$ in the support $\mathcal{X}$,
 }
\medskip
\\{$\bf(H_3)$}
{\it
 the regression function $m(\cdot)$ has continuous second order partial
derivatives  over  $\mathcal{X}$,
 }
\medskip
\\{$\bf(H_4)$}
%\begin{quote}
{\it the i.i.d. centered error regression terms $\varepsilon,
 \varepsilon_i$'s, have finite 6th moments, and are independent of the
 covariates $X,X_i$'s,
}
%\end{quote}
\medskip
\\{$\bf(H_5)$}
{\it
 the  probability  density function $f$ of $\varepsilon$
 has   bounded  continuous second order
 derivatives over  $\Rit$, and satisfies,
for  $h_p (e) = e^p f(e)$,
 $\sup_{e\in\Rit}|h_p^{(k)}(e)|<\infty$,
 $p\in[0,6]$, $k\in[0,2]$,
 }
\medskip
\\{$\bf(H_6)$}
{\it the density $\varphi$ of $(X, Y)$ has bounded continuous
second  order partial derivatives over $\Rit^d\times\Rit$,
 }
\vskip 0.2cm\noindent
{$\bf(H_7)$} {
 \it the Kernel functions $K_0$ and $K_1$
are symmetric, continuous over $\Rit^d$ with support in $[-1/2,
1/2]^d$ and $\int\!K_0(z)dz = 1$, $\int\!K_1 (z)dz = 1$,
 }
\medskip
\\{$\bf(H_8)$}
{\it
 the Kernel function $K_{2}$ has a compact support,  is three times
 continuously differentiable over
 $\Rit$, and satisfies $\int\!K_2(v)dv = 1$,
$\int\!v K_2(v)dv = 0$ and $\int\!|v^p K_2^{(\ell)}(v)| dv<\infty$
 for $p,\ell$ in $[0,3]$,
 }
\medskip
\\{$\bf(H_{9})$}
{\it
 the bandwidth $b_0$ decreases to $0$ and satisfies
 $\ln(1/b_0)/\ln(\ln n)\rightarrow\infty$ and
  $b_0^{d}/(nb_0^{2d})^p= O(b_0^{2p})$, $p\in[0,6]$,
 when $n\rightarrow\infty$,
}
\vskip 0.2cm\noindent
{$\bf(H_{10})$}
 {\it
the bandwidths $b_1$ and $h$ decrease to $0$ and are such that
$nb_1^{2d}\rightarrow\infty$ and
$n^{(d+8)}h^{7(d+4)}\rightarrow\infty$ when $n\rightarrow\infty$.
 }

\vskip 0.3cm \noindent
 Assumptions $(H_2)$, $(H_3)$, $(H_5)$ and $(H_6)$ impose that all the
functions
 to be estimated nonparametrically have two bounded derivatives.
 Consequently the conditions
$\int\!v K_j (v) dv = 0$, $j=0,1,2$, as assumed in $(H_7)$ and
$(H_8)$, represent
 standard conditions ensuring  that the bias of the resulting
 nonparametric estimators (\ref{mhat})
and (\ref{f2n}) are respectively of order $b_0^2$ and $b_0^2+h^2$.
 Assumption $(H_4)$ states independence between
the regression error terms
 and the covariates, which is the main condition for (\ref{Naive2}) to
 hold.
The differentiability of $K_2$ imposed in $(H_8)$ is more specific
to our
 two-steps estimation method. Assumption $(H_8)$ is used to expand the
 two-steps Kernel estimator
$\widehat{f}_{2n}$ in (\ref{fnhat}) around the unfeasible one
$\widetilde{f}_{2n}$ from
 (\ref{f2n}), using the derivatives of $K_2$
 up to third order and the differences
 $\widehat{m}_{in}(x)- m(x)$, $i\in[1,n]$, where
 $\widehat{m}_{in}(x)$ is a leave-one out
 version of the Kernel regression estimator (\ref{mhat}),
\begin{equation}
\widehat{m}_{in}(x)
 =
 \frac{\sum_{j=1\atop j\neq i}^nY_j
K_0\left(\frac{X_j-x}{b_0}\right)}
{\sum_{j=1\atop j\neq i}^n
K_0\left(\frac{X_j-x}{b_0}\right)}.
 \label{leave}
\end{equation}
Assumption $(H_9)$ is a standard condition to obtain  uniform
convergence
 of the regression estimator $\widehat{m}_n$ in (\ref{mhat}) (see for
 instance Einmahl and Mason, 2005), and also gives a similar consistency
 result for the leave-one-out estimator $\widehat{m}_{in}$.
  Assumption $(H_{10})$ is needed in the study of
  the difference between the
 feasible estimator $\widehat{f}_{2n}$ and the unfeasible
  estimator $\widetilde{f}_{2n}$.

\section{Main results}
 Our first main result establishes the order of the
difference $\widehat{f}_{2n}(\epsilon)-f(\epsilon)$. This is given
in the following subsection. Next, we shall give the optimal
bandwidths needed to estimate $f(\epsilon)$. We conclude this
section by proposing an asymptotic normality of the estimator
$\widehat{f}_{2n}(\epsilon)$.

\renewcommand{\thesubsection}{\thesection.\arabic{subsection}}

\subsection{Pointwise weak consistency}
In this subsection we deal the order of the difference
$\widehat{f}_{2n}(\epsilon)-f(\epsilon)$. We show that for $n$
large enough, the estimator $\widehat{f}_{2n}(\epsilon)$ is very
close to the theoretical density $f(\epsilon)$, as illustrated by
the following result.

\begin{theoreme}
Suppose that Assumptions $(H_1)-(H_{10})$ hold. Then for $n$ large
enough, we have
$$
\widehat{f}_{2n}(\epsilon)-f(\epsilon)
 =
 O_{\prob}
\biggl(
 AMSE(b_1,h)+RT_n(b_0,b_1,h)
 \biggr)^{1/2},
$$
where
$$
AMSE(b_1,h)
 =
 \esp_n
 \left[
\left(\widetilde{f}_{2n}(\epsilon)-f(\epsilon)\right)^2
 \right]
 =
 O_{\prob}
\left(b_1^4 +h^4+ \frac{1}{nb_1}\right),
$$
and
\begin{eqnarray*}
RT_n(b_0,b_1,h)
&=&
b_{0}^{4}
+
\left(b_{0}^d\vee b_{1}^d\right)
\left[
\frac{1}{nb_1^dh^3}
\left(
b_0^4
+
\frac{1}{nb_0^{d}}
\right)
+
\frac{1}{nb_0^d}
\right]
\\
&&
+
\left(b_{0}^d\vee b_{1}^d\right)
\left[
\frac{1}{nb_1^dh^5}
 \left(b_{0}^4
 +
 \frac{1}{nb_0^d}
 \right)^2
 +
\frac{1}{n^2b_0^{2d}h^3}
\right]
 \\
 &&
 +
\frac{1}{h^2}
\left(
b_0^4
+
\frac{1}{nb_0^d}
\right)^3
+
\frac{b_{0}^d\vee b_{1}^d}{h^7}
\left(
b_0^{4}
+
\frac{1}{nb_0^d}
\right)^3.
 \end{eqnarray*}
\label{Thm1}
\end{theoreme}

\noindent The result of Theorem \ref{Thm1} is based on the
evaluation of the difference between  $\widehat{f}_{2n}(\epsilon)$
and $\widetilde{f}_{2n}(\epsilon)$. This evaluation  gives an
indication about the impact of the estimation of  $m(\cdot)$ on
the nonparametric estimation of the regression error density.

\subsection{Optimal first-step and second-step bandwidths
 for the pointwise weak consistency}

 Our next result deals with the choice of the optimal bandwidth
$b_0$ used in the nonparametric estimation of the p.d.f of the
regression error term. We have the following theorem.

\refstepcounter{theorem}
\begin{theorem}
 Suppose that Assumptions
$(H_1)-(H_{10})$ are satisfied, and assume $b_0=b_1$. Define
$$
b_0^*
=
 b_0^*(h)
 =
\arg\min_{b_0}
 RT_n (b_0, b_0, h),
$$
where the minimization is performed over bandwidth $b_0$
fulfilling $(H_9)$. Then the optimal bandwidth $b_0^*$ satisfies
$$
b_0^*
\asymp
\max
\left\lbrace
\left(\frac{1}{n^2h^3}\right)^{\frac{1}{d+4}}
,
\left(\frac{1}{n^3h^7}\right)^{\frac{1}{2d+4}}
\right\rbrace,
$$
and we have
$$
RT_n(b_0^*,b_0^*, h)
\asymp
\frac{1}{n}
+
\max
\left\lbrace
\left(\frac{1}{n^2h^3}\right)^{\frac{4}{d+4}}
 ,
\left(\frac{1}{n^3h^7}\right)^{\frac{4}{2d+4}}
\right\rbrace.
$$
\label{Optimb0}
\end{theorem}
The  next theorem gives the conditions for which the estimator
$\widehat{f}_{2n}(\epsilon)$ reaches the optimal rate $n^{-2/5}$
when $b_0$ takes the value $b_0^*$.
 We prove that for
$d\leq 2$, the bandwidth that minimizes the term $AMSE(b_0^*,
h)+RT_n(b_0^*, b_0^*, h)$ has the same order as $n^{-1/5}$,
leading to the optimal order $n^{-2/5}$ for the term
$\left(AMSE(b_0^*, h)+RT_n(b_0^*,b_0^*, h)\right)^{1/2}$.

\begin{theorem}
Assume that $(H_1)-(H_{10})$ hold and set
$$
h^*
 =
\arg\min_{h}
\biggl(AMSE(b_0^*,h)+RT_n(b_0^*,b_0^*,h)\biggr),
$$
where $b_0^*=b_0^*(h)$ is defined as in Theorem \ref{Optimb0}.
 Then
\begin{enumerate}
\item For $d\leq 2$, the optimal bandwidth $h^*$ satisfies
$$
h^*
\asymp
\left(\frac{1}{n}\right)^{\frac{1}{5}},
$$
and  we have
$$
\biggl(
AMSE(b_0^*, h^*)
 +
 RT_n\left(b_0^*, h^*, h^*\right)
\biggr)^{\frac{1}{2}}
\asymp
\left(\frac{1}{n}\right)^{\frac{2}{5}}.
$$
\item For $d\geq 3$, $h^*$ satisfies
$$
h^*
\asymp
\left(\frac{1}{n}\right)^{\frac{3}{2d+11}},
$$
and  we have
$$
\biggl(
 AMSE(b_0^*, h^*)
 +
 RT_n(b_0^*, b_0^*, h^*)
\biggr)^{\frac{1}{2}}
\asymp
\left(\frac{1}{n}\right)^{\frac{6}{2d+11}}.
$$
\end{enumerate}
\label{Optimh}
\end{theorem}

\vskip 0.3cm\noindent
Theorem \ref{Optimh} follows from Theorem
\ref{Optimb0}, which  reveals that for $b_1$ proportional to
$n^{-1/5}$,  the bandwidth $b_0^*$ has the same order as
$$
\max
\left\lbrace
\left(\frac{1}{n}\right)^{\frac{7}{5(d+4)}}
,
\left(\frac{1}{n}\right)^{\frac{8}{5(2d+4)}}
\right\rbrace
=
\left(\frac{1}{n}\right)^{\frac{8}{5(2d+4)}}.
$$
 For $d\leq 2$, this order of $b_0^*$ is less
than the one of the optimal bandwidth $\widehat{b}_0$ obtained
for pointwise or mean square estimation of $m(\cdot)$ using a
nonparametric Kernel estimator. In fact, as seen in Chapter 3, the
optimal bandwidth $\widehat{b}_0$ for estimating $m(\cdot)$ is
obtained by minimizing the order of  the risk function
$$
r_n(b_0)
=
\esp
\left[
\int
\mathds{1}
\left(x\in\mathcal{X}\right)
\left(\widehat{m}_n(x)-m(x)\right)^2
 \widehat{g}_n^2(x)w(x)
 dx
 \right],
$$
which has the same order as $b_0^4+\left(1/(nb_0^d)\right)$,
leading to the optimal bandwidth $\widehat{b}_0=n^{-1/(d+4)}$. For
d=1, the optimal order of $b_0^*$ is $n^{-(1/5)\times(4/3)}$ which
goes to 0 slightly faster than $n^{-1/5}$, the optimal order of
the bandwidth for the mean square nonparametric estimation of
$m(\cdot)$. For $d=2$, the optimal order of $b_0^*$ is $n^{-1/5}$.
Again this order goes to 0 faster than the order $n^{-1/6}$ of the
optimal bandwidth for the nonparametric estimation of the
regression function with two covariates. But for $d\geq 3$, we
note that the order of $b_0^*$ goes to $0$ slowly than
$\widehat{b}_0$. Hence these sitauations reveal that the optimal
$\widehat{m}_n(\cdot)$ for estimating $f(\cdot)$ should have a
lower bias and a higher variance than the optimal Kernel
regression estimator of $m(\cdot)$.  This situation is the same as
the one noticed in Wang, Cai, Brown and Levine (2008) for the
estimation of the conditional variance function in a
heteroscedastic regression model. However these authors do not
investigate the order of the optimal bandwidth to be used for
estimating the regression function in their heteroscedastic setup.
Hence, as in Chapter 3, we conclude that an estimator of
$m(\cdot)$ with smaller bias should be preferred in our framework,
compared to the case where the regression function $m(\cdot)$ is
the parameter of interest.

\subsection{Asymptotic normality}
The aim of this subsection is to propose an asymptotic normality
of the estimator $\widehat{f}_{2n}(\epsilon)$. We have the
following result.

\begin{theorem}
Suppose that $b_0=b_1$ and assume
$$
{(\rm\bf{H}_{11}):}
\quad nb_0^{d+4}=O(1),
\quad nb_0^4h=o(1),
\quad nb_0^dh^3\rightarrow\infty,
$$
 when $n\rightarrow\infty$. Then under $(H_1)-(H_{10})$, we have
\begin{eqnarray*}
\sqrt{nh}
 \left(
\widehat{f}_{2n}(\epsilon)
 -
 \overline{f}_{2n}(\epsilon)
 \right)
\stackrel{d}{\rightarrow}
\mathcal{N}
\left(
0,
f(\epsilon)
\int
 K_2^2(v)
 dv
\right),
\end{eqnarray*}
where
 \begin{eqnarray*}
\overline{f}_{2n}(\epsilon)
 &=&
f(\epsilon) + \frac{b_0^2}{2}
 \int
 \mathds{1}
 \left(x\in\cal{X}\right)
 \frac{\partial^2 \varphi(x,\epsilon+m(x))}
 {\partial^2 x}
 dx
 \int
 z K_1(z)z^{\top}
 dz
 \\
 &&
 +
 \frac{h^2}{2}
 \int
 \mathds{1}
 \left(x\in\cal{X}\right)
 \frac{\partial^2 \varphi(x,\epsilon+m(x))}
 {\partial^2 y}
 dx
 \int
 v^2 K_2(v)
 dv
+
o\left(b_0^2+h^2\right).
\end{eqnarray*}
\label{Normalite2}
\end{theorem}

\noindent As seen in the comments of Theorem \ref{normalite} in
Chapter 3, we  can check that  for $d=1$, $h=h^*$ and
$b_1=b_0=b_0^*$, the conditions of  Assumption $(\rm\bf{H}_{11})$
are realizable with the bandwidths $b_0^*$ and $h^*$. But with
these bandwidths, the last constraint of $(\rm\bf{H}_{11})$ is not
satisfied for $d=2$, since for $b_0=b_0^*$ and $h=h^*$,
$nb_0^dh^3$ is bounded when $n$ goes to infinity.

\newpage
\section{ Conclusion}
 In this chapter, we investigated  the nonparametric Kernel estimation of the
p.d.f of the regression error using an integral method. The
difference between the feasible estimator which uses the estimated
regression function and the unfeasible one using the theoretical
regression function is studied. An optimal choice of the
first-step bandwidth used to estimate the regression function is
also established. Again, an asymptotic normality of the feasible
Kernel estimator and its rate-optimality are proposed. As in
Chapter 2, the contributions of the present chapter is the
analysis of the influence of the estimated regression function on
the regression errors p.d.f. Kernel estimator.

\vskip 0.3cm The strategy used here strategy is to use an approach
based on a two-steps procedure which, in a first step, integrates
a conditional p.d.f as  given in (\ref{fth}).  In a second step,
we build the Kernel estimator of $f(\epsilon)$ by estimating
nonparametrically the unknown functions in the integral terms of
(\ref{fth}). If this strategy can avoid
 the curse of dimensionality, a main aspect of our setup is
 to evaluate the impact of the estimation of
$m(\cdot)$ on the final integral Kernel estimator of $f(\cdot)$ in
the first nonparametric step, and to determine the optimal choice
of the first-step bandwidth $b_0$. For a such choice of $b_0$, our
results suggests that the optimal bandwidth to be used  should be
smaller than the optimal bandwidth for the mean square estimation
of $m(\cdot)$. This mean that the best choice for $b_0$ is the one
such that the estimator $\widehat{m}_n(\cdot)$ of the regression
has a lower bias and a higher variance than the optimal Kernel
regression of the estimation setup. With this choice of $b_0$, we
show that for $d\leq 2$, the  estimator
$\widehat{f}_{2n}(\epsilon)$ of $f(\epsilon)$ can reach the
optimal rate $n^{-2/5}$, which
 corresponds exactly to the  rate reached for the
Kernel density estimator of an univariate  variable. This reveals
that for $d\leq 2$, the integral Kernel estimator
$\widehat{f}_{2n}(\epsilon)$  is not affected by the curse of
dimensionality, since there is not a negative influence caused by
the  estimation of the optimal first-step bandwidth $b_0^{*}$.

\newpage
\setcounter{subsection}{0} \setcounter{equation}{0}
\renewcommand{\theequation}{\thesection.\arabic{equation}}

\section{Proofs section}
\subsection*{Proof of Theorem \ref{Thm1}}
The proof is a consequence of the two followings lemmas.
\begin{lemma}
Under $(H_1)-(H_{10})$, we have, when $n$ goes to infinity,
\begin{eqnarray*}
 \widehat{f}_{2n}(\epsilon)-\widetilde{f}_{2n}(\epsilon)
 &=&
O_{\prob}
\left[
b_{0}^{4}
+
\left(b_{0}^d\vee b_{1}^d\right)
\left(
\frac{1}{nb_1^dh^3}
\left(
b_0^4
+
\frac{1}{nb_0^{d}}
\right)
+
\frac{1}{nb_0^d}
\right)
\right]^{1/2}
\\
&&
+
O_{\prob}
\left[
\left(b_{0}^d\vee b_{1}^d\right)
\left(
\frac{1}{nb_1^dh^5}
 \left(b_{0}^4
 +
 \frac{1}{nb_0^d}
 \right)^2
 +
\frac{1}{n^2b_0^{2d}h^3}
\right)
 \right]^{1/2}
 \\
 &&
 +
O_{\prob}
\left[
\frac{1}{h^2}
\left(
b_0^4
+
\frac{1}{nb_0^d}
\right)^3
+
\frac{b_{0}^d\vee b_{1}^d}{h^7}
\left(
b_0^{4}
+
\frac{1}{nb_0^d}
\right)^3
 \right]^{1/2}.
\end{eqnarray*}
\label{Ordrefnf}
\end{lemma}

\begin{lemma}
If $(H_1)-(H_{10})$ hold, then
\begin{eqnarray*}
\widetilde{f}_{2n}(\epsilon)-f(\epsilon)
 =
 O_{\prob}
\left(b_1^4+h^4+\frac{1}{n h}\right)^{1/2}.
\end{eqnarray*}
\label{Bias}
\end{lemma}
Let now turn to the proof of  Theorem \ref{Thm1}. Using  Lemmas
\ref{Bias} and \ref{Ordrefnf}, we have
\begin{eqnarray*}
\lefteqn{
 \widehat{f}_{2n}(\epsilon)-f(\epsilon)
 =
 \left(\widetilde{f}_{2n}(\epsilon)-f(\epsilon)\right)
 +
 \widehat{f}_{2n}(\epsilon)-\widetilde{f}_{2n}(\epsilon)
 }
 \\
 &=&
O_{\prob}
\left(b_1^4+h^4+\frac{1}{n h}\right)^{1/2}
+
O_{\prob}
\left[
b_{0}^{4}
+
\left(b_{0}^d\vee b_{1}^d\right)
\left(
\frac{1}{nb_1^dh^3}
\left(
b_0^4
+
\frac{1}{nb_0^{d}}
\right)
+
\frac{1}{nb_0^d}
\right)
\right]^{1/2}
\\
&&
+
O_{\prob}
\left[
\left(b_{0}^d\vee b_{1}^d\right)
\left(
\frac{1}{nb_1^dh^5}
 \left(b_{0}^4
 +
 \frac{1}{nb_0^d}
 \right)^2
 +
\frac{1}{n^2b_0^{2d}h^3}
\right)
+\frac{1}{h^2}
\left(
b_0^4
+
\frac{1}{nb_0^d}
\right)^3
+
\frac{b_{0}^d\vee b_{1}^d}{h^7}
\left(
b_0^{4}
+
\frac{1}{nb_0^d}
\right)^3
\right]^{1/2},
\end{eqnarray*}
which yields the result of the Theorem.\eop

\vskip 0.3cm
 We now prove Lemmas \ref{Ordrefnf} and \ref{Bias}.

\subsection*{Proof of Lemma \ref{Ordrefnf}}

Let us introduce additional notations. Let $\widehat{m}_{in}(x)$
be as in (\ref{mchapi}) and define
\begin{eqnarray*}
S_{n}\left( x\right)
&=&
\frac{1}{nb_{1}^{d}h^{2}}
\sum_{i=1}^{n}
\left( \widehat{m}_{in}\left( x\right)
 -m\left(
x\right) \right)
K_{1}\left( \frac{ X_{i}-x}{b_{1}}\right)
 K_2^{(1)}
 \left(\frac{Y_i-\epsilon-m(x)}{h}
 \right),
\\
T_{n}\left( x\right)
&=&
\frac{1}{nb_{1}^{d}h^{3}}
\sum_{i=1}^{n}
\left(
\widehat{m}_{in}\left( x\right)
-
m\left( x\right)
\right)^{2}
K_{1}\left(\frac{X_{i}-x}{b_{1}}\right)
K_{2}^{(2)}\left(
\frac{Y_{i}-\epsilon-m(x)}{h} \right).
\end{eqnarray*}
 The proof of Lemma
\ref{Ordrefnf} is based on the following results.
\begin{lemma}
Define
$$
S_n
=
\int
\mathds{1}
\left(x\in\mathcal{X}\right)
S_n(x)
dx,
\quad
T_n
=
\int
\mathds{1}
\left(x\in\mathcal{X}\right)
T_n(x)
 dx.
$$
Then under  $(H_1)-(H_{10})$,
 we have
\begin{eqnarray*}
S_n
&=&
O_{\prob}
\left[
b_{0}^{4}
+
\left(b_{0}^d\vee b_{1}^d\right)
\left(
\frac{1}{nb_1^dh^3}
\left(
b_0^4
+
\frac{1}{nb_0^{d}}
\right)
+
\frac{1}{nb_0^d}
\right)
\right]^{1/2},
\\
T_n
&=&
O_{\prob}
\left[
\left(
b_{0}^{4}
+
\frac{1}{nb_0^{d}}
\right)^2
+
\left(b_{0}^d\vee b_{1}^d\right)
\left(
\frac{1}{nb_1^dh^5}
 \left(b_{0}^4
 +
 \frac{1}{nb_0^d}
 \right)^2
 +
 \frac{b_0^4}{nb_0^d}
 +
\frac{1}{n^2b_0^{2d}h^3}
\right)
 \right]^{1/2}.
\end{eqnarray*}
\label{OrdreST}
\end{lemma}

\begin{lemma}
Define
\begin{eqnarray*}
  R_n(x)
 =
 \frac{1}{nb_1^{d}h^4}
 \sum_{i=1}^n
 \left(
 \widehat{m}_{in}(x)-m(x)
 \right)^3
 K_1\left(
 \frac{X_i-x}{b_1}
 \right)
 \int_{0}^{1}
 (1-u)^2
 K_2^{(3)}
 \left(
 \frac{
 Y_i-\theta_{in}(x,u)}{h}
 \right)
 du,
\end{eqnarray*}
where $\theta_{in}(x,u)=\epsilon-m(x)
 -u\left(\widehat{m}_{in}(x)-m(x)\right)$, and set
 $$
R_n
=
\int
\mathds{1}
\left(x\in\mathcal{X}\right)
R_n(x)
dx.
$$
 If $(H_1)-(H_{10})$ hold, then
\begin{eqnarray*}
 R_n
 =
 O_{\prob}
\left[
\frac{1}{h^2}
\left(
b_0^4
+
\frac{1}{nb_0^d}
\right)^3
+
 \frac{b_{0}^d\vee b_{1}^d}{h^7}
\left(
b_0^{4}
+
\frac{1}{nb_0^d}
\right)^3
\right]^{1/2}.
\end{eqnarray*}
\label{OrdreR}
\end{lemma}

\begin{lemma}
Set
$$
P_n(x) = \frac{1}{nb_1^{d}h^2}
 \sum_{i=1}^n
 \left(
 \widehat{m}_n(x)-\widehat{m}_{in}(x)
 \right)
 K_1\left(
 \frac{X_i-x}{b_1}
 \right)
 \int_0^1
  K_2^{(1)}
 \left(\frac{Y_i-\widehat{\theta}_{in}(x,t)}{h}
 \right)
 dt,
 $$
  where $\widehat{\theta}_{in}(x,t)=\epsilon+\widehat{m}_{in}(x)
  +t\left(\widehat{m}_n(x)-\widehat{m}_{in}(x)\right)$, and define
 $$
 P_n
 =
\int
\mathds{1}
\left(x\in\mathcal{X}\right)
P_n(x)
dx.
$$
Then under $(H_1)-(H_{10})$, we have
 $$
 P_n
 =
 O_{\prob}
 \left(
 \frac{1}{n^2b_0^{2d}}
 +
 \frac{
 b_{0}^d\vee b_{1}^d}
 {n^2b_0^{2d} h^3}
 \right)^{1/2}.
$$
\label{OrdreP}
\end{lemma}

\noindent The proofs of these Lemmas are stated in Appendix B.

\newpage
 Let us now return to the proof of Lemma
\ref{Ordrefnf}. Observe that
\begin{eqnarray}
\nonumber
\lefteqn{
\widehat{\varphi}_n
\left(x,\epsilon+\widehat{m}_n(x)\right)
-
\widehat{\varphi}_n
\left(x,\epsilon+m(x)\right)
}
\\
 &=&
 \frac{1}{nb_1^{d}h}
 \sum_{i=1}^n
 K_1\left(
 \frac{X_i-x}{b_1}
 \right)
 \left[
 K_2\left(
 \frac{Y_i-\epsilon-\widehat{m}_n(x)}{h}
 \right)
 -
 K_2\left(
 \frac{Y_i-\epsilon-m(x)}{h}
 \right)
 \right],
\label{phichap}
\end{eqnarray}
where
\begin{eqnarray}
\nonumber
\lefteqn{
 K_2\left(
 \frac{Y_i-\epsilon-\widehat{m}_n(x)}{h}
 \right)
 -
 K_2\left(
 \frac{Y_i-\epsilon-m(x)}{h}
 \right)
 }
 \\\nonumber
 &=&
 K_2\left(
 \frac{Y_i-\epsilon-\widehat{m}_{in}(x)}{h}
 \right)
 -
  K_2
 \left(
 \frac{Y_i-\epsilon-m(x)}{h}
 \right)
 \\
 &&
 +
 \left[
 K_2\left(
 \frac{Y_i-\epsilon-\widehat{m}_n(x)}{h}
 \right)
 -
  K_2
 \left(
 \frac{Y_i-\epsilon-\widehat{m}_{in}(x)}{h}
 \right)
  \right].
 \label{phichap1}
\end{eqnarray}
 Since $K_2$ is three times continuously differentiable under $(H_8)$,
the Taylor's theorem with the integral remainder gives
\begin{eqnarray}
\nonumber
 \lefteqn{
 K_2\left(
 \frac{Y_i-\epsilon-\widehat{m}_{in}(x)}{h}
 \right)
 -
 K_2\left(
 \frac{Y_i-\epsilon-m(x)}{h}
 \right)
 }
 \\\nonumber
 &=&
 -
 \frac{1}{h}
 \left(
 \widehat{m}_{in}(x)-m(x)
 \right)
 K_2^{(1)}
 \left(\frac{Y_i-\epsilon-m(x)}{h}
 \right)
 \\\nonumber
 &&
 +\;
 \frac{1}{2h^2}
\left(\widehat{m}_{in}(x)-m(x)\right)^2
 K_2^{(2)}
 \left(\frac{Y_i-\epsilon-m(x)}{h}\right)
 \\\nonumber
 &&
 -\;
 \frac{1}{2h^3}
 \left(\widehat{m}_{in}(x)-m(x)\right)^3
 \int_{0}^{1}
 (1-u)^2
 K_2^{(3)}
 \left(
 \frac{
 Y_i-\epsilon-m(x)
 -u\left(\widehat{m}_{in}(x)-m(x)\right)}{h}
 \right)
 du.
 \\
 \label{TaylorK1}
\end{eqnarray}
Again, under $(H_8)$, we have
\begin{eqnarray*}
\lefteqn{
 K_2
 \left(
 \frac{Y_i-\epsilon-\widehat{m}_n(x)}{h}
 \right)
 -
 K_2
 \left(
 \frac{Y_i-\epsilon-\widehat{m}_{in}(x)}{h}
 \right)
 }
 \\
 &=&
 -\frac{1}{h}
 \left(
 \widehat{m}_n(x)-\widehat{m}_{in}(x)
 \right)
 \int_{0}^1
 K_2^{(1)}
 \left(\frac{Y_i-\epsilon-\widehat{m}_{in}(x)
 -t\left(\widehat{m}_n(x)-\widehat{m}_{in}(x)\right)}{h}
 \right)
 dt.
\end{eqnarray*}
Hence defining $S_n(x)$, $T_n(x)$, $R_n(x)$ and $P_n(x)$
respectively as in Lemmas \ref{OrdreST}, \ref{OrdreR} and
\ref{OrdreP},  the equality above, (\ref{TaylorK1}),
(\ref{phichap1})
 and (\ref{phichap}) give
 \begin{eqnarray*}
 \widehat{\varphi}_n
 \left(x,\epsilon+\widehat{m}_n(x)\right)
 -
 \widehat{\varphi}_n
 \left(x,\epsilon+m(x)\right)
 =
 -
 S_n(x)
 +
 \frac{T_n(x)}{2}
 -
 \frac{R_n(x)}{2}
 -
 P_n(x),
 \end{eqnarray*}
so that
\begin{eqnarray*}
\lefteqn{
 \widehat{f}_{2n}(\epsilon)-\widetilde{f}_{2n}(\epsilon)
 =
 -
 S_n
 +
 \frac{T_n}{2}
 -
 \frac{R_n}{2}
 -
 P_n
 }
 \\
 &=&
O_{\prob}
\left[
b_{0}^{4}
+
\left(b_{0}^d\vee b_{1}^d\right)
\left(
\frac{1}{nb_1^dh^3}
\left(
b_0^4
+
\frac{1}{nb_0^{d}}
\right)
+
\frac{1}{nb_0^d}
\right)
\right]^{1/2}
\\
&&
+
O_{\prob}
\left[
\left(
b_{0}^{4}
+
\frac{1}{nb_0^{d}}
\right)^2
+
\left(b_{0}^d\vee b_{1}^d\right)
\left(
\frac{1}{nb_1^dh^5}
 \left(b_{0}^4
 +
 \frac{1}{nb_0^d}
 \right)^2
 +
 \frac{b_0^4}{nb_0^d}
 +
\frac{1}{n^2b_0^{2d}h^3}
\right)
 \right]^{1/2}
 \\
 &&
 +
O_{\prob}
\left[
\frac{1}{h^2}
\left(
b_0^4
+
\frac{1}{nb_0^d}
\right)^3
+
 \frac{b_{0}^d\vee b_{1}^d}{h^7}
\left(
b_0^{4}
+
\frac{1}{nb_0^d}
\right)^3
+
\frac{1}{n^2b_0^{2d}}
 +
 \frac{
 b_{0}^d\vee b_{1}^d}
 {n^2b_0^{2d} h^3}
 \right]^{1/2}.
\end{eqnarray*}
Moreover,  since under $(H_9)$ $b_0$ goes to $0$ and that
$b_0^{d}/(n^pb_0^{2dp})= O(b_0^{2p})$, this gives for $p=1$,
\begin{eqnarray*}
\frac{b_0^4}{nb_0^d}
=
O\left(\frac{1}{nb_0^d}\right),
 \;
\frac{1}{n^2b_0^{2d}}
=
 O\left(b_0^4\right),
\;
\left(
b_{0}^{4}
+
\frac{1}{nb_0^{d}}
\right)^2
=
O\left(b_0^4\right).
\end{eqnarray*}
Hence it  follows that
\begin{eqnarray*}
 \widehat{f}_{2n}(\epsilon)-\widetilde{f}_{2n}(\epsilon)
 &=&
O_{\prob}
\left[
b_{0}^{4}
+
\left(b_{0}^d\vee b_{1}^d\right)
\left(
\frac{1}{nb_1^dh^3}
\left(
b_0^4
+
\frac{1}{nb_0^{d}}
\right)
+
\frac{1}{nb_0^d}
\right)
\right]^{1/2}
\\
&&
+
O_{\prob}
\left[
\left(b_{0}^d\vee b_{1}^d\right)
\left(
\frac{1}{nb_1^dh^5}
 \left(b_{0}^4
 +
 \frac{1}{nb_0^d}
 \right)^2
 +
\frac{1}{n^2b_0^{2d}h^3}
\right)
 \right]^{1/2}
 \\
 &&
 +
O_{\prob}
\left[
\frac{1}{h^2}
\left(
b_0^4
+
\frac{1}{nb_0^d}
\right)^3
+
\frac{b_{0}^d\vee b_{1}^d}{h^7}
\left(
b_0^{4}
+
\frac{1}{nb_0^d}
\right)^3
 \right]^{1/2},
\end{eqnarray*}
which ends the proof of the Lemma. \eop

\subsection*{Proof of Lemma \ref{Bias}}

Observe that
\begin{eqnarray}
\widetilde{f}_{2n}(\epsilon)-f(\epsilon)
 =
 \left(
 \widetilde{f}_{2n}(\epsilon)
 -
 \esp \widetilde{f}_{2n}(\epsilon)
 \right)
 +
 \left(
 \esp \widetilde{f}_{2n}(\epsilon)
 -
 f(\epsilon)
 \right).
 \label{Bias1}
\end{eqnarray}
For the first term in (\ref{Bias1}), the independence of the
$(X_i, Y_i)$'s gives
 \begin{eqnarray*}
 \lefteqn
  {
 \esp
 \left[
 \left(\widetilde{f}_{2n}(\epsilon)-\esp \widetilde{f}_{2n}(\epsilon)\right)^2
 \right]
 =
 \Var\left(\widetilde{f}_{2n}(\epsilon)\right)
  }
 \\
 &=&
 \Var
 \left[
  \frac{1}{nb_1^{d}h}
  \sum_{i=1}^n
 \int
 \mathds{1}\left(x\in\cal{X}\right)
 K_1
 \left(
 \frac{X_i-x}{b_1}
 \right)
 K_2
 \left(
 \frac{Y_i-\epsilon-m(x)}{h}
 \right)
 dx
 \right]
\\
&=&
 \frac{1}{(nb_1^{d}h)^2}
 \sum_{i=1}^n
 \Var
\left[
\int
\mathds{1} \left(x\in\cal{X}\right)
 K_1
 \left(
 \frac{X_i-x}{b_1}
 \right)
 K_2
 \left(
 \frac{Y_i-\epsilon-m(x)}{h}
 \right)
 dx
 \right]
\\
&\leq&
\frac{1}{(nb_1^{d}h)^2}
 \sum_{i=1}^n
  \esp
\left[
\int
 \mathds{1} \left(x\in\cal{X}\right)
 K_1
 \left(
 \frac{X_i-x}{b_1}
 \right)
 K_2
 \left(
 \frac{Y_i-\epsilon-m(x)}{h}
  \right)
  dx
 \right]^2.
\end{eqnarray*}
Moreover, note that by $(H_1)$, $(H_3)$ and  $(H_7)-(H_8)$, the
changes of variables $x=x_1+hz_1$,
$y_1=\epsilon+m(x_1+b_1z_1)+hv_1$ and the Cauchy-Schwarz
inequality give, since $\varphi(\cdot,\cdot)$ is bounded under
Assumption $(H_6)$,
\begin{eqnarray*}
\lefteqn
 {
 \sum_{i=1}^n
 \esp
\left[
\int
\mathds{1} \left(x\in\cal{X}\right)
 K_1
 \left(
 \frac{X_i-x}{b_1}
 \right)
 K_2
 \left(
 \frac{Y_i-\epsilon-m(x)}{h}
  \right)
  dx
  \right]^2
 }
 \\
 &=&
 n\int_{\Rit^d}
 dx_1
 \int_{\Rit}
 \left[
 \int
 \mathds{1}
 \left(x\in\cal{X}\right)
 K_1
 \left(
 \frac{x_1-x}{b_1}
 \right)
 K_2
 \left(
 \frac{y_1-\epsilon-m(x)}{h}
 \right)
 dx
 \right]^2
 \varphi\left(x_1, y_1\right)
  dy_1
 \\
 &=&
 n\int_{\Rit^d}
 dx_1
 \int_{\Rit}
 \left[
 b_1^d
 \int
 \mathds{1}
 \left(x_1+b_1z_1\in\cal{X}\right)
 K_1(z_1)
 K_2\left(
 \frac{y_1-\epsilon-m(x_1+b_1z_1)}{h}
 \right)
 dz_1
 \right]^2
 \varphi
 \left(x_1, y_1\right)
 dy_1
\\
 &\leq&
 Cnb_1^{2d}h
 \int_{\Rit^d}
 dz_1
 K_1^2(z_1)
 \int
 \mathds{1}
 \left(x_1+b_1z_1\in\cal{X}\right)
 dx_1
 \int_{\Rit}
 K_2^2(v_1)
 dv_1.
\end{eqnarray*}
Hence from the  two bounds above and the Tchebychev inequality, we
deduce
\begin{eqnarray}
\widetilde{f}_{2n}(\epsilon)-\esp \widetilde{f}_{2n}(\epsilon)
 =
 O_{\prob}
\left(\frac{1}{n h}\right)^{1/2}.
 \label{Bias2}
\end{eqnarray}
We now compute the order of the second term in (\ref{Bias1}).
Observe that
\begin{eqnarray}
 \nonumber
 \lefteqn{
 \esp \widetilde{f}_{2n}(\epsilon)
 =
 \esp
  \left[
  \frac{1}{nb_1^{d}h}
 \sum_{i=1}^n
 \int
 \mathds{1}
 \left(x\in\cal{X}\right)
 K_1\left(
 \frac{X_i-x}{b_1}
 \right)
 K_2\left(
 \frac{Y_i-\epsilon-m(x)}{h}
  \right)
  dx
  \right]
  }
  \\\nonumber
  &=&
 \frac{n}{nb_1^{d}h}
 \int
 \mathds{1}
 \left(x\in\cal{X}\right)
 \esp
 \left[
 K_1\left(
 \frac{X_1-x}{b_1}
 \right)
 K_2\left(
 \frac{Y_1-\epsilon-m(x)}{h}
  \right)
  \right]
  dx
  \\
  &=&
  \int
  \mathds{1}
 \left(x\in\cal{X}\right)
 \left[
  \int_{\Rit^d}
 dz_1
 \int_{\Rit}
 K_1(z_1) K_2(v_1)
 \varphi\left(x+b_1z_1,\epsilon+m(x)+hv_1\right)
  dv_1
\right]
 dx.
\label{Bias3}
\end{eqnarray}
By $(H_6)$, a second-order Taylor expansion yields, for  $z_1$
 and $v_1$ in the supports of $K_1$ and $K_2$, and $h$ and $b_1$ small
 enough,
\begin{eqnarray*}
 \varphi \left(x+b_1z_1,\epsilon+m(x)+hv_1\right)
 &=&
\varphi (x, \epsilon+m(x))
 +
b_1\frac{\partial\varphi(x, \epsilon+m(x))} {\partial x}
z_1^{\top}
+
h\frac{\partial\varphi (x, \epsilon+m(x))}
{\partial y} v_1
\\
&&
 +
\frac{b_1^2}{2}
 z_1\frac{\partial^2
 \varphi(x+\theta b_1z_1,\epsilon+m(x)+\theta b_1 v_1)}
 {\partial^2x} z_1^{\top}
 \\
 &&
+
b_1 h v_1
\frac{\partial^2 \varphi (x+\theta b_1 z_1,\epsilon+m(x)+\theta
b_1 v_1)}{\partial x\partial y} z_1^{\top}
\\
&&
 +
 \frac{h^2}{2}
 \frac{\partial^2 \varphi (x+\theta b_1 z_1,
\epsilon+m(x)+\theta b_1 v_1)} {\partial^2 y} v_1^2,
%\label{taylor}
\end{eqnarray*}
for some $\theta = \theta (x,\epsilon,b_1 z_1, h v_1)$ in $[0,1]$.
This gives, since $\int\! K_1 (z) dz = \int\! K_2 (v) dv =1$,
$\int\!z K_1(z) dz$ and that $\int\! v K_2 (v) dv$ vanishes under
$(H_7)-(H_8)$,
\begin{eqnarray*}
\lefteqn{
 \int_{\Rit^d}
 dz_1
 \int_{\Rit}
 K_1(z_1) K_2(v_1)
 \varphi\left(x+b_1z_1, \epsilon+m(x)+hv_1\right)
  dv_1
 }
\\
&&
-
\varphi(x, \epsilon+m(x))
-
\frac{b_1^2}{2} \frac{\partial^2
\varphi(x, \epsilon+m(x))} {\partial^2 x}
\int
z K_0(z)z^{\top}
dz
-
\frac{h^2}{2} \frac{\partial^2 \varphi (x, \epsilon+m(x))}
{\partial^2 y}
\int v^2 K_1 (v) dv
\\
& =&
\frac{b_1^2}{2}
\int \int z \left( \frac{\partial^2 \varphi
(x+\theta h_0 z, \epsilon+m(x)+\theta h_1 v)} {\partial^2 x}
-
\frac{\partial^2 \varphi (x, \epsilon+m(x))} {\partial^2 x}
\right) z^{\top}
K_1 (z) K_2 (v) dz dv
\\
&&
+
 b_1 h
\int \int v \left( \frac{\partial^2 \varphi(x+\theta b_1
z,\epsilon+m(x)+\theta b_1 v)} {\partial x\partial y}
-
\frac{\partial^2 \varphi (x, \epsilon+m(x))}
{\partial x \partial
y} \right) z^{\top}
K_1 (z) K_2 (v) dz dv
\\
&&
+
\frac{h^2}{2}
\int \int \left( \frac{\partial^2 \varphi
(x+\theta b_1 z, \epsilon+m(x)+\theta b_1 v)} {\partial^2 y}
-
\frac{\partial^2 \varphi (x, \epsilon+m(x))} {\partial^2 y}
\right)
v^2 K_1 (z) K_2 (v) dz dv.
\end{eqnarray*}
Hence by the Lebesgue Dominated Convergence Theorem, we have,
using (\ref{Bias3}) and (\ref{fth}),
\begin{eqnarray}
\nonumber \lefteqn{ \esp \widetilde{f}_{2n}(\epsilon) -
 \frac{b_1^2}{2}
 \int
 \mathds{1}
 \left(x\in\cal{X}\right)
 \frac{\partial^2 \varphi(x, \epsilon+m(x))}
 {\partial^2 x}
 dx
 \int
 z K_1(z)z^{\top}
  dz
}
\\\nonumber
&&
-
\frac{h^2}{2}
\int
 \mathds{1}
 \left(x\in\cal{X}\right)
\frac{\partial^2 \varphi (x,\epsilon+m(x))}
{\partial^2 y}
dx
\int
v^2 K_2(v) dv
\\\nonumber
&=&
 \int
 \mathds{1}
 \left(x\in\cal{X}\right)
 \varphi\left(x, \epsilon+m(x)\right)
 dx
 +
 o\left(b_1^2+h^2\right)
\\
&=&
 f(\epsilon)+ o\left(b_1^2+h^2\right),
 \label{Bias4}
\end{eqnarray}
so that
$$
\esp \widetilde{f}_{2n}(\epsilon)-f(\epsilon)
 =
 O\left(b_1^2+h^2\right).
$$
Finally, combining this result with (\ref{Bias2}) and
(\ref{Bias1}), we arrive at
$$
\widetilde{f}_{2n}(\epsilon)-f(\epsilon)
 =
 O_{\prob}
\left(b_1^4+h^4+\frac{1}{n h}\right)^{1/2}. \eop
$$

\subsection*{Proof of Theorem \ref{Optimb0}}
Observe that
\begin{eqnarray*}
RT_n(b_0,b_0,h)
&=&
b_{0}^{4}
+
\frac{1}{nh^3}
\left(
b_0^4
+
\frac{1}{nb_0^{d}}
\right)
+
\frac{1}{n}
+
\frac{1}{nh^5}
 \left(b_{0}^4
 +
 \frac{1}{nb_0^d}
 \right)^2
 \\
 &&
 +
\frac{1}{n^2b_0^{d}h^3}
 +
\frac{1}{h^2}
\left(
b_0^4
+
\frac{1}{nb_0^d}
\right)^3
+
\frac{b_{0}^d}{h^7}
\left(
b_0^{4}
+
\frac{1}{nb_0^d}
\right)^3.
 \end{eqnarray*}
and note that
$$
\left(\frac{1}{n^2h^3}\right)^{\frac{1}{d+4}}
 =
  \max
\left\lbrace
 \left(\frac{1}{n^2h^3}\right)^{\frac{1}{d+4}}
,
\left(\frac{1}{n^3h^7}\right)^{\frac{1}{2d+4}}
\right\rbrace
$$
if and only if $n^{4-d}h^{d+16}\rightarrow\infty$. To find the
order of $b_0^*$, we shall deal with the cases
$nb_0^{d+4}\rightarrow\infty$ and
 $nb_0^{d+4}=O(1)$.
\vskip 0.1cm\noindent First assume that
$nb_0^{d+4}\rightarrow\infty$. More precisely, we  suppose that
$b_0$ is in $\left[(u_n/n)^{1/(d+4)},\infty\right)$, where $u_n
\rightarrow\infty$. Since $1/(nb_0^d) = O(b_0^4)$ for all these
$b_0$, we have
$$
\left(b_0^4+\frac{1}{nb_0^d}\right)
\asymp
\left(b_0^4\right),
\quad
 \frac{1}{n^2b_0^dh^3}
 =
 O\left(\frac{b_0^4}{nh^3}\right).
$$
Hence the order of $b_0^*$ is computed  by minimizing the function
\begin{eqnarray*}
 b_0\rightarrow
 b_0^4
 +
 \frac{b_0^4}{nh^3}
 +
 \frac{1}{n}
 +
 \frac{1}{nh^5}
 \left(b_0^4\right)^2
 +
 \frac{1}{h^2}
 \left(b_0^4\right)^3
 +
 \frac{b_0^d}{h^7}
 \left(b_0^4\right)^3.
\end{eqnarray*}
Since this function is increasing with $b_0$, the minimum of $RT_n
(\cdot,\cdot,h)$ is achieved for $b_{0*}=(u_n/n)^{1/(d+4)}$. We
shall show later on that this choice of $b_{0*}$ is irrelevant
compared to the one arising when $nb_0^{d+4}=O(1)$.

\vskip 0.1cm Consider now the case $nb_0^{d+4}=O(1)$ i.e
$b_0^4=O\left(1/(nb_0^d)\right)$. This  gives, since $nb_0^{2d}$
diverges under $(H_9)$, using $b_0^{d}/(nb_0^{2d})^p=
O(b_0^{2p})$, $p=2$,
\begin{eqnarray*}
 \frac{1}{nh^3}
\left(
b_0^4
+
\frac{1}{nb_0^{d}}
\right)
\asymp
\frac{1}{n^2b_0^{d}h^3},
\quad
 \frac{1}{nh^5}
 \left(b_{0}^4
 +
 \frac{1}{nb_0^d}
 \right)^2
 =
 O\left(\frac{1}{n^3b_0^{2d}h^7}\right),
 \end{eqnarray*}
\begin{eqnarray*}
\frac{1}{h^2}
\left( b_0^4 + \frac{1}{nb_0^d} \right)^3
 =
 O\left(\frac{1}{n^2b_0^{d}h^3}\right)
\;\mbox{\rm and}\;\;
 \frac{b_0^d}{h^7}
 \left(b_0^4+\frac{1}{nb_0^d}\right)^3
 \asymp
 \left(\frac{1}{n^3b_0^{2d}h^7}\right).
 \end{eqnarray*}
Moreover if $nb_0^dh^4\rightarrow\infty$, we have
$$
 \frac{1}{n^3b_0^{2d}h^7}
 =
 O\left(\frac{1}{n^2b_0^{d}h^3}\right),
 \quad
 RT_n(b_0,b_0,h)
 =
 b_0^4+\frac{1}{n^2b_0^dh^3}+\frac{1}{n}.
 $$
Hence in this case, the order of $b_0^*$ is obtained by finding
the minimum of the function
$b_0^4+\left(1/n^2b_0^dh^3\right)+(1/n)$. The
 minimization of this function  gives  a solution $b_0$ such that
$$
b_0
\asymp
\left(\frac{1}{n^2h^3}\right)^{\frac{1}{d+4}},
\quad
RT_n(b_0,b_0,h)
\asymp
\frac{1}{n}
+
\left(\frac{1}{n^2h^3}\right)^{\frac{4}{d+4}}.
$$
This value  satisfies the constraints $nb_0^{d+4}=O(1)$ and
$nb_0^dh^4\rightarrow\infty$ when
$n^{4-d}h^{d+16}\rightarrow\infty$.

\vskip 0.1cm\noindent If now $nb_0^{d+4}=O(1)$ but
$nb_0^dh^4=O(1)$, we have
\begin{eqnarray*}
 \frac{1}{n^2b_0^{d}h^3}
 =
 O\left(\frac{1}{n^3b_0^{2d}h^7}\right),
 \quad
 RT_n(b_0,b_0,h)
 =
 b_0^4+\frac{1}{n^3b_0^{2d}h^7}+\frac{1}{n}.
\end{eqnarray*}
 In this case, the order of $b_0^*$ is achieved  by minimizing the function
$b_0^4+\left(1/n^3b_0^{2d}h^7\right)+(1/n)$, for which the
solution
 $b_0$ verifies
$$
b_0
\asymp
\left(\frac{1}{n^3h^7}\right)^{\frac{1}{2d+4}},
\quad
RT_n(b_0,b_0,h)
\asymp
\frac{1}{n}
+
\left(\frac{1}{n^3h^7}\right)^{\frac{4}{2d+4}}.
$$
This solution fulfills  the constraint $nb_0^dh^4=O(1)$ when
$n^{4-d}h^{d+16}=O(1)$. Hence we can conclude that for
$b_0^4=O\left(1/(nb_0^d)\right)$, the bandwidth $b_0^*$ satisfies
$$
b_0^*
\asymp
\max \left\lbrace
\left(\frac{1}{n^2h^3}\right)^{\frac{1}{d+4}}
,
\left(\frac{1}{n^3h^7}\right)^{\frac{1}{2d+4}}
\right\rbrace,
$$
which leads to
$$
RT_n\left(b_0^*,b_0^*, h\right)
\asymp
\frac{1}{n}
+
\max
\left\lbrace
\left(\frac{1}{n^2h^3}\right)^{\frac{4}{d+4}}
,
\left(\frac{1}{n^3h^7}\right)^{\frac{4}{2d+4}}
 \right\rbrace.
$$
We need now to compare the solution $b_0^*$ to the candidate
 $b_{0*}=(u_n/n)^{1/(d+4)}$ obtained when $nb_0^{d+4}\rightarrow\infty$.
  For this, we must  do a comparison between the
 orders of $RT_n(b_0^*,b_0^*,h)$  and $RT_n(b_{0*},b_{0*}, h)$. Since
 $RT_n(b_0,b_0,h)\geq b_0^4$, we have
 $RT_n(b_{0*},b_{0*},h)\geq(u_n/n)^{4/(d+4)}$, so that, for $n$
 large enough,
\begin{eqnarray*}
\frac{RT_n(b_0^*,b_0^*, h)}
{RT_n(b_{0*}, b_{0*},h)}
&\leq&
C\left[
\left(\frac{1}{n^2h^3}\right)^{\frac{1}{d+4}}
+
\left(\frac{1}{n^3h^7}\right)^{\frac{4}{2d+4}}
\right]
\left(\frac{n}{u_n}\right)^{\frac{4}{d+4}}
\\
&=&
o(1)
+
O\left(\frac{1}{u_n}\right)^{\frac{4}{d+4}}
\left(
\frac{1}{nh^{\frac{7(d+4)}{d+8}}}
\right)^{\frac{4(d+8)}{(2d+4)(d+4)}}
=
o(1),
\end{eqnarray*}
 using $u_n\rightarrow\infty$ and
$n^{(d+8)}h^{7(d+4)}\rightarrow\infty$ by $(H_{10})$. This shows
that $RT_n(b_0^*,b_0^*,h)\leq RT_n(b_{0*},b_{0^*},h)$ for $n$
large enough. This ends the proof of the Theorem, since $b_0^*$ is
the best candidate for the minimization of $RT_n(\cdot,\cdot, h)$.
\eop

\subsection*{Proof of Theorem \ref{Optimh}}
The proof is the same as the one of Theorem \ref{Optimbandw2} in
Chapter 3. \eop

\subsection*{Proof of Theorem \ref{Normalite2}}
The proof of the Theorem is based on the following Lemma.
\begin{lemma}
 Define
$$
\widetilde{\widetilde{f}}_{in}(\epsilon)
 =
 \frac{1}{b_1^{d}h}
 \int
 \mathds{1}
 \left(x\in\cal{X}\right)
 K_1\left(\frac{X_i-x}{b_1}\right)
 K_2\left(\frac{Y_i-\epsilon-m(x)}{b_1}\right)
 dx.
$$
Then, under $(H_1)$, $(H_6)-(H_8)$, we have, for  $b_1$ and $h$
and going to $0$ and for some constant $C>0$,
\begin{eqnarray*}
\esp \widetilde{\widetilde{f}}_{in}(\epsilon)
 &=&
 f(\epsilon)
 +
 \frac{b_1^2}{2}
 \int
 \mathds{1}
 \left(x\in\cal{X}\right)
 \frac{\partial^2 \varphi(x, \epsilon+m(x))}
 {\partial^2 x}
 dx
 \int
 z K_1(z)z^{\top} dz
 \\
 &&
 +
 \frac{b_1^2}{2}
 \int
 \mathds{1}
 \left(x\in\cal{X}\right)
 \frac{\partial^2 \varphi(x,\epsilon+m(x))}
 {\partial^2 y} dx
 \int
 v^2 K_2(v)
 dv
 +
 o\left(b_1^2+h^2\right),
\\
\Var \left(\widetilde{\widetilde{f}}_{in}(\epsilon)\right)
 &=&
\frac{f(\epsilon)}{h}
 \int
 K_2^2(v)
 dv
 +
o\left(\frac{1}{h}\right),
\\
\esp \left| \widetilde{\widetilde{f}}_{in}(\epsilon)
 -
 \esp
 \widetilde{\widetilde{f}}_{in}(\epsilon)
\right|^3
&\leq&
 \frac{C f(\epsilon)}{h^2}
 \int
 \int
 \left|
  K_1(z_1)
 K_2(v_1)
 \right|^3
 z_1 dv_1
 +
 o\left(\frac{1}{h^2}\right).
\end{eqnarray*}
\label{Momfin}
\end{lemma}
\noindent This Lemma is proved in Appendix B.

\vskip 0.3cm Let now turn to the proof of the Theorem
\ref{Normalite2}. Observe that
\begin{eqnarray}
\widehat{f}_{2n}(\epsilon) - \esp \widetilde{f}_{2n}(\epsilon) =
\left( \widetilde{f}_{2n}(\epsilon) - \esp
\widetilde{f}_{2n}(\epsilon) \right) + \left(
\widehat{f}_{2n}(\epsilon) - \widetilde{f}_{2n}(\epsilon) \right).
\label{Develop}
\end{eqnarray}
Let now $\widetilde{\widetilde{f}}_{in}(\epsilon)$ be as in Lemma
\ref{Momfin}, and note that $\widetilde{f}_{2n}(\epsilon)
=(1/n)\sum_{i=1}^n\widetilde{\widetilde{f}}_{in}(\epsilon)$. The
second and the third claims in Lemma \ref{Momfin} yield, since $h$
goes to $0$ under $(H_{10})$,
\begin{eqnarray*}
\frac { \sum_{i=1}^n \esp \left|
\widetilde{\widetilde{f}}_{in}(\epsilon)
 -
 \esp
 \widetilde{\widetilde{f}}_{in}(\epsilon)
\right|^3
}
{
\left(
\sum_{i=1}^n
 \Var\widetilde{\widetilde{f}}_{in}(\epsilon)
\right)^3}\leq
\frac{\frac{C nf(\epsilon)}{h^{2}}
 \displaystyle{\int}
 \displaystyle{\int}
 \left|
  K_1(z_1)
 K_2(v_1)
 \right|^3
 z_1 dv_1
 +
 o\left(\frac{1}{h^2}\right)
 }{
\left( \frac{nf(\epsilon)}{h}
\displaystyle{\int}
K_2^2(v)
dv
+
 o\left(\frac{n}{h}\right)
\right)^3
 }
 =O(h)= o(1).
\end{eqnarray*}
 Hence the Lyapounov Central Limit Theorem (Billingsley 1968, Theorem 7.3)
  gives, since $n h$
 diverges under $(H_{10})$,
$$
\frac{\widetilde{f}_{2n}(\epsilon)- \esp
\widetilde{f}_{2n}(\epsilon)}{\sqrt{\Var
\widetilde{f}_{2n}(\epsilon)}} =
 \frac{\widetilde{f}_{2n}(\epsilon)
 -
 \esp \widetilde{f}_{2n}(\epsilon)}
 {\sqrt{
 \frac{\Var\widetilde{\widetilde{f}}_{in}(\epsilon)}{n}}}
 \stackrel{d}{\rightarrow}
 \mathcal{N}\left(0,1\right),
$$
which yields, using the second result  in Lemma \ref{Momfin},
\begin{eqnarray}
\sqrt{nh} \left( \widetilde{f}_{2n}(\epsilon) - \esp
\widetilde{f}_{2n}(\epsilon) \right) \stackrel{d}{\rightarrow}
\mathcal{N} \left( 0, f(\epsilon) \int K_2^2(v) dv \right).
 \label{Develop1}
\end{eqnarray}
Observe now that  Lemma \ref{Ordrefnf} gives, for $b_1=b_0$,
\begin{eqnarray*}
\widehat{f}_{2n}(\epsilon)-\widetilde{f}_{2n}(\epsilon) &=&
O_{\prob} \left[ b_{0}^{4} + \frac{1}{nh^3} \left( b_0^4 +
\frac{1}{nb_0^{d}} \right) + \frac{1}{n} + \frac{1}{nh^5}
 \left(b_{0}^4
 +
 \frac{1}{nb_0^d}
 \right)^2
 \right]^{1/2}
 \\
 &&
 +
 \;O_{\prob}
 \left[
\frac{1}{n^2b_0^{d}h^3}
 +
\frac{1}{h^2}
\left(
b_0^4
+
\frac{1}{nb_0^d}
\right)^3
+
\frac{b_{0}^d}{h^7}
\left(
b_0^{4}
+
\frac{1}{nb_0^d}
\right)^3
\right]^{1/2}.
 \end{eqnarray*}
Moreover, since by Assumption $(\rm{H}_{11})$ we have
$nb_0^{d+4}=O(1)$, this  ensures that $nb_0^{2d}\rightarrow\infty$
under $(H_9)$, using $b_0^{d}/(nb_0^{2d})^p= O(b_0^{2p})$, $p=2$.
Therefore
\begin{eqnarray*}
 \frac{1}{nh^3}
\left(
b_0^4
+
\frac{1}{nb_0^{d}}
\right)
\asymp
\frac{1}{n^2b_0^{d}h^3},
\quad
 \frac{1}{nh^5}
 \left(b_{0}^4
 +
 \frac{1}{nb_0^d}
 \right)^2
 =
 O\left(\frac{1}{n^3b_0^{2d}h^7}\right),
 \end{eqnarray*}
\begin{eqnarray*}
\frac{1}{h^2}
\left( b_0^4 + \frac{1}{nb_0^d} \right)^3
 =
 O\left(\frac{1}{n^2b_0^{d}h^3}\right)
\;\mbox{\rm and}\;\;
 \frac{b_0^d}{h^7}
 \left(b_0^4+\frac{1}{nb_0^d}\right)^3
 \asymp
 \left(\frac{1}{n^3b_0^{2d}h^7}\right).
 \end{eqnarray*}
Hence, for $b_0$ and $h$ going to, it follows that
\begin{eqnarray*}
\sqrt{nh} \left( \widehat{f}_{2n}(\epsilon) -
\widetilde{f}_{2n}(\epsilon) \right) \asymp O_{\prob} \left[
 nh
 \left( b_0^4 +
\frac{1}{n^2b_0^dh^3}
+
\frac{1}{n}
+
 \frac{1}{n^3b_0^{2d}b_1^7}
\right)
\right]^{1/2}
=
o_{\prob}(1),
\end{eqnarray*}
since  $nb_0^4h=o(1)$ and $nb_0^{d}h^3\rightarrow\infty$ by
Assumption $(\rm{H}_{11})$. Hence from  (\ref{Develop1}) and
(\ref{Develop}), we deduce
$$
\sqrt{nh}
\left(
 \widehat{f}_{2n}(\epsilon)
 -
 \esp
\widetilde{f}_{2n}(\epsilon) \right) \stackrel{d}{\rightarrow}
\mathcal{N} \left(
 0,
f(\epsilon)
\int
K_2^2(v)
dv
\right).
$$
This proves  the Theorem, since the first result of Lemma
\ref{Momfin} gives for $b_1=b_0$,
\begin{eqnarray*}
\esp \widetilde{f}_{2n}(\epsilon) &=& f(\epsilon) +
\frac{b_0^2}{2}
 \int
 \mathds{1}
 \left(x\in\cal{X}\right)
 \frac{\partial^2 \varphi(x,\epsilon+m(x))}
 {\partial^2 x}
 dx
 \int
 z K_1(z)z^{\top}
 dz
 \\
 &&
 +
 \frac{h^2}{2}
 \int
 \mathds{1}
 \left(x\in\cal{X}\right)
 \frac{\partial^2 \varphi(x,\epsilon+m(x))}
 {\partial^2 y}
 dx
 \int
 v^2 K_2(v)
 dv
 +
 o\left(b_0^2+h^2\right)
 :=\overline{f}_{2n}(\epsilon).
 \eop
\end{eqnarray*}

\newpage
 \renewcommand{\theequation}{B.\arabic{equation}}
\renewcommand{\thesubsection}{B.\arabic{subsection}}
\setcounter{subsection}{0} \setcounter{equation}{0}

\begin{center}
\section*{Appendix B:  Proof of Lemmas \ref{OrdreST}-\ref{Momfin}}
\end{center}

\subsection*{Intermediate results for Lemmas \ref{OrdreST}-\ref{OrdreP}}

\begin{lemma}
If $(H_1)-(H_2)$, $(H_7)$ and $(H_9)$ are satisfied,
 we have
\begin{eqnarray*}
\sup_{x\in\mathcal{X}}
\left|
\widehat{g}_{n}(x) - g(x)
\right|
&=&
 O_{\prob}
 \left(
 b_0^4
 +
 \frac{\ln n}{nb_0^d}
 \right)^{1/2},
 \\
 \sup_{x\in\mathcal{X}}
 \left|
 \frac{1}
 {\widehat{g}_{n}(x)}
  -
 \frac{1}{g(x)}
 \right|
 &=&
 O_{\prob}
  \left(
 b_0^4
 +
 \frac{\ln n}{nb_0^d}
 \right)^{1/2}.
\end{eqnarray*}
\label{Ordregn}
\end{lemma}

\begin{lemma}
Let $\esp_{in}[\cdot]$ be the conditional mean given $\left(
X_1,\ldots,X_n,\varepsilon_k,k\neq i\right)$. Then if
$(H_1)-(H_5)$, $(H_8)$ and $(H_{10})$ hold, we have, for any
integer $i\in[1,n]$, $p\in[0,6]$ and $y\in\Rit$,
\begin{eqnarray}
\left|
 \esp_{in}
 \left[
 \varepsilon_i^p
 K_2^{(1)}
 \left(\frac{Y_i-y}{h}\right)
 \right]
 \right|
 \leq
  Ch^2,
 \quad
 \left|
 \esp_{in}
 \left[
 \varepsilon^p
  K_2^{(1)}
 \left(\frac{Y_i-y}{h}\right)^2
 \right]
 \right|
 \leq C h,
 \label{MomderfK1}
 \\
  \left|
 \esp_{in}
 \left[
 \varepsilon_i^p
 K_2^{(2)}
 \left(\frac{Y_i-y}{h}\right)
  \right]
  \right|
  \leq
  Ch^3,
 \quad
 \left|
 \esp_{in}
  \left[
 \varepsilon_i^p
 K_2^{(2)}
 \left(\frac{Y_i-y}{h}\right)^2
 \right]
 \right|
 \leq C h,
 \label{MomderfK2}
 \\
\left|
 \esp_{in}
 \left[
 \varepsilon_i^p
 K_2^{(3)}
\left(\frac{Y_i-y}{h}\right)
 \right]
 \right|
 \leq C h^3,
 \quad
 \left|
 \esp_{in}
 \left[
 \varepsilon_i^p
 K_2^{(3)}
\left(\frac{Y_i-y}{h}\right)^2
 \right]
 \right|
\leq C h,
 \label{MomderfK3}
\end{eqnarray}
for some constant $C>0$.
 \label{MomderfK}
\end{lemma}

Let $\esp_{n}[\cdot]$ and $\Var_n[\cdot]$ be respectively the
conditional mean and the conditional variance given $\left(
X_1,\ldots,X_n\right)$, and denote $b_{0}\vee b_{1}=\max \left(
b_{0},b_{1}\right)$. In the following, $S_n$ and $T_n$ are defined
as in Lemma \ref{OrdreST}. Then the following results are  used in
the proof of Lemmas \ref{OrdreST}, \ref{OrdreR} and \ref{OrdreP}.
\begin{lemma}
If $(H_1)-(H_{10})$ hold, then
\begin{equation*}
\mathbb{E}_{n}
\left[ S_{n}\right]
=
O_{\prob}
 \left( b_{0}^{2}\right),
\quad
\mathbb{E}_{n}
\left[ T_{n} \right]
=
O_{\mathbb{P}}
\left(
b_{0}^{4}
+
\frac{1}{nb_0^{d}}
\right).
\end{equation*}
\label{MeanU}
\end{lemma}

\begin{lemma}
Under  $(H_1)-(H_{10})$, we have
\begin{eqnarray*}
\Var_n
\left[ S_{n}\right]
&=&
O_{\prob}
\left(b_{0}^d\vee b_{1}^d\right)
\left[
\frac{1}{nb_1^dh^3}
\left(
 b_0^{4}
 +
 \frac{1}{nb_0^{d}}
 \right)
 +
 \frac{1}{nb_0^d}
\right],
\\
\Var_n
\left[ T_{n} \right]
&=&
O_{\prob}
\left(b_{0}^d\vee b_{1}^d\right)
\left[
\frac{1}{nb_1^dh^5}
 \left(b_{0}^4
 +
 \frac{1}{nb_0^d}
 \right)^2
 +
 \frac{b_0^4}{nb_0^d}
 +
\frac{1}{n^2b_0^{2d}h^3}
\right].
\end{eqnarray*}
\label{VarU}
\end{lemma}

\begin{lemma}
Define for all integer number $p$ in $[1,3]$,
\begin{equation*}
U_n(x)
=
U_{n}(x;p)
=
\frac{1}{nb_{1}^{d}h^{p+1}}
\sum_{i=1}^{n}
\left(\widehat{m}_{in}(x)-m(x)\right)^p
K_{1}\left(
\frac{X_{i}-x}{b_{1}} \right)
K_{2}^{(p)}\left(
\frac{Y_{i}-\epsilon-m(x)}{h}\right),
\end{equation*}
and assume that $(H_4)$ and $(H_7)$ hold. Consider $C$ large
enough and any $x_{1}$, $x_{2}$ in $\mathcal{X}$ with $\left\Vert
x_{2}-x_{1}\right\Vert \geq Cb_{0}\vee b_{1}$. Then $ U_{n}\left(
x_{1}\right) $ and $U_{n}\left( x_{2}\right) $ are independent
given $X_{1},\ldots ,X_{n}$. \label{IndepU}
\end{lemma}

\begin{lemma}
Set
\begin{eqnarray*}
\beta _{in}(x)
&=&
\frac{\sum_{1\leq j\neq
i\leq
n}\left(m(X_{j})-m( x)\right)
K_0\left(\frac{X_{j}-x
}{b_{0}}\right) }
{nb_{0}^{d}\widehat{g}_{n}\left( x\right)},
\\
\Sigma_{in}(x)
&=&
\frac{\sum_{1\leq j\neq i\leq
n}\varepsilon _{j}
K_0\left( \frac{X_{j}-x}{b_{0}}\right)}
{nb_{0}^{d}\widehat{g}_{n}\left( x\right)}.
\end{eqnarray*}
Then under $(H_1)-(H_5)$ and $(H_7)-(A_{9})$,
 we have, for all integers $p_1$ and $p_2$ in $[0,6]$,
\begin{eqnarray}
\sum_{i=1}^n
\int
\mathds{1}
\left(x\in\mathcal{X}\right)
\left\vert
\beta_{in}^{p_1}(x)
K_1^{p_2}
\left(\frac{X_i-x}{b_1}\right)
\right\vert
dx
&=&
O_{\prob}
 \left(nb_1^d\right)
\left( b_0^{2p_1} \right),
\label{Boundbetain}
\\
\sum_{i=1}^n
\int
\mathds{1}
\left(x\in\mathcal{X}\right)
\esp_n
\left\vert
\Sigma_{in}^{p_1}(x)
K_1^{p_2}
\left(\frac{X_i-x}{b_1}\right)
\right\vert
dx
&=&
O_{\prob}
\left(\frac{nb_1^d}{(nb_0^d)^{p_1/2}}\right).
\label{BoundSigmain}
\end{eqnarray}
\label{BoundEspmin}
\end{lemma}
\vskip 0.3cm\noindent
 The proof of  these lemmas are given in Appendix C.

\newpage
\subsection*{Proof of Lemma \ref{OrdreST}}

The proof follows directly from Lemmas \ref{MeanU} and \ref{VarU}.
Indeed, since the Tchebychev inequality, which ensures that
$$
A_n
=
O_{\prob}
\left(
\esp_n\left[A_n\right]
+
\Var_n^{1/2}
\left(A_n\right)
\right),
$$
Lemmas \ref{MeanU} and \ref{VarU} then give
\begin{eqnarray*}
S_n
&=&
O_{\prob}
\left[
b_{0}^{4}
+
\left(b_{0}^d\vee b_{1}^d\right)
\left(
\frac{1}{nb_1^dh^3}
\left(
b_0^4
+
\frac{1}{nb_0^{d}}
\right)
+
\frac{1}{nb_0^d}
\right)
\right]^{1/2},
\\
T_n
&=&
O_{\prob}
\left[
\left(
b_{0}^{4}
+
\frac{1}{nb_0^{d}}
\right)^2
+
\left(b_{0}^d\vee b_{1}^d\right)
\left(
\frac{1}{nb_1^dh^5}
 \left(b_{0}^4
 +
 \frac{1}{nb_0^d}
 \right)^2
 +
 \frac{b_0^4}{nb_0^d}
 +
\frac{1}{n^2b_0^{2d}h^3}
\right)
 \right]^{1/2},
\end{eqnarray*}
which proves Lemma \ref{OrdreST}. \eop

\subsection*{Proof of Lemma \ref{OrdreR}}

Set
$$
 R_n
  =
 \int
 \mathds{1}
 \left(x\in\mathcal{X}\right)
 R_n(x)
 dx.
$$
The proof of the Lemma proceeds by computing the conditional mean
and the conditional variance of $R_n$. For the conditional mean,
define
\begin{eqnarray*}
I_{in}(x)
 &=&
 \int_{0}^{1}
 (1-u)^2
  K_2^{(3)}
 \left(
 \frac{
 Y_i-\epsilon-m(x)
 -u\left(\widehat{m}_{in}(x)-m(x)\right)}{h}
\right)
 du,
 \\
 R_{in}(x)
 &=&
 \frac{1}{nb_1^{d}h^4}
 K_1\left(
 \frac{X_i-x}{b_1}
 \right)
\left(\widehat{m}_{in}(x) - m(x)\right)^3
I_{in}(x).
\end{eqnarray*}
This gives
\begin{eqnarray}
\esp_n\left[R_n\right]
=
\sum_{i=1}^n
\int
\mathds{1}
\left(x\in\mathcal{X}\right)
\esp_n
\left[R_{in}(x)\right]
dx,
\end{eqnarray}
where
\begin{eqnarray*}
 \esp_n
 \left[R_{in}(x)\right]
 =
 \frac{1}{nb_1^{d}h^4}
  K_1\left(
 \frac{X_i-x}{b_1}
 \right)
 \esp_{n}
 \left[
 \left(
 \widehat{m}_{in}(x) - m(x)
 \right)^3
 \esp_{in}
 \left[I_{in}(x)\right]
 \right].
\end{eqnarray*}
Moreover, since by Lemma \ref{MomderfK}-(\ref{MomderfK3}) we have
\begin{eqnarray*}
\left|
\esp_{in}
\left[I_{in}(x)\right]
\right|
 =
 \left|
 \int_{0}^{1}
 (1-u)^2
 \esp_{in}
 \left[
   K_2^{(3)}
 \left(
 \frac{
 Y_i-\epsilon-m(x)
 -u\left(\widehat{m}_{in}(x)-m(x)\right)}{h}
 \right)
 \right]
 du
 \right|
 \leq
 Ch^3,
\end{eqnarray*}
it then follows that, setting $p_1=3$ and $p_2=1$ Lemma
\ref{BoundEspmin},
\begin{eqnarray}
\nonumber
\lefteqn{
\left\vert
\esp_n
\left[R_n\right]
\right\vert
}
\\\nonumber
&\leq&
\frac{Ch^3}{nb_1^{d}h^4}
 \sum_{i=1}^n
 \int
 \mathds{1}
 (x\in \mathcal{X})
 \esp_{n}
 \left|
 \left(\widehat{m}_{in}(x)-m(x)\right)^3
   K_1\left(
 \frac{X_i-x}{b_1}
 \right)
 \right|
 dx
 \\\nonumber
 &\leq&
\frac{C}{nb_1^{d}h}
 \sum_{i=1}^n
 \int
 \mathds{1}
 (x\in \mathcal{X})
 \left\vert
 \beta_{in}^3(x)
  K_1\left(
 \frac{X_i-x}{b_1}
 \right)
 \right\vert
 dx
 \\\nonumber
 &&
 +
 \frac{C}{nb_1^{d}h}
 \sum_{i=1}^n
 \int
 \mathds{1}
 (x\in \mathcal{X})
 \esp_{n}
  \left|
 \Sigma_{in}^3(x)
   K_1\left(
 \frac{X_i-x}{b_1}
 \right)
 \right|
 dx
 \\
 &=&
O_{\prob}
\left[
\frac{1}{h^2}
\left(
b_0^4
+
\frac{1}{nb_0^d}
\right)^3
\right]^{1/2}.
\label{espRn}
\end{eqnarray}

 Consider  now the  conditional variance of
$R_n$. Let $C$ large enough and consider  $x_{1}$, $x_{2}$ in
$\mathcal{X}$ with $\left\Vert x_{2}-x_{1}\right\Vert \geq
Cb_{0}\vee b_{1}$.
 Then given $X_{1},\ldots ,X_{n}$ and under $(H_4)$,
  there exists two functions
 $\Phi_{1n}$ and $\Phi_{2n}$ such that
\begin{equation*}
R_{n} \left( x_{1}\right)
=
\Phi _{1n}
\left(\varepsilon_{i},i\in I_{1}\right)
\text{ and }
 R_{n}\left( x_{2}\right)
 =
\Phi_{2n}\left(\varepsilon _{i},i\in I_{2}\right),
\end{equation*}
with an empty $I_{1}\cap I_{2}$, since the Kernel functions $K_0$
and $K_1$ are compactly supported. Hence $R_{n}\left(
x_{1}\right)$ and $R_{n}\left( x_{2}\right) $ are independent
given $X_{1},\ldots ,X_{n}$, provided that $\left\Vert
x_{2}-x_{1}\right\Vert \geq Cb_{0}\vee b_{1}$, for $C$
sufficiently large.  Therefore
\begin{eqnarray}
\nonumber
\lefteqn{ \Var_n\left(R_{n}\right)}
\\\nonumber
&=&
\Var_n\left(\int\mathds{1}
\left( x\in \mathcal{X} \right)
R_{n}\left( x\right) dx\right)
 =
 \int \int
 \mathds{1}\left(\left(
x_{1},x_{2}\right)
\in \mathcal{X}^{2}\right)
\Cov_n\left(
R_{n}\left( x_{1}\right),
R_{n}\left( x_{2}\right)
\right)
dx_{1}dx_{2}
\\\nonumber
&\leq&
\int \int
\mathds{1}
\left( \left( x_{1},x_{2}\right)
\in \mathcal{X}
^{2},
\left\Vert x_{2}-x_{1}
\right\Vert
\leq Cb_{0}\vee
b_{1}\right)
\Var_n^{1/2}
\left( R_{n}\left( x_{1}\right)
\right)
\Var_n^{1/2}
\left( R_{n}\left( x_{2}\right)
\right)
dx_{1}dx_{2}
\\\nonumber
&\leq &
\frac{1}{2}
\int \int
\mathds{1}\left( \left(
x_{1},x_{2}\right)
\in \mathcal{X}^{2},
 \left\Vert
x_{2}-x_{1}
\right\Vert
\leq Cb_{0}\vee b_{1}\right)
\left\{
\Var_n\left( R_{n}\left( x_{1}\right)
\right)
+
\Var_n\left(
R_{n}\left( x_{2}\right)
\right) \right\} dx_{1}dx_{2}
\\
&\leq &
C\left( b_{0}\vee b_{1}\right)^{d}
\int
\mathds{1}\left(
x\in
\mathcal{X}\right)
\Var_n\left( R_{n}
\left( x\right) \right)
dx,
\label{VarR}
\end{eqnarray}
where
\begin{eqnarray}
\nonumber
\lefteqn{
\int
\mathds{1}
\left( x\in \mathcal{X} \right)
\Var_n\left(
R_{n}\left( x\right)\right)
dx
}
\\\nonumber
&=&
\sum_{i=1}^{n}
\int
\mathds{1}
\left( x\in \mathcal{X}\right)
\Var_n\left(
R_{in}\left( x\right)\right)
 dx
\\
&&
+
\sum_{1\leq i_{1}\neq i_{2}\leq n}
\int
\mathds{1}
\left(
x\in \mathcal{X}
\right)
\Cov_n
\left(
R_{i_{1}n}(x)
 ,
R_{i_{2}n}(x)
\right)
dx.
\label{VarR1}
\end{eqnarray}
For the conditional variances in (\ref{VarR1}), we have
\begin{eqnarray*}
 \Var_n
 \left(R_{in}(x)\right)
 \leq
 \frac{1}{(nb_1^{d}h^4)^2}
  K_0^2
 \left(
 \frac{X_i-x}{b_1}
 \right)
 \esp_n
 \left[
 \left(\widehat{m}_{in}(x)-m(x)\right)^6
  I_{in}^2(x)
  \right]
\end{eqnarray*}
with, applying Lemma \ref{MomderfK}-(\ref{MomderfK3}),
\begin{eqnarray*}
 \lefteqn{
 \esp_n
 \left[
 \left(\widehat{m}_{in}(x)-m(x)\right)^6
 I_{in}^2(x)
 \right]
 =
\esp_n
 \left[
 \left(\widehat{m}_{in}(x)-m(x)\right)^6
 \esp_{in}
 \left[
 I_{in}^2(x)
 \right]
\right]
 }
 \\
&\leq&
 C\esp_n
 \left[
 \left(\widehat{m}_{in}(x)-m(x)\right)^6
 \sup_{y\in\Rit}
 \esp_{in}
 \left[
  K_2^{(3)}
 \left(
 \frac{
 Y_i-y}{h}
 \right)^2
 \right]
 \right]
 du
 \\
 &\leq&
  Ch\esp_n
 \left[
 \left(
 \widehat{m}_{in}(x)-m(x)
 \right)^6
\right].
\end{eqnarray*}
 Hence from this result and Lemma \ref{BoundEspmin},  we
deduce
\begin{eqnarray}
\nonumber
\lefteqn{
\sum_{i=1}^{n}
\int
\mathds{1}
\left( x\in \mathcal{X} \right)
\Var_n\left(
R_{in}\left( x\right)\right)
 dx
}
\\\nonumber
&\leq&
\frac{C h}{(nb_1^{d}h^4)^2}
\sum_{i=1}^{n}
\int
\mathds{1}
\left( x\in \mathcal{X} \right)
\esp_n
 \left[
 \left(\widehat{m}_{in}(x)-m(x)\right)^6
\right]
 K_1^2
 \left(
 \frac{X_i-x}{b_1}
 \right)
 dx
\\\nonumber
&\leq&
\frac{C h}{(nb_1^{d}h^4)^2}
\sum_{i=1}^{n}
\int
\mathds{1}
\left( x\in \mathcal{X} \right)
\left(
 \beta_{in}^6(x)
 +
 \esp
 \left[\Sigma_{in}^6(x)\right]
 \right)
 K_1^2
 \left(
 \frac{X_i-x}{b_1}
 \right)
 dx
\\
&=&
\frac{O_{\prob}
 \left(nb_1^dh\right)}
 {(nb_1^{d}h^4)^2}
 \left(
 b_0^{4}
 +
 \frac{1}{nb_0^d}
\right)^3.
\label{VarR2}
\end{eqnarray}
Let now turn to the sum of the conditional covariances in
(\ref{VarR1}). We have
\begin{eqnarray*}
\left\vert
\Cov_n
\left(
R_{i_{1}n}(x)
,
R_{i_{2}n}(x)
\right)
\right\vert
\leq
\Var_n^{1/2}
\left(R_{i_{1}n}(x)\right)
\Var_n^{1/2}
\left(R_{i_{2}n}(x)\right),
\end{eqnarray*}
where
$$
\Var_n\left(R_{i_{1}n}(x)\right)
\leq
\frac{C h}{(nb_1^{d}h^4)^2}
\esp_n
 \left[
 \left(
 \widehat{m}_{i_{1}n}(x)-m(x)
 \right)^6
\right]
 K_1^2
 \left(
 \frac{X_{i_1}-x}{b_1}
 \right).
$$
Hence
\begin{eqnarray}
\nonumber
\lefteqn{
O_{\prob}
\left(\frac{(nb_1^{d}h^4)^2}{h}\right)
\sum_{1\leq i_{1}\neq i_{2}\leq n}
\int
\mathds{1}
\left(
x\in \mathcal{X}
\right)
\left\vert
\Cov_n \left(R_{i_{1}n}(x)
 ,
R_{i_{2}n}(x)
\right)
\right\vert
 dx }
\\\nonumber
&=&
\sum_{1\leq i_{1}\neq i_{2}\leq n}
\int
\mathds{1}
\left(
x\in \mathcal{X}
\right)
\esp_n^{1/2}
 \left[
 \left(
 \widehat{m}_{i_{1}n}(x)-m(x)
 \right)^6
\right]
\esp_n^{1/2}
 \left[
 \left(
 \widehat{m}_{i_{2}n}(x)-m(x)
 \right)^6
\right]
\\\nonumber
 &&
 \;\;\;\;\;\;\;
 \;\;\;\;\;\;\;
 \;\;\;\;\;\;\;
 \times
 \left\vert
 K_1
 \left(
 \frac{X_{i_1}-x}{b_1}
 \right)
 K_1
 \left(
 \frac{X_{i_2}-x}{b_1}
 \right)
 \right|
 dx
 \\\nonumber
&\leq&
\sum_{1\leq i_{1}\neq i_{2}\leq n}
\int
\mathds{1}
\left( x\in \mathcal{X} \right)
\esp_n
 \left[
 \left(
 \widehat{m}_{i_{1}n}(x)-m(x)
 \right)^6
\right]
 \left\vert
 K_1
 \left(
 \frac{X_{i_1}-x}{b_1}
 \right)
 K_1
 \left(
 \frac{X_{i_2}-x}{b_1}
 \right)
 \right|
 dx
 \\\nonumber
&&
+
\sum_{1\leq i_{1}\neq i_{2}\leq n}
\int
\mathds{1}
\left( x\in \mathcal{X} \right)
\esp_n
 \left[
 \left(
 \widehat{m}_{i_{2}n}(x)-m(x)
 \right)^6
\right]
 \left\vert
 K_1
 \left(
 \frac{X_{i_1}-x}{b_1}
 \right)
 K_1
 \left(
 \frac{X_{i_2}-x}{b_1}
 \right)
 \right|
 dx.
 \\
 \label{CovR}
\end{eqnarray}
Moreover, under $(H_7)$, the change of variable $x=u+b_1X_{i_2}$
and Lemma \ref{BoundEspmin} give
\begin{eqnarray*}
\lefteqn{
\sum_{1\leq i_{1}\neq i_{2}\leq n}
\int
\mathds{1}
\left( x\in \mathcal{X} \right)
\esp_n
 \left[
 \left(
 \widehat{m}_{i_{1}n}(x)-m(x)
 \right)^6
\right]
 \left\vert
 K_1
 \left(
 \frac{X_{i_1}-x}{b_1}
 \right)
 K_1
 \left(
 \frac{X_{i_2}-x}{b_1}
 \right)
 \right|
 dx
}
\\
&=&
b_1^d
\sum_{1\leq i_{1}\neq i_{2}\leq n}
\int
\mathds{1}
\left(u+b_1X_{i_2}
\in \mathcal{X}
\right)
\esp_n
 \left[
 \left(
 \widehat{m}_{i_{1}n}(u+b_1X_{i_2})
 -
 m(u+b_1X_{i_2})
 \right)^6
 \right]
 \\
 &&
 \;\;\;\;\;\;\;\;
 \;\;\;\;\;\;\;\;
 \;\;\;\;\;\;\;\;
 \times
 \left\vert
 K_1\left(u\right)
 K_1\left(
 \frac{X_{i_2}-u+b_1X_{i_2}}{b_1}
 \right)
 \right|
 du
 \\
 &=&
 O_{\prob}\left(nb_1^d\right)
\sum_{i=1}^n
\int
\mathds{1}
\left(x
\in \mathcal{X}
\right)
\esp_n
 \left[
 \left(
 \widehat{m}_{in}(x)
 -
 m(x)
 \right)^6
 \right]
 \left\vert
 K_1\left(
 \frac{X_i-x}{b_1}
 \right)
 \right|
 dx
 \\
&=&
 O_{\prob}\left(nb_1^d\right)
\sum_{i=1}^{n}
\int
\mathds{1}
\left(x\in \mathcal{X}\right)
\left(
 \beta_{in}^6(x)
 +
 \esp
 \left[\Sigma_{in}^6(x)\right]
 \right)
 \left\vert
 K_1
 \left(
 \frac{X_i-x}{b_1}
 \right)
 \right\vert
 dx
 \\
 &=&
 O_{\prob}\left(n^2b_1^{2d}\right)
 \left(
 b_0^{4}
+
\frac{1}{nb_0^d}
\right)^3.
\end{eqnarray*}
Therefore collecting this result and (\ref{CovR}), we arrive at
\begin{eqnarray*}
\sum_{1\leq i_{1}\neq i_{2}\leq n}
\int
\mathds{1}
\left(
x\in \mathcal{X}
\right)
\Cov_n
\left(R_{i_{1}n}(x)
 ,
R_{i_{2}n}(x)
\right)
dx
 =
 \frac{O_{\prob}
 \left(n^2b_1^{2d}h\right)}
 {(nb_1^{d}h^4)^2}
 \left(
 b_0^{4}
+
\frac{1}{nb_0^d}
\right)^3.
\end{eqnarray*}
Substituting this order and  (\ref{VarR2}) in (\ref{VarR1}), it
follows, since $nb_1^d\rightarrow\infty$ under $(H_{10})$,
\begin{eqnarray*}
\int
\mathds{1}
\left( x\in \mathcal{X}\right)
\Var_n\left(
R_{n}\left( x\right)\right)
dx
&=&
 O_{\prob}
\left[
\frac{nb_1^dh}
{(nb_1^{d}h^4)^2}
+
\frac{n^2b_1^{2d}h}
{(nb_1^{d}h^4)^2}
\right]
 \left(
 b_0^{4}
+
\frac{1}{nb_0^d}
\right)^3
\\
 &=&
 O_{\prob}
 \left(
 \frac{1}{h^7}
 \right)
 \left(
 b_0^{4}
+
\frac{1}{nb_0^d}
\right)^3.
\end{eqnarray*}
Hence  by (\ref{VarR}), (\ref{espRn}) and the Tchebychev
inequality, we have
\begin{eqnarray*}
R_n
=
O_{\prob}
\left[
\frac{1}{h^2}
\left(
b_0^4
+
\frac{1}{nb_0^d}
\right)^3
+
\frac{b_{0}^d\vee b_{1}^d}{h^7}
\left(
b_0^{4}
+
\frac{1}{nb_0^d}
\right)^3
\right]^{1/2},
\end{eqnarray*}
which proves the validity of the Lemma. \eop

\subsection*{Proof of Lemma \ref{OrdreP}}

Set
$$
P_n
=
\int
\mathds{1}
\left(x\in\mathcal{X}\right)
P_n(x) dx.
$$
The proof of the Lemma follows by computing the conditional mean
and the conditional variance of $P_n$. For the conditional mean,
define
\begin{eqnarray*}
\widehat{I}_{in}(x)
&=&
\int_0^1
 K_2^{(1)}
 \left(
 \frac{Y_i-\epsilon+\widehat{m}_{in}(x)
 -t\left(\widehat{m}_n(x)-\widehat{m}_{in}(x)\right)}
 {b_1}
\right)
 dt,
\\
P_{in}(x)
&=&
\frac{1}{nb_1^{d}h^2}
\left(\widehat{m}_n(x)-\widehat{m}_{in}(x)\right)
K_1\left(\frac{X_i-x}{b_1}\right)
\widehat{I}_{in}(x).
\end{eqnarray*}
Since
$$
\widehat{m}_n(x)-\widehat{m}_{in}(x)
 =
\frac{Y_i} {nb_0^d\widehat{g}_n(x)}
K_0\left(\frac{X_i-x}{b_0}\right),
$$
 and that $K_0$ is bounded under  $(H_7)$, Lemma \ref{Ordregn}
 gives
\begin{eqnarray}
\nonumber
&&
\esp_n \left[P_n\right]
=
\sum_{i=1}^n
\int
\mathds{1}
\left(x\in\mathcal{X}\right)
\esp_n
\left[P_{in}(x)\right]
dx
\\
&=&
O_{\prob}
\left(\frac{1}{nb_0^d}\right)
\left[
\frac{1}{nb_1^{d}h^2}
\sum_{i=1}^n
\int
\mathds{1}
\left(x\in\mathcal{X}\right)
\left\vert
 K_1\left(\frac{X_i-x}{b_1}\right)
\esp_n
\left[ Y_i
\widehat{I}_{in}(x)
\right]
\right\vert
 dx
 \right].
\label{OrdrePn1}
\end{eqnarray}
Moreover, observe that for any  $y\in\Rit$,
\begin{eqnarray*}
\esp_{in}
\left[
Y_i
K_2^{(1)}
\left(\frac{Y_i-y}{h}\right)
\right]
=
m(X_i)
\esp_{in}
\left[
K_2^{(1)}
\left(
\frac{Y_i-y}{h}
\right)
\right]
+
\esp_{in}
\left[
\varepsilon_i
K_2^{(1)}
\left( \frac{Y_i-y}{h} \right)
\right].
\end{eqnarray*}
Therefore, since $m(\cdot)$ is continuous on the compact support
$\mathcal{X}$ of the $X_i$'s, Lemma
\ref{MomderfK}-(\ref{MomderfK1}) yields
$$
\sup_{y\in\Rit}
\left|
\esp_{in}
\left[
Y_i
K_2^{(1)}
\left(\frac{Y_i-y}{h}\right)
\right]
\right|
\leq
Ch^2,
$$
uniformly for $i\in[1,n]$. Hence conditioning with respect to
$\left(X_1,\ldots,X_n,\varepsilon_k\right)$ yields that
\begin{eqnarray*}
\left\vert
\esp_n
\left[
Y_i
 \widehat{I}_{in}(x)
\right]
\right\vert
\leq
\left\vert
\sup_{y\in\Rit}
\int
\esp_{in}
\left[
Y_i
 K_2^{(1)}
 \left(
 \frac{Y_i-y}{b_1}
 \right)
 \right]
 dy
\right\vert
\leq C h^2,
\end{eqnarray*}
for all  $i$ and $x$. Combining this result with (\ref{OrdrePn1}),
we arrive at
\begin{eqnarray}
\nonumber
\esp_n
\left[P_n\right]
&=&
 O_{\prob}
\left(\frac{1}{nb_0^d}\right)
\left[
\frac{1}{nb_1^d}
\sum_{i=1}^n
\int
\mathds{1}
\left(x\in\mathcal{X}\right)
\left|
K_1\left(
\frac{X_i-x}{b_1}
\right)
\right|
 dx
\right]
\\
&=&
O_{\prob}
\left(\frac{1}{nb_0^d}\right).
\label{espPn}
\end{eqnarray}

\vskip 0.1cm Let now  consider the conditional variance of $P_n$.
Since
\begin{eqnarray*}
P_{in}(x)
=
\frac{1}{nb_1^{d}h^2}
\left[
\frac{Y_i}
{nb_0^d\widehat{g}_n(x)}
K_0\left(\frac{X_i-x}{b_0}\right)
\right]
K_1\left(\frac{X_i-x}{b_1}\right)
\widehat{I}_{in}(x),
\end{eqnarray*}
and that $K_0(\cdot)$ and $K_1(\cdot)$ have compact supports under
$(H_7)$ and $(H_8)$, it is shown that $P_{n}\left( x_{1}\right)$
and $P_{n}\left( x_{2}\right) $ are independent given
$X_{1},\ldots ,X_{n}$, provided that $\left\Vert
x_{2}-x_{1}\right\Vert \geq Cb_{0}\vee b_{1}$, for  $C$ large
enough. Hence arguing as for (\ref{VarR}) gives
\begin{eqnarray}
\Var_n\left(P_{n}\right)
\leq
C\left(b_{0}^d\vee b_{1}^d\right)
\int
\mathds{1}
\left(
x\in
\mathcal{X}\right)
\Var_n\left( P_{n}
\left( x\right)
\right)
dx,
\label{VarP}
\end{eqnarray}
where
\begin{eqnarray}
\nonumber
\lefteqn{
\int
\mathds{1}
\left( x\in \mathcal{X} \right)
\Var_n\left(
P_{n}\left( x\right) \right)
dx
}
\\\nonumber
&=&
\sum_{i=1}^{n}
\int
\mathds{1}
\left( x\in \mathcal{X} \right)
\Var_n\left(
P_{in}\left( x\right)\right)
dx
\\
&&
+
\sum_{1\leq i_{1}\neq i_{2}\leq n}
\int
\mathds{1}
\left(
x\in \mathcal{X}
\right)
\Cov_n
\left(P_{i_{1}n}(x)
 ,
P_{i_{2}n}(x)
\right)
dx.
\label{VarP1}
\end{eqnarray}
For the conditional variances in (\ref{VarP1}), first note that
\begin{eqnarray}
\nonumber
\Var_n\left(P_{in}(x)\right)
\leq
 \frac{1}{(nb_1^{d}h^2)^2}
 K_1^2\left(\frac{X_i-x}{b_1}\right)
 \left[
 \frac{1}
 {(nb_0^d)^2\widehat{g}_n^2(x)}
 K_0^2\left(\frac{X_i-x}{b_0}\right)
 \esp_n
 \left[
 \left(Y_i-m(x)\right)^2
 \widehat{I}_{in}^2(x)
 \right]
 \right].
 \\
\label{VarPn2}
\end{eqnarray}
Next, observe that for $X_i=z$ and $y\in\Rit$, and under
 $(H_1)$, $(H_3)-(H_5)$ and
$(H_7)$, we have
\begin{eqnarray*}
 \esp_n
 \left[
 Y_i^2
 K_2^{(1)}
\left(\frac{Y_i-y}{h}\right)^2
\right]
&=&
\int
\left(m(z)+e\right)^2
K_2^{(1)}
\left(\frac{m(z)+e-y}{h}\right)^2
f(e)de
\\
&\leq&
C h,
\end{eqnarray*}
 uniformly in $x$ and $i$. From this result
 and the Hölder inequality, we deduce
\begin{eqnarray*}
\lefteqn{
\esp_n
 \left[
 Y_i^2
 \widehat{I}_{in}^2(x)
 \right]
}
\\
 &\leq&
 \int_0^1
 \esp_n
 \left[
 Y_i^2
 K_2^{(1)}
 \left(
 \frac{Y_i-\epsilon+\widehat{m}_{in}(x)
 -t\left(\widehat{m}_n(x)-
 \widehat{m}_{in}(x)\right)}
 {h}\right)^2
 \right]
 dt
 \\
 &\leq&
 \sup_{y\in\Rit}
 \esp_n
\left[
Y_i^2
 K_2^{(1)}
\left(\frac{Y_i-y}{h}\right)^2
\right]
\\
&\leq&
 C h.
\end{eqnarray*}
Hence  by (\ref{VarPn2}) and Lemma \ref{Ordregn}, we have,  since
$K_0(\cdot)$ is bounded under $(H_7)$,
\begin{eqnarray*}
 \Var_n
 \left(
  P_{in}(x)
 \right)
 \leq
 \frac{C}{(nb_1^{d}h^2)^2}
 \times
 \frac{h}{(nb_0^{2d})^2\widehat{g}_n^2(x)}
 K_1^2\left(\frac{X_i-x}{b_1}\right),
\end{eqnarray*}
so that
\begin{eqnarray}
\nonumber
\lefteqn{
\sum_{i=1}^{n}
\int
\mathds{1}
\left( x\in \mathcal{X} \right)
\Var_n\left(P_{in}(x)\right)
dx
}
\\\nonumber
&=&
O_{\prob}
\left(\frac{h}{(nb_1^{d}h^2)^2}\right)
\left(\frac{1}{(nb_0^d)^2}\right)
\sum_{i=1}^{n}
\int
\mathds{1}
\left(x\in \mathcal{X} \right)
K_1^2\left(\frac{X_i-x}{b_1}\right)
dx
\\
&=&
O_{\prob}
\left(\frac{1}{nb_1^{d}h^3}\right)
\left(\frac{1}{n^2b_0^{2d}}\right).
\label{VarPn3}
\end{eqnarray}
 Let now consider the sum of the conditional covariances in
(\ref{VarP1}). We have, using the inequality above,
\begin{eqnarray*}
\left\vert
\Cov_n
\left(
P_{i_{1}n}(x)
,
P_{i_{2}n}(x)
\right)
\right\vert
 &\leq&
 \Var_n^{1/2}\left(P_{i_{1}n}(x)\right)
 \Var_n^{1/2}\left(P_{i_{2}n}(x)\right)
 \\
 &\leq&
 \frac{C}{(nb_1^{d}h^2)^2}
 \times
 \frac{h}{(nb_0^{2d})^2\widehat{g}_n^2(x)}
 \left\vert
 K_1\left(\frac{X_{i_1}-x}{b_1}\right)
 K_1\left(\frac{X_{i_2}-x}{b_1}\right)
 \right\vert.
\end{eqnarray*}
Hence from Lemma \ref{Ordregn}, we deduce
\begin{eqnarray*}
\lefteqn{
\sum_{1\leq i_{1}\neq i_{2}\leq n}
\int
\mathds{1}
\left(
x\in \mathcal{X}
\right)
\left\vert
\Cov_n
\left(P_{i_{1}n}(x)
 ,
P_{i_{2}n}(x)
\right)
\right\vert
dx
}
\\
&=&
O_{\prob}
\left(\frac{1}{(nb_1^{d}h^2)^2}\right)
\left(\frac{h}{(nb_0^d)^2}\right)
\sum_{1\leq i_{1}\neq i_{2}\leq n}
\int
\mathds{1}
\left(
x\in \mathcal{X}
\right)
\left\vert
 K_1\left(\frac{X_{i_1}-x}{b_1}\right)
 K_1\left(\frac{X_{i_2}-x}{b_1}\right)
 \right\vert
 dx
 \\
 &=&
 O_{\prob}
\left(\frac{1}{n^2b_0^{2d}h^3}\right).
\end{eqnarray*}
Substituting this order and (\ref{VarPn3}) in (\ref{VarP1}), and
using (\ref{VarP}), (\ref{espPn}) and the Tchebychev inequality,
we arrive at
\begin{eqnarray*}
 P_n
 &=&
 O_{\prob}
 \left[
 \frac{1}{nb_0^{d}}
 +
 \left(b_{0}^d\vee b_{1}^d\right)^{1/2}
 \left(
 \frac{1}{nb_1^{d}h^3}
 \left(\frac{1}{n^2b_0^{2d}}\right)
 +
 \frac{1}{n^2b_0^{2d}h^3}
 \right)^{1/2}
 \right]
 \\
 &=&
 O_{\prob}
 \left(
 \frac{1}{n^2b_0^{2d}}
 +
 \frac{b_{0}^d\vee b_{1}^d}
 {n^2b_0^{2d} h^3}
 \right)^{1/2}.
\end{eqnarray*}
This ends  the proof of the Lemma.\eop

\subsection*{Proof of Lemma \ref{Momfin}}

The first equality of the lemma is given by (\ref{Bias4}), since
$\widetilde{f}_{2n}(\epsilon)
=\left(1/n\right)
\sum_{i=1}^n\widetilde{\widetilde{f}}_{in}(\epsilon)$,
so that
\begin{eqnarray*}
\lefteqn{
\esp
\widetilde{\widetilde{f}}_{in}(\epsilon)
=
\esp
\widetilde{f}_{2n}(\epsilon)
}
 \\
 &=&
 f(\epsilon)
 +
 \frac{b_1^2}{2}
 \int
 \mathds{1}
 \left(x\in\cal{X}\right)
 \frac{\partial^2\varphi(x,\epsilon+m(x))}
 {\partial^2 x}
 dx
 \int
 z K_1(z)z^{\top}
 dz
 \\
 &&
 +
 \frac{h^2}{2}
 \int
 \mathds{1}
 \left(x\in\cal{X}\right)
 \frac{\partial^2 \varphi(x,\epsilon+m(x))}
 {\partial^2 y}
 dx
 \int
 v^2 K_2(v)
 dv
 +
 o\left(b_1^2+h^2\right).
\end{eqnarray*}

For the second result of the Lemma,  we have
\begin{eqnarray}
\nonumber \lefteqn{ \Var \left(
\widetilde{\widetilde{f}}_{in}(\epsilon) \right) = \esp
\left[\widehat{f}_{in}^2(\epsilon)\right] - \esp^2
\left[\widetilde{\widetilde{f}}_{in}(\epsilon)\right] }
\\
&=&
\frac{1}{b_1^{2d}h^2}
\esp
\left[
\left[
\int
\mathds{1}
\left(x\in\cal{X}\right)
 K_1
 \left(
 \frac{X_i-x}{b_1}
 \right)
 K_2
 \left(
 \frac{Y_i-\epsilon-m(x)}{b_1}
 \right)
 dx
 \right]^2
 \right]
 +
 O(1).
 \label{Varfin1}
\end{eqnarray}
Observe now that  the changes of variables $x=x_1+b_1z_1$ and
$y_1=\epsilon+m(x_1+b_1z_1)+b_1v_1$ give
\begin{eqnarray}
 \nonumber
 \lefteqn{
 \esp
\left[
\left[
\int
\mathds{1}
\left(x\in\cal{X}\right)
 K_1
 \left(
 \frac{X_i-x}{b_1}
 \right)
 K_2
 \left(
 \frac{Y_i-\epsilon-m(x)}{b_1}
 \right)
 dx
\right]^2
\right]
 }
 \\\nonumber
 &=&
 \int
 dx_1
 \int
 \left[
 \int
 \mathds{1}
 \left(x\in\cal{X}\right)
 K_1
 \left(
 \frac{x_1-x}{b_1}
 \right)
 K_2
 \left(
 \frac{y_1-\epsilon-m(x)}{b_1}
 \right)
 dx
 \right]^2
 \varphi\left(x_1, y_1\right)
  dy_1
 \\\nonumber
 &=&
 b_1^{2d}h
 \int
 dx_1
 \int
 \left[
  K_2(v_1)
 \int
 \mathds{1}
 \left(x_1+b_1z_1\in\cal{X}\right)
 K_1(z_1)
   dz_1
 \right]^2
 \varphi
 \left(x_1,\epsilon+m(x_1+b_1z_1)+b_1v_1\right)
 dv_1.
 \\
 \label{Varfin2}
\end{eqnarray}
Moreover,  note that under $(H_3)$ and $(H_6)$ we have
\begin{eqnarray*}
m\left(x_1+b_1z_1\right)
 &=&
 m(x_1)
+
b_1 z_1
\int_{0}^{1}
 m^{(1)}\left( x_1+tb_{1}z_1\right)
 dt,
 \\
\varphi
 \left(x_1,\epsilon+m(x_1+b_1z_1)+b_1v_1\right)
 &=&
 \varphi
 \left(x_1,\epsilon+m(x_1)\right)
 +
 b_1z_1\theta_n(x_1,z_1)
 \int_0^1
 \frac{\partial\varphi}{\partial y}
 \left(x_1, \bar{\theta}_n(u,x_1,z_1)\right)
 du,
\end{eqnarray*}
where
$$
\theta_n(x_1,z_1)
 =
 \int_{0}^{1}
 m^{(1)}\left( x_1+tb_{1}z_1\right)
 dt,
 \quad
\bar{\theta}_n(u,x_1,z_1)
=
\epsilon+m(x_1)+u\theta_n(x_1,z_1).
$$
Therefore
\begin{eqnarray}
\nonumber
\lefteqn{
\int
 dx_1
 \int
 \left[
  K_2(v_1)
 \int
 \mathds{1}
 \left(x_1+b_1z_1\in\cal{X}\right)
 K_1(z_1)
   dz_1
 \right]^2
 \varphi
 \left(x_1,\epsilon+m(x_1+b_1z_1)+b_1v_1\right)
 dv_1
 }
 \\\nonumber
 &=&
\int
 dx_1
 \int
 \left[
  K_2(v_1)
 \int
 \mathds{1}
 \left(x_1+b_1z_1\in\cal{X}\right)
 K_1(z_1)
   dz_1
 \right]^2
 \varphi
 \left(x_1,\epsilon+m(x_1)\right)
 dv_1
 +
 O(b_1)
\\\nonumber
&=&
 \int
 dx_1
 \int
 \left[
 \mathds{1}
 \left(x_1\in\cal{X}\right)
 K_2(v_1)
 \int
  K_1(z_1)
  dz_1
 \right]^2
 \varphi
 \left(x_1,\epsilon+m(x_1)\right)
 dv_1
 \\
 &&
 +
\int
 dx_1
 \int
 \delta_n\left(x_1, v_1\right)
 \varphi
 \left(x_1,\epsilon+m(x_1)\right)
 dv_1
 +
 O(b_1),
 \label{Varfin3}
\end{eqnarray}
where
\begin{eqnarray*}
\delta_n(x_1,v_1)
 =
 \left[
 K_2(v_1)
 \int
 \mathds{1}
 \left(x_1+b_1z_1\in\cal{X}\right)
 K_1(z_1)
  dz_1
 \right]^2
 -
 \left[
 \mathds{1}
 \left(x_1\in\cal{X}\right)
 K_2(v_1)
 \int
 K_1(z_1)
 dz_1
\right]^2.
\end{eqnarray*}
Applying the Lebesgue Dominated Convergence Theorem yields, for
$b_1$ going to $0$,
$$
 \int
 dx_1
 \int
 \delta_n\left(x_1, v_1\right)
 \varphi
 \left(x_1,\epsilon+m(x_1)\right)
 dv_1
=
o(1).
$$
Hence by (\ref{Varfin3}) and (\ref{fth}), we have, since
$\int\!K_1(z_1)dz_1=1$ under $(H_7)$,
\begin{eqnarray*}
\lefteqn{
\int
 dx_1
 \int
 \left[
 K_2(v_1)
 \int
 \mathds{1}
 \left(x_1+b_1z_1\in\cal{X}\right)
 K_1(z_1)
  dz_1
 \right]^2
 \varphi
 \left(x_1,\epsilon+m(x_1+b_1z_1)+b_1v_1\right)
 dv_1
 }
\\
&=&
 \int
 K_2^2(v_1)
 dv_1
 \int
 \mathds{1}
 \left(x_1\in\cal{X}\right)
 \varphi\left(x_1,\epsilon+m(x_1)\right)
 dx_1
 +
 o(1)
 \\
 &=&
 f(\epsilon)
\int
 K_2^2(v)
 dv
 +
 o(1).
\end{eqnarray*}
Combining this result with (\ref{Varfin2}) and (\ref{Varfin1}), we
arrive at
$$
\Var \left(\widetilde{\widetilde{f}}_{in}(\epsilon)\right)
 =
\frac{f(\epsilon)}{h}
 \int
 K_2^2(v)
 dv
+
o\left(\frac{1}{h}\right),
$$
which proves the second result of the lemma.

 \vskip 0.1cm
The last statement of Lemma is immediate. Indeed, the Triangular
and Convex inequalities and the Lebesgue Dominated Convergence
Theorem give, by (\ref{fth}),
\begin{eqnarray*}
\lefteqn{ \esp \left| \widetilde{\widetilde{f}}_{in}(\epsilon) -
\esp \widetilde{\widetilde{f}}_{in}(\epsilon) \right|^3 }
\\
&\leq&
\frac{C}{b_1^{3d}h^3}
\int
dx_1
\int
\left|
\int
\mathds{1}
\left(x\in\cal{X}\right)
 K_1
 \left(
 \frac{x_1-x}{b_1}
 \right)
 K_2
 \left(
 \frac{y_1-\epsilon-m(x)}{b_1}
 \right)
 dx
\right|^3
dy_1
 \\
 &=&
 \frac{Cb_1^{3d}h}{b_1^{3d}h^3}
 \int
 dx_1
 \int
 \left|
 \int
 \mathds{1}
 \left(x_1+b_1z_1\in\cal{X}\right)
 K_1(z_1)
 K_2(v_1)
  dz_1
 \right|^3
 \varphi
 \left(x_1,\epsilon+m(x_1+b_1z_1)+b_1v_1\right)
 dv_1
 \\
&=&
 \frac{C f(\epsilon)}{h^2}
 \int
 \int
 \left|
  K_1(z_1)
 K_2(v_1)
 \right|^3
 dz_1dv_1
 +
 o\left(\frac{1}{h^2}\right).
 \eop
\end{eqnarray*}

\newpage

\setcounter{equation}{0} \setcounter{subsection}{0}
\setcounter{lemma}{0}
\renewcommand{\theequation}{C.\arabic{equation}}
\renewcommand{\thesubsection}{C.\arabic{subsection}}
\begin{center}
\section*{Appendix C}
\end{center}

\subsection*{Proof of Lemma \ref{Ordregn}}
See the proof of Lemma \ref{Estig} in Chapter 3. \eop

\subsection*{Proof of Lemma \ref{MomderfK}}

For the first bound in (\ref{MomderfK1}), set $f_p(e)=e^pf(e)$,
and observe if $X_i=x$, we have by $(H_4)$ and the change of
variable $e=y-m(x)+hv$,
\begin{eqnarray}
\nonumber
\lefteqn{
 \esp_{in}
 \left[
\varepsilon_i^p
 K_2^{(1)}
 \left(\frac{Y_i-y}{h}\right)
 \right]
=
 \esp
 \left[
 \varepsilon_i^p
 K_2^{(1)}
 \left(\frac{\varepsilon_i+m(x)-y}{h}\right)
 \right]
 }
\\
&=& \int
 K_2^{(1)}
 \left(\frac{e+m(x)-y}{h}\right)
 f_p(e)
 de
 =
 h\int
 K_2^{(1)}(v)
 f_p\left(y-m(x)+h v\right)
 dv.
\label{firsteq}
\end{eqnarray}
Therefore, since $f_p$  has a bounded  continuous derivative under
$(A_{5})$ and that $\int\! K_2^{(1)}(v)dv =0$ under $(H_8)$, the
Taylor inequality gives
\begin{eqnarray*}
 \left|
\int
K_2^{(1)}
\left(\frac{e+m(x)-y}{h}\right)
 f_p(e) de
 \right|
 &=&
 h
 \left|
 \int
 K_2^{(1)}(v)
 \biggl[
 f_p\left(y-m(x)+h v\right)
 -
 f_p\left(y-m(x)\right)
 \biggr]
 \right|
 dv
 \\
 &\leq&
  h^2
 \sup_{u\in\Rit}
 |f_p^{(1)}(u)|
 \int
 |vK_2^{(1)}(v)|
 dv
 \\
&\leq&
 Ch^2,
\end{eqnarray*}
uniformly in $x\in\cal{X}$ and $y\in\Rit$. Hence from  this
inequality and (\ref{firsteq}), we deduce
\begin{eqnarray*}
\left|
\esp_{in}
 \left[
 \varepsilon_i^p
 K_2^{(1)}
 \left(\frac{Y_i-y}{h}\right)
 \right]
\right|
\leq
C h^2,
\end{eqnarray*}
for any $y\in\Rit$. This proves the first inequality in
(\ref{MomderfK1}).  The second bound of (\ref{MomderfK1}) is
immediate under $(H_5)$ and $(H_8)$, since for any $x$ in
$\cal{X}$, $\ell\in[1,3]$ and $y\in\Rit$,
\begin{eqnarray}
\nonumber
 \left|
 \esp_{in}
 \left[
 \varepsilon_i^p
 K_2^{(\ell)}
 \left(\frac{Y_i-y}{h}\right)^2
\mid X_i=x
 \right]
\right|
&=&
\left|
 h
 \int
 K_2^{(\ell)}(v)^2
 f_p\left((y-m(x)+h v\right)
 dv
\right|
\\\nonumber
&\leq&
 h\sup_{u\in\Rit}
 |f_p(u)|
 \int
 K_2^{(\ell)}(v)^2
 dv
 \\
&\leq&
 C h,
\label{firsteq1}
\end{eqnarray}
uniformly for $i$, $x$ and $y$. This proves (\ref{MomderfK1}).

 \vskip 0.1cm
 The proof of the second inequalities of
 (\ref{MomderfK2}) and (\ref{MomderfK3})  follows from  (\ref{firsteq1}).
 The first bounds in
(\ref{MomderfK2}) and (\ref{MomderfK3}) are proved simultaneously.
For  any integer $\ell$ in $\in\{2,3\}$ and $x\in\cal{X}$, we have
\begin{eqnarray}
\nonumber
 \esp_{in}
 \left[
 \varepsilon_i^p
 K_2^{(\ell)}
 \left(\frac{Y_i-y}{h}\right)
 \mid X_i=x
 \right]
 &=&
 \int
 K_2^{(2)}
 \left(\frac{e+m(x)-y}{h}\right)
 f_p(e)
 de
\\
&=&
 h
 \int
 K_2^{(2)}(v)
 f_p\left(y-m(x)+h v\right)
 dv.
 \label{firsteq2}
\end{eqnarray}
Under $(H_8)$,  the Kernel function $K_2(\cdot)$ is symmetric, has
a compact support and two
 continuous derivatives, with
$\int\! K_2^{(\ell)}(v)dv=0$ and
 $\int\!v K_2^{(\ell)}(v) dv=0$.
 Therefore, since  $f_p$ has a bounded continuous second
order derivative by $(H_5)$, the second order Taylor expansion
gives, for some $\theta=\theta(y,x,h v)$,
\begin{eqnarray*}
\lefteqn{
 \left|
 h
 \int
 K_2^{(\ell)}(v)
 f_p\left(y-m(x)+hv\right)
 dv
 \right|
  }
\\
&=&
 \left|
 h\int
 K_2^{(\ell)}(v)
\biggl[
 f_p\left(y-m(x)+h v\right)
 -
 f_p\left(y-m(x)\right)
 \biggr]
 dv
 \right|
\\
&=&
\left|
 h
 \int
 K_2^{(\ell)}(v)
\left[
 h v f_p^{(1)}\left(y-m(x)\right)
 +
 \frac{h^2v^2}{2}
 f_p^{(2)}\left(y-m(x)+\theta h v\right)
 \right]
 dv
 \right|
 \\
 &=&
\left|
\frac{h^3}{2}
 \int v^2
 K_2^{(\ell)}(v)
 f_p^{(2)}\left(y-m(x)+\theta hv\right)
 dv
 \right|
\\
 &\leq&
 Ch^3.
\end{eqnarray*}
Hence  from this bound and (\ref{firsteq2}), we deduce
$$
 \left|
 \esp_{in}
 \left[
\varepsilon_i^p
 K_2^{(\ell)}
 \left(\frac{Y-y}{h}\right)
 \right]
 \right|
\leq Ch^3,
$$
uniformly for $i$ and $y$. This  ends proof of the Lemma. \eop

\subsection*{Proof of Lemma \ref{MeanU}}

We have
\begin{equation*}
\mathbb{E}_{n}
\left[ S_{n}\right]
=
\int
\mathds{1} \left( x\in
\mathcal{X}\right)
\mathbb{E}_{n}
\left[S_{n}\left( x\right)
\right] dx,
\quad \mathbb{E}_{n}
\left[ T_{n}\right]
=
 \int
\mathds{1}
\left( x\in \mathcal{X}\right)
\mathbb{E}_{n}
\left[T_{n}\left( x\right) \right]
dx,
\end{equation*}
with
\begin{eqnarray*}
\mathbb{E}_{n}
\left[ S_{n}\left( x\right) \right]
&=&
\frac{1}{nb_{1}^{d}h^{2} }
\sum_{i=1}^{n}
\beta _{in}\left(
x\right)
K_{1}\left( \frac{X_{i}-x}{b_{1}} \right)
\mathbb{E}_{n}
\left[ K_{2}^{\left(1\right)}
\left(\frac{Y_{i}-\epsilon-m(x)
}{h}\right)
\right],
\\
\mathbb{E}_{n}
\left[T_{n}\left( x\right)\right]
&=&
\frac{1}{nb_{1}^{d}h^{3} }
\sum_{i=1}^{n} \left( \beta_{in}^{2}
\left( x\right)
+
\mathbb{E}_{n}
\left[
\Sigma _{in}^{2}
\left(
x\right)
\right]
\right)
K_{1}\left(\frac{X_{i}-x}{ b_{1}}\right)
\mathbb{E}_{n}
\left[
K_{2}^{\left( 2\right)
}\left( \frac{
Y_{i}-\epsilon-m(x) }{h}\right)
 \right].
\end{eqnarray*}
Observe first that under $(H_4)$,  Lemma
\ref{MomderfK}-(\ref{MomderfK1}) and Lemma \ref{Ordregn}  give
\begin{eqnarray*}
\lefteqn{ \sup_{x\in\mathcal{X}}
\left\vert
\frac{1}{nb_{1}^{d}h^{3}}
\sum_{i=1}^{n}
\mathbb{E}_{n}
\left[
\Sigma_{in}^{2}\left( x\right)
\right]
K_{1}\left(
\frac{X_{i}-x}{b_{1}}\right)
\mathbb{E}_{n}
\left[
K_{2}^{\left(
2\right) }\left(\frac{Y_{i}-\epsilon-m(x)}{h}\right)
\right]
\right\vert
}
\\
&\leq&
\frac{Ch^{3}}{h^{3}}
\sup_{x\in \mathcal{X}}
\left\vert
\frac{1}{ nb_{1}^{d}}
\sum_{i=1}^{n}\frac{\sum_{j=1}^{n}
K_0^{2}\left( \frac{X_{j}-x}{b_{0} }\right)}
{\left(nb_{0}^{d}
\widehat{g}_{n}\left(x\right)\right)^{2}}
K_{1}\left(\frac{X_{i}-x}{b_{1}}\right)
\right\vert
=
O_{\mathbb{P}}
\left( \frac{1}{nb_{0}^{d}}\right),
\end{eqnarray*}
and then
\begin{eqnarray*}
\left\vert
\int
\mathds{1}
\left( x\in \mathcal{X}\right)
\frac{1}{ nb_{1}^{d}h^{3}}
\sum_{i=1}^{n}
\mathbb{E}_{n}
\left[
\Sigma_{in}^{2}
\left(x\right)
\right]
K_{1}\left(\frac{X_{i}-x}{b_{1}}\right)
\mathbb{E}_{n}
\left[
K_{2}^{\left( 2\right)}
\left(\frac{Y_{i}-\epsilon-m(x)}{h}\right)
\right]
dx
\right\vert
=
O_{\mathbb{P}}
\left(\frac{1}{nb_{0}^{d}}\right).
\end{eqnarray*}
Consider now
\begin{equation*}
V_{n}\left( p\right)
=
\frac{1}{nb_{1}^{d}}
\int \mathds{1}
\left(
x\in \mathcal{X}\right)
\sum_{i=1}^{n}
\left\vert
\beta_{in}^p(x)
K_{1}\left(\frac{X_{i}-x}{b_{1}}\right)
\right\vert
 dx,
\end{equation*}
which is such that, using Lemma \ref{MomderfK}, the equality and
the bound above,
\begin{equation*}
\left\vert
\mathbb{E}_{n}
\left[ S_{n}\right]
\right\vert
\leq
CV_{n}
\left(1\right)
,
\quad
\left\vert
\mathbb{E}_{n}
\left[
T_{n}
 \right]
 \right\vert
 \leq CV_{n}\left( 2\right)
 +
O_{\mathbb{P}}
\left( \frac{1}{nb_{0}^{d}}\right).
\end{equation*}
Since Lemma \ref{BoundEspmin}-(\ref{Boundbetain}) ensures that
$V_n(p) = O_{\prob}\left(b_0^{2p}\right)$ for all integer number
$p\in[1,6]$, it then follows that
\begin{equation*}
\mathbb{E}_{n}
\left[ S_{n}\right]
 =
 O_{\prob}
 \left( b_{0}^{2}\right),
\quad
\mathbb{E}_{n}
\left[ T_{n} \right]
=
O_{\mathbb{P}}
\left(
b_{0}^{4}
+
\frac{1}{nb_{0}^{d}}
\right).
\end{equation*}
This proves the validity of the lemma. $\square$

\subsection*{Proof of Lemma \ref{VarU}}

Define $e_{in}(x)=\widehat{m}_{in}(x)-m(x)$, which is such that
\begin{equation*}
U_n(x)
=
U_{n}(x;p)
=
\frac{1}{nb_{1}^{d}h^{p+1}}
\sum_{i=1}^{n}
e_{in}^{p}\left( x\right)
K_{1}\left(
\frac{X_{i}-x}{b_{1}} \right)
K_{2}^{(p)} \left(
\frac{Y_{i}-\epsilon-m(x) }{h}\right).
\end{equation*}
Let
$$
U_{n}(p)
 =
 \int
 \mathds{1}
 \left( x\in \mathcal{X}\right)
 U_{n}\left(x\right)
dx,
$$
so that $S_n=U_n(1)$ and $T_n=U_n(2)$. Observe now that
  Lemma \ref{IndepU} gives
\begin{eqnarray}
\nonumber
 \lefteqn{
 \Var_n\left(U_{n}(p)\right)}
\\\nonumber
&=&
\Var_n\left(
\int\mathds{1}
\left( x\in \mathcal{X} \right)
U_{n}\left(x\right) dx\right)
 =
 \int \int
 \mathds{1}
 \left(\left(
x_{1},x_{2}\right)
\in \mathcal{X}^{2}\right)
\Cov_n\left(
U_{n}\left( x_{1}\right)
,
U_{n}\left( x_{2}\right)
\right)
dx_{1}dx_{2}
\\\nonumber
&=&
\int \int
 \mathds{1}
\left(
\left( x_{1},x_{2}\right)
\in \mathcal{X}^{2},
\left\Vert x_{2}-x_{1}
\right\Vert
\leq Cb_{0}\vee
b_{1}\right)
\Cov_n\left(
U_{n}\left( x_{1}\right)
,
U_{n}\left( x_{2}\right)
\right)
dx_{1}dx_{2}
\\\nonumber
&\leq&
\int \int
\mathds{1}
\left(
\left( x_{1},x_{2}\right)
\in \mathcal{X}^{2},
\left\Vert x_{2}-x_{1}
\right\Vert
\leq Cb_{0}\vee
b_{1}\right)
\Var_n^{1/2}
\left( U_{n}
\left( x_{1}\right)
\right)
\Var_n^{1/2}
\left(
U_{n}\left(
 x_{2}\right)
\right)
dx_{1}dx_{2}
\\\nonumber
&\leq &
\frac{1}{2}
\int \int
\mathds{1}
\left(
\left(
x_{1},x_{2}\right)
\in \mathcal{X}^{2},
 \left\Vert
x_{2}-x_{1}
\right\Vert
\leq Cb_{0}\vee b_{1}\right)
\left\{
\Var_n\left(
U_{n}\left( x_{1}\right)
\right)
+
\Var_n\left(
U_{n}\left( x_{2}\right)
\right)
\right\}
 dx_{1}dx_{2}
\\
&\leq &
C\left(b_{0}^d\vee b_{1}^d\right)
\int
\mathds{1}\left(
x\in
\mathcal{X}\right)
\Var_n\left( U_{n}
\left( x\right) \right)
dx,
\label{VarUb}
\end{eqnarray}
Moreover, we have
\begin{eqnarray}
\nonumber
\lefteqn{
\left( nb_{1}^{d}h^{p+1}\right)^{2}
\int
\mathds{1}
\left( x\in \mathcal{X} \right)
\Var_n\left(
U_{n}\left( x\right) \right)dx }
\\\nonumber
&=&
\sum_{i=1}^{n}
\int
\mathds{1}
\left( x\in \mathcal{X} \right)
K_{1}^{2}
\left(\frac{X_{i}-x}{b_{1}}\right)
 \Var_n\left(
W_{in}(x;p)\right)
 dx
\\\nonumber
&&
+
\sum_{1\leq i_{1}\neq i_{2}\leq n}
\int \mathds{1}
\left(
x\in \mathcal{X}
\right)
K_{1}
\left( \frac{X_{i_{1}}-x}{b_{1}}
\right)
K_{1}
\left( \frac{X_{i_{2}}-x}{b_{1}}\right)
\Cov_n \left(W_{i_{1}n}(x;p)
 ,
W_{i_{2}n}(x;p) \right) dx,
\\
\label{VarU1}
\end{eqnarray}
where
$$
W_{in}(x;p)
=
 e_{in}^{p}
 \left(x\right)
 K_{2}^{(p)} \left(\frac{
Y_{i}-\epsilon-m(x) }{h}\right).
$$
The first term in (\ref{VarU1}) yields, by Lemma \ref{MomderfK}
and Lemma \ref{BoundEspmin},
\begin{eqnarray}
\nonumber
\lefteqn{
\sum_{i=1}^{n}
\int\mathds{1}\left( x\in
\mathcal{X} \right)
 K_{1}^{2}
 \left(\frac{X_{i}-x}{b_{1}}\right)
\Var_n\left(W_{in}(x;p)\right)
dx
}
\\\nonumber
&\leq&
\sum_{i=1}^{n}
\int\mathds{1}
\left( x\in \mathcal{X}
\right)
K_{1}^{2}
\left(\frac{X_{i}-x}{b_{1}}\right)
\esp_n
\left[
e_{in}^{2p}\left( x\right)
 K_{2}^{(p)}
 \left( \frac{
Y_{i}-\epsilon-m(x) }{h}\right)^2
\right]
dx
\\\nonumber
&=&
\sum_{i=1}^{n}
\int
\mathds{1}
\left( x\in \mathcal{X} \right)
K_{1}^{2}
\left(\frac{X_{i}-x}{b_{1}}\right)
\esp_n
\left[
 e_{in}^{2p}\left( x\right)
\esp_{in}
\left[
K_{2}^{(p)}
\left(\frac{
Y_{i}-\epsilon-m(x)}{h}\right)^2
\right]
\right] dx
\\\nonumber
&\leq&
C h
\sum_{i=1}^{n}
\int
\mathds{1}
\left( x\in \mathcal{X}
\right)
\esp_n
\left[
e_{in}^{2p}\left( x\right)
\right]
K_{1}^{2}
\left(\frac{X_{i}-x}{b_{1}}\right)
dx
\\\nonumber
&\leq&
C h
\sum_{i=1}^{n}
\int
\mathds{1}
\left( x\in \mathcal{X}
\right)
\left(
\beta_{in}^{2p}\left(x\right)
+
\esp_n
\left[
\Sigma_{in}^{2p}(x)
\right]
\right)
K_{1}^{2}
\left(\frac{X_{i}-x}{b_{1}}
\right) dx
\\
&=&
 O_{\mathbb{P}}
 \left(nb_1^{d}h\right)
 \left(
 b_0^{4}
 +
 \frac{1}{nb_0^d}
 \right)^p.
 \label{VarU2}
\end{eqnarray}
For the sum of the conditional covariances in (\ref{VarU1}), set
$$
\widetilde{W}_n(p)
=
\sum_{1\leq i_{1}\neq
i_{2}\leq n}
\int
\mathds{1}
\left( x\in \mathcal{X} \right)
K_1\left( \frac{X_{i_{1}}-x}{b_{1}}\right)
K_1\left(\frac{X_{i_{2}}-x}{b_{1}}\right)
\Cov_n
\left(W_{i_{1}n}(x;p)
 ,
W_{i_{2}n}(x;p)
\right)
dx.
$$
We need to bound this term
 for $p\in[1,2]$. Since
$$
 W_{in}(x;p)
 =
 \left(
 \beta_{in}(x)
 +
 \Sigma_{in}(x)
 \right)^{p}
 \left(x\right)
 K_{2}^{(p)}
 \left(\frac{
 Y_{i}-\epsilon-m(x) }{h}
\right),
$$
the independence of the the $Y_j$'s gives, for any $i_1\neq i_2$,
\begin{eqnarray}
\nonumber
\lefteqn{
\Cov_n
\left(
W_{i_{1}n}(x;1)
,
W_{i_{2}n}(x;1)
\right)
}
\\\nonumber
&=&
\beta_{i_{1}n}(x)
\Cov_n
\left[
 K_{2}^{(1)}
 \left(
\frac{Y_{i_{1}}-\epsilon-m(x)}{h}
\right)
,
\Sigma_{i_{2}n}(x)
K_{2}^{(1)}
\left( \frac{Y_{i_{2}}-\epsilon-m(x)}{h} \right)
\right]
\\\nonumber
&&
+
\beta_{i_{2}n}(x)
\Cov_n
\left[
\Sigma_{i_{1}n}(x)
K_{2}^{(1)}
\left(\frac{Y_{i_{1}}-\epsilon-m(x) }{h} \right)
 ,
K_{2}^{(1)}
\left(\frac{Y_{i_{2}}-\epsilon-m(x)}{h}\right)
\right]
\\
&&
+
\Cov_n
\left[
\Sigma _{i_{1}n}(x)
K_{2}^{(1)}
\left(
\frac{Y_{i_{1}}-\epsilon-m(x)}{h}\right)
,
\Sigma_{i_{2}n}(x)
K_{2}^{(1)}
\left(\frac{Y_{i_{2}}-\epsilon-m(x)}{h}\right)
\right].
\label{CovU(1)}
\end{eqnarray}
Moreover, it is clear that the results of Lemma \ref{MomderfK}
remain valid with $\esp_n[\cdot]$, since
$\esp_n[A]=\esp_n[\esp_{in}[A]]$, where $\esp_{in}[\cdot]$
represents the conditional mean given
$\left(X_1,\ldots,X_n,\varepsilon_k,k\neq i\right)$. Therefore,
since $K_0(\cdot)$ is bounded under $(H_7)$, this yields, by
$(H_4)$ and  Lemma \ref{Ordregn},
\begin{eqnarray}
\nonumber
\lefteqn{
\left\vert
\beta_{i_{1}n}(x)
\Cov_n
\left[
K_{2}^{(1)}
\left(\frac{ Y_{i_{1}}-\epsilon-m(x)}{h}\right)
,
\Sigma_{i_{2}n}(x)
K_{2}^{(1)} \left(
\frac{Y_{i_{2}}-\epsilon-m(x)}{h}
\right)
\right]
\right\vert
}
\\\nonumber
&=&
\left\vert
\frac{\beta_{i_{1}n}(x)}
{nb_0^d\widehat{g}_n(x)}
K_0\left(\frac{X_{i_{2}}-x}{b_{0}}\right)
\esp_n
\left[
\varepsilon_{i_1}
K_{2}^{(1)}
\left(
\frac{Y_{i_{1}}-\epsilon-m(x)}{h}\right)
\right]
\esp_n
\left[
K_{2}^{(1)}
\left(\frac{Y_{i_{2}}-\epsilon-m(x)}{h}\right)
\right]
\right\vert.
\\
&=&
 O_{\prob}
 \left(\frac{h^4}{nb_0^d}\right)
 \left\vert
\beta_{i_{1}n}(x)
\right\vert,
\label{CovU1}
\end{eqnarray}
uniformly in $x$, $i_1$ and $i_2$.  We also have
\begin{eqnarray*}
\lefteqn{
\Cov_n
\left[
\Sigma_{i_{1}n}(x)
 K_{2}^{(1)}
 \left(
\frac{Y_{i_{1}}-\epsilon-m(x)}{h}
\right)
,
\Sigma _{i_{2}n}(x)
K_{2}^{(1)}
\left(\frac{Y_{i_{2}}-\epsilon-m(x)}{h}
\right)
\right]
}
\\
&=&
\frac{\esp
\left[\varepsilon^2\right]
}{(nb_0^d\widehat{g}_n(x))^2}
 \sum_{i_3=1\atop i_3\neq
i_1,i_2}^n
K_0^2\left(\frac{X_{i_{3}}-x}{b_{0}}\right)
\esp_n
\left[
 K_{2}^{(1)}
 \left( \frac{Y_{i_{1}}-\epsilon-m(x)}{h}
\right)
\right]
\esp_n
\left[
 K_{2}^{(1)}
 \left( \frac{Y_{i_{2}}-\epsilon-m(x)}{h}
\right)
\right]
\\
&&
+
\frac{1}{(nb_0^d\widehat{g}_n(x))^2}
 K_0\left(\frac{X_{i_{1}}-x}{b_{0}}\right)
 \esp_n
 \left[
 \varepsilon_{i_1}
 K_{2}^{(1)}
 \left( \frac{Y_{i_{1}}-\epsilon-m(x)}{h}
\right)
\right]
\\
&&\;\;\;\;\;\;\;\;\;\;\;\;\;\;\;\;\;\;\;\;\;
 \times
K_0\left(\frac{X_{i_{2}}-x}{b_{0}}\right)
 \esp_n
 \left[
 \varepsilon_{i_2}
 K_{2}^{(1)}
 \left( \frac{Y_{i_{2}}-\epsilon-m(x)}{h}
\right)
\right].
\end{eqnarray*}
This gives, by $(H_4)$, $(H_7)$, Lemma \ref{Ordregn} and Lemma
\ref{MomderfK},
\begin{eqnarray*}
 \lefteqn{
\left\vert
\Cov_n
\left[
\Sigma_{i_{1}n}(x)
K_{2}^{(1)} \left(
\frac{Y_{i_{1}}-\epsilon-m(x)}{h} \right)
,
 \Sigma _{i_{2}n}(x)
K_{2}^{(1)}
\left(\frac{Y_{i_{2}}-\epsilon-m(x)}{h} \right)
\right]
\right\vert
 }
\\
&=&
O_{\prob}
\left(\frac{h^{4}}{(nb_{0}^{d})^2}\right)
\sum_{i_3=1}^n
K^2\left(\frac{X_{i_{3}}-x}{b_{0}}\right)
+
O_{\prob}
\left(\frac{h^{4}}{(nb_{0}^{d})^2}\right)
=
O_{\prob}
\left(\frac{h^{4}}{nb_{0}^{d}}\right),
\end{eqnarray*}
 uniformly  for any  $i_{1}\neq i_{2}$,
 $x_1$ and $x_2$.
Collecting this result, (\ref{CovU1}) and (\ref{CovU(1)}), it
follows, using  Lemma \ref{Ordregn} and taking $p_1=1$ in Lemma
\ref{BoundEspmin}-(\ref{Boundbetain}),
\begin{eqnarray*}
\nonumber
\lefteqn{
O_{\prob}
\left(\frac{nb_{0}^{d}}{h^{4}}\right)
\widetilde{W}_n(1)
}
\\\nonumber
&=&
\sum_{1\leq i_{1}\neq i_{2}\leq n}
\int
\mathds{1}
\left( x\in
\mathcal{X} \right)
 \left\vert
 \beta_{i_{1}n}(x)
 K_{1}\left(
\frac{X_{i_{1}}-x}{b_{1}} \right)
K_{1}\left(
\frac{X_{i_{2}}-x}{b_{1}}\right)
\right\vert
dx
\\\nonumber
&&
+
\sum_{1\leq i_{1}\neq i_{2}\leq n}
 \int
 \mathds{1}
 \left(
x\in \mathcal{X} \right)
\left\vert
K_{1}\left(
\frac{X_{i_{1}}-x}{b_{1}}
\right)
K_{1}\left(
\frac{X_{i_{2}}-x}{b_{1}}\right)
\right\vert
dx
\\\nonumber
&\leq&
O_{\prob}\left(nb_1^d\right)
\sum_{i=1}^n
\int \mathds{1}
\left( x\in \mathcal{X} \right)
\left\vert
\beta_{in}(x)
K_{1}\left( \frac{X_{i}-x}{b_{1}}
\right)
\right\vert
 dx
+
O_{\prob}\left(n^2b_1^{2d}\right)
\\
&=&
O_{\prob}\left(n^2b_1^{2d}\right)
\left(
b_{0}^{2}
\right)
+
O_{\prob}\left(n^2b_1^{2d}\right)
=
O_{\prob}\left(n^2b_1^{2d}\right).
\end{eqnarray*}
Combining  this result with (\ref{VarU2}), (\ref{VarU1}) and
(\ref{VarUb}), we arrive at
\begin{eqnarray*}
\lefteqn{
\Var_n\left(S_n\right)
=
\Var_n\left(U_n(1)\right)
}
\\
&\leq&
O_{\prob}
\left(b_{0}^d\vee b_{1}^d\right)
\times
\frac{1}
{\left(nb_1^dh^2\right)^2}
 \left[
 nb_1^{d}h
 \left(
 b_0^{4}
 +
 \frac{1}{nb_{0}^{d}}
 \right)
 +
\widetilde{W}_n(1)
 \right]
\\
&=&
O_{\prob}
\left(b_{0}^d\vee b_{1}^d\right)
\times
\frac{1}
{\left(nb_1^dh^2\right)^2}
\left[
 nb_1^{d}h
 \left(
 b_0^{4}
 +
 \frac{1}{nb_0^d}
 \right)
 +
 \frac{nb_1^{2d}h^4}{b_0^d}
\right]
\\
&=&
O_{\prob}
\left(b_{0}^d\vee b_{1}^d\right)
\left[
\frac{1}{nb_1^dh^3}
\left(
 b_0^{4}
 +
 \frac{1}{nb_0^{d}}
 \right)
 +
 \frac{1}{nb_0^d}
\right].
\end{eqnarray*}
 This proves the first result of the Lemma.

 For the second, we also have  by (\ref{VarUb}),
 (\ref{VarU1}) and (\ref{VarU2}),
\begin{eqnarray*}
\Var_n\left(T_n\right)
=
\Var_n\left(U_n(2)\right)
=
\frac{O_{\prob}
\left(b_{0}^d\vee b_{1}^d\right)
}{\left(nb_1^dh^3\right)^2}
\left[
 nb_1^{d}h
 \left(
 b_{0}^4
 +
 \frac{1}{nb_{0}^{d}}
 \right)^2
 +
\widetilde{W}_n(2)
 \right].
\end{eqnarray*}
Hence the order of $\Var_n\left(T_n\right)$ follows from the
following result
\begin{equation}
\widetilde{W}_n(2)
=
O_{\prob}
\left[
\frac{h^6}{nb_0^d}
\left(n^2b_1^{2d}\right)
\left(
b_0^4
\right)
 +
\frac{h^3}{n^2b_0^{2d}}
\left(n^2b_1^{2d}\right)
\right].
\label{Wn2}
\end{equation}
Indeed, (\ref{Wn2}) and the  equality  before give
\begin{eqnarray*}
\lefteqn{
\Var_n\left(T_n\right)
}
\\
&=&
\frac{O_{\prob}
\left(b_{0}^d\vee b_{1}^d\right)}
{\left(nb_1^dh^3\right)^2}
\left[
nb_1^{d}h
 \left(b_{0}^4
  +
 \frac{1}{nb_{0}^{d}}
 \right)^2
 +
 \frac{h^6}{nb_0^d}
\left(n^2b_1^{2d}\right)
\left(
b_0^4
\right)
 +
\frac{h^3}{n^2b_0^{2d}}
\left(n^2b_1^{2d}\right)
\right]
\\
&=&
O_{\prob}
\left(b_{0}^d\vee b_{1}^d\right)
\left[
\frac{1}{nb_1^dh^5}
 \left(b_{0}^4
 +
 \frac{1}{nb_0^d}
 \right)^2
 +
 \frac{b_0^4}{nb_0^d}
 +
\frac{1}{n^2b_0^{2d}h^3}
\right].
\end{eqnarray*}
 This yields the  second result of the Lemma. We  now prove (\ref{Wn2}).
Observe that for $i_1\neq i_2$, we have
\begin{eqnarray}
\nonumber
\lefteqn{
\Cov_n
\left(
W_{i_{1}n}(x;2)
,
W_{i_{2}n}(x;2)
\right)
}
\\\nonumber
&=&
\Cov_n
\left[
\left(
\beta_{i_{1}n}(x)
+
\Sigma_{i_{1}n}(x)
\right)^2
 K_{2}^{(2)}
 \left(\frac{
Y_{i_1}-\epsilon-m(x)}{h}\right)
,
\left(
\beta_{i_{2}n}(x)
+
\Sigma_{i_{2}n}(x)
\right)^2
 K_{2}^{(2)}
 \left(\frac{
Y_{i_2}-\epsilon-m(x)}{h}\right)
\right]
\\\nonumber
&=&
\beta_{i_{1}n}^2(x)
\Cov_n
\left[
 K_{2}^{(2)}
 \left(
\frac{Y_{i_{1}}-\epsilon-m(x)}{h} \right)
,
\Sigma_{i_{2}n}^2(x)
K_{2}^{(2)}
\left( \frac{Y_{i_{2}}-\epsilon-m(x)}{h} \right)
\right]
\\\nonumber
&&
+
\beta_{i_{2}n}^2(x)
\Cov_n
\left[
\Sigma_{i_{1}n}^2(x)
K_{2}^{(1)}
\left(\frac{Y_{i_{1}}-\epsilon-m(x) }{h} \right)
 ,
K_{2}^{(2)}
\left( \frac{Y_{i_{2}}-\epsilon-m(x)}{h} \right)
\right]
\\\nonumber
&&
+
2\beta_{i_{1}n}^2(x)
\beta_{i_{2}n}(x)
\Cov_n
\left[
 K_{2}^{(2)}
 \left(
\frac{Y_{i_{1}}-\epsilon-m(x)}{h} \right)
,
\Sigma_{i_{2}n}(x)
K_{2}^{(2)}
\left( \frac{Y_{i_{2}}-\epsilon-m(x)}{h} \right)
\right]
\\\nonumber
&&
+
2\beta_{i_{1}n}(x)
\beta_{i_{2}n}^2(x)
\Cov_n
\left[
\Sigma_{i_{1}n}(x)
K_{2}^{(1)}
\left(\frac{Y_{i_{1}}-\epsilon-m(x) }{h} \right)
 ,
K_{2}^{(2)}
\left( \frac{Y_{i_{2}}-\epsilon-m(x)}{h} \right)
\right]
\\\nonumber
&&
+
2\beta_{i_{1}n}(x)
\Cov_n
\left[
\Sigma_{i_{1}n}(x)
K_{2}^{(1)}
\left(\frac{Y_{i_{1}}-\epsilon-m(x) }{h} \right)
 ,
\Sigma_{i_{2}n}^2(x)
K_{2}^{(2)}
\left( \frac{Y_{i_{2}}-\epsilon-m(x)}{h} \right)
\right]
\\\nonumber
&&
+
2\beta_{i_{2}n}(x)
\Cov_n
\left[
\Sigma_{i_{1}n}^2(x)
 K_{2}^{(2)}
 \left(
\frac{Y_{i_{1}}-\epsilon-m(x)}{h} \right)
,
\Sigma_{i_{1}n}(x)
K_{2}^{(2)}
\left( \frac{Y_{i_{2}}-\epsilon-m(x)}{h} \right)
\right]
\\\nonumber
&&
+
4\beta_{i_{1}n}(x)
\beta_{i_{2}n}(x)
\Cov_n
\left[
\Sigma_{i_{1}n}(x)
K_{2}^{(1)}
\left(\frac{Y_{i_{1}}-\epsilon-m(x) }{h} \right)
 ,
 \Sigma_{i_{2}n}(x)
K_{2}^{(2)}
\left( \frac{Y_{i_{2}}-\epsilon-m(x)}{h} \right)
\right]
\\
&&
+
\Cov_n
\left[
\Sigma _{i_{1}n}^2(x)
K_{2}^{(2)} \left(
\frac{Y_{i_{1}}-\epsilon-m(x)}{h} \right)
,
\Sigma_{i_{2}n}^2(x)
K_{2}^{(2)}
\left( \frac{Y_{i_{2}}-\epsilon-m(x)}{h} \right)
\right].
\label{CovU(2)}
\end{eqnarray}
The two-first terms in (\ref{CovU(2)}) are treated similarly,
since they are symmetric. Under $(H_4)$, we have, for any $i_1\neq
i_2$,
\begin{eqnarray*}
\lefteqn{
\beta_{i_{1}n}^2(x)
\Cov_n
\left[
 K_{2}^{(2)}
 \left(
\frac{Y_{i_{1}}-\epsilon-m(x)}{h} \right)
,
\Sigma_{i_{2}n}^2(x)
K_{2}^{(2)}
\left( \frac{Y_{i_{2}}-\epsilon-m(x)}{h} \right)
\right]
}
\\
&=&
\frac{\beta_{i_{1}n}^2(x)}
{\left(nb_0^d\widehat{g}_n(x)\right)^2}
\sum_{1\leq i_3\neq i_2\leq n}
K_0^2\left(\frac{X_{i_3}-x}{b_0}\right)
\Cov_n
\left[
 K_{2}^{(2)}
 \left(
\frac{Y_{i_{1}}-\epsilon-m(x)}{h} \right)
,
\varepsilon_{i_3}^2
K_{2}^{(2)}
\left( \frac{Y_{i_{2}}-\epsilon-m(x)}{h} \right)
\right]
\\
&=&
\frac{\beta_{i_{1}n}^2(x)}
{\left(nb_0^d\widehat{g}_n(x)\right)^2}
K_0^2\left(\frac{X_{i_1}-x}{b_0}\right)
\Cov_n
\left[
 K_{2}^{(2)}
 \left(
\frac{Y_{i_{1}}-\epsilon-m(x)}{h} \right)
,
\varepsilon_{i_1}^2
K_{2}^{(2)}
\left( \frac{Y_{i_{2}}-\epsilon-m(x)}{h} \right)
\right],
\end{eqnarray*}
with, using Lemma \ref{MomderfK},
\begin{eqnarray*}
\lefteqn{
\left\vert
\Cov_n
\left[
 K_{2}^{(2)}
 \left(
\frac{Y_{i_{1}}-\epsilon-m(x)}{h}\right)
,
\varepsilon_{i_1}^2
K_{2}^{(2)}
\left( \frac{Y_{i_{2}}-\epsilon-m(x)}{h}\right)
\right]
\right\vert
}
\\
&\leq&
\left|
 \esp_n
\left[
\varepsilon_{i_1}^2
 K_{2}^{(2)}
 \left(
\frac{Y_{i_{1}}-\epsilon-m(x)}{h} \right)
K_{2}^{(2)}
\left(\frac{Y_{i_{2}}-\epsilon-m(x)}{h}\right)
\right]
\right|
 \\
 &&
 +
 \left|
 \esp_n
 \left[
 \varepsilon_{i_1}^2
 K_{2}^{(2)}
\left(\frac{Y_{i_{2}}-\epsilon-m(x)}{h} \right)
\right]
\esp_n
\left[
 K_{2}^{(2)}
 \left(
\frac{Y_{i_{1}}-\epsilon-m(x)}{h}\right)
\right]
 \right|
\\
&\leq&
 C h^6.
\end{eqnarray*}
Therefore, since $K_0(\cdot)$ is bounded under $(H_7)$,  Lemma
\ref{Ordregn}  gives
\begin{eqnarray}
\nonumber
\lefteqn{
\left\vert
\beta_{i_{1}n}^2(x)
\Cov_n
\left[
 K_{2}^{(2)}
 \left(
\frac{Y_{i_{1}}-\epsilon-m(x)}{h} \right)
,
\Sigma_{i_{2}n}^2(x)
K_{2}^{(2)}
\left( \frac{Y_{i_{2}}-\epsilon-m(x)}{h}\right)
\right]
\right\vert
 }
\\
&\leq&
\frac{Ch^6\beta_{i_{1}n}^2(x)}
{\left(nb_0^d\widehat{g}_n(x)\right)^2}
K_0^2\left(\frac{X_{i_1}-x}{b_0}\right)
=
O_{\prob}
\left(\frac{h^6}{(nb_0^d)^2}\right)
 \beta_{i_{1}n}^2(x),
 \label{CovU21}
\end{eqnarray}
uniformly with respect to $i_1$, $i_2$ and $x$.

\vskip 0.1cm\noindent For the third and the fourth in
(\ref{CovU(2)}), we also have, uniformly for $i_1$, $i_2$ and $x$,
\begin{eqnarray}
\nonumber
\lefteqn{
\left\vert
\beta_{i_{1}n}(x)
\beta_{i_{2}n}^2(x)
\Cov_n
\left[
\Sigma_{i_{1}n}(x)
K_{2}^{(2)}
\left(\frac{Y_{i_{1}}-\epsilon-m(x) }{h}\right)
 ,
K_{2}^{(2)}
\left(\frac{Y_{i_{2}}-\epsilon-m(x)}{h}\right)
\right]
\right\vert
}
\\\nonumber
&=&
\left\vert
\frac{\beta_{i_{1}n}(x)\beta_{i_{2}n}^2(x)}
{nb_0^d\widehat{g}_n(x)}
K_0\left(\frac{X_{i_1}-x}{b_0}\right)
\esp_n
\left[
\varepsilon_{i_2}
K_{2}^{(2)}
\left( \frac{Y_{i_{1}}-\epsilon-m(x)}{h}\right)
\right]
\esp_n
\left[
 K_{2}^{(2)}
 \left(
\frac{Y_{i_{1}}-\epsilon-m(x)}{h} \right)
\right]
\right\vert
\\
&=&
O_{\prob}
\left(\frac{h^6}{nb_0^d}\right)
 \left\vert
 \beta_{i_{1}n}(x)
 \beta_{i_{2}n}^2(x)
 \right\vert.
 \label{CovU22}
\end{eqnarray}
 Further, note that
\begin{eqnarray*}
\lefteqn{
\beta_{i_{1}n}(x)
\Cov_n
\left[
\Sigma_{i_{1}n}(x)
K_{2}^{(2)}
\left(\frac{Y_{i_{1}}-\epsilon-m(x) }{h} \right)
 ,
\Sigma_{i_{2}n}^2(x)
K_{2}^{(2)}
\left(\frac{Y_{i_{2}}-\epsilon-m(x)}{h}\right)
\right]
}
\\
&=&
\beta_{i_{1}n}(x)
\esp_n
\left[
\Sigma_{i_{1}n}(x)
\Sigma_{i_{2}n}^2(x)
K_{2}^{(2)}
\left(\frac{Y_{i_{1}}-\epsilon-m(x)}{h}\right)
K_{2}^{(2)}
\left(\frac{Y_{i_{2}}-\epsilon-m(x)}{h}\right)
\right]
\\
&=&
\beta_{i_{1}n}(x)
\esp_n
\left[
\Sigma_{i_{2}n}^2(x)
K_{2}^{(2)}
\left(\frac{Y_{i_{1}}-\epsilon-m(x)}{h}\right)
\esp_{i_{2}n}
\left[
\Sigma_{i_{1}n}(x)
K_{2}^{(2)}
\left(\frac{Y_{i_{2}}-\epsilon-m(x)}{h}\right)
\right]
\right],
\end{eqnarray*}
where
\begin{eqnarray*}
\lefteqn{
\left\vert
\esp_{i_{2}n}
\left[
\Sigma_{i_{1}n}(x)
K_{2}^{(2)}
\left(\frac{Y_{i_{2}}-\epsilon-m(x)}{h}\right)
\right]
\right\vert
}
\\
&=&
\left\vert
\frac{1}{nb_0^d\widehat{g}_n(x)}
K\left(\frac{X_{i_2}-x}{b_0}\right)
\esp_n
\left[
\varepsilon_{i_2}
K_{2}^{(2)}\left(
\frac{Y_{i_{2}}-\epsilon-m(x)}{h}\right)
\right]
\right\vert
\leq
\frac{Ch^3}{nb_0^d\left|\widehat{g}_n(x)\right|}\;.
\end{eqnarray*}
  Therefore by $(H_7)$ and  Lemma \ref{Ordregn},
  we have, uniformly for $i_1$, $i_2$ and $x$,
\begin{eqnarray}
\nonumber
\lefteqn{
\left\vert
\beta_{i_{1}n}(x)
\Cov_n
\left[
\Sigma_{i_{1}n}(x)
K_{2}^{(2)}
\left(\frac{Y_{i_{1}}-\epsilon-m(x)}{h}\right)
 ,
\Sigma_{i_{2}n}^2(x)
K_{2}^{(2)}
\left(\frac{Y_{i_{2}}-\epsilon-m(x)}{h}\right)
\right]
\right\vert
}
\\\nonumber
&=&
O_{\prob}
\left(\frac{h^3}{nb_0^d}\right)
\left\vert
\beta_{i_{1}n}(x)
\right\vert
\esp_n
\left[
\Sigma_{i_{2}n}^2(x)
\right]
\\\nonumber
&\leq&
 O_{\prob}
 \left(\frac{h^3}{nb_0^d}\right)
 \left\vert
\beta_{i_{1}n}(x)
\right\vert
\times
\frac{\esp\left[\varepsilon^2\right]}
{\left(nb_0^d\widehat{g}_n(x)\right)^2}
\sum_{j=1}^n
K_0^2\left(\frac{X_j-x}{b_0}\right)
\\
&\leq&
O_{\prob}
\left(\frac{h^3}{n^2b_0^{2d}}\right)
\left\vert
\beta_{i_{1}n}(x)
\right\vert,
\label{CovU23}
\end{eqnarray}
We now treat the two last terms in (\ref{CovU(2)}). Observe that
 \begin{eqnarray*}
\lefteqn{
\Cov_n
\left[
\Sigma_{i_{1}n}(x)
 K_{2}^{(2)}
 \left(
\frac{Y_{i_{1}}-\epsilon-m(x)}{h}
\right)
,
\Sigma _{i_{2}n}(x)
K_{2}^{(2)}
\left(\frac{Y_{i_{2}}-\epsilon-m(x)}{h} \right)
\right]
}
\\
&=&
\frac{\esp
\left[\varepsilon^2\right]
}{(nb_0^d\widehat{g}_n(x))^2}
 \sum_{i_3=1\atop \i_3\neq
i_1,i_2}^n
K_0^2\left(\frac{X_{i_{3}}-x}{b_{0}}\right)
\esp_n
\left[
 K_{2}^{(2)}
 \left( \frac{Y_{i_{1}}-\epsilon-m(x)}{h}
\right)
\right]
\esp_n
\left[
 K_{2}^{(2)}
 \left( \frac{Y_{i_{2}}-\epsilon-m(x)}{h}
\right)
\right]
\\
&&
+
\frac{1}{(nb_0^d\widehat{g}_n(x))^2}
 K_0\left(\frac{X_{i_{1}}-x}{b_{0}}\right)
 \esp_n
 \left[
 \varepsilon_{i_1}
 K_{2}^{(2)}
 \left( \frac{Y_{i_{1}}-\epsilon-m(x)}{h}
\right)
\right]
\\
&&\;\;\;\;\;\;\;\;\;\;\;\;\;\;\;\;\;\;\;\;\;\;\;\;\;
 \times
K_0\left(\frac{X_{i_{2}}-x}{b_{0}}\right)
 \esp_n
 \left[
 \varepsilon_{i_2}
 K_{2}^{(2)}
 \left(
 \frac{Y_{i_{2}}-\epsilon-m(x)}{h}
\right)
\right].
\end{eqnarray*}
This gives, by  Lemma \ref{Ordregn}, Lemma \ref{MomderfK} and
uniformly with respect to $i_{1}\neq i_{2}$ and $x$,
\begin{eqnarray}
\nonumber
 \lefteqn{
\left\vert
\beta_{i_{1}n}(x)
\beta_{i_{2}n}(x)
\Cov_n
\left[
\Sigma_{i_{1}n}(x)
K_{2}^{(2)}
\left(
\frac{Y_{i_{1}}-\epsilon-m(x)}{h}
\right)
,
 \Sigma _{i_{2}n}(x)
K_{2}^{(2)}
\left(\frac{Y_{i_{2}}-\epsilon-m(x)}{h}
\right)
\right]
\right\vert
 }
\\\nonumber
&=&
O_{\prob}
\left[
\frac{h^{6}}
{(nb_{0}^{d})^2}
\sum_{i_3=1}^n
K_0^2
\left(\frac{X_{i_{3}}-x}{b_{0}}\right)
+
\frac{h^{6}}{(nb_{0}^{d})^2}
\right]
\left\vert
\beta_{i_{1}n}(x)
\beta_{i_{2}n}(x)
\right\vert
\\
&=&
O_{\prob}
\left(\frac{h^{6}}{nb_{0}^{d}}\right)
\left\vert
\beta_{i_{1}n}(x)
\beta_{i_{2}n}(x)
\right\vert.
\label{CovU23b}
\end{eqnarray}
Moreover,
\begin{eqnarray}
\nonumber
\lefteqn{
\left\vert
\Cov_n
\left[
\Sigma _{i_{1}n}^2(x)
K_{2}^{(2)} \left(
\frac{Y_{i_{1}}-\epsilon-m(x)}{h}\right)
,
\Sigma_{i_{2}n}^2(x)
K_{2}^{(2)}
\left( \frac{Y_{i_{2}}-\epsilon-m(x)}{h}\right)
\right]
\right\vert
}
\\\nonumber
&\leq&
\left\vert
\esp_n
\left[
\Sigma _{i_{1}n}^2(x)
\Sigma_{i_{2}n}^2(x)
K_{2}^{(2)} \left(
\frac{Y_{i_{1}}-\epsilon-m(x)}{h}\right)
K_{2}^{(2)}
\left( \frac{Y_{i_{2}}-\epsilon-m(x)}{h}\right)
\right]
\right\vert
\\\nonumber
&&
+
\left\vert
\esp_n
\left[
\Sigma _{i_{1}n}^2(x)
K_{2}^{(2)} \left(
\frac{Y_{i_{1}}-\epsilon-m(x)}{h}\right)
\right]
\esp_n
\left[
\Sigma _{i_{2}n}^2(x)
K_{2}^{(2)} \left(
\frac{Y_{i_{2}}-\epsilon-m(x)}{h} \right)
\right]
\right\vert,
\\
\label{CovU23bb}
\end{eqnarray}
with, using $(H_4)$, Lemma \ref{MomderfK} and Lemma \ref{Ordregn},
\begin{eqnarray}
\nonumber
\left\vert
\esp_n
\left[
\Sigma _{in}^2(x)
K_{2}^{(2)} \left(
\frac{Y_{i}-\epsilon-m(x)}{h} \right)
\right]
\right\vert
&=&
\left\vert
\esp_n\left[\Sigma _{in}^2(x)\right]
\esp_n
\left[
K_{2}^{(2)}
\left(
\frac{Y_{i}-\epsilon-m(x)}{h} \right)
\right]
\right\vert
\\\nonumber
&\leq&
\frac{Ch^3}
{\left(nb_0^d\widehat{g}_n(x)\right)^2}
\sum_{j=1}^n
K_0^2\left(\frac{X_j-x}{b_0}\right)
\\
&=&
O_{\prob}
\left(\frac{h^3}{nb_0^d}\right),
\label{CovU24}
\end{eqnarray}
uniformly for $i$ and $x$. Moreover, for the first term in Bound
(\ref{CovU(2)}), we have
\begin{eqnarray*}
\lefteqn{
\esp_n
\left[
\Sigma_{i_{1}n}^2(x)
\Sigma_{i_{2}n}^2(x)
K_{2}^{(2)}
\left(\frac{Y_{i_{1}}-\epsilon-m(x)}{h}\right)
K_{2}^{(2)}
\left(\frac{Y_{i_{2}}-\epsilon-m(x)}{h}\right)
\right]
}
\\
&=&
\esp_n
\left[
\Sigma_{i_{1}n}^2(x)
K_{2}^{(2)}
\left(\frac{Y_{i_{2}}-\epsilon-m(x)}{h}\right)
\esp_{i_{1}n}
\left[
\Sigma_{i_{2}n}^2(x)
K_{2}^{(2)}
\left(\frac{Y_{i_{1}}-\epsilon-m(x)}{h}\right)
\right]
\right],
\end{eqnarray*}
where
\begin{eqnarray*}
\lefteqn{
\left\vert
\esp_{i_{1}n}
\left[
\Sigma_{i_{2}n}^2(x)
K_{2}^{(2)}
\left(\frac{Y_{i_{1}}-\epsilon-m(x)}{h}\right)
\right]
\right\vert
}
\\
&=&
\left\vert
\frac{1}{\left(nb_0^d\widehat{g}_n(x)\right)^2}
\sum_{1\leq i_3\neq i_2\leq n}
K_0^2\left(\frac{X_{i_3}-x}{b_0}\right)
\esp_{i_{1}n}
\left[
\varepsilon_{i_3}
K_{2}^{(2)}\left(
\frac{Y_{i_{1}}-\epsilon-m(x)}{h}\right)
\right]
\right\vert
\\
&\leq&
\frac{Ch^3}{nb_0^d\left(\widehat{g}_n(x)\right)^2}\;.
\end{eqnarray*}
 Therefore, since $K_2^{(2)}$ is bounded under $(H_7)$, it
 follows, by Lemma \ref{Ordregn},
 and uniformly with respect to  $x$, $i_1$ and
 $i_2$,
\begin{eqnarray*}
\lefteqn{
\left\vert
\esp_n
\left[
\Sigma_{i_{1}n}^2(x)
\Sigma_{i_{2}n}^2(x)
K_{2}^{(2)}
\left(\frac{Y_{i_{1}}-\epsilon-m(x)}{h}\right)
K_{2}^{(2)}
\left(\frac{Y_{i_{2}}-\epsilon-m(x)}{h}\right)
\right]
\right\vert
}
\\
&=&
O_{\prob}
\left(\frac{h^3}{nb_0^d}\right)
\esp_n
\left[
\Sigma_{i_{2}n}^2(x)
\right]
=
O_{\prob}
\left(\frac{h^3}{n^2b_0^{2d}}\right),
\end{eqnarray*}
Hence from (\ref{CovU24}) and (\ref{CovU23bb}),
 we deduce
$$
\left\vert
\Cov_n
\left[
\Sigma _{i_{1}n}^2(x)
K_{2}^{(2)}
\left(
\frac{Y_{i_{1}}-\epsilon-m(x)}{h}\right)
,
\Sigma_{i_{2}n}^2(x)
K_{2}^{(2)}
\left(\frac{Y_{i_{2}}-\epsilon-m(x)}{h}\right)
\right]
\right\vert
=
O_{\prob}
\left(\frac{h^3}{n^2b_0^{2d}}\right),
$$
 uniformly in $x$, $i_1$ and $i_2$. Collecting this result,
(\ref{CovU23})-(\ref{CovU23b}) and (\ref{CovU21})-(\ref{CovU22}),
it follows then by Equality (\ref{CovU(2)}),
\begin{eqnarray}
\nonumber
\lefteqn{
\left\vert\widetilde{W}_n(2)\right\vert
=
\left\vert
\sum_{1\leq i_{1}\neq i_{2}\leq n} \int \mathds{1}
\left( x\in \mathcal{X} \right) K_1\left(
\frac{X_{i_{1}}-x}{b_{1}}\right)
K_1\left(\frac{X_{i_{2}}-x}{b_{1}}\right) \Cov_n
\left(W_{i_{1}n}(x;2)
 ,
W_{i_{2}n}(x;2) \right) dx
\right\vert }
\\\nonumber
&=&
O_{\prob}
\left(\frac{h^6}{(nb_0^d)^2}\right)
\sum_{1\leq
i_{1}\neq i_{2}\leq n}
 \int
 \mathds{1}
 \left(
x\in \mathcal{X} \right)
\left\vert
\beta_{i_{1}n}^2(x)
K_{1}\left(
\frac{X_{i_{1}}-x}{b_{1}}
\right)
K_{1}\left(
\frac{X_{i_{2}}-x}{b_{1}}\right)
\right\vert
dx
\\\nonumber
&&
+
O_{\prob}
\left(\frac{h^6}{nb_0^d}\right)
\sum_{1\leq i_{1}\neq i_{2}\leq n}
 \int
 \mathds{1}
 \left(
x\in \mathcal{X} \right)
\left\vert
\beta_{i_{1}n}(x)
 \beta_{i_{2}n}^2(x)
K_{1}\left(
\frac{X_{i_{1}}-x}{b_{1}}
\right)
K_{1}\left(
\frac{X_{i_{2}}-x}{b_{1}}\right)
\right\vert
dx
\\\nonumber
&&
+
O_{\prob}
\left(\frac{h^3}{n^2b_0^{2d}}\right)
\sum_{1\leq i_{1}\neq i_{2}\leq n}
 \int
 \mathds{1}
 \left(
x\in \mathcal{X} \right)
\left\vert
\beta_{i_{1}n}(x)
K_{1}\left(
\frac{X_{i_{1}}-x}{b_{1}}
\right)
K_{1}\left(
\frac{X_{i_{2}}-x}{b_{1}}\right)
\right\vert
dx
\\\nonumber
&&
+
O_{\prob}
\left(\frac{h^6}{nb_0^d}\right)
\sum_{1\leq i_{1}\neq i_{2}\leq n}
 \int
 \mathds{1}
 \left(
x\in \mathcal{X} \right)
\left\vert
\beta_{i_{1}n}(x)
 \beta_{i_{2}n}(x)
K_{1}\left(
\frac{X_{i_{1}}-x}{b_{1}}
\right)
K_{1}\left(
\frac{X_{i_{2}}-x}{b_{1}}\right)
\right\vert
dx
\\
&&
+
O_{\prob}
\left(\frac{h^3}{n^2b_0^{2d}}\right)
\sum_{1\leq i_{1}\neq i_{2}\leq n}
 \int
 \mathds{1}
 \left(
x\in \mathcal{X} \right)
\left\vert
K_{1}\left(
\frac{X_{i_{1}}-x}{b_{1}}
\right)
K_{1}\left(
\frac{X_{i_{2}}-x}{b_{1}}\right)
\right\vert
dx.
\label{CovU25}
\end{eqnarray}
Moreover, note that for any integers $p_1$ and $p_2$ in $[0,2]$,
\begin{eqnarray*}
\lefteqn{
\sum_{1\leq i_{1}\neq i_{2}\leq n}
 \int
 \mathds{1}
 \left(
x\in \mathcal{X} \right)
\left\vert
\beta_{i_{1}n}^{p_1}(x)
\beta_{i_{2}n}^{p_2}(x)
K_{1}\left(
\frac{X_{i_{1}}-x}{b_{1}}
\right)
K_{1}\left(
\frac{X_{i_{2}}-x}{b_{1}}\right)
\right\vert
dx
}
\\
&\leq&
 \sum_{1\leq i_{1}\neq i_{2}\leq n}
 \int
 \mathds{1}
 \left(
x\in \mathcal{X} \right)
\left\vert
\beta_{i_{1}n}^{p_1+p_2}(x)
K_{1}\left(
\frac{X_{i_{1}}-x}{b_{1}}
\right)
K_{1}\left(
\frac{X_{i_{2}}-x}{b_{1}}\right)
\right\vert
dx
\\
&&
 +
 \sum_{1\leq i_{1}\neq i_{2}\leq n}
 \int
 \mathds{1}
 \left(
x\in \mathcal{X} \right)
\left\vert
\beta_{i_{2}n}^{p_1+p_2}(x)
K_{1}\left(
\frac{X_{i_{1}}-x}{b_{1}}
\right)
K_{1}\left(
\frac{X_{i_{2}}-x}{b_{1}}\right)
\right\vert
dx.
\end{eqnarray*}
Since $(H_7)$ and Lemma \ref{BoundEspmin}-(\ref{Boundbetain})
give, for $p=p_1+p_2$,
\begin{eqnarray*}
\lefteqn{
\sum_{1\leq i_{1}\neq i_{2}\leq n}
 \int
 \mathds{1}
 \left(
x\in \mathcal{X} \right)
\left\vert
\beta_{i_{1}n}^{p_1+p_2}(x)
K_{1}\left(
\frac{X_{i_{1}}-x}{b_{1}}
\right)
K_{1}\left(
\frac{X_{i_{2}}-x}{b_{1}}\right)
\right\vert
dx
}
\\
&=&
b_1^d
\sum_{1\leq i_{1}\neq i_{2}\leq n}
\int
\mathds{1}
\left(u+b_1X_{i_2}
\in \mathcal{X}
\right)
 \left\vert
\beta_{i_{1}n}^{p_1+p_2}
\left(u+b_1X_{i_2}\right)
 K_1\left(u\right)
 K_1\left(
 \frac{X_{i_2}-u-b_1X_{i_2}}{b_1}
 \right)
 \right\vert
 du
\\
&=&
O_{\prob}\left(nb_1^d\right)
\sum_{i=1}^n
\int
\mathds{1}
\left( x\in \mathcal{X} \right)
\left\vert
\beta_{in}^{p_1+p_2}(x)
K_{1}\left( \frac{X_{i}-x}{b_{1}}\right)
\right\vert
 dx
\\
&=&
O_{\prob}
\left(n^2b_1^{2d}\right)
\left(
b_0^{2p}
\right),
\end{eqnarray*}
it the follows, by this result, the bound above and
(\ref{CovU25}),
\begin{eqnarray*}
\widetilde{W}_n(2)
=
O_{\prob}
\left[
\frac{h^6}{nb_0^d}
\left(n^2b_1^{2d}\right)
\left(
b_0^4
\right)
 +
\frac{h^3}{n^2b_0^{2d}}
\left(n^2b_1^{2d}\right)
\right].
\end{eqnarray*}
This proves (\ref{Wn2}) and then completes the proof of the
Lemma.\eop

\subsection*{Proof of Lemma \ref{IndepU}}

The lemma follows directly from the fact that given $X_{1},\ldots
,X_{n}$, we have $U_{n} \left( x_{1}\right) = \Phi _{1n} \left(
\varepsilon _{i},i\in I_{1}\right)$ and $U_{n}\left( x_{2}\right)
= \Phi_{2n}\left( \varepsilon _{i},i\in I_{2}\right)$,  with an
empty $I_{1}\cap I_{2}$, since the Kernel functions are compactly
supported and $\left\Vert x_{2}-x_{1}\right\Vert \geq Cb_{0}\vee
b_{1}$ for a sufficiently large $C$. $\square $

\subsection*{Proof of Lemma \ref{BoundEspmin}}

Define
\begin{equation*}
V_{n}
=
\frac{1}{nb_{1}^{d}}
\int
\mathds{1}
\left(
x\in \mathcal{X}\right)
\sum_{i=1}^{n}
\left\vert
\beta_{in}^{p_1}\left( x\right)
K_{1}^{p_2}
\left(\frac{X_{i}-x}{b_{1}}\right)
\right\vert
 dx,
\end{equation*}
and
\begin{equation*}
\Delta_{j}\left( x\right)
=
 \left( m\left( X_{j}\right)-m\left(
x\right) \right)
K_0\left( \frac{X_{j}-x}{b_{0}}\right),
\end{equation*}
which is such that, using Lemma \ref{Ordregn},
\begin{eqnarray*}
\lefteqn{
\left\vert
\beta _{in}
\left( x\right)
\right\vert
=
\left\vert
\frac{\sum_{1\leq j\neq i\leq n}
\Delta _{j}\left(
x\right) }
{ nb_{0}^{d}
\widehat{g}_{n}\left( x\right)}
\right\vert
}
\\
&\leq &
\sup_{x\in \mathcal{X}}
\left|
\frac{1}
{\widehat{g}_{n}\left( x\right)}
\right|
\times
\frac{1}{nb_{0}^{d}}
\left(
\left\vert
\sum_{1\leq j\neq i\leq n}
\left( \Delta _{j}\left( x\right)
-
\mathbb{E}
\left[ \Delta_{j}\left( x\right)
\right] \right)
\right\vert
+
\left\vert
\sum_{1\leq j\neq i\leq n}
\mathbb{E}
\left[
\Delta_{j}\left(
x\right)
\right]
\right\vert
\right)
\\
&\leq &
\frac{O_{\mathbb{P}}\left( 1\right)}
{nb_{0}^{d}} \left(
\left\vert
\sum_{1\leq j\neq i\leq n}
\left( \Delta_{j}\left(
x\right)
 -
 \mathbb{E}
 \left[
\Delta_{j}\left( x\right)
\right]
\right)
\right\vert
+
\left\vert
\sum_{1\leq j\neq i\leq n}
\mathbb{E}
\left[
\Delta_{j}\left(
x\right)
\right]
\right\vert
\right),
\end{eqnarray*}
uniformly in $x$. This gives using the Markov Inequality which
ensures that $ A_{n}=O_{\prob}\left(
\mathbb{E}\left|A_{n}\right|\right)$,
\begin{eqnarray}
\nonumber
\lefteqn{
\left\vert
 V_{n}
 \right\vert
 }
\\\nonumber
&\leq&
\frac{O_{\mathbb{P}}\left(1\right)
}{nb_{1}^{d}}\frac{1}
{\left( nb_{0}^{d}\right)^{p_1}}
\sum_{i=1}^{n}
\int
\mathds{1}
\left( x\in \mathcal{X} \right)
 \left\lbrace
\left\vert
\sum_{1\leq j\neq i\leq n}
\left(
\Delta_{j}\left(x\right)
 -
 \mathbb{E}
 \left[
\Delta_{j}
\left( x\right)
\right]
\right)
\right\vert
+
\left\vert
\sum_{1\leq j\neq i\leq n}
\mathbb{E}
\left[
\Delta_{j}\left(
x\right)
\right]
\right\vert
\right\rbrace^{p_1}
 \\\nonumber
 &&
 \;\;\;\;\;\;\;
 \;\;\;\;\;\;\;\;
 \;\;\;\;\;\;\;\;\;\;
 \;\;\;\;\;\;\;\;\;\;\;
\;\;\;\;\;\;\;\;\;\;\;\;\;\;
\times
\mathbb{E}
\left|
 K_{1}^{p_2}
 \left(\frac{X_{i}-x}{b_{1}}\right)
 \right|
  dx
\\
&\leq &
\frac{O_{\mathbb{P}}
\left( 1\right)}
{\left(
nb_{0}^{d}\right)^{p_1}}
\int
\mathds{1}
\left( x\in
\mathcal{X}\right)
\left\{
\mathbb{E}
\left[
\left\vert
\sum_{j=2}^{n}
\left(
\Delta _{j}
\left( x\right)
-
\mathbb{E}
\left[
\Delta_{j}
\left( x\right)
\right]
\right)
\right\vert^{p_1}
\right]
+
\left\vert
\sum_{j=2}^{n}
\mathbb{E}
\left[
\Delta_{j}
\left( x\right)
\right]
\right\vert^{p_1}
\right\}
dx.
\label{Vnp}
\end{eqnarray}
We bound the two resulting integrals in (\ref{Vnp}). For the
first, the  Marcinkiewicz-Zygmund inequality (see e.g Chow and
Teicher, 2003, p. 386), the H\"{o}lder and the  Minkowski
inequalities give
\begin{eqnarray}
\nonumber
\lefteqn{
\int \mathds{1}
\left( x\in \mathcal{X}\right)
\mathbb{E}
\left[
\left\vert
\sum_{j=2}^{n}
\left(
\Delta
_{j}\left( x\right)
- \mathbb{E}
\left[
\Delta _{j}\left( x\right)
\right]
\right)
\right\vert^{p_1}
\right] dx }
\\\nonumber
&\leq &
\int
\mathds{1}
\left( x\in \mathcal{X}\right)
\mathbb{E}^{1/2}
\left[
\left\vert
\sum_{j=2}^{n}
\left( \Delta
_{j}\left( x\right)
-\mathbb{E}
\left[
\Delta _{j}\left( x\right)
\right] \right)
\right\vert^{2p_1}
\right] dx
\\\nonumber
&\leq &
 \int
 \mathds{1}
 \left( x\in \mathcal{X}\right)
\mathbb{E}^{1/2}
\left[
\left\vert
\sum_{j=2}^{n}
\left(
\Delta
_{j}\left( x\right)
-
 \mathbb{E}
 \left[
 \Delta _{j}
 \left( x\right)
\right]
\right)^{2}
\right\vert^{p_1}
\right]
dx
\\\nonumber
&=&
\int
\mathds{1}
\left( x\in \mathcal{X}\right)
\left\{
\mathbb{E}^{1/p_1}
\left[
\left\vert
\sum_{j=2}^{n}
\left( \Delta _{j}\left( x\right)
-
\mathbb{E}
\left[
\Delta_{j}
\left( x\right)
\right]
\right)^{2}
\right\vert
^{p}\right]
\right\}^{p_{1}/2}
dx
\\\nonumber
&\leq&
\int
\mathds{1}
\left( x\in \mathcal{X}\right)
\left\{
\sum_{j=2}^{n}
\mathbb{E}^{1/p_{1}}
\left[
\left\vert
\left( \Delta
_{j}\left( x\right)
-
 \mathbb{E}
 \left[ \Delta _{j}
 \left( x\right)
\right]
\right)^{2}
\right\vert^{p_1}
\right]
\right\}^{p_{1}/2}
dx
\\\nonumber
&\leq &
C \int
\mathds{1}
\left( x\in \mathcal{X}\right)
\left\{
\sum_{j=2}^{n}
\mathbb{E}^{1/p_{1}}
\left[
\Delta _{j}^{2p_1}
\left(
x\right)
\right]
\right\}^{p_{1}/2}dx
\\\nonumber
&=&
C\int
\mathds{1}
\left( x\in \mathcal{X}\right)
\left\{
\sum_{j=2}^{n}
\left[
\int
\left(\left( m(z)-m(x)\right)
K_0\left(\frac{z-x}{b_{0}}\right)
\right)^{2p_1}
g\left( z\right)
 dz
 \right]^{1/p_{1}}
 \right\}^{p_{1}/2}
 dx
 \\\nonumber
&=&
C\int
\mathds{1}
\left( x\in \mathcal{X}\right)
\left\{
\sum_{j=2}^{n}
\left[ b_{0}^{d}
\int
\left(\left( m(x+b_{0}u)-m(x)\right)
K_0\left(
u\right)
\right)^{2p}
g\left( x+b_{0}u\right)
du \right]^{1/p_{1}}
\right\}^{p_{1}/2}
dx
\\
&\leq &
C\left\{
n\left[ b_{0}^{d}b_{0}^{2p_{1}}\right]^{1/p_{1}}
\right\}^{p_{1}/2}
=
 O\left(
 \left( n^{p_1}b_{0}^{d}\right)^{1/2}b_{0}^{p_1}
 \right).
\label{Vnp1}
\end{eqnarray}

For the second resulting integral in (\ref{Vnp}), we have, since
the  $\Delta_j(x)$'s are identically distributed,
\begin{eqnarray*}
\lefteqn{
\int
\mathds{1}
\left( x\in \mathcal{X}\right)
\left\vert
\sum_{j=2}^{n}
\mathbb{E}
\left[
\Delta_{j}
\left(
x\right)
\right]
\right\vert^{p_1}
dx
\leq
n^{p_1}
\int
\mathds{1}
\left( x\in \mathcal{X}\right)
\left\vert
\mathbb{E}
\left[
\Delta_{1}
\left( x\right)
\right]
\right\vert^{p_1}
dx
}
\\
&\leq&
 n^{p_1}
 \int
 \mathds{1}
 \left( x\in \mathcal{X}\right)
 \left\vert
b_{0}^{d}
\int
\left(
 m\left( x+b_{0}u\right) -m\left( x\right)
\right)
g\left( x+b_{0}u\right)
K_0\left( u\right)
du\right\vert^{p_1}
dx
\\
&\leq &
C n^{p_1}
\left[
\left( b_{0}^{d}
\times b_{0}^{2}\right)^{p_1}
\right]
=
O\left( nb_{0}^{d+2}\right)^{p_1},
\end{eqnarray*}
using $\int\!uK_0(u)du=0$ and the fact that expect for those $x$
at a distance $O(b_{0})$ of the boundaries of $\mathcal{X}$, we
have for all $u$ in the support of $K_0\left( \cdot \right) $,
\begin{eqnarray*}
\lefteqn{
\left( m\left( x+b_{0}u\right)
-
 m\left( x\right) \right)
g\left( x+b_{0}u\right)
}
\\
&=&
b_{0}\left( m^{\left( 1\right)}
\left( x\right) u^{T}+b_{0}u
\int_{0}^{1}
\left( 1-t\right)
m^{\left( 2\right)}
\left(
x+tb_{0}u\right)
dt
u^{T}\right)
\left( g\left( x\right)
+
 b_{0}
\int_{0}^{1}g^{\left( 1\right)
} \left( x+tb_{0}u\right)
dtu^{T}\right).
\end{eqnarray*}
Substituting the order  in the bound above and (\ref{Vnp1}) in
(\ref{Vnp}), we obtain
\begin{eqnarray*}
 V_{n}
=
\frac{O_{\mathbb{P}}
\left( 1\right)}
{ \left(
nb_{0}^{d}\right)^{p_1}}
\left[
\left(n^{p_1}b_{0}^{d}\right)^{1/2}b_{0}^{p_1}
+
\left(nb_{0}^{d+2}\right)^{p_1}
\right]
=
O_{\prob}
\left(
b_0^{2p_{1}}
\right),
\end{eqnarray*}
since under $(H_9)$, we have  $b_0^{d}/(nb_0^{2d})^p=
O(b_0^{2p})$, for all $p$ in $[0,6]$. This proves
(\ref{Boundbetain}).

\vskip 0.3cm  Let now turn to (\ref{BoundSigmain}). The
H\"{o}lder, the Marcinkiewicz-Zygmund and the Minkowski
inequalities give
\begin{eqnarray*}
\lefteqn{
\sum_{i=1}^{n}
\int\mathds{1}
\left( x\in \mathcal{X}
\right)
\esp_n
\left|
\Sigma_{in}^{p_1}(x)
K_{1}^{p_2}
\left(\frac{X_{i}-x}{b_{1}}\right)
\right|
dx
}
\\
&=&
\left\vert
\sum_{i=1}^{n}
\int
\frac{\mathds{1}
 \left( x\in \mathcal{X}
\right)}
{\left|nb_{0}^{d}\widehat{g}_{n}
\left(
x\right)\right|^{p_1}}
\mathbb{E}_{n}
\left[
\left\vert
\sum_{1\leq
j\neq i\leq n}\varepsilon _{j}
 K_0\left(
\frac{X_{j}-x}{b_{0}}\right)
\right\vert ^{p_1}
\right]
\left\vert
K_{1}^{p_2}
\left(\frac{X_{i}-x}{b_{1}}\right)
\right\vert
dx
\right\vert
\\
&\leq&
\sum_{i=1}^{n}
\int
\frac{\mathds{1}
 \left( x\in \mathcal{X}
\right)}
{\left|nb_{0}^{d}\widehat{g}_{n}
\left(
x\right)\right|^{p_1}}
\mathbb{E}_{n}^{1/2}
\left[
\left\vert
\sum_{1\leq
j\neq i\leq n}\varepsilon _{j}
 K_0\left(
\frac{X_{j}-x}{b_{0}}\right)
\right\vert^{2p_{1}}
\right]
\left\vert
K_{1}^{p_2}
\left(\frac{X_{i}-x}{b_{1}}\right)
\right\vert
 dx
\\
&\leq &
C\sum_{i=1}^{n}
\int
\frac{\mathds{1}
 \left( x\in \mathcal{X}
\right) }
{\left| nb_{0}^{d}\widehat{g}_{n}
\left(
x\right)
\right|^{p_1}}
\left\{
\mathbb{E}_{n}^{1/p_1}
\left[
\left\vert
\sum_{1\leq j\neq i\leq n}\varepsilon _{j}^{2}
K_0^{2}\left(\frac{X_{j}-x}{ b_{0}}\right)
\right\vert^{p_1}
\right]
\right\}^{p_{1}/2}
\left\vert
 K_{1}^{p_2}
 \left(\frac{X_{i}-x}{b_{1}}\right)
 \right\vert
 dx
\\
&\leq &
 C\sum_{i=1}^{n}
 \int
 \frac{\mathds{1}
 \left( x\in \mathcal{X}
\right) }
{\left|
 nb_{0}^{d}\widehat{g}_{n}
\left(
x\right)\right|^{p_1}}
\left\{
\sum_{1\leq j\neq i\leq
n}\mathbb{E}_{n}^{1/p_{1}}
\left[
\left\vert
\varepsilon _{j}^{2}
K_0^{2}
\left( \frac{X_{j}-x }{b_{0}}\right)
\right\vert^{p_1}
\right]
\right\}^{p_{1}/2}
\left\vert
 K_{1}^{p_2}
 \left(\frac{X_{i}-x}{b_{1}}\right)
 \right\vert
dx
\\
&\leq &
C\sum_{i=1}^{n}
\int
\frac{\mathds{1}
 \left( x\in \mathcal{X}
\right)}
 {(nb_{0}^{d})^{p_{1}/2}
\left|\widehat{g}_{n}
\left( x\right)\right|^{p_1}}
\left\{
\frac{1}{nb_0^d}
\sum_{1\leq j\neq i\leq n}
\left\vert
K_0\left(
\frac{X_{j}-x}{b_{0}} \right)
\right\vert
\right\}^{p_{1}/2}
\left\vert
 K_{1}^{p_2}
 \left(\frac{X_{i}-x}{b_{1}}\right)
 \right\vert
dx.
\end{eqnarray*}
It then follows from Lemma \ref{Ordregn} that
\begin{eqnarray*}
\sum_{i=1}^{n}
\int
\mathds{1}
\left( x\in \mathcal{X} \right)
\esp_n
\left|
\Sigma_{in}^{p_1}(x)
 K_{1}^{p_2}
\left(\frac{X_{i}-x}{b_{1}}\right)
 \right|
 dx
 &=&
 O_{\mathbb{P}}
\left( \frac{nb_{1}^{d}}
{\left( nb_{0}^{d}\right)^{p_{1}/2}}\right)
\sup_{x\in \mathcal{X}}\frac{1}{nb_{1}^{d}}
\sum_{i=1}^{n}
\left\vert
 K_{1}^{p_2}
 \left(\frac{X_{i}-x}{b_{1}}\right)
 \right\vert
\\
 &=&
O_{\mathbb{P}}
\left( \frac{nb_{1}^{d}}{\left(
nb_{0}^{d}\right)^{p_{1}/2}}\right).
\end{eqnarray*}
This proves (\ref{BoundSigmain}) and completes the proof of the
Lemma.\eop

\chapter[Simulation study] {Simulation study}

\setcounter{subsection}{0} \setcounter{equation}{0}
\renewcommand{\theequation}{\thesection.\arabic{equation}}

\renewcommand{\thefootnote}{\arabic{footnote}}
\setcounter{footnote}{1}
\setlength{\baselineskip}{.26in}

{\bf Abstract:} In this chapter we present our numerical results.
We analyze and compare the performances of the Kernel density
estimator $\widehat{f}_{1n}$, based on the estimated residuals,
and the ones of the integral Kernel estimator $\widehat{f}_{2n}$.
This comparison is made in the univariate case with a quadratic
model, as described in the next section.  The chapter is organized
as follows. Section 5.1 is devoted to the description of our
simulation framework. Section 5.2 investigates the global study
for the estimators $\widehat{f}_{1n}$ and $\widehat{f}_{2n}$. We
compare in that section the performances of these estimators in
the sense of the Average Integrated Squared Error (AISE). Section
5.3 deals with the pointwise study of our two Kernel estimators,
and compare their  Average Squared Error (ASE), while section 5.4
investigates their asymptotic normality.

\section{Description of our simulation framework}

Let us consider  the following quadratic model
\begin{equation}
Y=3X^2+2X+1+\varepsilon,
\label{quadmodel}
\end{equation}
where $\varepsilon\sim\mbox N(0,1)$ and $X\sim\mbox U[-1,1]$. For
our numerical study, we generate $T=100$ independent samples
$\left(X_{k1},\varepsilon_{k1}\right),\left(X_{k2},\varepsilon_{k2}\right),
\ldots,\left(X_{kn},\varepsilon_{kn}\right)$, $k=1,\ldots,T$, of
size $n=200$, from the model (\ref{quadmodel}).
% Let
%$\widehat{f}_{1}(\epsilon)$ denotes the Kernel density estimator
%of $f(\epsilon)$ based on the estimated residuals, and let
%$\widehat{f}_{2}(\epsilon)$ be the Kernel estimation of
%$f(\epsilon)$ based on the integral transformation. Recall that
%\begin{eqnarray*}
%\widehat{f}_{1n}(\epsilon)
% &=&
%\frac{1} {b_1\sum_{i=1}^n
%\mathds{1}
%\left(X_i\in\mathcal{X}_0\right)}
%\sum_{i=1}^n \mathds{1}
%\left(
%X_i\in \mathcal{X}_0\right)
%K_1\left(
%\frac{\widehat{\varepsilon}_i-\epsilon}{b_1}
%\right),
%%\label{fnchap1}
%\\
%\widehat{f}_{2n}(\epsilon)
%&=&
%\int
%\mathds{1}
% \left(x\in\mathcal{X}\right)
%\widehat{\varphi}_n
%\left(x,\epsilon+\widehat{m}_n(x)\right) dx,
%%\label{fnchap2}
%\end{eqnarray*}
%where the estimators $\widehat{\varepsilon}_i$,
%$\widehat{\varphi}_n(\cdot)$ and $\widehat{m}_n(\cdot)$ are
%defined as in Chapters 3 and 4.
 Define, for any integer $i\in[1,200]$ and
any integer $k\in[1,100]$,
$$
Y_{ki}
 =
 3X_{ki}^2+2X_{ki}+1+\varepsilon_{ki}.
$$
We denote by $\widehat{f^*}_{1k}(\epsilon)$ and
$\widehat{f^*}_{2k}(\epsilon)$ the simulated versions of the
estimators $\widehat{f}_{jn}(\epsilon)$ ($j=1,2$) based on $k^{\rm
th}$ sample
$\left(X_{k1},\varepsilon_{k1}\right),\left(X_{k2},\varepsilon_{k2}\right),
\ldots,\left(X_{kn},\varepsilon_{kn}\right)$.
%\begin{eqnarray*}
%\widehat{f}_{1k}(\epsilon)
% &=&
%\frac{1}{nb_1}
%\sum_{i=1}^n
%K_1\left(\frac{\widehat{\varepsilon}_{ki}-\epsilon}{b_1}\right),
%\\
%\widehat{f}_{2k}(\epsilon)
%&=&
%\int_{-1}^1 \widehat{\varphi}_n
%\left(x,\epsilon+\widehat{m}_{k}(x)\right)
%dx,
%\end{eqnarray*}
%where $\widehat{\varepsilon}_{ki}
%=Y_{ki}-\widehat{m}_{in}(X_{ki})$,
%$$
%\widehat{m}_{k}(x)
%=
% \frac{ \sum_{j=1}^{n}Y_{jk}
%K\left(\frac{X_{jk}-x}{b_{0}}\right)}
%{\sum_{j=1}^{n}
%K\left(\frac{X_{jk}-x}{b_{0}}\right)},
%\quad
%\widehat{m}_{in}(X_{ki})
%=
%\frac{\sum_{j=1\atop j\neq
%i}^nY_{kj}
% K_0\left(\frac{X_{kj}-X_{ki}}{b_0}\right)}
%{\sum_{j=1\atop j\neq i}^n
%K_0\left(\frac{X_{kj}-X_{ki}}{b_0}\right)}.
%
Hence the estimators  $\widehat{f}_{jn}(\epsilon)$ are
approximated by
%the Average Kernel Density Estimators
\begin{equation}
\overline{\widehat{f}}_{jn}(\epsilon)
 =
\frac{1}{T}\sum_{k=1}^T
\widehat{f^*}_{jk}(\epsilon),
\quad j=1,2.
\label{simuldensity}
\end{equation}
For the estimator $\widehat{f}_{1n}(\epsilon)$, we do not make a
truncation and consider $\mathcal{X}_0=[-1,1]$ in the estimator
$\widehat{f}_{1n}(\epsilon)$. We also denote by
$\widetilde{f}_{1n}(\epsilon)$ the Kernel estimator of
$f(\epsilon)$ based on the true residuals, and by
$\widetilde{f}_{2n}(\epsilon)$ the integral Kernel estimator of
$f(\epsilon)$ based on the true regression function. That is,
\begin{eqnarray*}
\widetilde{f}_{1n}(\epsilon)
 &=&
\frac{1}{nb_1}
\sum_{i=1}^n
K_1\left(\frac{\varepsilon_i-\epsilon}{b_1}\right),
\\
\widetilde{f}_{2n}(\epsilon)
&=&
\int_{-1}^1
\widehat{\varphi}_n
\left(x,\epsilon+m(x)\right)
dx,
\end{eqnarray*}
where $m(x)=3x^2+2x+1$, and $\widehat{\varphi_n}$ is defined as in
Chapter 4. Hence we can approximate these estimators by
%the the simulated Average Kernel Density Estimators
\begin{equation}
\overline{\widetilde{f}}_{jn}(\epsilon)
 =
\frac{1}{T}\sum_{k=1}^T
\widetilde{f^*}_{jk}(\epsilon),
\quad
j=1,2,
\label{simuldensityb}
\end{equation}
where $\widetilde{f^*}_{jk}(\epsilon)$ is the Kernel version
 of $\widetilde{f}_{jn}(\epsilon)$ based on the $k^{\rm
th}$ generated sample.

\vskip 0.1cm
 For the choice of the Kernels functions
$K_{\ell}$, $\ell=0,1,2$, we consider the Epanechnikov Kernel
function
$$
K(x)= K_0(x) =K_1(x)
=
\frac{3}{4}
\left(1-x^2\right)
\mathds{1}
\left(|x|\leq 1\right),
$$
and the  the biquadratic or biweight Kernel function
$$
K_2(x)
=
 \frac{15}{16}
 \left(1-x^2\right)^2
 \mathds{1}
\left(|x|\leq 1\right).
$$

Recall that the numerical value of $\widetilde{f}_{2n}(\epsilon)$
is approximated by the Riemann sum
$$
S_n(\epsilon)
=
\sum_{j=1}^p
\widehat{\varphi}_n
\left(x_j,\epsilon+m(x_j)\right)
\left(x_j-x_{j-1}\right),
$$
where $\{x_0,x_1,\ldots,x_p\}$ is a set of points such that
$-1=x_0<x_1<\ldots<x_p=1$. In our setup, the sequence $(x_j)$ is
chosen such that $p=100$ and
$$
x_j=-1+\frac{2j}{p},
\quad j=1,\ldots,p.
$$

\section{Global study}

\begin{figure}[H]
\caption{The  AISE of the  Kernel estimator $\widehat{f}_{1n}$
based on the estimated residuals.}
 \centering
\includegraphics[scale=0.45]{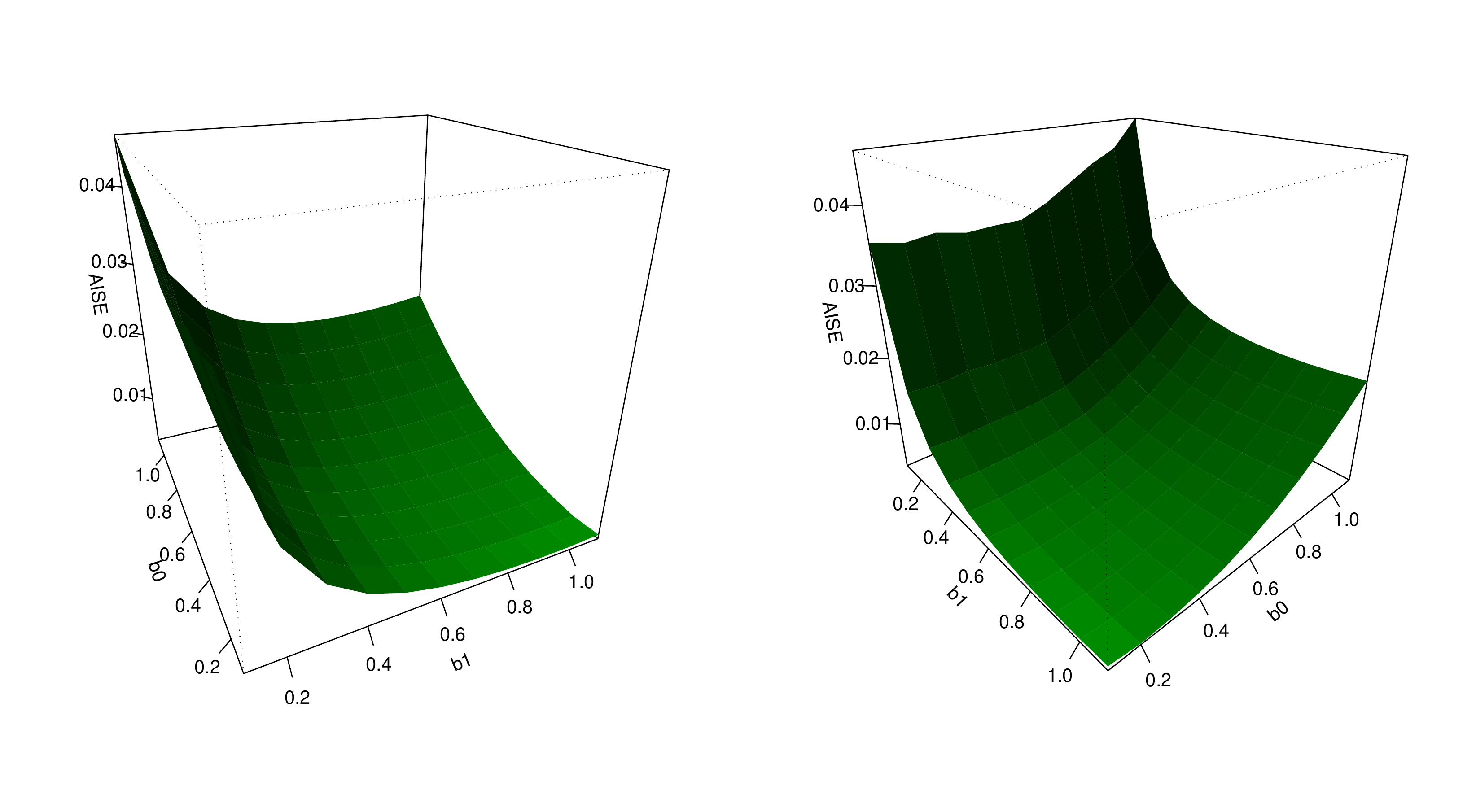}
\label{Aise1}
\end{figure}

\begin{figure}[H]
\caption{The AISE of the integral Kernel estimator
 $\widehat{f}_{2n}$.}
 \centering
\includegraphics[scale=0.45]{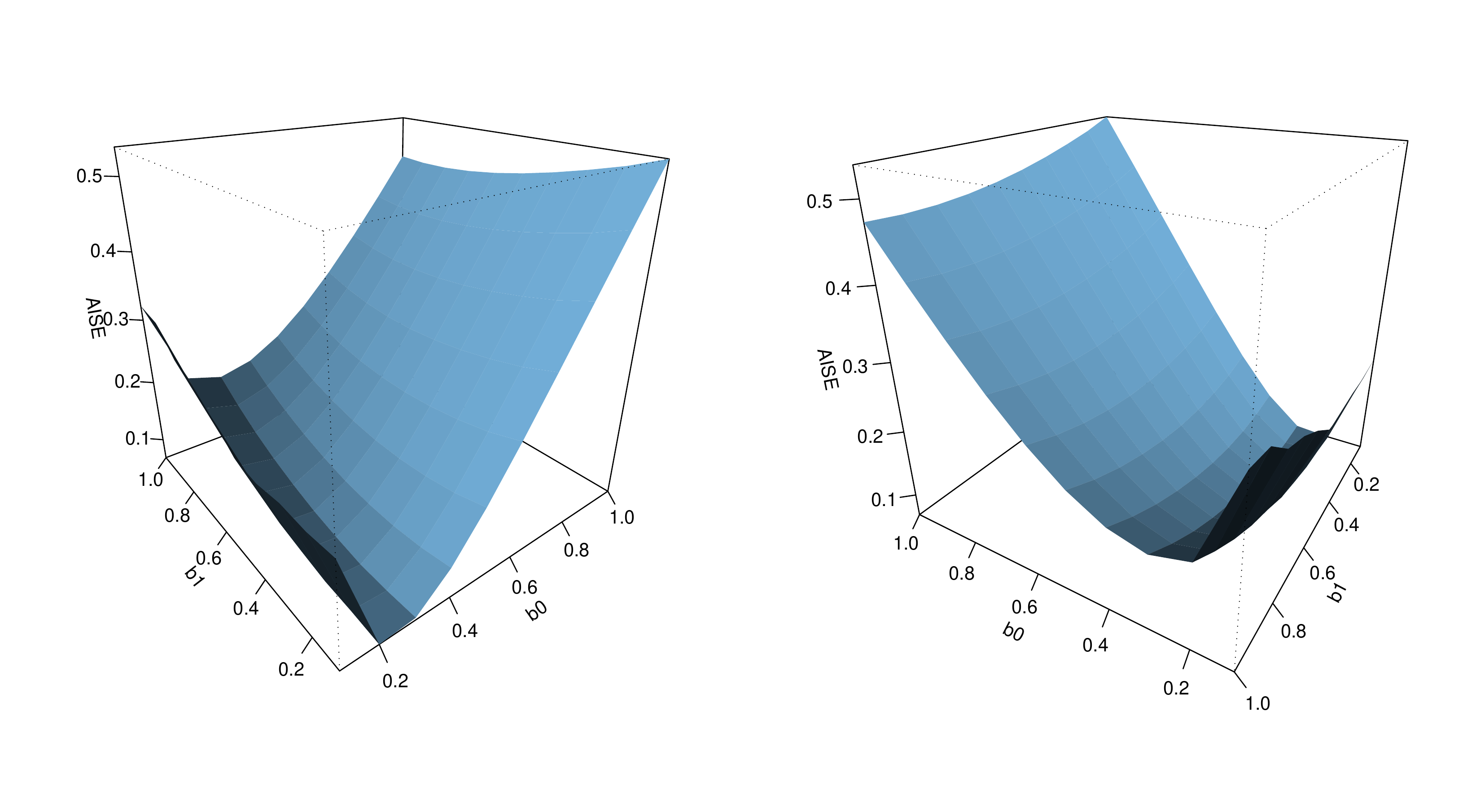}
\label{Aise2}
\end{figure}

In the nonparametric density estimation, it is known that a proper
choice of the bandwidths is crucial for the precision of the
estimator.  In our simulations setup, we first find the simulated
optimal bandwidths for the estimators $\widehat{f}_{jn}$, $j=1,2$.
To that aim, we need to apply the Mean Itegrated Square Error
(MISE) criterion which consists to minimize the quantities
$$
{\rm MISE}(\widehat{f}_{jn})
  =
 \esp
\left[
\int_{-A}^{A}
\left(\widehat{f}_{jn}(t)-f(t)\right)^2
 dt
 \right],
$$
where $[-A,A]$ is a set that contains all the simulated residuals.
In this subsection, we suppose that $[-A,A]=[-5,5]$. For the sake
of simplicity, we assume that $h=b_1$ for the estimator
$\widehat{f}_{2n}(\epsilon)$. Using the $T$ generated samples, we
can approximate
 the MISE of the estimators $\widehat{f}_{jn}$ by the simulated
Average Integrated Square Error (AISE) defined as follows:
\begin{eqnarray*}
{\rm AISE}(\widehat{f}_{jn})
= {\rm
AISE}(\widehat{f}_{jn})(b_1,b_0)
 =
 \frac{1}{T}\sum_{k=1}^T
 \int
\left(\widehat{f^*}_{jk}(t)-f(t)\right)^2
 dt.
\end{eqnarray*}
Now for each $j$, we denote by $(\widehat{b}_{1j},
\widehat{b}_{0j})$  the optimal bandwidths that minimize the above
AISE.
 These bandwidths are simulated from $T=100$ other independent
 samples of size $n=200$ generated from the model
 (\ref{quadmodel}),
 and different of the samples
 $\left(X_{k1},\varepsilon_{k1}\right),\left(X_{k2},\varepsilon_{k2}\right),
\ldots,\left(X_{kn},\varepsilon_{kn}\right)$.
%\begin{eqnarray*}
%(\widehat{b}_{1j},\widehat{b}_{0j})
%\approx \arg\min_{(b_1,b_0)}
%\left[
%\frac{1}{T}\sum_{k=1}^T\sum_{p=1}^L
%\left(\widehat{f^*}_{jk}(u_p)-f(u_p)\right)^2
%\left(u_p-u_{p-1}\right)
% \right],
%\end{eqnarray*}
%where   $\{u_0,u_1,\ldots,u_L\}$ is a partition of the interval
%$[-5,5]$ such that
%$$
%u_j=-5+\frac{10j}{L},
%\quad j=1,\ldots,L.
%$$
% We assume here that $L=100$, and that the bandwidths $b_0$ and $b_1$
%  vary on $[0.1,1.1]$ in the set $\{h_j=0.1+(0.01)\times
%j, 1\leq j\leq 100\}$. Similarly, using the same partition
%$\{u_0,u_1,\ldots,u_L\}$, the optimal bandwidths $\widetilde{b}_j$
%of the
% estimators $\widetilde{f}_{jn}$  satisfy
%\begin{eqnarray*}
%\widetilde{b}_j
%\approx \arg\min_{b_1}
% \left[
%\frac{1}{T}\sum_{k=1}^T\sum_{p=1}^L
%\left(\widetilde{f^*}_{jk}(u_p)-f(u_p)\right)^2
%\left(u_p-u_{p-1}\right)
% \right].
%\end{eqnarray*}
\vskip 0.2cm
 In  Figures \ref{Aise1} and \ref{Aise2}, we plot the
AISE of the estimators $\widehat{f}_{1n}$ and $\widehat{f}_{2n}$
when $b_0$ and $b_1$ vary on $[0.1,1.1]$ in the set
$\{h_j=0.1+(0.01)\times j, 1\leq j\leq 100\}$.
 The first plot shows that the optimal bandwidths
$(\widehat{b}_{11},\widehat{b}_{01})$ for the Kernel estimator
$\widehat{f}_{1n}$ would be achieved when the couple
$(\widehat{b}_{11}, \widehat{b}_{01})$ is very close to $(1,0.2)$,
while the second plot reveals that $(\widehat{b}_{12},
\widehat{b}_{02})$ should be achieved at the neighborhood of
$(0.2,0.2)$.

\vskip 0.1cm These graphical results about the bandwidths
$(\widehat{b}_{1j}, \widehat{b}_{0j})$ are confirmed by the
numerical results of Table \ref{Tab1}, in which we give the
optimal bandwidths for estimators $\widehat{f}_{jn}$ and
$\widetilde{f}_{jn}$, and their corresponding AISE. For the
 $(\widehat{b}_{1j},\widehat{b}_{0j})$, we observe that
$\widehat{b}_{01}$ is approximately as small as
$\widehat{b}_{02}$, while $\widehat{b}_{11}$ and
$\widehat{b}_{12}$ are clearly different, same as
$\widetilde{b}_1$ and $\widetilde{b}_2$. The results of Table
\ref{Tab1} also reveal that ${\rm
AISE}(\widehat{f}_{1n})(\widehat{b}_{11}, \widehat{b}_{01})<{\rm
AISE}(\widehat{f}_{2n})(\widehat{b}_{12}, \widehat{b}_{02})$,
${\rm AISE}(\widehat{f}_{2n})$ being approximately twice as big as
${\rm AISE}(\widehat{f}_{1n})$. This would suggest that for a
judicious choice of the bandwidths $(b_0,b_1)$, the AISE of the
estimator $\widehat{f}_{1n}$ is smaller than the one of
$\widehat{f}_{2n}$. Consequently $\widehat{f}_{1n}$ should be
preferred to $\widehat{f}_{2n}$ for the estimation of p.d.f of the
residuals.
 Moreover, Table \ref{Tab1} shows that
 ${\rm AISE}(\widehat{f}_{1n})(\widehat{b}_{11},
\widehat{b}_{01})\approx {\rm
AISE}(\widetilde{f}_{1n})(\widetilde{b}_1)$ and that  ${\rm
AISE}(\widehat{f}_{2n})(\widehat{b}_{12}, \widehat{b}_{02})<{\rm
AISE}(\widetilde{f}_{2n})(\widetilde{b}_2)$.
%In a semiparametric
%context, Müller, Schick and
% Wefelmeyer (2004) have shown that for the estimation of
% linear functionals
% of the error distribution, the estimators that use the
%nonparametric residuals may have a smaller asymptotic variance
%compared to estimators based upon the true errors. A raison that
%may explain this amazing situation is the fact that the estimators
%$\widehat{f}_{jn}$ built from  the estimated regression function
%better use the fact that the residuals $\varepsilon_i$ have mean
%zero, contrarily to the estimators $\widetilde{f}_{jn}$ based on
%the true regression function.
\newpage

\newpage

\begin{table}[H]
\caption{\label{Tab1} The optimal bandwidths $(\widehat{b}_{1j},
\widehat{b}_{0j})$,  $\widetilde{b}_j$ and their corresponding
AISE when $b_0$ and $b_1$ vary on $[0.1,1.1]$.}
\small{
\begin{center}
\tabcolsep=3pt
\begin{tabular}{|c|c|c|c|c|c|c|c|c|c|}
 \hline
 \multicolumn{3}{|c|}{$\widehat{f}_{1n}$} &
\multicolumn{2}{|c|}{$\widetilde{f}_{1n}$} &
\multicolumn{3}{|c|}{$\widehat{f}_{2n}$} &
\multicolumn{2}{|c|}{$\widetilde{f}_{2n}$}
\\
\cline{1-10}
$\widehat{b}_{11}$ & $\widehat{b}_{01}$ & ${\rm
AISE}(\widehat{b}_{11}, \widehat{b}_{01})$& $\widetilde{b}_1$ &
${\rm AISE}(\widetilde{b}_1)$& $\widehat{b}_{12}$ &
$\widehat{b}_{02}$  & ${\rm
AISE}(\widehat{b}_{12},\widehat{b}_{02})$& $\widetilde{b}_2$ &
${\rm AISE}(\widetilde{b}_2)$
\\
\hline
  0.95 & 0.19 & 0.003141035 & 1.01 & {\bf 0.003083492} &
 0.24  & 0.17 & {\bf 0.006217096} & 0.22& 0.006406112
 \\
\hline
\end{tabular}
\end{center}
}
\end{table}

\noindent Table \ref{Tab1} shows that the optimal first-step
bandwidths for the estimators $\widehat{f}_{1n}$ and
$\widehat{f}_{2n}$ would be very small, as recommended in Wang,
Brown, Cai and Levine (2008).

\begin{figure}[H]
  \caption{From top to bottom, the  $0.95\rm th$ confidence band of
  $\widehat{f}_{1n}$, the p.d.f of N(0,1), the Average Kernel estimator
  $\overline{\widehat{f}}_{1n}$ and
 the  $0.05\rm th$ confidence band of $\widehat{f}_{1n}$  when
 $b_1=\widehat{b}_{11}=0.95$ and $b_0=\widehat{b}_{01}=0.19$.}
\centering
\includegraphics[scale=0.4]{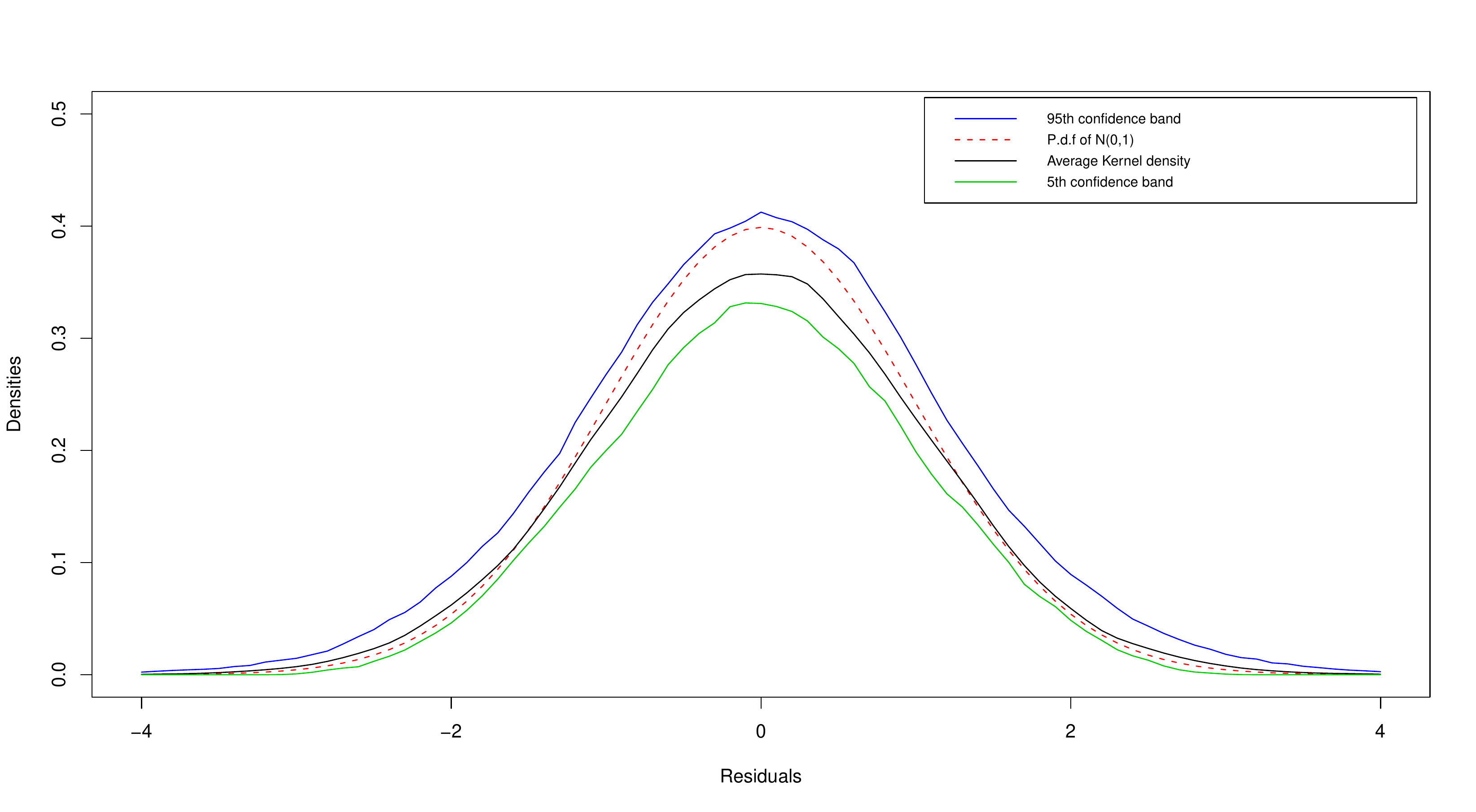}
 \label{Fig1}
\end{figure}

\begin{figure}[H]
\caption{From top to bottom, the  $0.95\rm th$ confidence band of
  $\widehat{f}_{2n}$, the p.d.f of N(0,1), the Average Kernel estimator
  $\overline{\widehat{f}}_{2n}$ and
the  $0.05\rm th$ confidence band of $\widehat{f}_{2n}$ when
 $b_1=\widehat{b}_{12}=0.24$ and $b_0=\widehat{b}_{02}=0.12$.}
\centering
\includegraphics[scale=0.4]{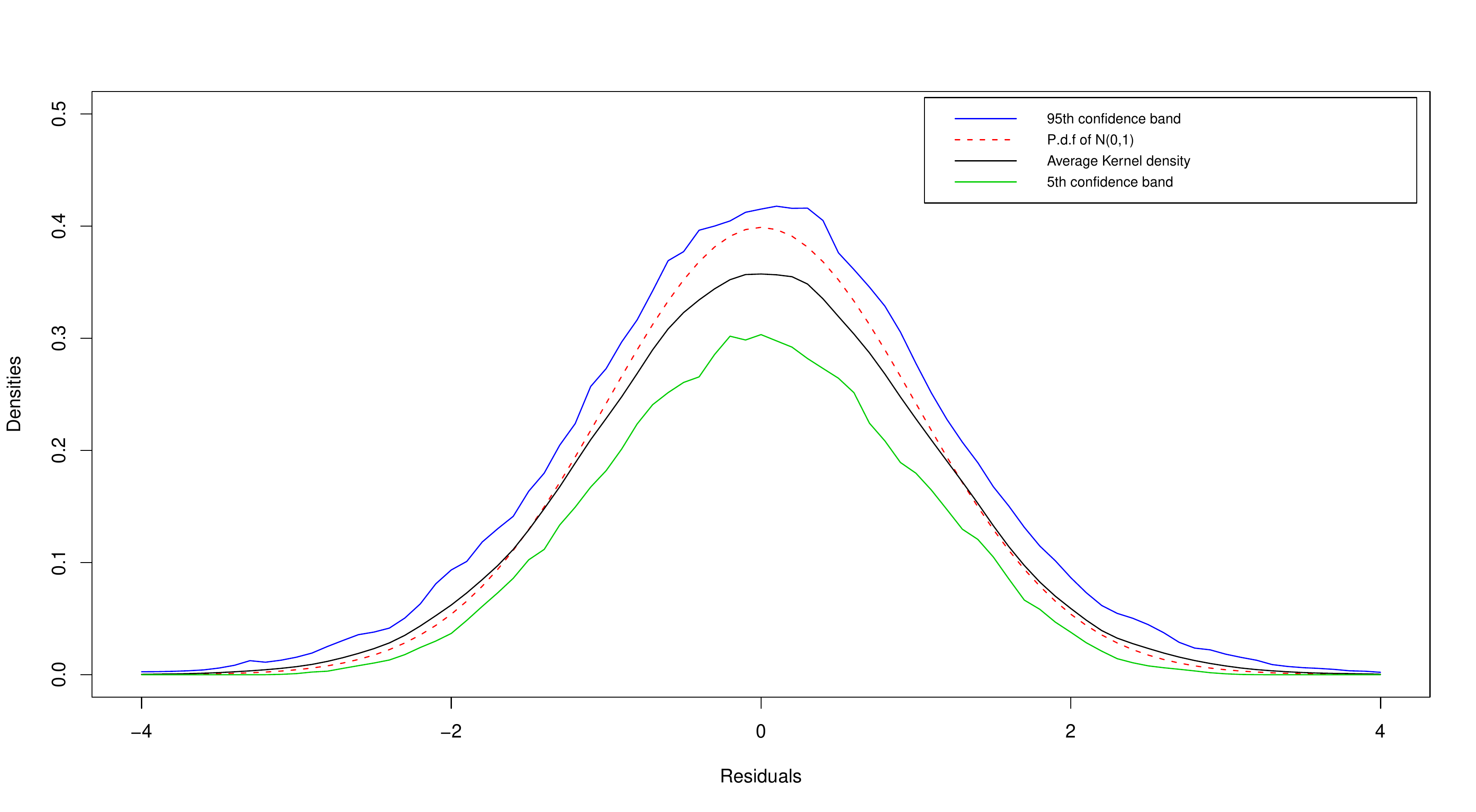}
 \label{Fig2}
\end{figure}

 We now define for $\alpha=0.05$ and $\alpha=0.95$,
 the $\alpha\rm th$ confidence band
 $\widehat{f}_{jn}(\cdot,\alpha)$ of the
 estimator $\widehat{f}_{jn}(\cdot)$ as follows.
 For each $j$ and any
$\epsilon\in[-5,5]$, we consider the $T$ ordered values
$\widehat{f^*}_{j,(k)}(\epsilon)$ of the
$\widehat{f^*}_{jk}(\epsilon)$'s such that
$\widehat{f^*}_{j,(1)}(\epsilon) \leq
\widehat{f^*}_{j,(2)}(\epsilon) \leq \ldots \leq
\widehat{f^*}_{j,(T)}(\epsilon)$. Hence the function
$\widehat{f}_{jn}(\alpha,\cdot)$ is defined as
$$
\widehat{f}_{jn}(\epsilon,\alpha) =
 \widehat{f^*}_{j,(\alpha T)}(\epsilon),
 \quad
 \epsilon\in[-5,5].
 $$
Using the optimal bandwidths $(\widehat{b}_{1j},\widehat{b}_{0j})$
described above, we represent in  Figures \ref{Fig1} and
\ref{Fig2} the Average Kernel estimators
$\overline{\widehat{f}}_{jn}$, the p.d.f of $N(0,1)$, the $0.95\rm
th$ and the $0.05\rm th$ confidence bands of the estimators
$\widehat{f}_{jn}$, $j=1,2$. These plots can be useful for having
a general idea about the confidence interval  of the density $f$.
For example, we  see that for $\epsilon$ varying in the
neighborhood of $0$,  we have $\widehat{f}_{jn}(\epsilon,0.05)
<f(\epsilon)<\widehat{f}_{jn}(\epsilon,0.95)$.

\vskip 0.1cm\noindent In each of the Figures \ref{Fig1} and
\ref{Fig2}, the bias of the estimated density is quite important
around the inflexion point $\epsilon=0$, but the true density
function remains in the good confidence interval. We also notice
that the graphics plotted in Figure \ref{Fig2} are less smooth
than the ones represented in Figure \ref{Fig1}. This may explain
the fact that ${\rm
AISE}(\widehat{f}_{1n})(\widehat{b}_{11},\widehat{b}_{01}) <{\rm
AISE}(\widehat{f}_{2n})(\widehat{b}_{12},\widehat{b}_{02})$.

 \section{Pointwise study}

In this section, we are interested in the pointwise study of the
estimators $\widehat{f}_{jn}(\epsilon)$ and
$\widetilde{f}_{jn}(\epsilon)$. First, we compare the Average
Square Errors (ASE) of these estimators  at the points
$\epsilon=-1,0,1$. In a second time, a comparison of the bias and
variances of these estimators is established, and next their
asymptotic normality is investigated.

\renewcommand{\thesubsection}{\thesection.\arabic{subsection}}
\subsection{Comparison of the ASE}

Let $(\widehat{f}_{jn}(\epsilon),\widehat{f^*}_{jk}(\epsilon))$
and
 $(\widetilde{f}_{jn}(\epsilon), \widetilde{f^*}_{jk}(\epsilon))$
  be as in the previous subsection. We
compare  the pointwise ASE of the estimators
$\widehat{f}_{jn}(\epsilon)$ to the ones of the estimators
$\widetilde{f}_{jn}(\epsilon)$. These ASE are defined as
\begin{eqnarray*}
{\rm ASE}(\widehat{f}_{jn})(\epsilon)
 &=&
 %{\rm ASE}(\widehat{f}_{jn})(\epsilon,b_1,b_0)
 %=
 \frac{1}{T}\sum_{k=1}^T
 \left(\widehat{f^*}_{jk}(\epsilon)-f(\epsilon)\right)^2,
\\
{\rm ASE}(\widetilde{f}_{jn})(\epsilon)
 &=&
 %{\rm ASE}(\widetilde{f}_{jn})(\epsilon,b_1)
% =
 \frac{1}{T}\sum_{k=1}^T
 \left(\widetilde{f^*}_{jk}(\epsilon)-f(\epsilon)\right)^2.
 \end{eqnarray*}
  The comparison of the ASE is done at
the points $\epsilon=-1,0,1$, using respectively the pointwise
optimal bandwidths
\begin{eqnarray*}
%(\widehat{b}_{1j},\widehat{b}_{0j})
%&=&
(\widehat{b}_{1j}(\epsilon),\widehat{b}_{0j}(\epsilon))
 =
\arg\min_{(b_1,b_0)}
{\rm ASE} (\widehat{f}_{jn})(\epsilon),
\quad
%\widetilde{b}_j
%&=&
\widetilde{b}_j(\epsilon)
 =
\arg\min_{(b_1,b_0)}
{\rm ASE} (\widetilde{f}_{jn})(\epsilon).
\end{eqnarray*}
 As in the global study, these
bandwidths are based upon $T=100$ new independent
 samples of size $n=200$ generated from the model (\ref{quadmodel}),
 and different of the samples that are used for
 computing  ${\rm ASE}(\widehat{f}_{jn})(\epsilon,b_1,b_0)$ and ${\rm
 ASE}(\widetilde{f}_{jn})(\epsilon,b_1)$.
In this  section, the minimizations of the ASE are performed for
$b_1$ and $b_0$ varying on $[0.1,3]$, in the set
$\{h_j=0.1+(0.01)\times j, 1\leq j\leq 290\}$.
%Table \ref{Bandw1}
% groups the simulated values of the  bandwidths
%$(\widehat{b}_{1j}(\epsilon),\widehat{b}_{0j}(\epsilon))$ and
%$\widetilde{b}_j(\epsilon)$ when $T=100$ and $n=200$. With the
%help of these bandwidths, we compute the
% ASE of the estimators
%$\widehat{f}_j(\epsilon)$ and $\widetilde{f}_j(\epsilon)$. These
%ASE are based upon  $T=100$ new independent samples
%$\left(x_{k1},e_{k1}\right),\left(x_{k2},e_{k2}\right),
%\ldots,\left(x_{kn},e_{kn}\right)$ ($k=1,\ldots,T$) of size
%$n=200$, generated from the model (\ref{quadmodel}).
For $j=1,2$ and $\epsilon=-1,0,1$, the optimal values of the ${\rm
ASE}(\widehat{f}_{jn})(\epsilon)$ and ${\rm
ASE}(\widetilde{f}_{jn})(\epsilon)$ are gathered in Tables
\ref{Bandw2} and \ref{Bandw3}. These values show that for any
$\epsilon=-1,0,1$,
\begin{eqnarray}
{\rm ASE}(\widehat{f}_{1n})
(\epsilon,\widehat{b}_{11},\widehat{b}_{01})
 <
 {\rm ASE}(\widehat{f}_{2n})
(\epsilon,\widehat{b}_{12},\widehat{b}_{02}).
 \label{Ase}
\end{eqnarray}
 This fact parallels the results of the
Global study in which we saw that for an optimal choice of the
bandwidth, the  AISE of the estimator $\widehat{f}_{1n}$ is
smaller
 than the one of the estimator $\widehat{f}_{2n}$. Consequently the
 pointwise estimator $\widehat{f}_{1n}(\epsilon)$ should also be preferred to
 the estimator $\widehat{f}_{2n}(\epsilon)$ for the nonparametric
 Kernel estimation of $f(\epsilon)$.

%\begin{table}[H]
%\caption{\label{Bandw1} Values of the  optimal bandwidths
%$(\widehat{b}_{1j}(\epsilon),\widehat{b}_{0j}(\epsilon))$ and
%$\widetilde{b}_j(\epsilon)$ when $b_0,b_1\in[0.1,3]$.}
%\small{
%\begin{center}
%\tabcolsep=3pt
%\begin{tabular}{|c|c|c|c|c|c|c|c|c|c|c|c|}
% \hline
% \multicolumn{4}{|c|}{$\epsilon=-1$}&
%\multicolumn{4}{|c|}{$\epsilon=0$} &
%\multicolumn{4}{|c|}{$\epsilon=1$}
%\\
%\cline{1-12} $(\widehat{b}_{11},\widehat{b}_{01})$ &
%$\widetilde{b}_1$ & $(\widehat{b}_{12},\widehat{b}_{02})$ &
%$\widetilde{b}_2$ & $(\widehat{b}_{11},\widehat{b}_{01})$ &
%$\widetilde{b}_1$ & $(\widehat{b}_{12},\widehat{b}_{02})$ &
%$\widetilde{b}_2$ & $(\widehat{b}_{11},\widehat{b}_{01})$ &
%$\widetilde{b}_1$ & $(\widehat{b}_{12},\widehat{b}_{02})$ &
%$\widetilde{b}_2$
%\\
%\hline
%  (1.38,2.85) & 1.95 & (0.29,0.10) & 0.33 &
%  (0.84,0.23) & 0.88 & (0.19,0.10)  & 0.18 &
%  (1.78,0.31) & 1.87 & (0.29,0.1) & 0.19
% \\
%\hline
%\end{tabular}
%\end{center}
%}
%\end{table}

\begin{table}[H]
\caption{\label{Bandw2} ASE of $\widehat{f}_{1n}(\epsilon)$ and
$\widetilde{f}_{1n}(\epsilon)$ using the bandwidths
$(\widehat{b}_{11}(\epsilon),\widehat{b}_{01}(\epsilon))$ and
$\widetilde{b}_1(\epsilon)$.}
\small{
\begin{center}
\tabcolsep=3pt
\begin{tabular}{|c|c|c|c|c|c|}
 \hline
 \multicolumn{2}{|c|}{$\epsilon=-1$}&
\multicolumn{2}{|c|}{$\epsilon=0$} &
\multicolumn{2}{|c|}{$\epsilon=1$}
\\
\cline{1-6} ${\rm ASE}(\widehat{f}_{1n})$ & ${\rm
ASE}(\widetilde{f}_{1n})$ & ${\rm ASE}(\widehat{f}_{1n})$ & ${\rm
ASE}(\widetilde{f}_{1n})$ & ${\rm ASE}(\widehat{f}_{1n})$ & ${\rm
ASE}(\widetilde{f}_{1n})$
\\
\hline
   {\bf 0.00020536762} & 0.00023502221 & 0.0015443395 &
   {\bf 0.0011523854} & {\bf 0.00013607338} & 0.00028682107
 \\
\hline
\end{tabular}
\end{center}
}
\end{table}

\begin{table}[H]
\caption{\label{Bandw3} ASE of $\widehat{f}_{2n}(\epsilon)$ and
$\widetilde{f}_{2n}(\epsilon)$ based on
 the bandwidths
$(\widehat{b}_{12}(\epsilon),\widehat{b}_{02}(\epsilon))$ and
$\widetilde{b}_2(\epsilon)$.}
 \small{
\begin{center}
\tabcolsep=3pt
\begin{tabular}{|c|c|c|c|c|c|}
 \hline
 \multicolumn{2}{|c|}{$\epsilon=-1$}&
\multicolumn{2}{|c|}{$\epsilon=0$} &
\multicolumn{2}{|c|}{$\epsilon=1$}
\\
\cline{1-6} ${\rm ASE}(\widehat{f}_{2n})$ & ${\rm
ASE}(\widetilde{f}_{2n})$ & ${\rm ASE}(\widehat{f}_{2n})$ & ${\rm
ASE}(\widetilde{f}_{2n})$ & ${\rm ASE}(\widehat{f}_{2n})$ & ${\rm
ASE}(\widetilde{f}_{2n})$
\\
\hline
  {\bf 0.00086381543} &0.0009047470 & 0.0026912672 &
  {\bf 0.0025553406} & {\bf 0.00092628778} & 0.0014950721
 \\
\hline
\end{tabular}
\end{center}
}
\end{table}

\noindent From Tables \ref{Bandw2} and \ref{Bandw3}, we also
notice that for $j=1,2$,
\begin{eqnarray*}
{\rm ASE}(\widehat{f}_{jn})(0)
\approx {\rm
ASE}(\widetilde{f}_{jn})(0),
 \quad
 {\rm ASE}(\widehat{f}_{jn})(\epsilon)
 <
 {\rm ASE}(\widetilde{f}_{jn})(\epsilon),
 \quad\epsilon=-1,1.
\end{eqnarray*}
For the estimation of linear functionals
 of the error distribution in
 a semiparametric context, Müller, Schick and
 Wefelmeyer (2004) have shown that the estimators using the
estimated residuals may have a smaller asymptotic variance
compared to estimators that are based on the true errors.
 A reason that may explain this  effect is that the
 estimators $\widetilde{f}_{jn}(\epsilon)$ do not use the fact that the
 residuals $\varepsilon_i$ have mean zero, contrarily to the estimators
$\widehat{f}_{jn}(\epsilon)$. Note however that the improvement of
$\widehat{f}_{1n}(\epsilon)$ on $\widetilde{f}_{1n}(\epsilon)$ is
much more clear-cut in the pointwise setup than in the global one.

Nevertheless, we observe that for the estimator $\widehat{f}_{1n}$
based on the estimated residuals, the values of the ASE are quite
different at the points $\epsilon=-1$ and $\epsilon=1$. We then
attempt to explain this situation by analyzing the behavior of the
error terms  around these points.
 Define, for
any integers $k\in[1,T]$ and $i\in[1,n]$,
$$
\widehat{\delta}_{ki}(\epsilon)
=
\left(\widehat{\varepsilon}_{ki}-\varepsilon_{ki}\right)
\mathds{1}
\left(
\left|\widehat{\varepsilon}_{ki}-\epsilon\right|
\leq
\widehat{b}_{11}(\epsilon)
\right),
$$
where $\widehat{\varepsilon}_{ki}$ is the Kernel empirical version
of $\varepsilon_{ki}$ based on the  optimal first-step
 bandwidth $\widehat{b}_{01}(\epsilon)$ for the estimator
$\widehat{f}_{1n}(\epsilon)$. We then define the empirical mean
$\overline{\delta}(\epsilon)$ and the empirical variance
$\sigma_{\bar{\delta}}^2(\epsilon)$ of the
$\widehat{\delta}_{ki}(\epsilon)$'s as
$$
\overline{\delta}(\epsilon)
=
\frac{1}{nT} \sum_{k=1}^T
\sum_{i=1}^n \widehat{\delta}_{ki}(\epsilon),
\quad
\sigma_{\bar{\delta}}^2(\epsilon)
=
\frac{1}{nT} \sum_{k=1}^T
\sum_{i=1}^n \left(
\widehat{\delta}_{ki}(\epsilon)-\overline{\delta}(\epsilon)
\right)^2.
$$

\begin{table}[H]
\caption{\label{Varesid} Values of the empirical means
$\overline{\delta}(\epsilon)$ and the empirical variances
$\sigma_{\bar{\delta}}^2(\epsilon)$ for $\epsilon=-1,1$.}
\begin{center}
\tabcolsep=3pt
\begin{tabular}{|c|c|c|c|}
 \hline
 \multicolumn{2}{|c|}{$\epsilon=-1$}&
\multicolumn{2}{|c|}{$\epsilon=1$}
\\
\cline{1-4}
 $\overline{\delta}(\epsilon)$ &
 $\sigma_{\bar{\delta}}^2(\epsilon)$ &
  $\overline{\delta}(\epsilon)$ &
 $\sigma_{\bar{\delta}}^2(\epsilon)$
\\
\hline
  -0.3911658 & 0.3331359 & 0.03403744 & 0.04659867
 \\
\hline
\end{tabular}
\end{center}
\end{table}

\noindent In Table \ref{Varesid}, we evaluate the quantities
$\overline{\delta}(\epsilon)$ and
$\sigma_{\bar{\delta}}^2(\epsilon)$, using the bandwidths
$b_0=\widehat{b}_{01}(\epsilon)$  and
$b_1=\widehat{b}_{11}(\epsilon)$, $\epsilon=-1,1$. We observe that
the variables $\widehat{\delta}_{ki}(-1)$ have a lower empirical
bias and a higher empirical variance than the data
$\widehat{\delta}_{ki}(1)$. Hence around the point $\epsilon=-1$,
the error percentage for the estimation of the
 true residuals $\varepsilon_{ki}$ by the nonparametric residuals
 $\widehat{\varepsilon}_{ki}$ is more important  than around the point
 $\epsilon=1$. This may explain the difference of the ASE at the
 points $\epsilon=-1,1$ for the estimator $\widehat{f}_{1n}$,
  as seen in Table \ref{Bandw2}.

\subsection{Comparison of the bias and variances}

In this subsection, we suppose that the estimators
 $\widehat{f^*}_{jk}(\epsilon)$
and $\widetilde{f^*}_{jk}(\epsilon)$ ($j=1,2)$ defined in the
previous subsection are respectively based upon the optimal
bandwidths
$(\widehat{b}_{1j}(\epsilon),\widehat{b}_{0j}(\epsilon))$ and
$\widetilde{b}_{j}(\epsilon)$. For each  $j$, let
$\widehat{B}_{jn}(\epsilon)$ and $\widetilde{B}_{jn}(\epsilon)$ be
respectively the empirical
 bias  of the estimated densities $\widehat{f}_{jn}(\epsilon)$ and
 $\widetilde{f}_{jn}(\epsilon)$.
These estimated quantities are defined as
\begin{eqnarray*}
\widehat{B}_{jn}(\epsilon)
=
\frac{1}{T}\sum_{k=1}^T
\left(\widehat{f^*}_{jk}(\epsilon)-f(\epsilon)\right),
\quad
\widetilde{B}_{jn}(\epsilon)
=
\frac{1}{T}\sum_{k=1}^T
\left(\widetilde{f^*}_{jk}(\epsilon)-f(\epsilon)\right).
\end{eqnarray*}
 The simulated values of the bias
$\widehat{B}_{jn}(\epsilon)$ and $\widetilde{B}_{jn}(\epsilon)$ at
the points $\epsilon=-1,0,1$ are represented in Table \ref{Tab2}
and \ref{Tab2b}.

\begin{table}[H]
\caption{\label{Tab2} Optimal values of the  bias
$\widehat{B}_{1n}(\epsilon)$ and $\widetilde{B}_{1n}(\epsilon)$.}
 \small{
\begin{center}
\tabcolsep=3pt
\begin{tabular}{|c|c|c|c|c|c|}
 \hline
 \multicolumn{2}{|c|}{$\epsilon=-1$}&
\multicolumn{2}{|c|}{$\epsilon=0$} &
\multicolumn{2}{|c|}{$\epsilon=1$}
\\
\cline{1-6} $\widehat{B}_{1n}(\epsilon)$
&$\widetilde{B}_{1n}(\epsilon)$& $\widehat{B}_{1n}(\epsilon)$ &
$\widetilde{B}_{1n}(\epsilon)$&
 $\widehat{B}_{1n}(\epsilon)$ & $\widetilde{B}_{1n}(\epsilon)$
\\
\hline
  {\bf -0.005290204} & -0.008126554 &
  -0.02450615 & {\bf -0.01726208} &
  {\bf -0.005446647} &  -0.008392434
 \\
\hline
\end{tabular}
\end{center}
}
\end{table}

\begin{table}[H]
\caption{\label{Tab2b} Optimal values of the bias
$\widehat{B}_{2n}(\epsilon)$ and $\widetilde{B}_{2n}(\epsilon)$.}
 \small{
\begin{center}
\tabcolsep=3pt
\begin{tabular}{|c|c|c|c|c|c|}
 \hline
 \multicolumn{2}{|c|}{$\epsilon=-1$}&
\multicolumn{2}{|c|}{$\epsilon=0$} &
\multicolumn{2}{|c|}{$\epsilon=1$}
\\
\cline{1-6}
$\widehat{B}_{2n}(\epsilon)$
&$\widetilde{B}_{2n}(\epsilon)$&
$\widehat{B}_{2n}(\epsilon)$ &
$\widetilde{B}_{2n}(\epsilon)$&
 $\widehat{B}_{2n}(\epsilon)$ &
  $\widetilde{B}_{2n}(\epsilon)$
 \\
\hline
  {\bf -0.01615247} & -0.01661985 &
   {\bf -0.02447482} & -0.02612883 &
  -0.01803243 &  {\bf -0.01262712}
 \\
\hline
\end{tabular}
\end{center}
}
\end{table}

\noindent Table \ref{Tab2} reveals that
$|\widehat{B}_{1n}(\epsilon)|<|\widehat{B}_{2n}(\epsilon)|$ for
$\epsilon=-1$ and $\epsilon=1$, and that
$|\widehat{B}_{1n}(0)|\approx|\widehat{B}_{2n}(0)|$. This
indicates that at the points $\epsilon=-1$ and $\epsilon=1$, the
estimator $\widehat{f}_{1n}(\epsilon)$ would be less biased than
the estimator $\widehat{f}_{2n}(\epsilon)$.

\vskip 0.1cm\noindent
 Moreover for $\epsilon=-1$ and $\epsilon=1$,
 the estimator $\widehat{f}_{1n}(\epsilon)$ is much less biased than
 the estimator $\widetilde{f}_{1n}(\epsilon)$. Consequently, there is a positive
 influence of the bandwidth $\widehat{b}_{01}(\epsilon)$ on the
 bias of $\widehat{f}_{1n}(\epsilon)$. But this situation contrasts with
 the one observed at $\epsilon=0$, for which
  $\widehat{f}_{1n}(\epsilon)$ is more biased
 than $\widetilde{f}_{1n}(\epsilon)$.
\vskip 0.1cm \noindent
For $\widehat{f}_{2n}(\epsilon)$ and
$\widetilde{f}_{2n}(\epsilon)$, we note that the bias of these
estimators are very close at the points $\epsilon=-1,0,1$.
 This means  that the estimation of the
 regression function has a negligible impact on the
 bias of the estimator $\widehat{f}_{2n}(\epsilon)$.

\vskip 0.3cm
 Now, let  $\widehat{V}_{jn}(\epsilon)$ and
$\widetilde{V}_{jn}(\epsilon)$ be the estimated variances of
$\widehat{f}_{jn}(\epsilon)$ and $\widetilde{f}_{jn}(\epsilon)$
defined as
\begin{eqnarray*}
\widehat{V}_{jn}(\epsilon)
=
 \frac{1}{T}\sum_{k=1}^T
\left(\widehat{f^*}_{jk}(\epsilon) -
\widehat{\mu}_{jn}(\epsilon)\right)^2,
\quad
\widetilde{V}_{jn}(\epsilon)
 =
 \frac{1}{T}\sum_{k=1}^T
\left(\widetilde{f^*}_{jk}(\epsilon)
-
\widetilde{\mu}_{jn}(\epsilon)\right)^2,
\end{eqnarray*}
where
$$
\widehat{\mu}_{jn}(\epsilon)
=
 \frac{1}{T}\sum_{k=1}^T
\widehat{f^*}_{jk}(\epsilon),
 \quad
\widetilde{\mu}_{jn}(\epsilon)
=
\frac{1}{T}\sum_{k=1}^T
\widetilde{f^*}_{jk}(\epsilon).
$$
The simulated values of these empirical parameters are gathered in
Tables \ref{Tab3} and \ref{Tab3b}.

\begin{table}[H]
\caption{\label{Tab3} Optimal values of  variances
$\widehat{V}_{1n}(\epsilon)$ and $\widetilde{V}_{1n}(\epsilon)$.}
 \small{
\begin{center}
\tabcolsep=3pt
\begin{tabular}{|c|c|c|c|c|c|}
 \hline
 \multicolumn{2}{|c|}{$\epsilon=-1$}&
 \multicolumn{2}{|c|}{$\epsilon=0$} &
 \multicolumn{2}{|c|}{$\epsilon=1$}
\\
\cline{1-6} $\widehat{V}_{1n}(\epsilon)$ &
$\widetilde{V}_{1n}(\epsilon)$& $\widehat{V}_{1n}(\epsilon)$ &
$\widetilde{V}_{1n}(\epsilon)$& $\widehat{V}_{1n}(\epsilon)$ &
$\widetilde{V}_{1n}(\epsilon)$
\\
\hline
   0.0001773813 & {\bf 0.0001689813} &
  0.0009437874 &  {\bf 0.0008544052} &
  {\bf 0.0001064074} &  0.000216388
 \\
\hline
\end{tabular}
\end{center}
}
\end{table}

\begin{table}[H]
\caption{\label{Tab3b} Optimal values of the variances
$\widehat{V}_{2n}(\epsilon)$ and $\widetilde{V}_{2n}(\epsilon)$.}
 \small{
\begin{center}
\tabcolsep=3pt
\begin{tabular}{|c|c|c|c|c|c|}
 \hline
 \multicolumn{2}{|c|}{$\epsilon=-1$}&
 \multicolumn{2}{|c|}{$\epsilon=0$} &
 \multicolumn{2}{|c|}{$\epsilon=1$}
\\
\cline{1-6} $\widehat{V}_{2n}(\epsilon)$ &
$\widetilde{V}_{2n}(\epsilon)$& $\widehat{V}_{2n}(\epsilon)$ &
$\widetilde{V}_{2n}(\epsilon)$& $\widehat{V}_{2n}(\epsilon)$ &
$\widetilde{V}_{2n}(\epsilon)$
\\
\hline
   {\bf 0.0006029132} &  0.0006285278 &
   0.00209225 & {\bf 0.001872625} &
  {\bf 0.0006011193} & 0.001335628
 \\
\hline
\end{tabular}
\end{center}
}
\end{table}

\noindent
 From Table \ref{Tab3}, we notice that
$\widehat{V}_{1n}(\epsilon)<\widehat{V}_{2n}(\epsilon)$ for
$\epsilon=-1,0,1$. Consequently the estimator
$\widehat{f}_{1n}(\epsilon)$ should be preferred to
$\widehat{f}_{2n}(\epsilon)$, since the latter estimator is less
efficient than the first one.

\vskip 0.1cm\noindent
Moreover, we observe that
$\widehat{V}_{1n}(\epsilon)$ is  much less than
$\widetilde{V}_{1n}(\epsilon)$ at $\epsilon=1$, and slightly equal
to $\widetilde{V}_{1n}(\epsilon)$ when $\epsilon=-1$ and
$\epsilon=0$.
 This means that the estimation of the residuals
  may have a positive influence on the final estimator
$\widehat{f}_{1n}(\epsilon)$.

 \vskip 0.1cm\noindent
 For the variances $\widehat{V}_{2n}(\epsilon)$
 and $\widetilde{V}_{2n}(\epsilon)$, it is seen that
 the first variance is much less than the
 latter one at $\epsilon=1$, and very close to
$\widetilde{V}_{2n}(\epsilon)$ when  $\epsilon=-1$ and
$\epsilon=0$. Hence the estimation of the regression function $m$
may have a positive effect  on the estimator
$\widehat{f}_{2n}(\epsilon)$.

\vskip 0.1cm \textcolor{red}{\bf In conclusion, we note that at
the points $\epsilon=-1,0,1$, the estimator
$\widehat{f}_{1n}(\epsilon)$ dominates the estimator
$\widehat{f}_{2n}(\epsilon)$ for the ASE, the bias and the
variance}. As in the Global study,  this suggests that the first
estimator should be preferred to the second one when we are
interested in their Pointwise study.

\newpage

\subsection{Asymptotic normality}

We  examine here the asymptotic normality of the estimators
$\widehat{f}_{jn}(\epsilon)$, for $j=1,2$ and $\epsilon=-1,0,1$.
To that aim,  we introduce the standardized variables
$$
\widehat{Z}_{jn}(\epsilon)
 =
\frac{\sqrt{n\widehat{b}_{1j}(\epsilon)}
\left(\widehat{f}_{jn}(\epsilon)-f(\epsilon)\right)}
{\sqrt{f(\epsilon)\int\! K_1^2(v) dv}}, \quad
\widehat{Z^*}_{jk}(\epsilon)
 =
\frac{\sqrt{n\widehat{b}_{1j}(\epsilon)}
\left(\widehat{f^*}_{jk}(\epsilon)-f(\epsilon)\right)}
{\sqrt{f(\epsilon)\int\! K_1^2(v) dv}}, \quad k=1,\ldots,T,
$$
where the $\widehat{f^*}_{jk}(\epsilon)$'s are defined as in
previous subsection, for the evaluation of the bias and variances.
The empirical mean $\overline{\widehat{\mu}}_{jn}(\epsilon)$ and
the empirical variance
$\overline{\widehat{\sigma}}_{jn}^2(\epsilon)$ of the
$\widehat{Z}_{jn}(\epsilon)$'s are such that
$$
\overline{\widehat{\mu}}_{jn}(\epsilon)
=
\frac{1}{T} \sum_{k=1}^T
\widehat{Z^*}_{jk}(\epsilon),
\quad
\overline{\widehat{\sigma}}_{jn}^2(\epsilon)
=
\frac{1}{T}\sum_{k=1}^T
\left(
\widehat{Z^*}_{jk}(\epsilon)
-
\overline{\widehat{\mu}}_{jn}(\epsilon)
\right)^2.
$$

\textcolor{red}{\large\bf Are the data
$\widehat{Z}_{jn}(\epsilon)$ normal distributed?}

\vskip 0.1cm\noindent For each $j$ and $\epsilon$, we wish to test
the hypothesis
\textcolor{blue}{\bf
$$
{\rm {H}_{0j}(\epsilon):}
\;\widehat{Z}_{jn}(\epsilon)
\sim
N\left(\mu_j(\epsilon),\sigma_j^2(\epsilon)\right)
\;{\rm versus}\;\;
 {\rm
{H}_{1j}(\epsilon):} \;
\widehat{Z}_{jn}(\epsilon)
\not\sim
N\left(\mu_j(\epsilon),\sigma_j^2(\epsilon)\right),
$$
} where the parameters $\mu_j(\epsilon)$ and
$\sigma_j^2(\epsilon)$ are unknown and have to  be estimated. The
normality of the data $\widehat{Z}_{jn}(\epsilon)$ can
 be tested by an analytical method such as the Lilliefors
method for the Kolmogorov-Smirnov test. Let us  perform this
Lilliefors test that the data $\widehat{Z}_{jn}(\epsilon)$ come
from the normal distribution. For this, we denote by $\widehat{\rm
KS}_j(\epsilon)$ and $\widehat{p}_j(\epsilon)$ respectively as the
Kolmogorov-Smirnov statistic and the $p$-value of the above test.
With the Lilliefors's method, the evaluation of the $p$-values
$\widehat{p}_j(\epsilon)$ and the statistics $\widehat{\rm
KS}_j(\epsilon)$ accounts for the estimations of $\mu_j(\epsilon)$
and $\sigma_j^2(\epsilon)$. For the characteristics and the
properties of the KS or Lilliefors's test, see Massey (1951),
Shorack and Wellner (1986), Dallal and Wilkinson (1986), Lehmann
and Romano (1998), and Thode (2002). In Table \ref{KS} we have
gathered the numerical values of the $\widehat{\rm
KS}_j(\epsilon)$'s and the $\widehat{p}_j(\epsilon)$'s.

\begin{table}[H]
\caption{\label{KS} Values of the statistics $\widehat{\rm
KS}_j(\epsilon)$ and the p-values $\widehat{p}_j(\epsilon)$ of the
$\widehat{Z}_j(\epsilon)$'s.}
 \small{
\begin{center}
\tabcolsep=3pt
\begin{tabular}{|c|c|c|c|c|c|c|c|c|c|c|c|}
 \hline
 \multicolumn{4}{|c|}{$\epsilon=-1$}&
\multicolumn{4}{|c|}{$\epsilon=0$}&
\multicolumn{4}{|c|}{$\epsilon=1$}
\\
\cline{1-12}
$\widehat{\rm KS}_1(\epsilon)$ &
$\widehat{p}_1(\epsilon)$ & $\widehat{\rm KS}_2(\epsilon)$ &
$\widehat{p}_2(\epsilon)$ & $\widehat{\rm KS}_1(\epsilon)$ &
$\widehat{p}_1(\epsilon)$ & $\widehat{\rm KS}_2(\epsilon)$ &
$\widehat{p}_2(\epsilon)$ & $\widehat{\rm KS}_1(\epsilon)$ &
$\widehat{p}_1(\epsilon)$ & $\widehat{\rm KS}_2(\epsilon)$ &
$\widehat{p}_2(\epsilon)$
\\
\hline
 0.0506 & 0.9598 & 0.0427 & 0.9933 &
 0.0713 & 0.6891 & 0.0518 & 0.9516 &
  0.0746 & 0.6347 & 0.0882 &  0.4176
\\
\hline
\end{tabular}
\end{center}
}
\end{table}
\noindent  The results of Table \ref{KS}  show that the hypothesis
on the normality of the data is accepted, since
$\widehat{p}_j(\epsilon)>0.05=\alpha$ (a default value of the
level of significance). \textcolor{red}{\bf Hence according to the
Lilliefors method, we can accept the fact that the data
$\widehat{Z}_{jn}(\epsilon)$ come from a normal distribution.}

\begin{figure}[H]
\caption{From left to right: normal Q-Q plot of the data
$Z_{1n}(-1)$ and $Z_{2n}(-1)$.} \centering
\includegraphics[scale=0.45]{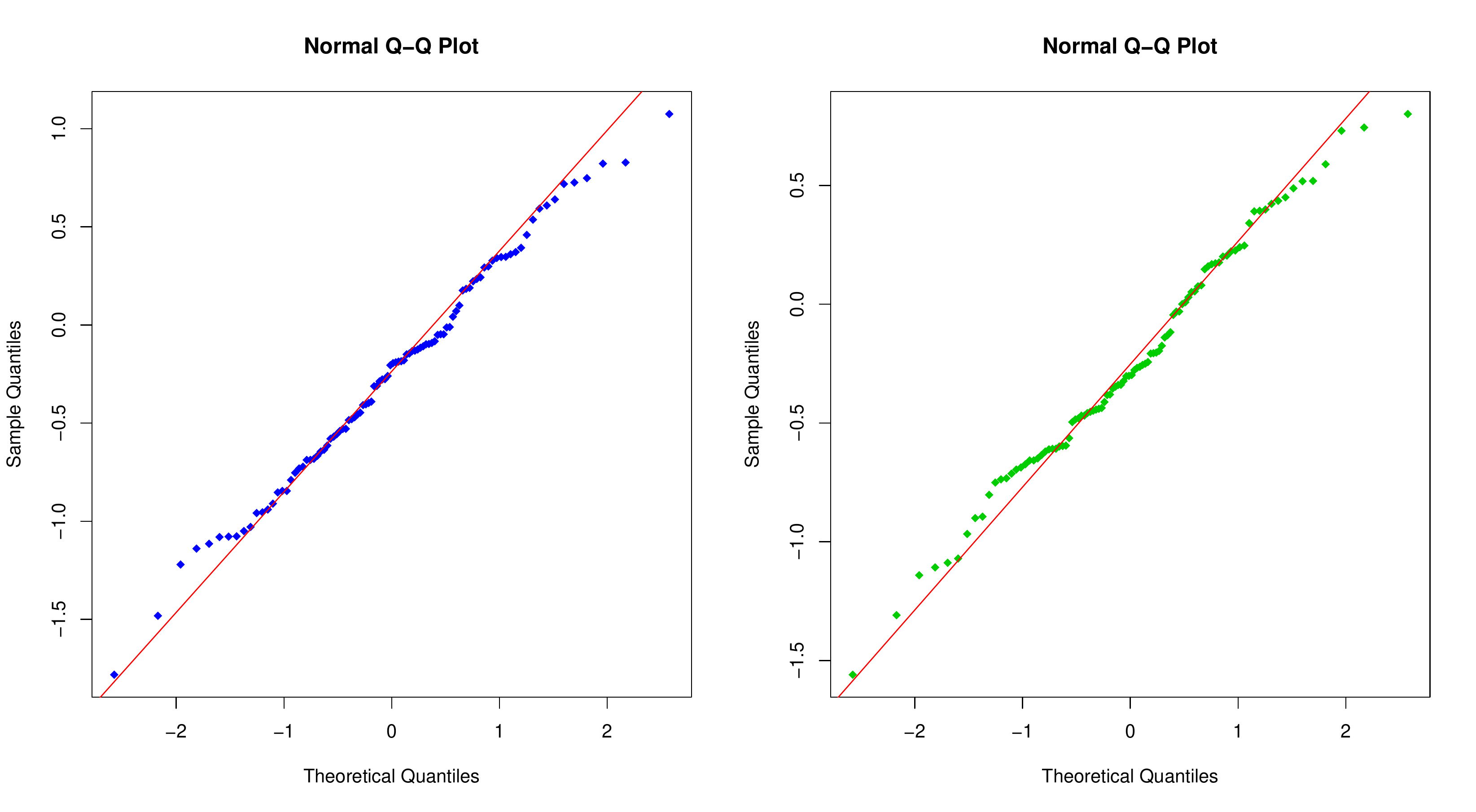}
 \label{QQmoins1}
\end{figure}

\begin{figure}[H]
\caption{Normal Q-Q plot of the data $Z_{1n}(0)$ and $Z_{2n}(0)$.}
\centering
\includegraphics[scale=0.45]{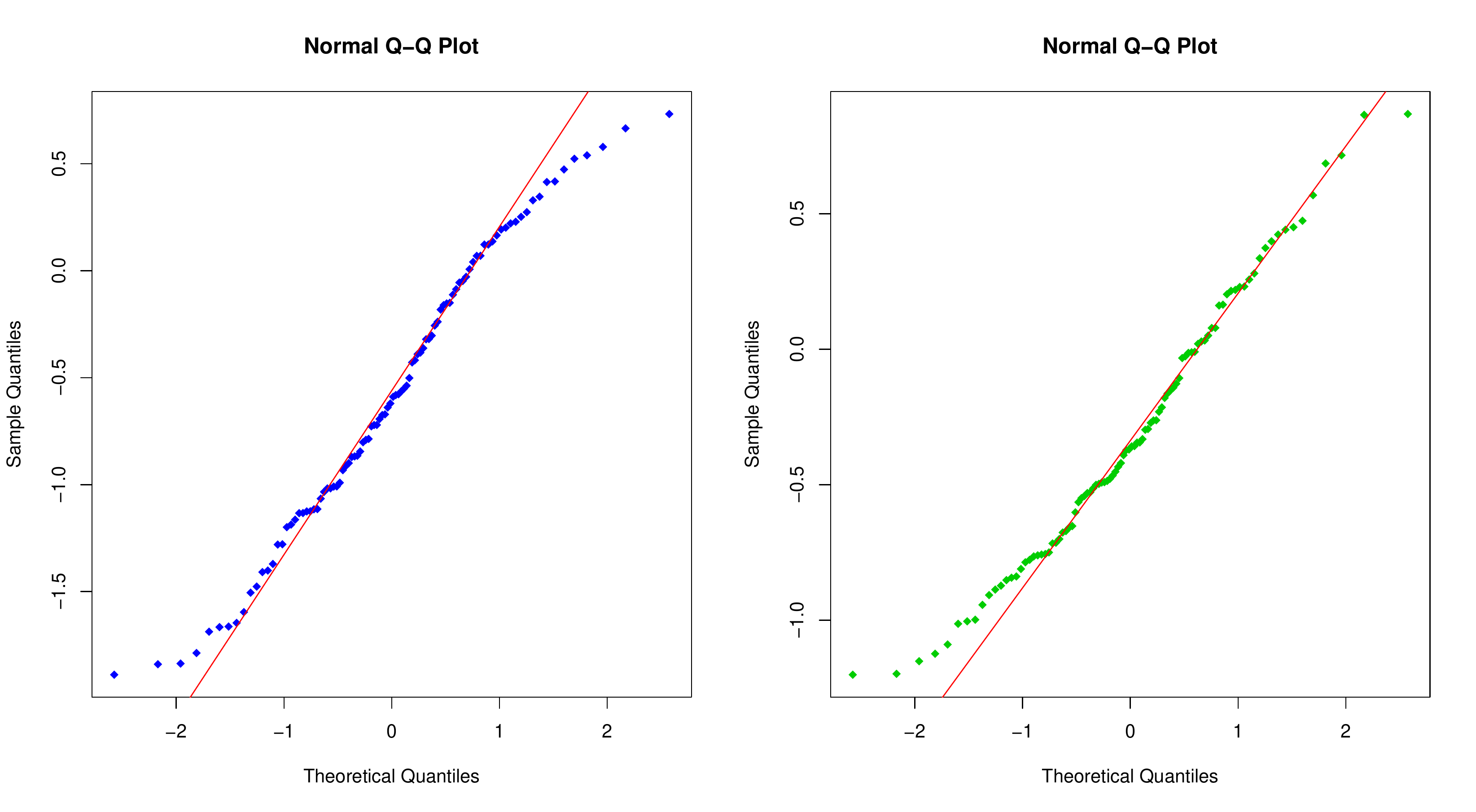}
 \label{QQzero}
\end{figure}

\begin{figure}[H]
\caption{Normal Q-Q plot of the data $Z_{1n}(1)$ and $Z_{2n}(1)$.}
\centering
\includegraphics[scale=0.45]{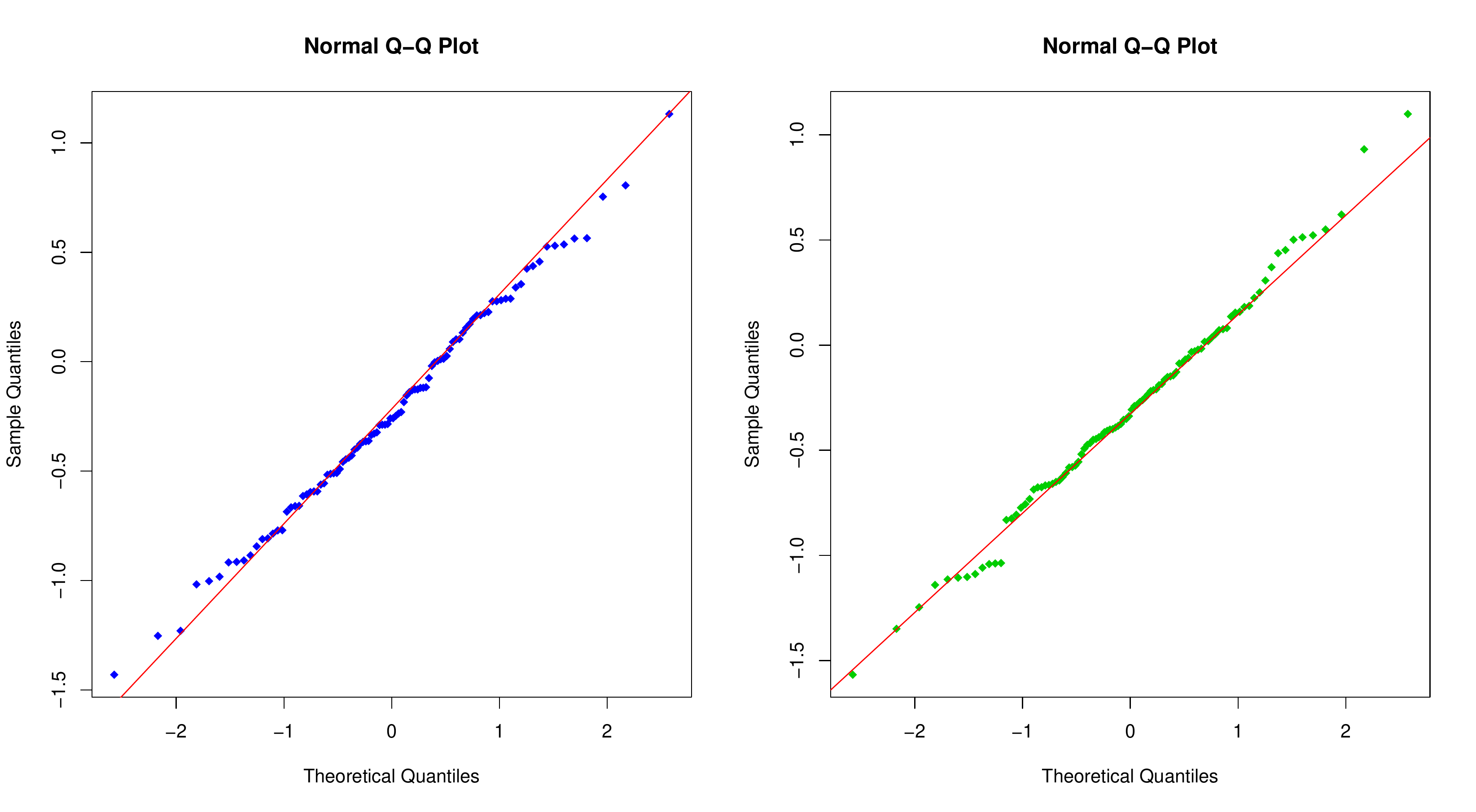}
 \label{QQun}
\end{figure}

\vskip 0.3cm Beside the Lilliefors test, there exists a graphical
method for investigating the normality of the data. This method is
the normal Q-Q plots of the variables
$\widehat{Z}_{jn}(\epsilon)$. The Q-Q plot provides a graphical
way to determine the level of
 normality.  If the data fall exactly along a reference line
(called the Henry's line), then the hypothesis on their normality
can be receivable. If the empirical data deviate widely from this
line, the data are non-normal. In Figures \ref{QQmoins1},
\ref{QQzero} and \ref{QQun}, we represent the normal Q-Q plots of
the data $\widehat{Z}_{1n}(\epsilon)$ and
$\widehat{Z}_{2n}(\epsilon)$ for $\epsilon=-1,0,1$.

\noindent From the above figures, the hypothesis ${\bf \rm
{H}_{0j}(\epsilon)}$ can be receivable. However, we can have some
doubts about the symmetry of the $\widehat{Z}_{jn}(\epsilon)$'s,
since we notice that they deviate slightly from the Henry's line
at the tails of the distribution. This augurs that the
distribution of the variables $\widehat{Z}_{jn}(\epsilon)$ should
 not be symmetric.

\vskip 0.1cm \noindent

\textcolor{red}{\large\bf Do the data $\widehat{Z}_{jn}(\epsilon)$
come from the standard normal $N(0,1)$?}

\vskip 0.1cm\noindent Since the data $\widehat{Z}_{jn}(\epsilon)$
are standardized variables, we can wonder if they come from the
normal distribution $N(0,1)$. To give some elements of answer to
this question, we first compute the empirical bias and variances
of the $\widehat{Z}_{jn}(\epsilon)$'s. These quantities are
grouped in Table \ref{Tab8}.

 \begin{table}[H]
\caption{\label{Tab8} Values of the empirical means and variances
of the $\widehat{Z}_{jn}(\epsilon)$'s.}
 \small{
\begin{center}
\tabcolsep=3pt
\begin{tabular}{|c|c|c|c|c|c|c|c|c|c|c|c|}
 \hline
 \multicolumn{4}{|c|}{$\epsilon=-1$}&
\multicolumn{4}{|c|}{$\epsilon=0$} &
\multicolumn{4}{|c|}{$\epsilon=1$}
\\
\cline{1-12}
 $\overline{\widehat{\mu}}_{1n}(\epsilon)$ &
$\overline{\widehat{\sigma}}_{1n}^2(\epsilon)$ &
$\overline{\widehat{\mu}}_{2n}(\epsilon)$ &
$\overline{\widehat{\sigma}}_{2n}^2(\epsilon)$ &
$\overline{\widehat{\mu}}_{1n}(\epsilon)$ &
$\overline{\widehat{\sigma}}_{1n}^2(\epsilon)$ &
$\overline{\widehat{\mu}}_{2n}(\epsilon)$ &
$\overline{\widehat{\sigma}}_{2n}^2(\epsilon)$ &
$\overline{\widehat{\mu}}_{1n}(\epsilon)$ &
$\overline{\widehat{\sigma}}_{1n}^2(\epsilon)$ &
$\overline{\widehat{\mu}}_{2n}(\epsilon)$ &
$\overline{\widehat{\sigma}}_{2n}^2(\epsilon)$
\\
\hline
  -0.1459 & 0.4900 & -0.3105 & 0.4098 &
  -0.5798 & 0.6389 & -0.3124 & 0.4380 &
   -0.2491 & 0.4465 & -0.3097 & 0.4323
 \\
\hline
\end{tabular}
\end{center}
}
\end{table}
\noindent This table shows that the empirical paremeters
$\overline{\widehat{\mu}}_{jn}(\epsilon)$ and
$\overline{\widehat{\sigma}}_{1n}^2(\epsilon)$ are clearly
different to $0$ and $1$, which correspond respectively to the
theoretical mean and variance of the normal $N(0,1)$.

\vskip 0.1cm
 We now evaluate the empirical quantiles of the variables
$\widehat{Z}_{jn}(\epsilon)$. For each $j$, we consider the
ordered values $\widehat{Z^*}_{j,(k)}(\epsilon)$ of the
$\widehat{Z^*}_{jk}(\epsilon)$'s  such that
$\widehat{Z^*}_{j,(1)}(\epsilon) \leq
\widehat{Z^*}_{j,(2)}(\epsilon) \leq \ldots \leq
\widehat{Z^*}_{j,(T)}(\epsilon)$.
 Hence for any $\alpha\in[0,1]$, the $\alpha^{\rm th}$ empirical
 quantiles
 of the $\widehat{Z}_{jn}(\epsilon)$'s are defined as
$\widehat{Z}_{jn}(\epsilon,\alpha)=\widehat{Z^*}_{j,(\alpha
T)}(\epsilon)$. In Table \ref{Tab6}, we give the simulated values
of these quantiles when $\alpha=0.05$ and $\alpha=0.95$.

\begin{table}[H]
\caption{\label{Tab6} Values of the  quantiles
 $a_j(\epsilon)=\widehat{Z^*}_{j,(0.05\times T)}(\epsilon)$ and
 $c_j(\epsilon)=\widehat{Z^*}_{j,(0.95\times T)}(\epsilon)$.}
 \small{
\begin{center}
\tabcolsep=3pt
\begin{tabular}{|c|c|c|c|c|c|c|c|c|c|c|c|}
 \hline
 \multicolumn{4}{|c|}{$\epsilon=-1$} &
\multicolumn{4}{|c|}{$\epsilon=0$} &
\multicolumn{4}{|c|}{$\epsilon=1$}
\\
\cline{1-12}
$a_1(\epsilon)$ & $c_1(\epsilon)$& $a_2(\epsilon)$ &
$c_2(\epsilon)$& $a_1(\epsilon)$ & $c_1(\epsilon)$ &
$a_2(\epsilon)$ & $c_2(\epsilon)$& $a_1(\epsilon)$ &
$c_1(\epsilon)$& $a_2(\epsilon)$ & $c_2(\epsilon)$
\\
\hline
 -0.9719 & 0.6654 & -0.9790 & 0.4006 &
 -1.6874 & 0.4734 & -1.0893  & 0.4738 &
 -1.1377 & 0.5389 &  -1.1109 & 0.5071
\\
\hline
\end{tabular}
\end{center}
}
\end{table}
\noindent From Table \ref{Tab6}, we note first that the quantities
$\widehat{Z}_{jn}(\epsilon,0.05)$ and
$\widehat{Z}_{jn}(\epsilon,0.95)$ are globally clearly different
to $-1.64$ and $1.64$, the corresponding theoretical quantiles of
the normal distribution $N(0,1)$. Next, the values of these
quantiles show that the variables $\widehat{Z}_{jn}(\epsilon)$ are
globally asymmetric. This suggests that the
$\widehat{Z}_{jn}(\epsilon)$'s are not distributed according to
the normal variable $N(0,1)$.

\vskip 0.3cm We now estimate the confidence intervals for the
theoretical quantiles of the $\widehat{Z}_{jn}(\epsilon)$'s. This
estimation is done under  $\bf \rm {H}_{0j}(\epsilon)$, and based
on the following result. For $\alpha\in]0,1[$ and
$T\rightarrow\infty$,
$$
\frac{\widehat{Z^*}_{j,(\alpha
T)(\epsilon)}-Q_{j,\alpha}(\epsilon)} {\sqrt{\alpha(\alpha-1)/
\left(T f^2\left(Q_{j,\alpha}(\epsilon)\right)\right)}}
 \stackrel{d}{\longrightarrow} N(0,1),
$$
where $Q_{j,\alpha}(\epsilon)$ is the theoretical $\alpha^{\rm
th}$ quantile of the variable $\widehat{Z}_{jn}(\epsilon)$, and
$f(\cdot)$ the p.d.f of the normal variable $N(0,1)$. This result
can be found, for example, in Tassi (1985). A consequence of a
such result is that an asymptotic confidence interval for the
$Q_{j,\alpha}(\epsilon)$'s, with a level of confidence $1-\alpha$,
is given by
$$
\widehat{I}_{j,\alpha}(\epsilon)
=
\left[ \widehat{Z^*}_{j,(\alpha
T)}(\epsilon)-\frac{q_{\alpha/2}\sqrt{\alpha(1-\alpha)}}{\sqrt{T}
f\left(\widehat{Z^*}_{j,(\alpha T)}(\epsilon)\right)} ,
\widehat{Z^*}_{j,(\alpha
T)}(\epsilon)+\frac{q_{\alpha/2}\sqrt{\alpha(1-\alpha)}}{\sqrt{T}
f\left(\widehat{Z^*}_{j,(\alpha T)}(\epsilon)\right)} \right],
$$
where $q_{\alpha/2}$ denotes the $(1 -\alpha/2)$ quantile of the
standard normal distribution. In Tables \ref{Tab9} and
\ref{Tab10}, we give the estimations of the
$\widehat{I}_{j,\alpha}(\epsilon)$ when $\alpha=0.05$ and
$\alpha=0.95$, with $T=100$ and $n=200$. As seen above in Table
\ref{Tab6}, the results of Tables \ref{Tab9} and \ref{Tab10} also
reveal that the quantiles $Q_{j,0.05}(\epsilon)$ and
$Q_{j,0.95}(\epsilon)$ should be respectively quite different to
$-1.64$ and $1.64$.

\begin{table}[H]
\caption{\label{Tab9} Confidence intervals of the theoretical
quantiles $Q_{j,\alpha}(\epsilon)$ when $\alpha=0.05$.}
 \small{
\begin{center}
\tabcolsep=3pt
\begin{tabular}{|c|c|c|c|c|c|}
 \hline
 \multicolumn{2}{|c|}{$\epsilon=-1$}&
\multicolumn{2}{|c|}{$\epsilon=0$} &
\multicolumn{2}{|c|}{$\epsilon=1$}
\\
\cline{1-6}
$\widehat{I}_{1,\alpha}(\epsilon)$ &
$\widehat{I}_{2,\alpha}(\epsilon)$ &
$\widehat{I}_{1,\alpha}(\epsilon)$ &
$\widehat{I}_{2,\alpha}(\epsilon)$ &
$\widehat{I}_{1,\alpha}(\epsilon)$ &
$\widehat{I}_{2,\alpha}(\epsilon)$
\\
\hline
$\left[-1.143, -0.800\right]$ & $\left[-1.152,
-0.806\right]$ & $\left[-2.132, -1.243\right]$ & $\left[-1.283,
-0.896\right]$ & $\left[-1.342, -0.933\right]$ & $\left[-1.309,
-0.913\right]$
 \\
\hline
\end{tabular}
\end{center}
}
\end{table}

\begin{table}[H]
\caption{\label{Tab10} Confidence intervals of the theoretical
quantiles $Q_{j,\alpha}(\epsilon)$ when $\alpha=0.95$.}
 \small{
\begin{center}
\tabcolsep=3pt
\begin{tabular}{|c|c|c|c|c|c|}
 \hline
 \multicolumn{2}{|c|}{$\epsilon=-1$}&
\multicolumn{2}{|c|}{$\epsilon=0$} &
\multicolumn{2}{|c|}{$\epsilon=1$}
\\
\cline{1-6}
 $\widehat{I}_{1,\alpha}(\epsilon)$ &
$\widehat{I}_{2,\alpha}(\epsilon)$ &
$\widehat{I}_{1,\alpha}(\epsilon)$ &
$\widehat{I}_{2,\alpha}(\epsilon)$ &
$\widehat{I}_{1,\alpha}(\epsilon)$ &
$\widehat{I}_{2,\alpha}(\epsilon)$
\\
\hline
$\left[0.5318, 0.7990\right]$ & $\left[0.3257, 0.5172\right]$ &
$\left[0.3536, 0.5934\right]$ & $\left[0.3549, 0.5941\right]$ &
$\left[0.4159, 0.6635\right]$ & $\left[0.3852, 0.6297\right]$
 \\
\hline
\end{tabular}
\end{center}
}
\end{table}
The above results indicates that the data
$\widehat{Z}_{jn}(\epsilon)$ are not distributed according to the
standard normal $N(0,1)$. We then attempt to verify if this
situation is due to the influence of estimated first-step
bandwidths $\widehat{b}_{0j}(\epsilon)$ on the variables
$\widehat{Z}_{jn}(\epsilon)$. For this,
 we  test the hypothesis
$\widetilde{Z}_{jn}(\epsilon)\sim N\left(0, 1\right)$ versus
$\widetilde{Z}_{jn}(\epsilon)\not\sim N\left(0, 1\right)$, where
$$
\widetilde{Z}_{jn}(\epsilon)
 =
\frac{\sqrt{n\widetilde{b}_j(\epsilon)}
\left(\widetilde{f}_{jn}(\epsilon)-f(\epsilon)\right)}
{\sqrt{f(\epsilon)\int\! K_1^2(v) dv}},
$$
 and $\widetilde{f}_{jn}(\epsilon)$ being at in the beginning of the
subsection, in the comparison of the bias and variances. To
perform the test, we consider $T$ independent replications
$$
\widetilde{Z^*}_{jk}(\epsilon)
 =
\frac{\sqrt{n\widetilde{b}_j(\epsilon)}
\left(\widetilde{f^*}_{jk}(\epsilon)-f(\epsilon)\right)}
{\sqrt{f(\epsilon)\int\! K_1^2(v) dv}}, \quad k=1,\ldots,T,
$$
of the variables $\widetilde{Z}_{jn}(\epsilon)$. We denote by
$\widetilde{D}_{jn}(\epsilon)$ the Kolmogorov-Smirnov statistic
associated with the test. For our goodness of fit test, the null
hypothesis is rejected with a level of significance $\alpha$ if
$\sqrt{T}\widetilde{D}_{jn}(\epsilon)>K_{\alpha}$, where
$K_{\alpha}$ satisfies
$$
\prob\left(
\sqrt{T}
\widetilde{D}_{jn}(\epsilon)\leq
K_{\alpha}\right)
=
\prob\left(\widetilde{D}_{jn}(\epsilon)\leq
\frac{K_{\alpha}}{\sqrt{T}}\right)
=
1-\alpha.
$$
 The tables of critical values of
the goodness of fit test to the standard normal variable can be
found in the statistic literature. See, for example, Smirnov
(1948), Miller (1956), Gibbons and Chakraborti (2003). Some of the
results for the asymptotic approximations of the critical value
$K_{\alpha}$ based on the ration $C_{\alpha}=K_{\alpha}/\sqrt{T}$
are:

\begin{table}[H]
\begin{center}
\begin{tabular}{|c|c|c|c|c|c|}
 \hline
 $\prob\left(\widetilde{D}_j(\epsilon)>C_{\alpha}\right)$ &
  $0.20$ & $0.15$ & $0.10$ & $0.05$ & $0.01$
 \\
 \hline
  $K_{\alpha}$ & $1.07$ & $1.14$ &
  $1.22$ & $1.36$ & $1.63$
 \\
\hline
\end{tabular}
\end{center}
\end{table}
\noindent In Table \ref{Critics},  we give the  values of the
statistics $\sqrt{T}\widetilde{D}_{jn}(\epsilon)$ for $T=100$,
$n=200$, $j=1,2$ and $\epsilon=-1,0,1$. The results obtained here
show that for the level $\alpha=0.05$, the null hypothesis
$\widetilde{Z}_{jn}(\epsilon)\sim N(0,1)$ is rejected, since
$\sqrt{T}\widetilde{D}_{jn}(\epsilon)>K_{\alpha}=1.36$, for all
$j$ and $\epsilon$.

\begin{table}[H]
\caption{\label{Critics} Values of the statistics
$\sqrt{T}\widetilde{D}_{jn}(\epsilon)$ for $T=100$, $j=1,2$ and
$\epsilon=-1,0,1$.}
 \small{
\begin{center}
\tabcolsep=3pt
\begin{tabular}{|c|c|c|c|c|c|}
 \hline
 \multicolumn{2}{|c|}{$\epsilon=-1$}&
\multicolumn{2}{|c|}{$\epsilon=0$} &
\multicolumn{2}{|c|}{$\epsilon=1$}
\\
\cline{1-6}
 $\sqrt{T}\widetilde{D}_{1n}(\epsilon)$ &
$\sqrt{T}\widetilde{D}_{2n}(\epsilon)$ &
$\sqrt{T}\widetilde{D}_{1n}(\epsilon)$ &
$\sqrt{T}\widetilde{D}_{2n}(\epsilon)$ &
$\sqrt{T}\widetilde{D}_{1n}(\epsilon)$ &
$\sqrt{T}\widetilde{D}_{2n}(\epsilon)$
 \\
  \hline
  3.159609 & 3.354464 &
  2.780215 & 2.676096 &
  2.465744 & 1.890398
  \\
 \hline
 \end{tabular}
\end{center}
}
\end{table}

We now attempt to explain the non-validity of the hypothesis
$\widetilde{Z}_{jn}(\epsilon)\sim N(0,1)$ by computing the
empirical mean $\overline{\widetilde{\mu}}_{jn}(\epsilon)$ and the
empirical variance
$\overline{\widetilde{\sigma}}_{jn}^2(\epsilon)$ of the data
$\widetilde{Z}_{jn}(\epsilon)$.
% These empirical
%parameters are grouped in the following table.

\begin{table}[H]
\caption{\label{MVZT} Values of the empirical means and variances
of the $\widetilde{Z}_{jn}(\epsilon)$'s.}
 \small{
\begin{center}
\tabcolsep=3pt
\begin{tabular}{|c|c|c|c|c|c|c|c|c|c|c|c|}
 \hline
 \multicolumn{4}{|c|}{$\epsilon=-1$}&
\multicolumn{4}{|c|}{$\epsilon=0$} &
\multicolumn{4}{|c|}{$\epsilon=1$}
\\
\cline{1-12}
 $\overline{\widetilde{\mu}}_{1n}(\epsilon)$ &
$\overline{\widetilde{\sigma}}_{1n}^2(\epsilon)$ &
$\overline{\widetilde{\mu}}_{2n}(\epsilon)$ &
$\overline{\widetilde{\sigma}}_{2n}^2(\epsilon)$ &
$\overline{\widetilde{\mu}}_{1n}(\epsilon)$ &
$\overline{\widetilde{\sigma}}_{1n}^2(\epsilon)$ &
$\overline{\widetilde{\mu}}_{2n}(\epsilon)$ &
$\overline{\widetilde{\sigma}}_{2n}^2(\epsilon)$ &
$\overline{\widetilde{\mu}}_{1n}(\epsilon)$ &
$\overline{\widetilde{\sigma}}_{1n}^2(\epsilon)$ &
$\overline{\widetilde{\mu}}_{2n}(\epsilon)$ &
$\overline{\widetilde{\sigma}}_{2n}^2(\epsilon)$
\\
\hline
  -0.4402 & 0.5322 & -0.4044 & 0.4163 &
  -0.4300 & 0.5743 & -0.2880 & 0.4685 &
   -0.3274 & 0.6124 & -0.0773 & 0.5473
 \\
\hline
\end{tabular}
\end{center}
}
\end{table}
\noindent Table \ref{MVZT} shows that the estimated quantities
$\overline{\widetilde{\mu}}_{jn}(\epsilon)$ and
$\overline{\widetilde{\sigma}}_{jn}^2(\epsilon)$ are clearly
different to $0$ and $1$. This should explain the rejection of the
hypothesis $\widetilde{Z}_{jn}(\epsilon)\sim N(0, 1)$, as seen
above. Hence the results of our simulation study reveal that with
the optimal step bandwitdhs $(\widehat{b}_{0j}(\epsilon),
\widehat{b}_{1j}(\epsilon))$ and $\widetilde{b}_j(\epsilon)$, the
variables $\widehat{Z}_{jn}(\epsilon)$ and
$\widetilde{Z}_{jn}(\epsilon)$ are not distributed according to
the standard distribution $N(0,1)$. However, the impact of the
estimated optimal first-step bandwidths
$\widehat{b}_{0j}(\epsilon)$ on the asymptotic normality of the
variables $\widehat{Z}_{jn}(\epsilon)$ may not be so important as
augured by the results obtained with the data
$\widehat{Z}_{jn}(\epsilon)$.

\newpage

\section{Conclusion}
The aim of this subsection was to analyze and compare the
performances of the Kernel density estimator $\widehat{f}_{1n}$,
based on the estimated residuals, and the ones of the integral
Kernel estimator $\widehat{f}_{2n}$. Several aspects have been
noticed in our simulation study.
 First, in the global framework, our numerical results
show that the estimator $\widehat{f}_{1n}$ should be
 preferred to the estimator $\widehat{f}_{2n}$.
The reason is  that the optimal AISE of the latter estimator is
much more higher than the one of the first estimator. For the
evaluation of  the bandwidths
$(\widehat{b}_{1j},\widehat{b}_{0j})$ that minimize the AISE of
the estimators $\widehat{f}_{jn}$ $(j=1,2)$, our numerical results
indicates that $\widehat{b}_{01}$ is much smaller than
$\widehat{b}_{11}$, and that $\widehat{b}_{02}$ is approximately
as small as $\widehat{b}_{12}$.

\vskip 0.3cm
 Next, for the pointwise study which is made at the points
 $\epsilon=-1,0,1$,
we observe that
 $\widehat{f}_{1n}(\epsilon)$ dominates
$\widehat{f}_{2n}(\epsilon)$ for $\epsilon=-1$ and $\epsilon=1$ as
well as for the ASE, the bias and the variance.
 Further,  the ASE of the estimators
$\widehat{f}_{jn}(\epsilon)$ are nearly the same
 as the  ones of the estimators $\widetilde{f}_{jn}(\epsilon)$ for
$\epsilon=0$, and lower than the ASE of
  $\widetilde{f}_{jn}(\epsilon)$ for $\epsilon=-1$ and $\epsilon=1$.
In a semiparametric context, Müller, Schick and
 Wefelmeyer (2004) have shown that for the estimation of linear functionals
 of the error distribution, the estimators that use the
estimated residuals may have a smaller asymptotic variance
compared to the estimators based upon the true errors. Some of our
simulation results suggest that a similar conclusion may hold when
estimating the p.d.f. of regression residuals. In fact,
 for $\epsilon=1$,
the variances of the estimators $\widetilde{f}_{jn}(\epsilon)$ are
higher than the ones of the estimators
$\widehat{f}_{jn}(\epsilon)$. This shows that the estimation of
the first-step bandwidth $b_0$ may have a positive influence when
estimating $f(\epsilon)$.

\vskip 0.3cm
 The study of the asymptotic normality of the
standardized variables $\widehat{Z}_{jn}(\epsilon)$ and
$\widetilde{Z}_{jn}(\epsilon)$, based on the  density estimators
$\widehat{f}_{jn}(\epsilon)$ and $\widetilde{f}_{jn}(\epsilon)$,
reveals that the data $\widehat{Z}_{jn}(\epsilon)$ and
$\widetilde{Z}_{jn}(\epsilon)$ are normal, but are not distributed
according to the standard normal variable. This means that the
normal approximation of these variables by the normal $N(0,1)$ is
not satisfying for a small size of the samples ($n=200$ in our
framework). Therefore, it will be interesting, in a future works,
to use the boostrap method for obtaining an alternative
approximation of the considered variables. This will be one of the
main aspects of the perspectives of our future researches, as
illustrated at the end of this thesis.

\chapter[Appendix] {Appendix}

\setcounter{subsection}{0} \setcounter{equation}{0}
\renewcommand{\theequation}{\thesection.\arabic{equation}}

\renewcommand{\thefootnote}{\arabic{footnote}}
\setcounter{footnote}{1} \setlength{\baselineskip}{.26in}

{\bf Abstract:} This chapter contains some  results which have an
interest themselves and are used in  Chapter 3 and Chapter 4. We
begin with the  Lyapounov Central Limit Theorem for triangular
arrays which is used, for example, in the proof of Proposition
\ref{Prop}  and  Theorem \ref{normalite}.
  We also recall Theorem 1 and Theorem 2 in Einmahl and Mason (2005).
  These results are need in the validation of Lemma \ref{Estig}.
  We conclude by Theorem 2 in Whitlle (1960) and the
Marcinkievicz-Zygmund inequality (see e.g Chow and Teicher 2003,
p. 386) which are very useful for proving  Lemma
\ref{BoundEspmchap}  and Lemma \ref{BoundEspmin}.

\section{Lyapounov's Central Limit Theorem}

For each integer $n\geq 1$, let $\left\lbrace
X_{1n},X_{2n},\ldots,X_{nn}\right\rbrace$ be a collection of
random variables such that $X_{1n},X_{2n},\ldots,X_{nn}$ are
independent. Then $\left\lbrace
X_{1n},X_{2n},\ldots,X_{nn}\right\rbrace$ is called a triangular
array of independent variables.

\begin{theorem} (Lyapounov's Theorem)
\\ For all integer $n\geq 1$, assume that the variables $X_{in}$,
$1\leq i\leq n$, are independent with $\esp\left[X_{in}\right]=0$
for all $i$. Let
$\alpha_n=\sqrt{\sum_{i=1}^n\Var\left(X_{in}\right)}$. If there
exists $\delta>0$ such that
$$
\lim_{n\rightarrow\infty} \alpha_n^{-(2+\alpha)}
\sum_{i=1}^n\esp\left[|X_{in}|^{2+\alpha}\right]=0,
$$
then
$$
\frac{X_{1n}+X_{2n}+\ldots+X_{nn}}
{\sqrt{\sum_{i=1}^n\Var\left(X_{in}\right)}}
\stackrel{d}{\longrightarrow} N(0,1)
$$
when $n\rightarrow\infty$.
\end{theorem}
\vskip 0.1cm\noindent
This result can be found, for example, in
Billingsley (1968, Theorem 7.3).

\section{Uniform in bandwidth consistency of kernel-type function estimators}

In this section, we give two results concerning the uniform in
bandwidth consistency of kernel-type estimators, such that the
density estimator and the regression function estimator. The
results proposed here are established in Einmahl and Mason (2005).
They are one of the keys of our main results in Chapter 3 and
Chapter 4.

\vskip 0.1cm
The first result we give concerns the Kernel density estimator.
 Let $X_1,X_2\ldots,X_n$ be i.i.d $\Rit^d$, $d\geq 1$,
valued random variables and
assume that the common distribution function of the variables has a Lebesgue
 density function, which we denote by $f$. The Kernel density estimator of $f$
  based upon the sample $X_1,X_2,\ldots,X_n$, a Kernel function $K$ and a
  bandwidth $0<h=h(n)<1$ is defined as
  $$
  \widehat{f}_{n,h}(x)
  =
  \frac{1}{nh}\sum_{i=1}^n
  K\left(\frac{x-X_i}{h^{1/d}}\right),
  \quad x\in\Rit^d.
  $$
For any  function $G$ defined and bounded on $\Rit^d$, we denote
by $\|G\|_{\infty}$ the uniform norm of $G$ such that
$$
\|G\|_{\infty}
=
\sup_{x\in\Rit^d}\left|G(x)\right|.
$$
The following theorem is proposed Einmahl and Mason (2005, p. 1382).
\begin{theorem}
Assume that the Kernel function $K$ is symmetric, continuous over $\Rit^d$ with
support contained in $[-1/2,1/2]^d$ and $\int\! K(x)dx=1$.
If the density function $f$
is continuous and bounded on its support, then we have for any $C>0$,
with probability $1$,
$$
\limsup_{n\rightarrow\infty}
\sup_{C\left(\ln(n)/n\right)\leq h\leq 1}
\|\widehat{f}_{n,h}-\esp\widehat{f}_{n,h}\|_{\infty}
=
O\left(
\frac{\sqrt{\ln\left(1/h\right)
\vee\ln\left(\ln n\right)}}{nh}
\right).
$$
\label{EM1}
\end{theorem}

\noindent{\bf Remark:} Choosing a sequence $h=h(n)$ satisfying
$(nh/\ln n)\rightarrow\infty$ and $\ln(1/h)/\ln\left(\ln
n\right)\rightarrow\infty$, one obtains, with probabilty 1,
$$
 \|\widehat{f}_{n,h}-\esp\widehat{f}_{n,h}\|_{\infty} =
O\left(\sqrt{\left(\ln\left(1/h\right)\right)/(nh)}\right),
$$
which is Theorem 1 of Giné and Guillou (2005).

\vskip 0.3cm The other kinds of kernel-type estimators treated by
Einmahl and Mason is the regression Kernel estimators. For the
illustration, consider i.i.d $(d+1)$-dimensional random vectors
$(X,Y),(X_1,Y_1),(X_2,Y_2),\ldots,(X_n,Y_n)$, where the
$Y$-variables are one-dimensional. We assume that $X$ has a
marginal Lebesgue density function $f$ and that the regression
function
$$
m(x)=\esp\left[Y\mid X=x\right],
\quad x\in\Rit^d.
$$
exists. Let $\widehat{m}_{n,h}(x)$ be the Nadaraya-Watson estimator
of $m(x)$ with bandwidth
$0<h<1$, that is,
$$
\widehat{m}_{n,h}(x)
=
\frac{\sum_{i=1}^n Y_i K\left((x-X_i)/h^{1/d}\right)}
{\sum_{i=1}^n K\left((x-X_i)/h^{1/d}\right)}\;.
$$
With the above setup, we have the following uniform in bandwidth result.
 Let $K$ and $h$ be as in
the previous section, and set
$$
\overline{r}(x,h)
=
h^{-1}\esp\left[Y K\left(\frac{x-X}{h^{1/d}}\right)\right],
\quad
\overline{f}(x,h)
=
h^{-1}\esp\left[K\left(\frac{x-X}{h^{1/d}}\right)\right].
$$
For any subset $I$ of $\Rit^d$, let $I^{\epsilon}$ denote its
closed $\epsilon$-neighborhood with respect
to the maximum norm $|\cdot|_{+}$ on $\Rit^d$, that is,
$|x|_{+}=\max_{1\leq i\leq n}|x_i|$, $x\in\Rit^d$.
Set further for any function $\psi:\Rit^d\rightarrow\Rit$,
$\|\psi\|_{I}=\sup_{x\in\Rit}\left|\psi(x)\right|$.

\begin{theorem} {\rm (Einmahl and Mason 2005, p. 1384)}
\\ Let $I$ be a compact subset of $\Rit^d$ of $\Rit^d$ and assume
that the Kernel function $K$ satisfies
the condition of Theorem \ref{EM1}. Suppose further that there
 exists an $\epsilon>0$ so that
$f$ is continuous and strictly positive on $J:=I^{\epsilon}$. If
we assume  that for some $p>2$,
$$
\sup_{z\in J}\esp\left(|Y|^p\mid X=z\right)
:=\alpha<\infty,
$$
we have for any $C>0$ and $b_n\searrow 0$ with $\gamma=\gamma(p)=1-2/p$,
$$
\limsup_{n\rightarrow\infty}
\sup_{C\left(\ln(n)/n\right)^{\gamma}\leq h\leq b_n}
\|\widehat{m}_{n,h}-\overline{r}(\cdot,h)/\overline{f}(\cdot,h)\|_{I}
=
O\left(
\frac{\sqrt{\ln\left(1/h\right)
\vee\ln\left(\ln n\right)}}{nh}
\right),
$$
almost surely.
\end{theorem}

\section{Bounds for the moments of linear forms in independent
variables}

The aim of this section is to propose absolute moments of linear
forms in independent statistical variables. The first result we
give here is
 established by Whitlle (1960, Theorem 2).
Consider the linear form $L=\sum_{j=1}^n a_j\zeta_j$, where the
$\zeta_j$'s are
 assumed to be independent mean-zero random variables, but not necessarily
 to be distributed identically. In what follows, we shall write
  $$
 C(p)
 =
 \frac{2^{p/2}}{\sqrt{\pi}}
 \int_{-\infty}^{+\infty}
 |x|^p e^{-x^2} dx,
 $$
and $\gamma_j(p)=\left(\esp\left|\zeta_j\right|^p\right)^{1/p}$,
$p>0$, provided that
 these quantities exist.

\begin{theorem} {\rm (Whittle, 1960)}
 \\Then the following inequality is valid
 $$
 \esp\left(|L|^p\right)
 \leq 2^pC(p)\left(\sum_{j=1}^p\gamma_j^2(p)\right)^{p/2},
 $$
provided that $p\geq 2$ and the right-hand member exists.
Moreover, if all the $\zeta_j$ have symmetric distributions, then
the right-hand member may be divided by $2^p$.
\end{theorem}

The second result we give is the Marcinkiewicz-Zygmund Inequality
(See Chow and Teicher 2003, p.386). For any $p\geq 1$, let
$\|\cdot\|_p$ denotes the $L^p$-norm, that is,
$\left\|X\right\|_p=\left(\esp|X|^p\right)^{1/p}$ for any random
variable $X$ such that $\esp\left(|X|^p\right)<\infty$.

\begin{theorem} {\rm Marcinkiewicz-Zygmund Inequality}
\\ If $\{X_n, n\geq 1\}$ are independent random variables
with $\esp[X_n]=0$ for all $n$, then for every $p\geq 1$, there
exist positive constant
 $A_p$ and $B_p$ depending only
upon $p$ for which
$$
A_p\left\|\left(\sum_{j=1}X_j^2\right)^{1/2}\right\|_{p}
\leq\left\|\sum_{j=1}X_j\right\|_{p}\leq
B_p\left\|\left(\sum_{j=1}X_j^2\right)^{1/2}\right\|_{p}.
$$
\end{theorem}

\noindent
The proof of this Theorem can be found, for example, in
Chow and Teicher (2003, p. 386).

\chapter*{Perspectives}
\addcontentsline{toc}{chapter}{Perspectives}
\section*{Abstract}
\noindent In this section, we sketch some perspectives for
possible future researches. First, we  have seen in
  our simulation study  that the
estimator of $f(\epsilon)$ introduced in Chapter 3 would be
preferred to the one proposed in Chapter 4. However, it would be
very interesting to compare the theoretical bias of the two
estimators for determining the estimator that have to be used in a
given context.

\vskip 0.3cm
 Our numerical results also reveal  a curious
situation: the estimator $\widehat{f}_{jn}(\epsilon)$ ($j=1,2$) is
sometimes more efficient than the estimator
$\widetilde{f}_{jn}(\epsilon)$ when we are interested in their
pointwise study. This situation comes from the evaluation of the
second order of $\widehat{f}_{jn}(\epsilon)$, that is
$\widehat{f}_{jn}(\epsilon)-\widetilde{f}_{jn}(\epsilon)$, which
possibly allows to improve the performances of
$\widehat{f}_{jn}(\epsilon)$. This curious siuation makes one to
think that the term
$\widehat{f}_{jn}(\epsilon)-\widetilde{f}_{jn}(\epsilon)$ is worth
thinking about and deserved further consideration. We shall also
attempt to obtain the uniform weak consistency for the difference
 $\widehat{f}_{jn}(\epsilon)-\esp_n\widehat{f}_{jn}(\epsilon)$.

\vskip 0.3cm All the results proposed in this thesis are obtained
in the case of a homoscedastic model. Then another axis for future
researches will concern the extension of our results in a
heteroscedastic framework, when the variance function depends upon
the explanatory variable.

\section*{Résumé}
\noindent Dans cette partie, nous donnons une esquisse des
perspectives de recherche pour nos futurs travaux. D'abord, les
résultats de nos simulations numériques montrent que l'estimateur
de $f(\epsilon)$ introduit au Chapitre 3 devrait être préferé à
celui défini au Chapitre 4. Cependant, il serait intéressant de
comparer de façon théorique les biais des deux estimateurs. Ce
sera l'un des problèmes sur lesquels nous nous pencherons dans nos
recherches ultérieures.

\vskip 0.3cm Les résultats de nos simulations montrent également
un point assez curieux: l'estimateur $\widehat{f}_{jn}(\epsilon)$
($j=1,2$) est parfois plus efficace que l'estimateur
$\widetilde{f}_{jn}(\epsilon)$ lorsqu'on les étudie
ponctuellement. Cette situation est due au second ordre de
$\widehat{f}_{jn}(\epsilon)$, c'est à dire
$\widehat{f}_{jn}(\epsilon)-\widetilde{f}_{jn}(\epsilon)$, qui
permet éventuellement d'améliorer les performances de
$\widehat{f}_{jn}(\epsilon)$. Ce deuxième ordre mériterait d'être
étudié de façon plus poussée. Nous tenterons aussi d'obtenir des
résultats de consistance uniforme pour la différence entre
$\widehat{f}_{jn}(\epsilon)$ et
$\esp_n\widehat{f}_{jn}(\epsilon)$.

\vskip 0.3cm Tous les résultats proposés dans cette thèse ont été
obtenus dans un modèle de régression homoscédastique. Un autre axe
de recherche pour nos futurs travaux sera de voir si des résultats
comparables peuvent être obtenus dans le cas du modèle
hétéroscédastique, où l'erreur du modèle dépend de la variable
explicative.

\end{document}